\documentclass[11pt]{amsart}
\usepackage[linkcolor=blue, citecolor = blue]{hyperref}
\hypersetup{colorlinks=true}
\usepackage{amsmath, amsfonts, amsthm, amssymb} 
\usepackage{tikz}
\usepackage{ytableau}
\usepackage{caption}
\usepackage{subcaption}
\usepackage{graphicx} 
\usepackage[makeroom]{cancel}
\usepackage[margin=1.0in]{geometry}
\usepackage{enumitem}

\usetikzlibrary{arrows.meta}

\definecolor{red 1}{HTML}{B8548D} 
\definecolor{red 2}{HTML}{CCE0EB} 
\definecolor{red 3}{HTML}{EAADAC} 
\definecolor{red 4}{HTML}{C98181} 
\definecolor{red 5}{HTML}{B8545B} 
\definecolor{red 6}{HTML}{69A5C3} 
\definecolor{red 7}{HTML}{7F8CCC} 
\definecolor{red 8}{HTML}{80BF80} 

\definecolor{sblack}{RGB}{0,0,0}
\definecolor{sorange}{RGB}{230,159,0}
\definecolor{sskyblue}{RGB}{86,180,233}
\definecolor{sbluishgreen}{RGB}{0,158,115}
\definecolor{syellow}{RGB}{240,228,66}
\definecolor{sblue}{RGB}{0,114,178}
\definecolor{svermilion}{RGB}{213,94,0}
\definecolor{sredpurple}{RGB}{204,121,167}

\newtheorem*{theorem*}{Theorem}
\newtheorem{theorem}{Theorem}[section]
\newtheorem{lemma}[theorem]{Lemma}
\newtheorem{proposition}[theorem]{Proposition}
\newtheorem*{corollary*}{Corollary}
\newtheorem{corollary}[theorem]{Corollary}

\theoremstyle{definition}
\newtheorem{remark}[theorem]{Remark}
\newtheorem{definition}[theorem]{Definition}
\newtheorem{example}[theorem]{Example}

\usetikzlibrary{decorations.markings}
\usetikzlibrary{decorations.pathmorphing}
\usetikzlibrary{matrix,decorations.pathreplacing}

\pgfkeys{tikz/mymatrixenv/.style={decoration=brace,every left delimiter/.style={xshift=4pt},every right delimiter/.style={xshift=-4pt}}}
\pgfkeys{tikz/mymatrix/.style={matrix of math nodes,left delimiter=[,right delimiter={]},inner sep=1pt,row sep=0em,column sep=0em,nodes={inner sep=6pt}}}

\title[Flagged Grothendieck enumerators and ribbon decomposition determinants]
{Flagged Grothendieck enumerators and ribbon decomposition determinants}
\author[Alibek Adilzhan \and Damir Yeliussizov]{Alibek Adilzhan \and Damir Yeliussizov}
\address{Kazakh-British Technical University, Almaty, Kazakhstan}
\email{\href{mailto:alibek.adilzhanov@gmail.com}{alibek.adilzhanov@gmail.com}}
\email{\href{mailto:yeldamir@gmail.com}{yeldamir@gmail.com}}

\begin{document}
\tikzset{
    -->-/.style={decoration={
    markings,
    mark=at position .9 with {\arrow{>}}},postaction={decorate}},
    ->-/.style={decoration={
    markings,
    mark=at position .7 with {\arrow{>}}},postaction={decorate}}
}

\begin{abstract}
We establish Hamel--Goulden-type ribbon decomposition determinantal formulas for a new class of symmetric functions called flagged Grothendieck enumerators. In particular, these functions generalize both canonical refined stable Grothendieck polynomials and their duals. The proof uses Talaska's formula via flows on planar lattice graphs with systems of paths and cycles. 
\end{abstract}

\maketitle

\setcounter{tocdepth}{1}
\tableofcontents

\section{Introduction}
In this paper, we develop Hamel--Goulden-type ribbon decomposition determinantal formulas for a new class of symmetric functions which we call {\it flagged Grothendieck enumerators}. In particular, these functions generalize both stable Grothendieck polynomials and their duals. 

\medskip

Stable Grothendieck polynomials are $K$-theoretic extensions of Schur polynomials which were originally studied by Fomin and Kirillov \cite{fk1, fk2} as a special case of Grothendieck polynomials due to Lascoux and Sch\"utzenberger \cite{ls0} (the latter are $K$-theoretic extensions of Schubert polynomials). These symmetric functions possess a number of remarkable properties and they were extensively studied from algebraic, combinatorial, geometric, representation theoretic, and  probabilistic perspectives, see \cite{fk1, fk2, lenart, buch, lp, motegi, gal,  yel17, yel19, hs, yel20, cp, pech, gp, ims, ims2} for some references. 

\medskip

We shall extend one of their most general known versions given by {\it canonical refined stable Grothendieck polynomials} 
and their duals 
studied in \cite{hwang1, hwang2}, which are equipped with two extra sets of refined parameters that lead to nice symmetries and duality. These functions generalize (dual) refined stable Grothendieck polynomials studied in \cite{ggl, cp}, and canonical versions studied in \cite{yel17}. In fact, we introduce more general {\it flagged Grothendieck enumerators} indexed by skew signatures and establish Hamel--Goulden-type formulas for them. Notably, these functions generalize both families of canonical refined stable Grothendieck polynomials and their duals (via the Hall inner product), as well other some other quite large families of symmetric functions. 

\medskip

While stable Grothendieck polynomials have some known Jacobi--Trudi and bialternant type determinantal formulas, see \cite{hwang1, hwang2} and references therein, no other types of determinantal formulas (such as e.g. Giambelli or Lascoux--Pragacz) were known for them. 
For Schur functions, the most general unifying determinantal formulas are due to Hamel--Goulden \cite{hg}, which arise via various ribbon decompositions of the underlying shape. Such formulas remarkably unify Jacobi--Trudi, Giambelli, and Lascoux--Pragacz \cite{ls, lasp} formulas. 

\medskip
 
Hamel--Goulden-type formulas for stable Grothendieck polynomials can be useful for several reasons: (i) algebraic, as they can reveal some structural properties; (ii) combinatorial, as they provide interesting lattice and tableaux constructions; (iii) computational, as they can be computationally efficient in certain cases. 

\subsection{Flagged Grothendieck enumerators}
We introduce a new family of symmetric functions called {\it flagged Grothendieck enumerators} denoted by $\mathsf{G}^{\mathbf{a}, \mathbf{b}}_{\lambda / \mu}(\mathbf{x}; \mathbf{U}, \mathbf{V})$ for variable set $\mathbf{x} = (x_i)$, parameter sets of matrices $\mathbf{U} = (u_{ij}), \mathbf{V} = (v_{ij})$, and defined as weighted sums via {\it trio tableaux} which are triples of semistandard tableaux with certain conditions, see \textsection\,\ref{sec:fg} for definitions. These power series are symmetric in $\mathbf{x}$, they are indexed  by {\it skew signatures} $\lambda/\mu$ (see definitions in \textsection\,\ref{sec:2}) and have flag conditions $\mathbf{a}, \mathbf{b}$ imposed on the parameter variable sets $\mathbf{U}, \mathbf{V}$. 

\subsubsection{Some special cases} Flagged Grothendieck enumerators generalize the following known families of symmetric functions. 

\begin{itemize}

\item 
{\bf Canonical refined stable Grothendieck polynomials.} 
For partitions $\lambda, \mu$ of length at most $n$ and flags given by $a_i = i$, $b_i = \lambda_i$, flagged Grothendieck enumerators specialize to {\it canonical refined stable Grothendieck polynomials} $\{G_{\lambda/\mu}(\mathbf{x} ; \boldsymbol{\alpha}, \boldsymbol{\beta}) \}$ (studied in \cite{hwang1, hwang2})
$$
G_{\lambda/\mu}(\mathbf{x} ; \boldsymbol{\alpha}, \boldsymbol{\beta}) = C \cdot 
\mathsf{G}^{\mathbf{a}, \mathbf{b}}_{\lambda/\mu}(\mathbf{x} ; \mathbf{U}, \mathbf{V} ) |_{G\text{-spec.}}, 
\qquad C = \prod_{i \in [n]} \prod_{j } (1 - \beta_i x_j),
$$
where by $G$-specialization we mean the following parameter substitutions: $u_{i,j} \to \alpha_i - \beta_{j}, v_{i,j} \to \beta_{i} - \alpha_{j}$. 

\item 
{\bf Dual canonical refined stable Grothendieck polynomials.} 
For partitions $\lambda, \mu$ of length at most $n$ and flags given by $a_i = i-1$, $b_i = \mu_i$, flagged Grothendieck enumerators specialize to {\it dual canonical refined stable Grothendieck polynomials} $\{g_{\lambda/\mu}(\mathbf{x} ; \boldsymbol{\alpha}, \boldsymbol{\beta}) \}$ (studied in \cite{hwang1, hwang2})
$$
g_{\lambda/\mu}(\mathbf{x}_m; \boldsymbol{\alpha}, \boldsymbol{\beta}) = \mathsf{G}^{\mathbf{a}, \mathbf{b}}_{\mu/\lambda}(\varnothing ; \mathbf{U}, \mathbf{V} ) |_{g-\text{spec.}},
$$
where by $g$-specialization we mean the following parameter substitutions: $u_{i,j} = \alpha_i - \beta_j$, $v_{i,j} = z_i - \alpha_{j-m}$ and $z_i = x_i$ for $i \in [m]$, $z_i = \beta_{i - m}$ for $i > m$. 


\item 
{\bf Factorial supersymmetric Schur functions. }
For partitions $\lambda, \mu$ of length at most $n$ and flags given by $a_i = m$, $b_i = m + \mu_i + 1 - j$, flagged Grothendieck enumerators specialize to {\it factorial supersymmetric Schur functions} $s_{\lambda/\mu}(\mathbf{x} / \mathbf{y} | \mathbf{z})$ (studied in \cite{mol98}) 
$$
s_{\lambda/\mu}(\mathbf{x}_m / \mathbf{y}_{m'} \,|\, \mathbf{z}) = \mathsf{G}^{\mathbf{a},  \mathbf{b}}_{\mu/\lambda}(\varnothing; \mathbf{U}, \mathbf{V}) |_{f\text{-spec.}},
$$
where by $f$-specialization we mean the following parameter substitution: $u_{i,j} = y_j + z_i,$ $v_{i,j} = x_i - z_j$. 

\item
{\bf Macdonald's ninth variation of Schur functions. }
For partitions $\lambda, \mu$ of length at most $n$ and flags given by $a_i = m$, $b_i = m + \mu_i + 1 - i$, flagged Grothendieck enumerators specialize to {\it Macdondald's ninth variation of Schur functions} $s_{\lambda/\mu}(\mathbf{X})$ with $\mathbf{X} = (x_{i,j})_{i \in [m], j \in \mathbb{Z}}$ (studied in \cite{macdonald1, fkingb})
$$
s_{\lambda/\mu}(\mathbf{X}) = \mathsf{G}^{\mathbf{a},  \mathbf{b}}_{\mu/\lambda}(\varnothing; \mathbf{0}, \mathbf{V})_{s\text{-spec.}},
$$
where by $s$-specialization we mean the following parameter substitution: $v_{i,j} = x_{i, j - i}$.
\end{itemize}

Let us note that each of the above families is quite large and covers a number of other known families of symmetric functions. 

\subsubsection{Duality with column flagged functions} There are also {\it column} flagged Grothendieck enumerators $\mathsf{G}^{\mathbf{a}', \mathbf{b}'}_{\lambda'/\mu'}$ (given in \textsection\,\ref{sec:dual}) which satisfy the following duality under the standard ring involution $\omega$ of the ring of symmetric functions extended to completion,
$$
\omega_{\mathbf{x}}(\mathsf{G}^{\mathbf{a}, \mathbf{b}}_{\lambda/\mu}(\mathbf{x}; \mathbf{U}, \mathbf{V})) = \mathsf{G}^{\mathbf{a}', \mathbf{b}'}_{\lambda'/\mu'}(\mathbf{x}; \mathbf{U}^{T}, \mathbf{V}^{T}),
$$
where $\mathbf{U}^T$ denotes the transpose matrix, see Theorem~\ref{thm:omega}.

\subsection{Main result}
Our main result 
establishes Hamel--Goulden-type ribbon decomposition determinantal formulas for flagged Grothendieck enumerators 
and hence for canonical refined stable Grothendieck functions in particular. 

\begin{theorem}[Hamel--Goulden formulas for flagged Grothendieck numerators, cf. Theorem~\ref{thm:HG-sfG}]
    \label{thm:HG-sfG-intro}
    Let $\Theta = (\theta_1, \ldots, \theta_k)$ be an outer ribbon decomposition of skew signature $\lambda/\mu$ (see \textsection\,\ref{sec:rd}). Let $\mathbf{a}, \mathbf{b}$ be $\mathsf{G}$-flags and $\mathbf{a}^{ij}, \mathbf{b}^{ij}$ be induced ribbon $\mathsf{G}$-flags for the ribbons $\theta_i \# \theta_j$ (see \textsection\,\ref{sec:det}). Then the following formula holds:  
    \begin{align*}
        \mathsf{G}^{\mathbf{a}, \mathbf{b}}_{\lambda/\mu}(\mathbf{x} ; \mathbf{U} , \mathbf{V}) 
        = (-1)^{\mathsf{hp}^{\Theta}_{\lambda/\mu}} \cdot \mathsf{G}^{\mathbf{a}^{\varepsilon}, \mathbf{b}^{\varepsilon}}_{\varepsilon}(\mathbf{x} ;  \mathbf{U} , \mathbf{V})^{-k+1} \cdot \det \left[ (-1)^{\mathsf{hp}^{\Theta}_{\theta_i \# \theta_j}} \cdot \mathsf{G}^{\mathbf{a}^{ij}, \mathbf{b}^{ij}}_{\theta_i \# \theta_j}(\mathbf{x} ; \mathbf{U} , \mathbf{V})\right]_{1 \le i,j \le k}. 
    \end{align*}
    Here $\#$ is an extended signed Hamel--Goulden operation, 
    $\varepsilon$ is an empty ribbon contour which depends on $\Theta$, 
    and $\mathsf{hp}^{\Theta}_{\lambda/\mu}$ is a hook-pipe counting function (see \textsection\,\ref{sec:rd}). 
\end{theorem}

Let us emphasize that unlike the ordinary Hamel--Goulden formulas, our formula contains the nontrivial factor $\mathsf{G}^{\mathbf{a}^{\varepsilon}, \mathbf{b}^{\varepsilon}}_{\varepsilon}(\mathbf{x} ;  \mathbf{y} , \mathbf{z})^{-k+1}$ whose enumerator is computed on empty ribbon contour $\varepsilon$ depending on $\Theta$. 

In particular, as special cases we obtain both new and known Jacobi--Trudi-type formulas and new Giambelli-type formulas for canonical refined stable Grothendieck functions. 

We also note that this theorem generalizes Hamel--Goulden-type ribbon decomposition determinants for 
dual canonical refined stable Grothendieck polynomials, which we studied in \cite{ay25} and whose proofs are much simpler than for the above case of flagged Grothendieck enumerators.


\subsection{Key tools and highlights of the proof}
Most of the text is devoted to the proof of Theorem~\ref{thm:HG-sfG-intro}, which is after some preparations given in \textsection\,\ref{sec:proof} (a more detailed plan of the proof is outlined in \textsection\,\ref{sec:proofplan}). 
Our proof involves rather unusual tools and technical components, compared to other similar problems on symmetric functions.  
Let us highlight some of them: 
\begin{itemize}
\item {\bf Flows on planar lattice graphs with cycles.} We use Talaska's \cite{tal} generalization of the Lindstr\"om--Gessel--Viennot (LGV) lemma \cite{lgv} via flows on planar directed graphs with systems of paths and cycles. 
Let us note that such constructions have not arised in other known symmetric function determinants, e.g. known examples of Hamel--Goulden formulas (and their variations) are proved combinatorially using the usual LGV lemma, see e.g. \cite{hg, cll, cyy, fkingb}, including the case of dual stable Grothendieck polynomials studied by us in \cite{ay25}. 

\item {\bf Signed skew signatures. }
Our main functions are indexed by signed skew signatures $\lambda/\mu$ which generalize integer partitions, where we do {\it not} require the usual condition $\lambda \supseteq \mu$. In this approach we introduce the sign reciprocity $\lambda/\mu = - \mu/\lambda$ which is useful for obtaining e.g. a specialization to 
dual stable Grothendieck polynomials and some other families as shown above. 

\item {\bf Extended ribbon decompositions and operations.}
The usual Hamel--Goulden formulas are written via the known $\#$ operation for combining ribbons \cite{hg, cyy}. 
We need to extend this operation to work with skew signatures and signed shapes, see \textsection\,\ref{sec:rd} for details.  

\item {\bf Non-identity permutations.} 
One important technical challenge in establishing the needed weight-preserving bijection lies in path systems corresponding to non-identity permutations, which also contribute nontrivial monomial terms to our functions. This is carried out in detail in \textsection\,\ref{sec:proof}. 
\end{itemize}








\subsection{Outline of the paper}
In \textsection\,\ref{sec:2} we give preliminary new definitions on signed skew shapes used throughout the text. 
In \textsection\,\ref{sec:fg} we introduce flagged Grothendieck polynomials via trio tableaux and show how they specialize to canonical refined stable Grothendieck polynomials. 
In \textsection\,\ref{sec:rd} we introduce extended signed ribbon decompositions and defined extended Hamel--Goulden $\#$ operations that will be used the main determinantal formulas.
In \textsection\,\ref{sec:det} we define induced ribbon flags and state main Hamel--Goulden-type formulas for flagged Grothendieck enumerators. 
In \textsection\,\ref{sec:tal} we recall Talaska's formula and its main corollary used in the proof for our formulas. 
In \textsection\,\ref{sec:proof} we give proof of the main formula stated earlier. 
In \textsection\,\ref{sec:dual} we define column flagged Grothendieck enumerators, show duality with row flagged functions under the standard ring involution, and give corresponding Hamel--Goulden-type formulas for these functions. 
In \textsection\,\ref{sec:spec} we discuss some special cases of our formulas such as Jacobi--Trudi and skew Giambelli type determinants. 
In \textsection\,\ref{sec:spec1} we discuss some known families of symmetric functions which can be produced as specializations of flagged Grothendieck enumerators. 


\section{Signatures, null strips, and signed skew shapes}\label{sec:2}
We denote $[n] := \{1, \ldots, n\}$. 

We use the following generalized version of integer partitions. 
\begin{definition}[Signatures]
A {\it signature} of length $n$ is an integer sequence $\lambda = (\lambda_1 \ge \cdots \ge \lambda_n) \in \mathbb{Z}^n$. 
We denote by $\mathsf{Sign}_n$ the set of signatures of length $n$. 
For signatures $\lambda = (\lambda_i), \mu = (\mu_i) \in \mathsf{Sign}_n$ we define the following {\it skew diagram} 
given by the set 
$$
D(\lambda, \mu) := \{(i,j) : 1 \le i \le n,  \min(\mu_i, \lambda_i) < j \le \max(\mu_i, \lambda_i)\} \subset \mathbb{Z}^2.
$$
We write $\mu \subseteq \lambda$ when $\mu_i \le \lambda_i$ for all $i \in [n]$. 

An element $\gamma = (i, j) \in D(\lambda, \mu)$ is called a {\it cell} (or a box). 
The {\it content} of cell $\gamma = (i,j)$ is defined by $c(\gamma) = j - i$. For $\gamma = (i,j)$ we also write $\mathrm{row}(\gamma) = i$ and $\mathrm{col}(\gamma) = j$. 




A signature $\lambda \in \mathbb{Z}_{\ge 0}^n$ is called a {\it partition} and we denote by $\mathsf{Par}_n \subset \mathsf{Sign}_n$ for the corresponding set of partitions. 

\end{definition}

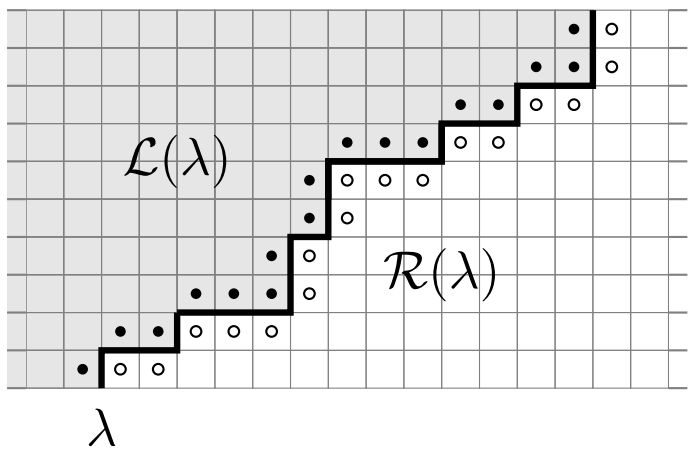
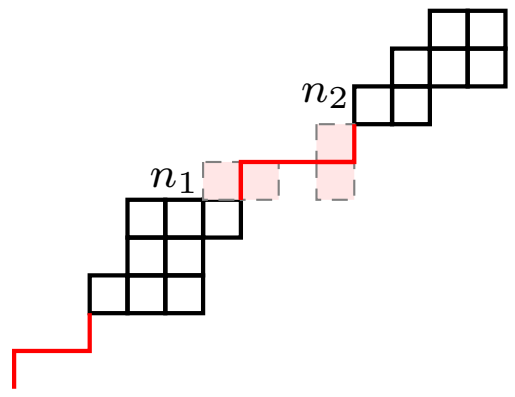
\begin{figure}
    \centering
    \begin{subfigure}[t]{0.45\textwidth}
    \centering
    \begin{tikzpicture}[y=-1cm, scale=0.25, font=\small]

        \foreach \x in {-4,...,10} {
            \draw[gray, very thin, fill=black!10] (\x, 1) rectangle ++(1, 1);
        }
        \foreach \x in {-4,...,10} {
            \draw[gray, very thin, fill=black!10] (\x, 2) rectangle ++(1, 1);
        }
        \foreach \x in {-4,...,8} {
            \draw[gray, very thin, fill=black!10] (\x, 3) rectangle ++(1, 1);
        }
        \foreach \x in {-4,...,6} {
            \draw[gray, very thin, fill=black!10] (\x, 4) rectangle ++(1, 1);
        }
        \foreach \x in {-4,...,3} {
            \draw[gray, very thin, fill=black!10] (\x, 5) rectangle ++(1, 1);
        }
        \foreach \x in {-4,...,3} {
            \draw[gray, very thin, fill=black!10] (\x, 6) rectangle ++(1, 1);
        }
        \foreach \x in {-4,...,2} {
            \draw[gray, very thin, fill=black!10] (\x, 7) rectangle ++(1, 1);
        }
        \foreach \x in {-4,...,2} {
            \draw[gray, very thin, fill=black!10] (\x, 8) rectangle ++(1, 1);
        }
        \foreach \x in {-4,...,-1} {
            \draw[gray, very thin, fill=black!10] (\x, 9) rectangle ++(1, 1);
        }
        \foreach \x in {-4,...,-3} {
            \draw[gray, very thin, fill=black!10] (\x, 10) rectangle ++(1, 1);
        }

        \foreach \x in {10.5, ..., 10.5} {
            \fill (\x, 1.5) circle (4pt) node {};
        }
        \foreach \x in {9.5,...,10.5} {
            \fill (\x, 2.5) circle (4pt) node {};
        }
        \foreach \x in {7.5,...,8.5} {
            \fill (\x, 3.5) circle (4pt) node {};
        }
        \foreach \x in {4.5,...,6.5} {
            \fill (\x, 4.5) circle (4pt) node {};
        }
        \foreach \x in {3.5,...,3.5} {
            \fill (\x, 5.5) circle (4pt) node {};
        }
        \foreach \x in {3.5,...,3.5} {
            \fill (\x, 6.5) circle (4pt) node {};
        }
        \foreach \x in {2.5,...,2.5} {
            \fill (\x, 7.5) circle (4pt) node {};
        }
        \foreach \x in {0.5,...,2.5} {
            \fill (\x, 8.5) circle (4pt) node {};
        }
        \foreach \x in {-1.5,...,-0.5} {
            \fill (\x, 9.5) circle (4pt) node {};
        }
        \foreach \x in {-2.5,...,-2.5} {
            \fill (\x, 10.5) circle (4pt) node {};
        }

        \foreach \x in {11.5, ..., 11.5} {
            \filldraw[fill=white, draw=black] (\x, 1.5) circle (4pt) node {};
        }
        \foreach \x in {11.5,...,11.5} {
            \filldraw[fill=white, draw=black] (\x, 2.5) circle (4pt) node {};
        }
        \foreach \x in {9.5,...,10.5} {
            \filldraw[fill=white, draw=black] (\x, 3.5) circle (4pt) node {};
        }
        \foreach \x in {7.5,...,8.5} {
            \filldraw[fill=white, draw=black] (\x, 4.5) circle (4pt) node {};
        }
        \foreach \x in {4.5,...,6.5} {
            \filldraw[fill=white, draw=black] (\x, 5.5) circle (4pt) node {};
        }
        \foreach \x in {4.5,...,4.5} {
            \filldraw[fill=white, draw=black] (\x, 6.5) circle (4pt) node {};
        }
        \foreach \x in {3.5,...,3.5} {
            \filldraw[fill=white, draw=black] (\x, 7.5) circle (4pt) node {};
        }
        \foreach \x in {3.5,...,3.5} {
            \filldraw[fill=white, draw=black] (\x, 8.5) circle (4pt) node {};
        }
        \foreach \x in {0.5,...,2.5} {
            \filldraw[fill=white, draw=black] (\x, 9.5) circle (4pt) node {};
        }
        \foreach \x in {-1.5,...,-0.5} {
            \filldraw[fill=white, draw=black] (\x, 10.5) circle (4pt) node {};
        }



        \foreach \y in {1, ..., 10} {
            \foreach \x in {-4,...,12} {
                \draw[gray, very thin] (\x, \y) rectangle ++(1, 1);
            }
        }

        \foreach \y in {1, ..., 10} {
            \fill[black!10] (-4.5, \y) rectangle ++(0.5, 1);
        }

        \draw[very thick] (4, 6) -- (4, 5) -- (7, 5) -- (7, 4);
        \draw[very thick] (0, 9) -- (3, 9) -- (3, 7) -- (4, 7) -- (4, 5) -- (7, 5) -- (7, 4) -- (9, 4) -- (9, 3) -- (11, 3) -- (11, 1);
        \draw[very thick] (0, 9) -- (0, 10) -- (-2, 10) -- (-2, 11);

        \foreach \x in {1,...,11} {
            \draw[gray, thin] (-4, \x) -- (-4.5, \x);
            }
        \foreach \x in {1,...,11} {
            \draw[gray, thin] (13, \x) -- (13.5, \x);
        }



        \fill (-2, 11) node[below] {$\lambda$};
        \fill (0, 5) node {$\mathcal{L}(\lambda)$};
        \fill (7, 8) node {$\mathcal{R}(\lambda)$};

    \end{tikzpicture}
    \subcaption{The bold line corresponds to a signature $\lambda$; the gray cells belong to $\mathcal{L}(\lambda)$, and the white cells belong to $\mathcal{R}(\lambda)$. The black (resp. white) dots are the border cells of $\mathcal{L}(\lambda)$ (resp. $\mathcal{R}(\lambda)$).}
    \label{fig:l-and-r-of-gpar-1-a}
    \end{subfigure}
    \hspace{0.5cm}
    \begin{subfigure}[t]{0.45\textwidth}
    \centering
    \begin{tikzpicture}[y=-1cm, scale=0.25, font=\scriptsize]
        \foreach \x in {9,...,10} {
            \draw[thick] (\x, 1) rectangle ++(1, 1);
        }
        \foreach \x in {8,...,10} {
            \draw[thick] (\x, 2) rectangle ++(1, 1);
        }
        \foreach \x in {7,...,8} {
            \draw[thick] (\x, 3) rectangle ++(1, 1);
        }
        \foreach \x in {1,...,3} {
            \draw[thick] (\x, 6) rectangle ++(1, 1);
        }
        \foreach \x in {1,...,2} {
            \draw[thick] (\x, 7) rectangle ++(1, 1);
        }
        \foreach \x in {0,...,2} {
            \draw[thick] (\x, 8) rectangle ++(1, 1);
        }

        \foreach \x in {3,...,4} {
            \draw[thin, gray, densely dashed, fill=red!10] (\x, 5) rectangle ++(1, 1);
        }

        \foreach \x in {4,...,5} {
            \draw[thin, gray, densely dashed, fill=red!10] (6, \x) rectangle ++(1, 1);
        }

        \draw[thick, red] (4, 6) -- (4, 5) -- (7, 5) -- (7, 4);
        \draw[thick, red] (0, 9) -- (0, 10) -- (-2, 10) -- (-2, 11);
        \draw[thick, red] (4, 6) -- (4, 5);
        \draw[thick, red] (6, 5) -- (7, 5);




        \fill (2.25, 5.5) node {$n_1$};
        \fill (6.25, 3.25) node {$n_2$};

    \end{tikzpicture}
    \subcaption{An example of a skew shape $\lambda/\mu$. Null strips are shown by red lines; a vertical and a horizontal  null-strips  $n_1, n_2$ are shown by the pairs of cells they are defined by.}
    \label{fig:l-and-r-of-gpar-1-b}
    \end{subfigure}
    \caption{Signatures and shapes.}
    \label{fig:l-and-r-of-gpar-1}
\end{figure}

\begin{definition}
For $\lambda \in \mathsf{Sign}_n$ we define the sets $\mathcal{L}(\lambda), \mathcal{R}(\lambda) \subset \mathbb{Z}^2$ given by 
$$
\mathcal{L}(\lambda) := \{(i,j) : i \in [n], j \le \lambda_i \}, \quad \mathcal{R}(\lambda) := \{(i,j) : i \in [n], j > \lambda_i \}. 
$$
A cell $(i,j) \in \mathcal{L}(\lambda)$ (resp. in $\mathcal{R}(\lambda)$) is a {\it border cell} if $(i,j+1) \not\in \mathcal{L}(\lambda)$ or $(i+1,j) \not\in \mathcal{L}(\lambda)$ (resp. if $(i,j-1) \not\in \mathcal{R}(\lambda)$ or $(i-1,j) \not\in \mathcal{R}(\lambda)$). See Fig.~\ref{fig:l-and-r-of-gpar-1-a}.
\end{definition}

\begin{definition}[Null strips]
Let $\lambda, \mu \in \mathsf{Sign}_n$ where $\lambda \subseteq \mu$ or $\mu \subseteq \lambda$. 
A {\it vertical (resp. horizontal) null strip} is a pair of cells $((i,j), (i,j+1))$ (resp. $((i,j), (i+1,j))$) such that $(i,j)$ is a border cell of $\mathcal{L}(\mu)$ and $(i,j+1)$ is a border cell of $\mathcal{R}(\lambda)$ (resp. $(i+1,j)$ is a border cell of $\mathcal{R}(\lambda)$). We denote by $N(\lambda, \mu)$ the set of null strips (vertical or horizontal) defined by $\mathcal{L}(\mu)$ and $\mathcal{R}(\lambda)$. See Fig.~\ref{fig:l-and-r-of-gpar-1-b}. 
\end{definition}

\begin{remark}
Let $((i,j), (i,j+1))$ be a vertical strip. It follows that  $j = \lambda_i = \mu_i$ and thus geometrically it can be viewed as the right edge of the cell $(i, \lambda_i)$. Similarly, for a horizontal strip $((i,j), (i+1,j))$ we have $\mu_i > \lambda_{i+1}$ and it be viewed as bottom edge of $(i,j)$. See Fig.~\ref{fig:l-and-r-of-gpar-1-b}. 
\end{remark}

\begin{definition}[Signed skew shapes]
For $\lambda, \mu \in \mathsf{Sign}_n$  where $\lambda \subseteq \mu$ or $\mu \subseteq \lambda$, we define the {\it skew shape} 
$$
\lambda/\mu := (D(\lambda, \mu), N(\lambda, \mu))
$$
as the pair of the corresponding skew diagram and the set of null strips.  
When $\lambda \subseteq \mu$ (resp. $\mu \subseteq \lambda$) then the corresponding skew shape, skew diagram, and null strips are called {\it positive} (resp. {\it negative}). We also write 
$$
-\lambda/\mu := \mu/\lambda.
$$
We shall also denote 
$$
\mathcal{L}(\lambda/\mu) := \mathcal{L}(\mu), \qquad \mathcal{R}(\lambda/\mu) := \mathcal{R}(\lambda). 
$$
\end{definition}


\begin{proposition}
A skew shape $\lambda/\mu$ is fully defined by the sets $\mathcal{L}(\mu)$ and $\mathcal{R}(\lambda)$. 
\end{proposition}


\section{Flagged Grothendieck enumerators}\label{sec:fg}
For a set of variables $\mathbf{x} = (x_i)$ we denote $\mathbf{x}_m = (x_i)_{i \le m}$ and $\overline{\mathbf{x}} = (-x_i)$. 

\subsection{$\mathsf{G}$-flags}
Let us denote $\mathbb{Z}_{-} := \mathbb{Z} \cup \{ -\infty\}$ so that $-\infty < i$ for any $i \in \mathbb{Z}$. 

\begin{definition}[Row $\mathsf{G}$-flags]
    The vectors $\mathbf{a} = (a_i) \in \mathbb{Z}_{-}^{n}, \mathbf{b} = (b_i) \in \mathbb{Z}_{-}^{n}$ are called (row) {\it $\mathsf{G}$-flags} if they satisfy the following conditions:
    \begin{align*}
        a_i \le a_{i + 1}, \qquad b_i \ge b_{i + 1}, \qquad i \ge 1.
    \end{align*}

\end{definition}

\begin{remark}
$\mathsf{G}$-flags slightly differ from usual flags defined for flagged Schur functions, see \cite{wachs, macdonald1}. 
\end{remark}

\subsection{Flagged trio tableaux}


Let us recall that a {\it semistandard} tableau is a filling of a shape with entries weakly increasing along rows (from left to right) and strictly increasing along columns (from top to bottom). Similarly, a {\it reverse semistandard} tableau is a filling of a shape with entries weakly decreasing along rows (from left to right) and strictly decreasing along columns (from top to bottom).

\begin{figure}
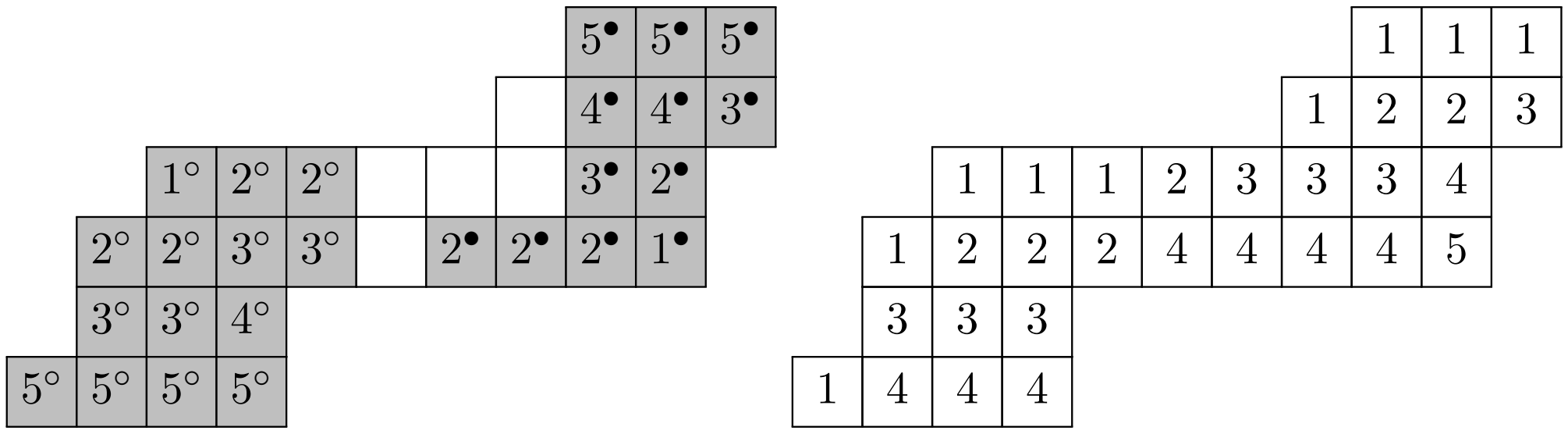

    \centering
    \begin{ytableau}
        \none & \none & \none & \none & \none & \none & \none & \none & *(lightgray) 5^{\bullet} & *(lightgray) 5^{\bullet} & *(lightgray) 5^{\bullet} \\
        \none & \none & \none & \none & \none & \none & \none & *(white) & *(lightgray) 4^{\bullet} & *(lightgray) 4^{\bullet} & *(lightgray) 3^{\bullet} \\
        \none & \none & *(lightgray) 1^{\circ} & *(lightgray) 2^{\circ} & *(lightgray) 2^{\circ} & *(white)  & *(white)  & *(white)  & *(lightgray) 3^{\bullet} & *(lightgray) 2^{\bullet} \\
        \none & *(lightgray) 2^{\circ} & *(lightgray) 2^{\circ} & *(lightgray) 3^{\circ} & *(lightgray) 3^{\circ} & *(white) & *(lightgray) 2^{\bullet} & *(lightgray) 2^{\bullet} & *(lightgray) 2^{\bullet} & *(lightgray) 1^{\bullet} \\
        \none & *(lightgray) 3^{\circ} & *(lightgray) 3^{\circ} & *(lightgray) 4^{\circ} \\ 
        *(lightgray) 5^{\circ} & *(lightgray) 5^{\circ} & *(lightgray) 5^{\circ} & *(lightgray) 5^{\circ}
    \end{ytableau}
    \centering
    \begin{ytableau}
        \none & \none & \none & \none & \none & \none & \none & \none & 1 &  1 &  1 \\
        \none & \none & \none & \none & \none & \none & \none &  1 &  2 &  2 &  3 \\
        \none & \none &  1 &  1 &  1 &  2  &  3 &  3  &  3 &  4 \\
        \none &  1 &  2 &  2 &  2 &  4 &  4 &  4 &  4 &  5 \\
        \none &  3 &  3 &  3 \\ 
         1 &  4 &  4 &  4
    \end{ytableau}
    \caption{An example of a trio tableau. On the left figure: the tableaux $T^{\circ}$ and $T^{\bullet}$. On the right figure: the tableau $T$. Here we have $\lambda = (8,8,8,6,4,4)$, $\mu = (8,7,5,5,4,4)$ and $\tau = (11,11,10,10,4,4)$, $\nu = (8,7,2,1,1,0)$.} 
    \label{fig:trio}
\end{figure}

\begin{definition}[Flagged trio tableaux]
    Let $\lambda, \mu \in \mathsf{Sign}_n$ and $m$ be nonnegative integer. Let $\mathbf{a} = (a_i), \mathbf{b} = (b_i)$ be  $\mathsf{G}$-flags for skew shape $\lambda/\mu$. 
    
    A {\it flagged trio tableau} for shape $\lambda/\mu$ is a triple $\mathbf{T} = (T^{\circ}, T, T^{\bullet})$ of tableaux such that 
    
    (i) (Left tableau) $T^{\circ}$ is a semistandard tableau of some negative shape $\nu/\mu$ with $\nu \subseteq 
    \mu$ and $\nu \in \mathsf{Sign}_n$, filled with entries from the set 
    $\mathbb{Z}^{\circ} := \{ i^{\circ} \,|\,  i \in \mathbb{Z} \}$
    such that the elements in $i$-th row are at most $a_i^{\circ}$.
    
    (ii) (Right tableau) $T^{\bullet}$ is a reverse semistandard tableau of some negative shape $\lambda/\tau$ with $\lambda \subseteq \tau$ and $\tau \in \mathsf{Sign}_n$, filled with entries from the set 
    $\mathbb{Z}^{\bullet} := \{ i^{\bullet} \,|\, i \in \mathbb{Z}\}$
    such that the elements in $i$-th row are at least $b_i^{\bullet}$. 
    
    (iii) (Central tableau) $T$ is a semistandard tableau of shape $\tau/\nu$ with $\nu \subseteq \tau$ filled with entries from $[m]$. 
    
    See an example of trio tableau in Fig.~\ref{fig:trio}.
    
    The {\it weight} of $\mathbf{T}$ is defined as 
    $
    \mathrm{wt}(\mathbf{T}) := \mathrm{wt}(T^{\circ}) \cdot \mathrm{wt}(T) \cdot \mathrm{wt}(T^{\bullet}),     
    $
    and 
    \begin{align*}
    \mathrm{wt}(T^{\circ}) := \prod_{(i,j) \in D(\nu, \mu)}  \mathrm{wt}(T^{\circ}_{i,j}), \quad 
    \mathrm{wt}(T) := \prod_{(i,j) \in D(\tau, \nu)} \mathrm{wt}(T_{i,j}), \quad 
    \mathrm{wt}(T^{\bullet}) := \prod_{(i,j) \in D(\lambda, \tau)} \mathrm{wt}(T^{\bullet}_{i,j}),  
    \end{align*}
    where denoting $T^{\circ}_{i,j} = a^{\circ}$, $T^{}_{i,j} = a^{}$, or $T^{\bullet}_{i,j} = a^{\bullet}$ and $c = c(i,j)$ we define the cell weights as
    $$
    \mathrm{wt}(T^{\circ}_{i,j}) := v_{a, a + c}, \quad  
    \mathrm{wt}(T^{}_{i,j}) := x_{a}, \quad 
    \mathrm{wt}(T^{\bullet}_{i,j}) := u_{a, a - c},      
    $$
    for sets of variables $\mathbf{x}_m = (x_1,\ldots, x_m)$, $\mathbf{U} = (u_{i,j})_{i,j \in \mathbb{Z}}$, $\mathbf{V} = (v_{i,j})_{i,j \in \mathbb{Z}}$. 
    
    We denote by $TT_{\lambda/\mu}(\mathbf{a}, \mathbf{b}, m)$ the set of flagged trio tableaux of shape $\lambda/\mu$ with $\mathsf{G}$-flags $\mathbf{a}, \mathbf{b}$. 
    We also denote by $TT^{\circ}_{\nu/\mu}(\mathbf{a})$, $TT^{\bullet}_{\lambda/\tau}(\mathbf{b})$, and $TT^{}_{\tau/\nu}(m)$ the corresponding sets of left, right, and central tableaux, respectively.  
\end{definition}

\subsection{Flagged Grothendieck enumerators}

\begin{definition}[Flagged Grothendieck enumerators]
Let $m \in \mathbb{Z}_{\ge 0}$ and $\mathbf{x}_m = (x_1,\ldots, x_m)$, $\mathbf{U} = (u_{i,j})_{i,j \in \mathbb{Z}}$, $\mathbf{V} = (v_{i,j})_{i,j \in \mathbb{Z}}$ be sets of variables. 
    Let $\lambda/\mu$ be skew shape of signatures $\lambda, \mu \in \mathsf{Sign}_n$ with  $\mathsf{G}$-flags $\mathbf{a}, \mathbf{b}$.
    We define the {\it flagged Grothendieck enumerators} as follows:  
    \begin{align*}
        \mathsf{G}^{\mathbf{a}, \mathbf{b}}_{\lambda / \mu}(\mathbf{x}_m ; \mathbf{U} , \mathbf{V}) 
        := 
        \sum_{\mathbf{T}\, \in\, TT_{\lambda / \mu}(\mathbf{a}, \mathbf{b}, m) } \mathrm{wt}(\mathbf{T}),
    \end{align*}
    i.e. as the weighted sum over flagged trio tableaux. 
\end{definition}

It is not difficult to see that flagged Grothendieck enumerators satisfy the following basic stability property 
$$
\mathsf{G}^{\mathbf{a}, \mathbf{b}}_{\lambda / \mu}(\mathbf{x}_m ; \mathbf{U} , \mathbf{V}) = 
\mathsf{G}^{\mathbf{a}, \mathbf{b}}_{\lambda / \mu}(\mathbf{x}_{m+1} ; \mathbf{U} , \mathbf{V}) |_{x_{m+1} = 0}.
$$
and hence we can consider them in infinitely many variables $\mathbf{x}$. 

The following proposition 
can easily be proved using the Bender--Knuth involutions performed on the central tableau of a trio tableau. 

\begin{proposition} 
The flagged Grothendieck enumerators $\mathsf{G}^{\mathbf{a}, \mathbf{b}}_{\lambda/\mu}(\mathbf{x} ; \mathbf{U} , \mathbf{V})$ 
are symmetric power series in $\mathbf{x}$ variables. 
\end{proposition}


\subsection{Specialization to canonical refined stable Grothendieck polynomials}
For sets of variables $\mathbf{x} = (x_1, x_2, \ldots)$, $\boldsymbol{\alpha} = (\alpha_1, \alpha_2, \ldots)$, $\boldsymbol{\beta} = (\beta_1, \beta_2, \ldots)$, we consider the {\it canonical refined stable Grothendieck polynomials} $G_{\lambda/\mu}(\mathbf{x}; \boldsymbol{\alpha}, \boldsymbol{\beta})$ indexed by skew partitions which were studied in \cite{hwang1, hwang2}.
These functions are symmetric in $\mathbf{x}$, we have 
$$G_{\lambda/\mu}(\mathbf{x}; \boldsymbol{\alpha}, \boldsymbol{\beta}) = s_{\lambda/\mu}(\mathbf{x}) + \text{higher degree terms},$$  
which is an infinite series with growing degrees, where $s_{\lambda/\mu}(\mathbf{x})$ are skew Schur functions. 
The following duality 
$$
\omega (G_{\lambda/\mu}(\mathbf{x}; \boldsymbol{\alpha}, \boldsymbol{\beta})) = G_{\lambda'/\mu'}(\mathbf{x}; \boldsymbol{\beta}, \boldsymbol{\alpha})
$$
holds for the action of the standard involution $\omega$ of the ring of symmetric functions extended to completion. 
These functions can be defined via different tableaux formulas, see \cite{hwang1, hwang2}; and we shall adopt for them trio tableaux. 

\subsubsection{$G$-specialization}

Let us denote by {\it $G$-specialization} the following substitution of variables in flagged Grothendieck enumerators: 
$$
\mathsf{G}^{\mathbf{a}, \mathbf{b}}_{\lambda/\mu}(\mathbf{x} ; \boldsymbol{\alpha}, \boldsymbol{\beta}) := \mathsf{G}^{\mathbf{a}, \mathbf{b}}_{\lambda/\mu}(\mathbf{x} ; \mathbf{U} , \mathbf{V}) |_{G\text{-spec.}}
$$
where we specialize $u_{i,j} = \alpha_i - \beta_j$, $v_{i,j} = \beta_i - \alpha_j$ and set $\alpha_i = \beta_i = 0$ for $i \le 0$. 

\medskip

The following proposition shows how flagged Grothendieck enumerators generalize {canonical refined stable Grothendieck polynomials}. 

\begin{proposition}\label{prop:gspec}
Let $\lambda, \mu \in \mathsf{Par}_n$ with $\mu \subseteq \lambda$ and $\mathbf{a} = (a_i), \mathbf{b} = (b_i)$ be $\mathsf{G}$-flags given by $a_i = i$, $b_i = \lambda_i$.  
 Then we have 
$$
G_{\lambda/\mu}(\mathbf{x}_m ; \boldsymbol{\alpha}, \boldsymbol{\beta}) = C \cdot 
\mathsf{G}^{\mathbf{a}, \mathbf{b}}_{\lambda/\mu}(\mathbf{x}_m ; \boldsymbol{\alpha}, \boldsymbol{\beta}), 
\qquad C = \prod_{i \in [n]} \prod_{j \in [m]} (1 - \beta_i x_j).
$$
\end{proposition}

\begin{proof}
It follows from \cite[Thm.~6.4]{hwang1} that $G_{\lambda/\mu}(\mathbf{x}_m ; \boldsymbol{\alpha}, \boldsymbol{\beta})/C$ can be expressed as weighted sums over trio tableaux with the given concrete flags. 
\end{proof}

\section{Extended signed ribbon decompositions}\label{sec:rd}

We assume that  $\lambda/\mu$ with $\lambda, \mu \in \mathsf{Sign}_n$ is a given skew signature shape throughout. 

\subsection{Ribbons and decompositions} 
A {\it ribbon} $\theta$ is a connected skew diagram which does not contain $2 \times 2$ box. If $\theta$ is positive, its bottom-left cell is called the {\it head}, and its top-right cell is called the {\it tail}. If $\theta$ is negative, the directions are reversed, i.e. the bottom-left cell is called the tail, and the top-right cell is called the head. 

A {\it null ribbon} is a single null strip $(\alpha, \beta)$. If it is positive, the head is $\alpha$ and the tail is $\beta$. If it is negative, the head is $\beta$ and the tail is $\alpha$. 

\begin{figure}
    \centering
    \begin{subfigure}[t]{0.45\textwidth}
    \centering
    \begin{tikzpicture}[y=-1cm, scale=0.25, font=\scriptsize, baseline={([yshift=-12mm]current bounding box.center)}]



  
        \foreach \x in {6,...,9} {
            \draw[thick] (\x, 0) rectangle ++(1, 1);
        }
        \foreach \x in {6,...,6} {
            \draw[thick] (\x, 1) rectangle ++(1, 1);
        }
        \foreach \x in {6,...,6} {
            \draw[thick] (\x, 2) rectangle ++(1, 1);
        }
        \foreach \x in {3,...,6} {
            \draw[thick] (\x, 3) rectangle ++(1, 1);
        }
        \foreach \x in {-1,...,3} {
            \draw[thick] (\x, 4) rectangle ++(1, 1);
        }
        
        \foreach \x in {4,...,5} {
            \draw[thick] (-1, \x) -- (-1.3, \x);
        }
        \foreach \x in {0,...,1} {
            \draw[thick] (10, \x) -- (10.3, \x);
        }




    \end{tikzpicture}
    \subcaption{Infinite ribbon $\hat{\theta}$.\\}
    \label{fig:tiling-1-a}
    \end{subfigure}
    \begin{subfigure}[t]{0.45\textwidth}
    \centering
    \begin{tikzpicture}[y=-1cm, scale=0.25, font=\tiny]



  
        \foreach \x in {6,...,9} {
            \draw[thick, fill=red] (\x, 0) rectangle ++(1, 1);
        }
        \foreach \x in {6,...,6} {
            \draw[thick, fill=red] (\x, 1) rectangle ++(1, 1);
        }
        \foreach \x in {6,...,6} {
            \draw[thick, fill=red] (\x, 2) rectangle ++(1, 1);
        }
        \foreach \x in {3,...,6} {
            \draw[thick, fill=red] (\x, 3) rectangle ++(1, 1);
        }
        \foreach \x in {-1,...,3} {
            \draw[thick, fill=red] (\x, 4) rectangle ++(1, 1);
        }

        \foreach \x in {7,...,10} {
            \draw[thick, fill=orange] (\x, 1) rectangle ++(1, 1);
        }
        \foreach \x in {7,...,7} {
            \draw[thick, fill=orange] (\x, 2) rectangle ++(1, 1);
        }
        \foreach \x in {7,...,7} {
            \draw[thick, fill=orange] (\x, 3) rectangle ++(1, 1);
        }
        \foreach \x in {4,...,7} {
            \draw[thick, fill=orange] (\x, 4) rectangle ++(1, 1);
        }
        \foreach \x in {0,...,4} {
            \draw[thick, fill=orange] (\x, 5) rectangle ++(1, 1);
        }

        \foreach \x in {8,...,11} {
            \draw[thick, fill=green] (\x, 2) rectangle ++(1, 1);
        }
        \foreach \x in {8,...,8} {
            \draw[thick, fill=green] (\x, 3) rectangle ++(1, 1);
        }
        \foreach \x in {8,...,8} {
            \draw[thick, fill=green] (\x, 4) rectangle ++(1, 1);
        }
        \foreach \x in {5,...,8} {
            \draw[thick, fill=green] (\x, 5) rectangle ++(1, 1);
        }
        \foreach \x in {1,...,5} {
            \draw[thick, fill=green] (\x, 6) rectangle ++(1, 1);
        }

        \foreach \x in {9,...,12} {
            \draw[thick, fill=blue!60] (\x, 3) rectangle ++(1, 1);
        }
        \foreach \x in {9,...,9} {
            \draw[thick, fill=blue!60] (\x, 4) rectangle ++(1, 1);
        }
        \foreach \x in {9,...,9} {
            \draw[thick, fill=blue!60] (\x, 5) rectangle ++(1, 1);
        }
        \foreach \x in {6,...,9} {
            \draw[thick, fill=blue!60] (\x, 6) rectangle ++(1, 1);
        }
        \foreach \x in {2,...,6} {
            \draw[thick, fill=blue!60] (\x, 7) rectangle ++(1, 1);
        }
        
        \foreach \x in {4,...,5} {
            \draw[thick] (-1, \x) -- (-1.3, \x);
        }
        \foreach \x in {0,...,2} {
            \draw[thick] (10, \x) -- (10.3, \x);
        }

        \foreach \x in {5,...,6} {
            \draw[thick] (-0, \x) -- (-0.3, \x);
        }
        \foreach \x in {1,...,2} {
            \draw[thick] (11, \x) -- (11.3, \x);
        }

        \foreach \x in {6,...,7} {
            \draw[thick] (1, \x) -- (0.7, \x);
        }
        \foreach \x in {2,...,3} {
            \draw[thick] (12, \x) -- (12.3, \x);
        }

        \foreach \x in {7,...,8} {
            \draw[thick] (2, \x) -- (1.7, \x);
        }
        \foreach \x in {3,...,4} {
            \draw[thick] (13, \x) -- (13.3, \x);
        }



        \fill[black] (1.5, 7.5) node[left] {$\hat{\theta}_4$};
        \fill[black] (0.5, 6.5) node[left] {$\hat{\theta}_3$};
        \fill[black] (-0.5, 5.5) node[left] {$\hat{\theta}_2$};
        \fill[black] (-1.5, 4.5) node[left] {$\hat{\theta}_1$};

        \fill (2.5, 1) node[scale=1.5] {$\boldsymbol{\ddots}$};
        \fill (9.5, 8) node[scale=1.5] {$\boldsymbol{\ddots}$};

    \end{tikzpicture}
    \subcaption{Tiling the $\mathbb{Z}^2$ with tiles.}
    \label{fig:tiling-1-b}
    \end{subfigure} \\[1cm]
    \begin{subfigure}[t]{0.45\textwidth}
    \centering
    \begin{tikzpicture}[y=-1cm, scale=0.25, font=\scriptsize]
        \foreach \x in {8,...,10} {
            \draw[thick] (\x, 1) rectangle ++(1, 1);
        }
        \foreach \x in {7,...,10} {
            \draw[thick] (\x, 2) rectangle ++(1, 1);
        }
        \foreach \x in {7,...,9} {
            \draw[thick] (\x, 3) rectangle ++(1, 1);
        }
        \foreach \x in {1,...,3} {
            \draw[thick] (\x, 6) rectangle ++(1, 1);
        }
        \foreach \x in {0,...,2} {
            \draw[thick] (\x, 7) rectangle ++(1, 1);
        }
        \foreach \x in {0,...,2} {
            \draw[thick] (\x, 8) rectangle ++(1, 1);
        }
        
        \draw[thick] (4, 6) -- (4, 5) -- (7, 5) -- (7, 4);



        \draw [red, thick, densely dashed] (-3, 6.5) -- (3.5, 6.5) -- (3.5, 5.5) -- (6.5, 5.5) -- (6.5, 2.5) -- (13.5, 2.5);


    \end{tikzpicture}
    \subcaption{Intersection of $\lambda/\mu$ and $\hat{\theta}_1$}
    \label{fig:tiling-1-c}
    \end{subfigure}
    \begin{subfigure}[t]{0.45\textwidth}
    \centering
    \begin{tikzpicture}[y=-1cm, scale=0.25, font=\scriptsize]
        \foreach \x in {8,...,10} {
            \draw[thick] (\x, 1) rectangle ++(1, 1);
        }
        \foreach \x in {7,...,10} {
            \draw[thick] (\x, 2) rectangle ++(1, 1);
        }
        \foreach \x in {7,...,9} {
            \draw[thick, fill=red!60] (\x, 3) rectangle ++(1, 1);
        }
        \foreach \x in {1,...,3} {
            \draw[thick, fill=red!60] (\x, 6) rectangle ++(1, 1);
        }
        \foreach \x in {0,...,2} {
            \draw[thick] (\x, 7) rectangle ++(1, 1);
        }
        \foreach \x in {0,...,2} {
            \draw[thick] (\x, 8) rectangle ++(1, 1);
        }

        \foreach \x in {3,...,4} {
            \draw[thin, gray, densely dashed, fill=red!10] (\x, 5) rectangle ++(1, 1);
        }

        \foreach \x in {4,...,5} {
            \draw[thin, gray, densely dashed, fill=red!10] (6, \x) rectangle ++(1, 1);
        }

        \draw[thick] (4, 6) -- (4, 5) -- (7, 5) -- (7, 4);
        \draw[thick, red] (4, 6) -- (4, 5);
        \draw[thick, red] (6, 5) -- (7, 5);





    \end{tikzpicture}
    \subcaption{Decomposition ribbons and null strips obtained by the intersection.}
    \label{fig:tiling-1-d}
    \end{subfigure}
    \\[1cm]
            \begin{subfigure}[t]{0.45\textwidth}
    \centering
    \begin{tikzpicture}[y=-1cm, scale=0.25, font=\scriptsize]
        \foreach \x in {8,...,10} {
            \draw[thick, fill=orange!50] (\x, 1) rectangle ++(1, 1);
        }
        \foreach \x in {7,...,10} {
            \draw[thick, fill=red!50] (\x, 2) rectangle ++(1, 1);
        }
        \foreach \x in {7,...,9} {
            \draw[thick, fill=green!50] (\x, 3) rectangle ++(1, 1);
        }
        \foreach \x in {1,...,3} {
            \draw[thick, fill=red!50] (\x, 6) rectangle ++(1, 1);
        }
        \foreach \x in {0,...,2} {
            \draw[thick, fill=green!50] (\x, 7) rectangle ++(1, 1);
        }
        \foreach \x in {0,...,2} {
            \draw[thick, fill=blue!50] (\x, 8) rectangle ++(1, 1);
        }
        
        \draw[thick] (4, 6) -- (4, 5) -- (7, 5) -- (7, 4);


        \foreach \x in {3,...,4} {
            \draw[thin, gray, densely dashed, fill=red!10] (\x, 5) rectangle ++(1, 1);
        }

        \foreach \x in {4,...,5} {
            \draw[thin, gray, densely dashed, fill=red!10] (6, \x) rectangle ++(1, 1);
        }


        \draw [orange, thick, dashed] (-3, 5.5) -- (2.5, 5.5) -- (2.5, 4.5) -- (5.5, 4.5) -- (5.5, 1.5) -- (13.5, 1.5);
        \draw [red, thick, dashed] (-3, 6.5) -- (3.5, 6.5) -- (3.5, 5.5) -- (6.5, 5.5) -- (6.5, 2.5) -- (13.5, 2.5);
        \draw [green, thick, dashed] (-3, 7.5) -- (4.5, 7.5) -- (4.5, 6.5) -- (7.5, 6.5) -- (7.5, 3.5) -- (13.5, 3.5);
        \draw [blue, thick, dashed] (-3, 8.5) -- (5.5, 8.5) -- (5.5, 7.5) -- (8.5, 7.5) -- (8.5, 4.5) -- (13.5, 4.5);


    \end{tikzpicture}
    \subcaption{Ribbons and null-strips in outer decomposition.}
    \label{fig:tiling-1-e}
    \end{subfigure}
    
    \caption{Examples of ribbons, null strips, and outer decompositions.}
    \label{fig:tiling-1}
\end{figure}
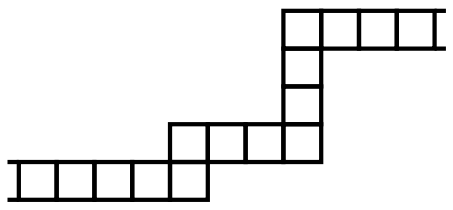
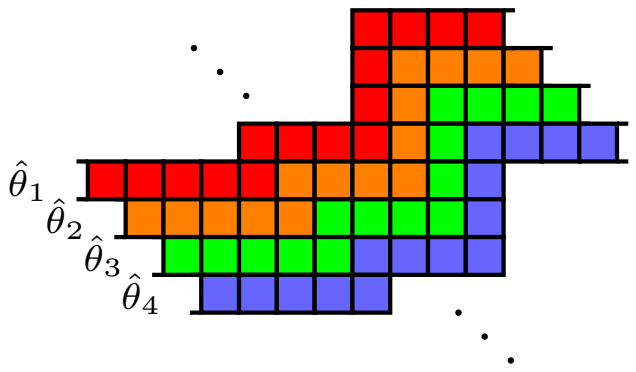
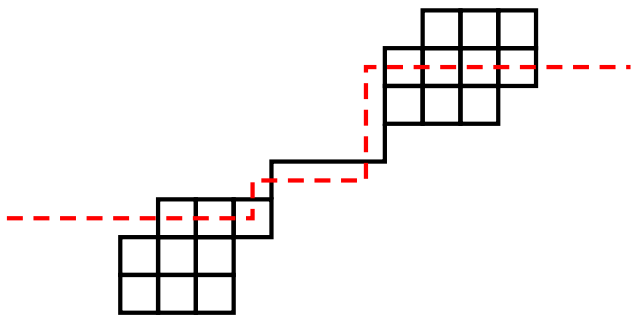
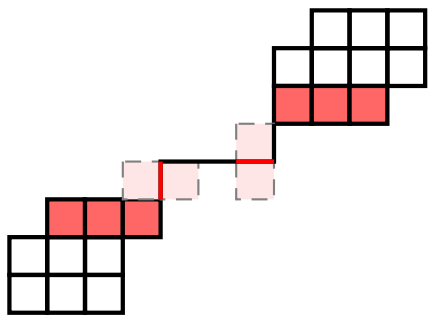
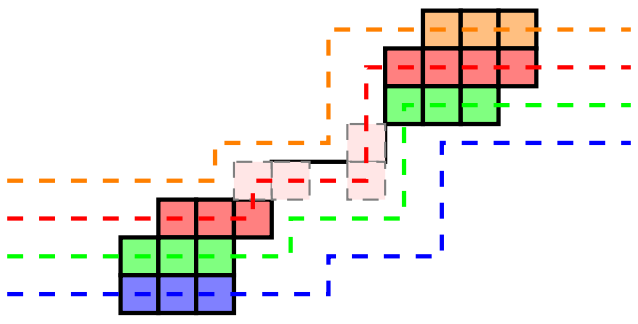

An {\it infinite ribbon} $\hat{\theta}$ is a (positive) ribbon with the head $(n, -\infty)$ and tail $(1, \infty)$ for some $n \ge 1$, see Fig.~\ref{fig:tiling-1-a}. 
It is not hard to see that for every content $c \in \mathbb{Z}$ there is unique cell $\alpha$ with content $c$ in $\hat\theta$.
For $k \in \mathbb{Z}$ we define the translated copy $\hat\theta_k$ of $\hat\theta$ given by 
$$
\hat\theta_k := \{(i+k,j+k) : (i,j) \in \hat\theta \}.
$$
    It is not hard to see that any cell $\alpha \in \mathbb{Z}^2$ belongs to some unique copy $\hat{\theta}_k$, and thus the set of all copies of $\hat{\theta}$ {\it tiles} the set $\mathbb{Z}^2$ (see Fig.~\ref{fig:tiling-1-b}), i.e. 
    \begin{equation*}
        \bigcup_{k \in \mathbb{Z}} \hat{\theta}_k = \mathbb{Z}^2.
    \end{equation*}
    We call each copy $\hat{\theta}_k$ a {\it tile}.

\begin{remark}
The idea of infinite ribbons appears in \cite{np}. 
\end{remark}

\begin{definition}[Outer decompositions]
    Let $\lambda/\mu = ({D}(\lambda, \mu), {N}(\lambda, \mu))$ be given skew shape. Let $\hat{\theta}_k$ be a tile. Let $\rho \subseteq {D}(\lambda, \mu) \cap \hat{\theta}_k$ be a ribbon whose bottom-left cell belongs to the bottom or left perimeter of ${D}(\lambda, \mu)$, and whose top-right cell belongs to the top or right perimeter of ${D}(\lambda, \mu)$. Then $\rho$ is called a {\it decomposition ribbon}. Consider the following sets:
    \begin{align*}
        &A = \{\rho \subseteq {D}(\lambda, \mu) \cap \hat{\theta}_k\ |\ \rho \ \text{is a decompoisition ribbon}\}, \\ 
        &B = \{(\alpha, \beta) \in {N}(\lambda, \mu)\ |\ \alpha, \beta \in \hat{\theta}_k\}.
    \end{align*}
    We define the intersection of $\lambda/\mu$ with $\hat{\theta}_k$ as follows (see Fig.~\ref{fig:tiling-1-c}, \ref{fig:tiling-1-d}): 
    \[
        \lambda/\mu \cap \hat{\theta}_k := A \cup B.
    \]
    Then the {\it outer decomposition} $\Theta$ of $\lambda/\mu$ is defined as the set of all decomposition ribbons and null strips shared by $\lambda/\mu$ and all tiles $\hat{\theta}_k$ for $k \in \mathbb{Z}$ (see Fig.~\ref{fig:tiling-1-e}), i.e.: 
    \begin{equation*}
        \Theta = \bigcup_{k \in \mathbb{Z}} \lambda/\mu \cap \hat{\theta}_k.
    \end{equation*}
\end{definition}

\begin{remark}
We treat the null-strips in $\Theta$ as ribbons. Throughout the text by ribbon in outer decomposition we shall mean either a ribbon or a null-strip.
\end{remark}


\subsection{Pipe functions}

It will be useful for us to describe ribbon decompositions using pipe functions defined below. 

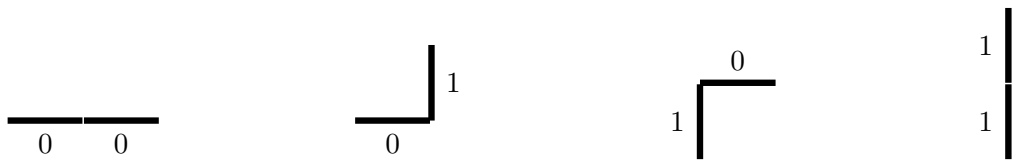
\begin{figure}
    \centering
    \begin{tikzpicture}[scale=0.5, transform shape]
        \draw[very thick] (1, 1) -- (1.99, 1);
        \node at (1.5, 0.7) {$0$};
        \draw[very thick] (2.01, 1) -- (3, 1);
        \node at (2.5, 0.7) {$0$};
    \end{tikzpicture}
    \hspace{1cm}
    \begin{tikzpicture}[scale=0.5, transform shape]
        \draw[very thick] (1, 1) -- (1.99, 1);
        \node at (1.5, 0.7) {$0$};
        \draw[very thick] (2.01, 1) -- (2.01, 2);
        \node at (2.3, 1.5) {$1$};
    \end{tikzpicture}
    \hspace{1cm}
    \begin{tikzpicture}[scale=0.5, transform shape]
        \draw[very thick] (1, 1) -- (1, 1.99);
        \node at (0.7, 1.5) {$1$};
        \draw[very thick] (1, 2.01) -- (2, 2.01);
        \node at (1.5, 2.3) {$0$};
    \end{tikzpicture}
    \hspace{1cm}
    \begin{tikzpicture}[scale=0.5, transform shape]
        \draw[very thick] (1, 1) -- (1, 1.99);
        \node at (0.7, 1.5) {$1$};
        \draw[very thick] (1, 2.01) -- (1, 3);
        \node at (0.7, 2.5) {$1$};
    \end{tikzpicture}
    \caption{The pipes corresponding to the pairs $(0, 0), (0, 1), (1, 0), (1,1)$.}
    \label{fig:possible-pipes}
\end{figure}

Let us define the set 
$$\mathcal{P} := \{(0, 0), (0, 1), (1, 0), (1, 1)\},$$ where $0$ and $1$ mean horizontal and vertical unit directions, respectively. To each pair in $\mathcal{P}$ we match a {\it pipe}, whose directions are determined by corresponding vectors,  
see Fig.~\ref{fig:possible-pipes}. 
We call a pipe $p \in \mathcal{P}$ {\it hook pipe} if $p = (0, 1)$ or $p = (1, 0)$. We denote the set of hook pipes as $\mathcal{P}_H := \{(0, 1), (1, 0)\}$. 


Note that for each $c \in \mathbb{Z}$ there is only one cell with content $c$ in $\hat{\theta}$. Thus, the cell with content $c + 1$ stays either on the right, or on top of the cell with content $c$. It turns out that this data (location of the cell of content $c + 1$ for all $c \in \mathbb{Z}$) is crucial for constructing the lattice (which will be defined in later sections), and we define {\it pipe functions} which store that data. 

\begin{definition}[Pipe functions]
We define {\it pipe function} 
$\pi : \mathbb{Z} \rightarrow \mathcal{P} 
$ 
which assigns to cells (of a shape) 
with content $c$ some pipe $\pi(c) \in \mathcal{P}$, according to the following restrictions: 
    $$\pi(c - 1)^2 = \pi(c)^1 \text{ for all } c \in \mathbb{Z}, \text{ where } \pi(c) = (\pi(c)^1, \pi(c)^2).$$ 
It is not hard to see that any pipe function $\pi$ induces some infinite ribbon $\hat{\theta} = \hat{\theta}(\pi)$. 
For simplicity, we shall also write $\pi(\alpha) := \pi(c(\alpha))$ for a cell $\alpha$ throughout the text.

    For $\lambda, \mu \in \mathsf{Sign}_n$, consider the contents $a = \mathrm{min}(\lambda_n, \mu_n) - n$ and $b = \mathrm{max}(\lambda_1, \mu_1)$. We say that pipe function $\pi$ and shape $\lambda/\mu$ are {\it compatible} if 
    \begin{enumerate}
        \item $\pi(a) = (0, \cdot)$ and $\pi(b) = (\cdot, 0)$. 
        \item $\pi(c) = (0, 0)$ for all $c < a$ and $c > b$.
    \end{enumerate}
    In this text we consider only pipe functions that are compatible with the given shape. 

\end{definition}

\begin{figure}
    \centering
    \[
    \begin{array}{cccccccccc}
        \dots & -4 & -3 & -2 & -1 & 0 & 1 & 2 & 3 & \dots \\ 
        \begin{tikzpicture}[scale=0.3, transform shape] 
        \draw[red, very thick] (1, 1) -- (1.99, 1);
        \draw[red, very thick] (2.01, 1) -- (3, 1);
    \end{tikzpicture} & 
    \begin{tikzpicture}[scale=0.3, transform shape] 
        \draw[red, very thick] (1, 1) -- (1.99, 1);
        \draw[red, very thick] (2.01, 1) -- (3, 1);
    \end{tikzpicture} & 
    \begin{tikzpicture}[scale=0.3, transform shape] 
        \draw[red, very thick] (1, 1) -- (1.99, 1);
        \draw[red, very thick] (2.01, 1) -- (2.01, 2);
    \end{tikzpicture} & 
    \begin{tikzpicture}[scale=0.3, transform shape] 
        \draw[red, very thick] (1, 1) -- (1, 1.99);
        \draw[red, very thick] (1, 2.01) -- (2, 2.01);
    \end{tikzpicture} &
    \begin{tikzpicture}[scale=0.3, transform shape] 
        \draw[red, very thick] (1, 1) -- (1.99, 1);
        \draw[red, very thick] (2.01, 1) -- (3, 1);
    \end{tikzpicture} & 
    \begin{tikzpicture}[scale=0.3, transform shape] 
        \draw[red, very thick] (1, 1) -- (1.99, 1);
        \draw[red, very thick] (2.01, 1) -- (2.01, 2);
    \end{tikzpicture} & 
    \begin{tikzpicture}[scale=0.3, transform shape] 
        \draw[red, very thick] (1, 1) -- (1, 1.99);
        \draw[red, very thick] (1, 2.01) -- (1, 3);
    \end{tikzpicture} &
    \begin{tikzpicture}[scale=0.3, transform shape] 
        \draw[red, very thick] (1, 1) -- (1, 1.99);
        \draw[red, very thick] (1, 2.01) -- (2, 2.01);
    \end{tikzpicture} &
    \begin{tikzpicture}[scale=0.3, transform shape] 
        \draw[red, very thick] (1, 1) -- (1.99, 1);
        \draw[red, very thick] (2.01, 1) -- (3, 1);
    \end{tikzpicture} & 
    \begin{tikzpicture}[scale=0.3, transform shape] 
        \draw[red, very thick] (1, 1) -- (1.99, 1);
        \draw[red, very thick] (2.01, 1) -- (3, 1);
    \end{tikzpicture}
    \end{array}
    \]
    \begin{tikzpicture}[y=-1cm, scale=0.50, transform shape]
	\foreach \x in {1,...,2} {
        \draw[very thick] (\x, 1) rectangle ++(1, 1);
    }
    \foreach \x in {0,...,2} {
        \draw[very thick] (\x, 2) rectangle ++(1, 1);
    }
    \foreach \x in {0,...,1} {
        \draw[very thick] (\x, 3) rectangle ++(1, 1);
    }

    \foreach \x in {-2,...,0} {
        \draw[gray, dashed] (\x, 1) rectangle ++(1, 1);
    }
    \foreach \x in {-2,...,-1} {
        \draw[gray, dashed] (\x, 2) rectangle ++(1, 1);
    }
    \foreach \x in {-2,...,-1} {
        \draw[gray, dashed] (\x, 3) rectangle ++(1, 1);
    }


    \foreach \x in {3,...,4} {
        \draw[gray, dashed] (\x, 1) rectangle ++(1, 1);
    }
    \foreach \x in {3,...,4} {
        \draw[gray, dashed] (\x, 2) rectangle ++(1, 1);
    }
    \foreach \x in {2,...,4} {
        \draw[gray, dashed] (\x, 3) rectangle ++(1, 1);
    }

    \foreach \x in {-1, ..., 5} {
        \foreach \y in {1, ..., 3}
        \pgfmathsetmacro{\z}{int(\x - \y)}
        \node at (\x - 0.5, \y + 0.5) {$\z$};
    }

    \foreach \x in {1,...,4} {
        \draw[gray, dashed] (-2, \x) -- (-2.3, \x);
    }
    \foreach \x in {1,...,4} {
        \draw[gray, dashed] (5, \x) -- (5.3, \x);
    }

    \end{tikzpicture}
    \hspace{0.5cm}
    \begin{tikzpicture}[y=-1cm, scale=0.35, transform shape]

    \draw[red, very thick] (-0.5, -0.5) -- (1, -0.5) -- (2, -0.5) -- (2.5, -0.5) -- (2.5, -1) -- (2.5, -1.5) -- (3, -1.5) -- (4, -1.5) -- (4.5, -1.5) -- (4.5, -2) -- (4.5, -3) -- (4.5, -3.5) -- (5, -3.5) -- (6, -3.5) -- (7.5, -3.5);
	\foreach \x in {4,...,6} {
        \draw[very thick] (\x, -4) rectangle ++(1, 1);
    }
    \foreach \x in {4,...,4} {
        \draw[very thick] (\x, -3) rectangle ++(1, 1);
    }
    \foreach \x in {2,...,4} {
        \draw[very thick] (\x, -2) rectangle ++(1, 1);
    }
    \foreach \x in {0,...,2} {
        \draw[very thick] (\x, -1) rectangle ++(1, 1);
    }




    \foreach \x in {-1,...,0} {
        \draw[very thick] (0, \x) -- (-0.3, \x);
    }
    \foreach \x in {-4,...,-3} {
        \draw[very thick] (7, \x) -- (7.3, \x);
    }




    


    \end{tikzpicture}
    \caption{A skew shape with contents, and a pipe function $\pi$ (above) with the corresponding infinite ribbon.}
    \label{fig:pipe-func-ex}
\end{figure}
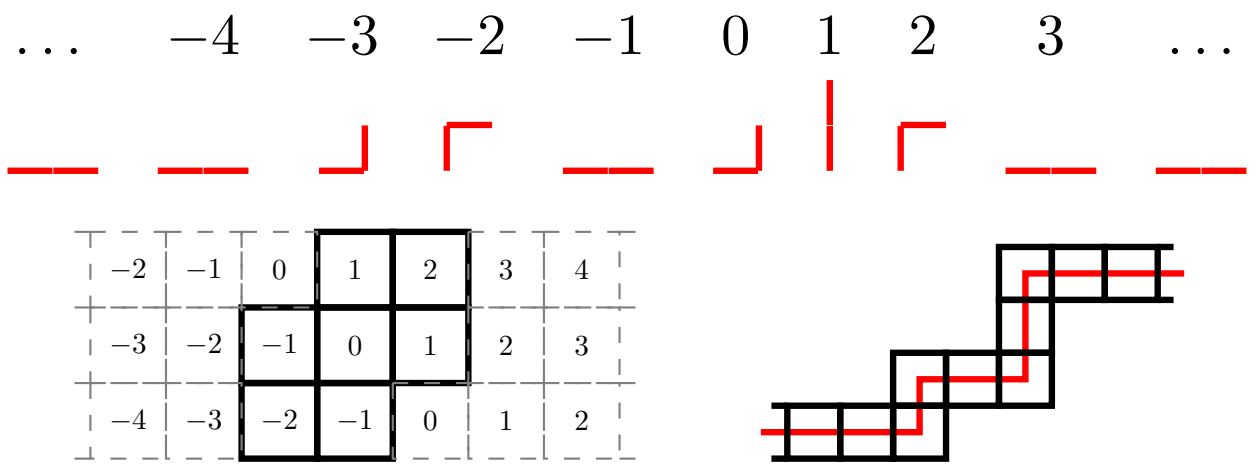

\begin{example}
    Let $\lambda = (3,3,2),\ \mu = (1,0,0)$. An example of a pipe function $\pi$ is shown in Fig.~\ref{fig:pipe-func-ex} and it is defined as follows 
        \[
    \begin{array}{rclrclrcl}
        -4 & \mapsto & (0, 0) \quad & -1 & \mapsto & (0, 0) \quad & 2 & \mapsto & (1, 0) \\ 
        -3 & \mapsto & (0, 1) \quad & 0 & \mapsto & (0, 1) \quad & 3 & \mapsto & (0, 0) \\ 
        -2 & \mapsto & (1, 0) \quad & 1 & \mapsto & (1, 1) \quad & c & \mapsto & (0, 0), 
    \end{array}
    \]
    for all $c \in \mathbb{Z} \setminus [-2, 3]$. 
\end{example}

\begin{definition}[Hook-pipe counting function] 
For a given  skew shape $\lambda/\mu$ and pipe function $\pi$ we denote by $\mathsf{hp}^{\pi}_{\lambda/\mu}$ the number of hook pipes in the negative cells of the diagram $D(\lambda,\mu)$. More formally, 
$$
\mathsf{hp}^{\pi}_{\lambda/\mu} := 
            \begin{cases}
                0, \quad & \text{if ${D}(\lambda, \mu)$ is positive}, \\ 
                |\{\alpha \in {D}(\lambda, \mu)\ |\ \pi(\alpha) \in \mathcal{P}_H\}|, \quad & \text{if ${D}(\lambda, \mu)$ is negative}.
            \end{cases}
$$
\end{definition}

\subsection{Extended Hamel--Goulden operation}
There is a well-known Hamel--Goulden operation $\#$ for combining pairs of ribbons in $\Theta = (\theta_i)$ (\cite{hg}, \cite{cyy}) which results in certain ribbons $\theta_i \# \theta_j$. We shall refine that $\#$ operation (on signed shapes) so that it will produce shapes for ribbons with the same cells but also with null strips.

Informally, we define $\theta_i \# \theta_j$ as a natural extension of sub-ribbon of $\hat{\theta}$ (whoose cell contents lie in interval $[p, q]$) into a shape compatible with $\pi$. 
The formal definition is the following.

\begin{definition}
    Let $\hat{\theta} = \hat\theta(\pi)$ be an infinite ribbon. Let $a, b \in \mathbb{Z}$ with $a \le b + 1$. We define $[a, b]_{\hat{\theta}}$ as a sub-ribbon of $\hat{\theta}$ which contains the cells with contents from the interval $[a, b]$; or as a null strip $(\beta, \alpha)$ if $a = b + 1$ where $\alpha, \beta \in \hat{\theta}$ with $c(\beta) = b, c(\alpha) = a$. 

    We then construct certain skew shape $\mathrm{sh}[a, b]_{\hat{\theta}}$ with the diagram (or null strip) $[a, b]_{\hat{\theta}}$. 

    First, construct ribbon $\theta$ from $\hat{\theta}$ as follows: 
    
    (i) if $\pi(a) = (1, \cdot)$ delete the last row of $\hat{\theta}$.

    (ii) if $\pi(b) = (\cdot, 1)$ delete the first row of $\hat{\theta}$, then shift each cell $(i, j) \in \hat{\theta}$ to $(i - 1, j - 1)$.

    Let $n$ be number of rows in the resulted ribbon $\theta$. Let $\xi, \kappa \in \mathsf{Sign}_n$ with $\kappa \subseteq \xi$, whose diagram of shape $\xi/\kappa$ is $\theta$. 
    Let 
    \[
    A = \{\alpha \in \theta\ |\ c(\alpha) < a\}, \quad B = \{\alpha \in \theta\ |\ c(\alpha) > b\}.
    \]
    We construct the set $\mathcal{L}(a, b)$ from ${\mathcal{L}(\kappa)}$ by adding the elements of $A$ if $\pi(a) = (0, \cdot)$ and the elements of $B$ if $\pi(b) = (\cdot, 1)$. Similarly, we construct the set $\mathcal{R}(a, b)$ from ${\mathcal{R}(\kappa)}$ by adding the elements of $A$ if $\pi(a) = (1, \cdot)$ and the elements of $B$ if $\pi(b) = (\cdot, 0)$. 
    
    Then define 
    \[
    \mathrm{sh}[a, b]_{\hat{\theta}} := (\mathcal{L}(a, b), \mathcal{R}(a, b)). 
    \] 
    We note that $\mathrm{sh}[a, b]_{\hat{\theta}}$ is a positive shape compatible with $\pi$. 

\end{definition}


\begin{definition}[Extended $\#$-operation]
    Let $\hat{\theta} = \hat\theta(\pi)$ be an infinite ribbon, $\lambda/\mu$ be given skew shape and $\Theta = (\theta_i)$ be an outer ribbon decomposition of $\lambda/\mu$. Let $\delta_i, \gamma_i$ be the head and the tail of $\theta_i$. 
    
    (i) If $\lambda/\mu$ is positive, define 
    \begin{equation*}
        \theta_i \# \theta_j := 
        \begin{cases}
            \mathrm{sh}[c(\delta_i), c(\gamma_j)]_{\hat{\theta}},\quad &\text{if}\ c(\delta_i) \le c(\gamma_j) + 1 \\
            -\mathrm{sh}[c(\gamma_j) + 1, c(\delta_i) - 1]_{\hat{\theta}},\quad &\text{otherwise}.
        \end{cases}
    \end{equation*}

    (ii) If $\lambda/\mu$ is negative, define 
    \begin{equation*}
        \theta_i \# \theta_j := 
        \begin{cases}
            -\mathrm{sh}[c(\gamma_j), c(\delta_i)]_{\hat{\theta}},\quad &\text{if}\ c(\gamma_j) \le c(\delta_j) + 1 \\
            \mathrm{sh}[c(\delta_i) + 1, c(\gamma_j) - 1]_{\hat{\theta}},\quad &\text{otherwise}.
        \end{cases}
    \end{equation*}
\end{definition}

Next we define {\it empty shape} as a shape $\kappa/\kappa$ compatible with $\pi$ for some $\kappa \in \mathrm{Sign}_n$. The formal definition is given below. 
\begin{definition}[Empty shapes]
    Let $\hat{\theta} = \hat\theta(\pi)$ be an infinite ribbon. Let $\theta$ be the ribbon obtained from $\hat{\theta}$ by deleting its last row. Let $n$ be number of rows in $\theta$; and $\xi, \kappa \in \mathsf{Sign}_n$ such that the diagram of $\xi/\kappa$ is $\theta$. We define the following empty shape:
    \[
        \varepsilon(\hat{\theta}) := \kappa/\kappa.
    \]
 Note that empty shape $\varepsilon$ is compatible with the pipe function $\pi$ as well.
\end{definition}

\begin{example}
    Let $\hat{\theta}$ be the infinite ribbon drawn in Fig.~\ref{fig:sub-ribbons-of-infinite-ribbon-1-a}. The cells with contents $a, b$ are highlighted by dots. It follows that $\pi(a) = (1, \cdot), \pi(b) = (\cdot, 0)$, thus $\theta$ is obtained by deleting the last row from $\hat{\theta}$ (see Fig.~\ref{fig:sub-ribbons-of-infinite-ribbon-1-b}); then by definition, the cells of $\theta$ with content larger than $b$ and less than $a$ are added to $\mathcal{R}(\xi), \mathcal{L}(\eta)$ (resp.); this defines the shape $\mathrm{sh}[a, b]_{\hat{\theta}}$ depicted in Fig.~\ref{fig:sub-ribbons-of-infinite-ribbon-1-c}. The empty ribbon obtained from $\hat{\theta}$ is depicted in Fig.~\ref{fig:sub-ribbons-of-infinite-ribbon-1-d}.
\end{example}

\begin{proposition}
    Let $\pi$ be pipe function compatible with the shape $\lambda/\mu$. Let $\hat{\theta} = \hat{\theta}(\pi)$ be infinite ribbon, and $\Theta = (\theta_1, \ldots, \theta_k)$ be outer ribbon decomposition of $\lambda/\mu$ induced by $\hat{\theta}$. Then $\pi$ is compatible with the shapes $\theta_i \# \theta_j$ ($i, j \in [k]$) and empty shape $\varepsilon(\hat{\theta})$. 
\end{proposition}

\begin{figure}
    \centering
    \begin{subfigure}[t]{0.45\textwidth}
    \centering
    \begin{tikzpicture}[y=-1cm, scale=0.25, font=\scriptsize]
        \foreach \x in {7,...,10} {
            \draw[thick] (\x, 1) rectangle ++(1, 1);
        }
        \foreach \x in {7,...,7} {
            \draw[thick] (\x, 2) rectangle ++(1, 1);
        }
        \foreach \x in {4,...,7} {
            \draw[thick] (\x, 3) rectangle ++(1, 1);
        }
        \foreach \x in {4,...,4} {
            \draw[thick] (\x, 4) rectangle ++(1, 1);
        }
        \foreach \x in {1,...,4} {
            \draw[thick] (\x, 5) rectangle ++(1, 1);
        }
        \foreach \x in {1,...,1} {
            \draw[thick] (\x, 6) rectangle ++(1, 1);
        }
        \foreach \x in {-1,...,1} {
            \draw[thick] (\x, 7) rectangle ++(1, 1);
        }
        \foreach \x in {-3,...,-1} {
            \draw[thick] (\x, 8) rectangle ++(1, 1);
        }
        \foreach \x in {-3,...,-3} {
            \draw[thick] (\x, 9) rectangle ++(1, 1);
        }
        \foreach \x in {-5,...,-3} {
            \draw[thick] (\x, 10) rectangle ++(1, 1);
        }






        \foreach \x in {10,...,11} {
            \draw[thick] (-5, \x) -- (-5.3, \x);
            }
        \foreach \x in {1,...,2} {
            \draw[thick] (11, \x) -- (11.3, \x);
        }

        \fill (1.5, 6.5) circle (4pt) node[left=1.5pt] {$a$};
        \fill (5.5, 3.5) circle (4pt) node[above=1.5pt] {$b$};


        \fill[black] (-5, 10.5) node[left=1.5pt] {$\hat{\theta}_0$};

    \end{tikzpicture}
    \subcaption{Infinite ribbon $\hat{\theta}$.\\}
    \label{fig:sub-ribbons-of-infinite-ribbon-1-a}
    \end{subfigure}
    \begin{subfigure}[t]{0.45\textwidth}
    \centering
    \begin{tikzpicture}[y=-1cm, scale=0.25, font=\scriptsize]
        \foreach \x in {7,...,10} {
            \draw[thick] (\x, 1) rectangle ++(1, 1);
        }
        \foreach \x in {7,...,7} {
            \draw[thick] (\x, 2) rectangle ++(1, 1);
        }
        \foreach \x in {4,...,7} {
            \draw[thick] (\x, 3) rectangle ++(1, 1);
        }
        \foreach \x in {4,...,4} {
            \draw[thick] (\x, 4) rectangle ++(1, 1);
        }
        \foreach \x in {1,...,4} {
            \draw[thick] (\x, 5) rectangle ++(1, 1);
        }
        \foreach \x in {1,...,1} {
            \draw[thick] (\x, 6) rectangle ++(1, 1);
        }
        \foreach \x in {-1,...,1} {
            \draw[thick] (\x, 7) rectangle ++(1, 1);
        }
        \foreach \x in {-3,...,-1} {
            \draw[thick] (\x, 8) rectangle ++(1, 1);
        }
        \foreach \x in {-3,...,-3} {
            \draw[thick] (\x, 9) rectangle ++(1, 1);
        }
        \foreach \x in {-5,...,-3} {
            \draw[thin, gray, densely dashed] (\x, 10) rectangle ++(1, 1);
        }






        \foreach \x in {10,...,11} {
            \draw[thin, gray, densely dashed] (-5, \x) -- (-5.3, \x);
            }
        \foreach \x in {1,...,2} {
            \draw[thick] (11, \x) -- (11.3, \x);
        }

        \fill (1.5, 6.5) circle (4pt) node[left=1.5pt] {$a$};
        \fill (5.5, 3.5) circle (4pt) node[above=1.5pt] {$b$};


        \fill[black] (12, 1.5) node {$\theta$};
        \fill[gray] (-5, 10.5) node[left=1.5pt] {$\hat{\theta}_0$};

    \end{tikzpicture}
    \subcaption{Ribbon $\theta$ obtained from $\hat{\theta}$ by deleting last row.}
    \label{fig:sub-ribbons-of-infinite-ribbon-1-b}
    \end{subfigure}
    \\[1cm]
    \begin{subfigure}[t]{0.45\textwidth}
    \centering
    \begin{tikzpicture}[y=-1cm, scale=0.25, font=\scriptsize]
        \foreach \x in {7,...,10} {
            \draw[thin, gray, densely dashed] (\x, 1) rectangle ++(1, 1);
        }
        \foreach \x in {7,...,7} {
            \draw[thin, gray, densely dashed] (\x, 2) rectangle ++(1, 1);
        }
        \foreach \x in {4,...,7} {
            \draw[thin, gray, densely dashed] (\x, 3) rectangle ++(1, 1);
        }
        \foreach \x in {4,...,5} {
            \draw[thick] (\x, 3) rectangle ++(1, 1);
        }
        \foreach \x in {4,...,4} {
            \draw[thick] (\x, 4) rectangle ++(1, 1);
        }
        \foreach \x in {1,...,4} {
            \draw[thick] (\x, 5) rectangle ++(1, 1);
        }
        \foreach \x in {1,...,1} {
            \draw[thick] (\x, 6) rectangle ++(1, 1);
        }
        \foreach \x in {-1,...,1} {
            \draw[thin, gray, densely dashed] (\x, 7) rectangle ++(1, 1);
        }
        \foreach \x in {-3,...,-1} {
            \draw[thin, gray, densely dashed] (\x, 8) rectangle ++(1, 1);
        }
        \foreach \x in {-3,...,-3} {
            \draw[thin, gray, densely dashed] (\x, 9) rectangle ++(1, 1);
        }






        \foreach \x in {1,...,2} {
            \draw[thin, gray, densely dashed] (11, \x) -- (11.3, \x);
        }

        \draw[thick] (6, 3) -- (7, 3) -- (7, 1);
        \draw[thick] (1, 7) -- (-1, 7) -- (-1, 8) -- (-3, 8) -- (-3, 10);

        \fill (1.5, 6.5) circle (4pt) node[left=1.5pt] {$a$};
        \fill (5.5, 3.5) circle (4pt) node[above=1.5pt] {$b$};


        \fill[gray] (12, 1.5) node {$\theta$};

    \end{tikzpicture}
    \subcaption{Shape $\mathrm{sh}[a, b]_{\hat{\theta}}$.}
    \label{fig:sub-ribbons-of-infinite-ribbon-1-c}
    \end{subfigure}
    \begin{subfigure}[t]{0.45\textwidth}
    \centering
    \begin{tikzpicture}[y=-1cm, scale=0.25, font=\scriptsize]
        \foreach \x in {7,...,10} {
            \draw[thin, gray, densely dashed] (\x, 1) rectangle ++(1, 1);
        }
        \foreach \x in {7,...,7} {
            \draw[thin, gray, densely dashed] (\x, 2) rectangle ++(1, 1);
        }
        \foreach \x in {4,...,7} {
            \draw[thin, gray, densely dashed] (\x, 3) rectangle ++(1, 1);
        }
        \foreach \x in {4,...,4} {
            \draw[thin, gray, densely dashed] (\x, 4) rectangle ++(1, 1);
        }
        \foreach \x in {1,...,4} {
            \draw[thin, gray, densely dashed] (\x, 5) rectangle ++(1, 1);
        }
        \foreach \x in {1,...,1} {
            \draw[thin, gray, densely dashed] (\x, 6) rectangle ++(1, 1);
        }
        \foreach \x in {-1,...,1} {
            \draw[thin, gray, densely dashed] (\x, 7) rectangle ++(1, 1);
        }
        \foreach \x in {-3,...,-1} {
            \draw[thin, gray, densely dashed] (\x, 8) rectangle ++(1, 1);
        }
        \foreach \x in {-3,...,-3} {
            \draw[thin, gray, densely dashed] (\x, 9) rectangle ++(1, 1);
        }
       
        \foreach \x in {1,...,2} {
            \draw[thin, gray, densely dashed] (11, \x) -- (11.3, \x);
        }

        \draw[thick] (6, 3) -- (7, 3) -- (7, 1);
        \draw[thick] (6, 3) -- (4, 3) -- (4, 5) -- (1, 5) -- (1, 7);
        \draw[thick] (1, 7) -- (-1, 7) -- (-1, 8) -- (-3, 8) -- (-3, 10);

        \fill[gray] (12, 1.5) node {$\theta$};
    \end{tikzpicture}
    \subcaption{Empty-shape $\varepsilon(\hat{\theta})$.}
    \label{fig:sub-ribbons-of-infinite-ribbon-1-d}
    \end{subfigure}
    \caption{Examples of ribbons and shapes.}
    \label{fig:sub-ribbons-of-infinite-ribbon-1}
\end{figure}
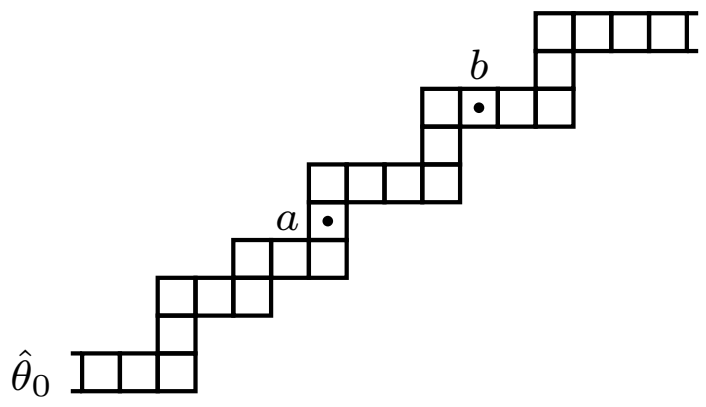
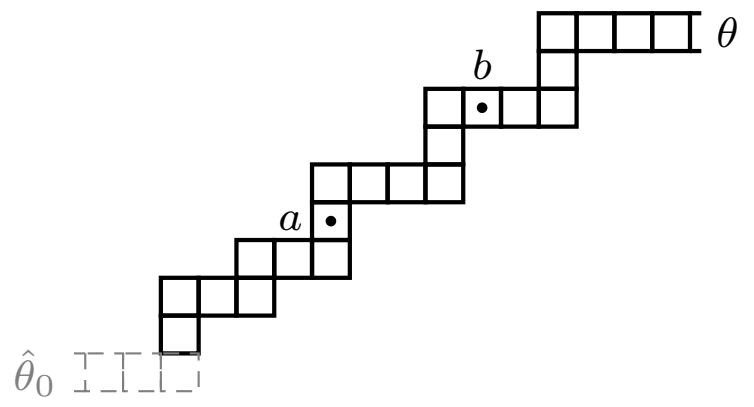
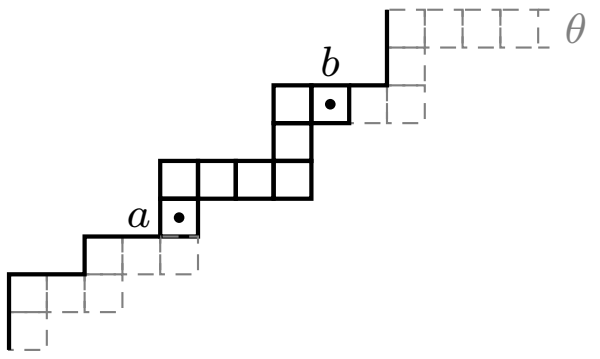
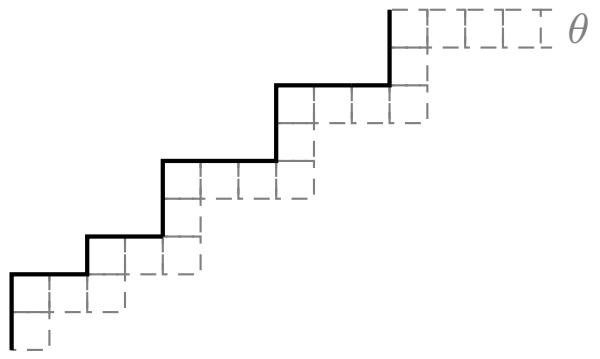

\begin{figure}
    \centering
    \begin{tikzpicture}[y=-1cm, scale=0.4, transform shape]
	\foreach \x in {3,...,5} {
        \draw[very thick] (\x, 1) rectangle ++(1, 1);
    }
    \foreach \x in {2,...,5} {
        \draw[very thick] (\x, 2) rectangle ++(1, 1);
    }
    \foreach \x in {2,...,4} {
        \draw[very thick] (\x, 3) rectangle ++(1, 1);
    }
    \foreach \x in {1,...,4} {
        \draw[very thick] (\x, 4) rectangle ++(1, 1);
    }
    \foreach \x in {1,...,2} {
        \draw[very thick] (\x, 5) rectangle ++(1, 1);
    }

    \foreach \x in {-4,...,2} {
        \draw[gray, dashed] (\x, 1) rectangle ++(1, 1);
    }
    \foreach \x in {-4,...,1} {
        \draw[gray, dashed] (\x, 2) rectangle ++(1, 1);
    }
    \foreach \x in {-4,...,1} {
        \draw[gray, dashed] (\x, 3) rectangle ++(1, 1);
    }
    \foreach \x in {-4,...,0} {
        \draw[gray, dashed] (\x, 4) rectangle ++(1, 1);
    }
    \foreach \x in {-4,...,0} {
        \draw[gray, dashed] (\x, 5) rectangle ++(1, 1);
    }


    \foreach \x in {6,...,11} {
        \draw[gray, dashed] (\x, 1) rectangle ++(1, 1);
    }
    \foreach \x in {6,...,11} {
        \draw[gray, dashed] (\x, 2) rectangle ++(1, 1);
    }
    \foreach \x in {5,...,11} {
        \draw[gray, dashed] (\x, 3) rectangle ++(1, 1);
    }
    \foreach \x in {5,...,11} {
        \draw[gray, dashed] (\x, 4) rectangle ++(1, 1);
    }
    \foreach \x in {3,...,11} {
        \draw[gray, dashed] (\x, 5) rectangle ++(1, 1);
    }

    \foreach \x in {1,...,6} {
        \draw[gray] (-4, \x) -- (-4.3, \x);
    }
    \foreach \x in {1,...,6} {
        \draw[gray] (12, \x) -- (12.3, \x);
    }

	\draw [red, thick] (1.5, 6) -- (1.5, 4.5) -- (3.5, 4.5) -- (3.5, 3.5) -- (5, 3.5);
    \draw [red, thick] (5.5, 3) -- (5.5, 1);
    \draw [red, thick, densely dashed] (5, 3.5) -- (5.5, 3.5) -- (5.5, 3);

    \draw [red, thick, densely dashed] (-4.5, 5.5) -- (0.5, 5.5) -- (0.5, 3.5) -- (2, 3.5);
    \draw [red, thick] (2, 3.5) -- (2.5, 3.5) -- (2.5, 2.5) -- (4.5, 2.5) -- (4.5, 1);

    \draw [red, thick, densely dashed] (-4.5, 4.5) -- (-0.5, 4.5) -- (-0.5, 2.5) -- (1.5, 2.5) -- (1.5, 1.5) -- (3, 1.5);
    \draw [red, thick] (3, 1.5) -- (3.5, 1.5) -- (3.5, 1);
    \draw [red, thick, densely dashed] (-4.5, 3.5) -- (-1.5, 3.5) -- (-1.5, 1.5) -- (0.5, 1.5) -- (0.5, 1);
    \draw [red, thick, densely dashed] (-4.5, 2.5) -- (-2.5, 2.5) -- (-2.5, 1);
    \draw [red, thick, densely dashed] (-4.5, 1.5) -- (-3.5, 1.5) -- (-3.5, 1);

    \draw [red, thick] (2.5, 6) -- (2.5, 5.5) -- (3, 5.5);
    \draw [red, thick] (4.5, 5) -- (4.5, 4.5) -- (5, 4.5);
    
    \draw [red, thick, densely dashed] (2.5, 6) -- (2.5, 5.5) -- (4.5, 5.5) -- (4.5, 4.5) -- (6.5, 4.5) -- (6.5, 1.5) -- (12.5, 1.5);
    \draw [red, thick, densely dashed] (5.5, 6) -- (5.5, 5.5) -- (7.5, 5.5) -- (7.5, 2.5) -- (12.5, 2.5);
    \draw [red, thick, densely dashed] (8.5, 6) -- (8.5, 3.5) -- (12.5, 3.5);
    \draw [red, thick, densely dashed] (9.5, 6) -- (9.5, 4.5) -- (12.5, 4.5);
    \draw [red, thick, densely dashed] (10.5, 6) -- (10.5, 5.5) -- (12.5, 5.5);

    \node at (-5, 3.5) {$\ldots$};
    \node at (13, 3.5) {$\ldots$};

    \end{tikzpicture}
\qquad
    \begin{tikzpicture}[y=-1cm, scale=0.4, transform shape]
	\foreach \x in {3,...,3} {
        \draw[very thick] (\x, 1) rectangle ++(1, 1);
    }
    \foreach \x in {1,...,3} {
        \draw[very thick] (\x, 2) rectangle ++(1, 1);
    }
    \foreach \x in {1,...,1} {
        \draw[very thick] (\x, 3) rectangle ++(1, 1);
    }

    \draw [very thick] (4, 1) -- (4, 0);
    \draw [very thick] (0, 6) -- (0, 4) -- (1, 4);

    \foreach \x in {-4,...,0} {
        \draw[gray, dashed] (\x, 2) rectangle ++(1, 1);
    }
    \foreach \x in {-4,...,0} {
        \draw[gray, dashed] (\x, 3) rectangle ++(1, 1);
    }
    \foreach \x in {-4,...,-1} {
        \draw[gray, dashed] (\x, 4) rectangle ++(1, 1);
    }
    \foreach \x in {-4,...,-1} {
        \draw[gray, dashed] (\x, 5) rectangle ++(1, 1);
    }

    \foreach \x in {4,...,11} {
        \draw[gray, dashed] (\x, 0) rectangle ++(1, 1);
    }
    \foreach \x in {4,...,11} {
        \draw[gray, dashed] (\x, 1) rectangle ++(1, 1);
    }
    \foreach \x in {4,...,11} {
        \draw[gray, dashed] (\x, 2) rectangle ++(1, 1);
    }
    \foreach \x in {2,...,11} {
        \draw[gray, dashed] (\x, 3) rectangle ++(1, 1);
    }

    \draw [red, thick, densely dashed] (-4.5, 5.5) -- (-0.5, 5.5) -- (-0.5, 3.5) -- (1.5, 3.5) -- (1.5, 2.5) -- (3.5, 2.5) -- (3.5, 1);
    \draw [red, thick, densely dashed] (-4.5, 4.5) -- (-1.5, 4.5) -- (-1.5, 2.5) -- (0.5, 2.5) -- (0.5, 2);
    \draw [red, thick, densely dashed] (-4.5, 3.5) -- (-2.5, 3.5) -- (-2.5, 2);
    \draw [red, thick, densely dashed] (-4.5, 2.5) -- (-3.5, 2.5) -- (-3.5, 2);

    \draw [red, thick, densely dashed] (2.5, 4) -- (2.5, 3.5) -- (4.5, 3.5) -- (4.5, 0.5) -- (12.5, 0.5);
    \draw [red, thick, densely dashed] (5.5, 4) -- (5.5, 1.5) -- (12.5, 1.5);
    \draw [red, thick, densely dashed] (6.5, 4) -- (6.5, 2.5) -- (12.5, 2.5);
    \draw [red, thick, densely dashed] (7.5, 4) -- (7.5, 3.5) -- (12.5, 3.5);

    \foreach \x in {2,...,6} {
        \draw[gray] (-4, \x) -- (-4.3, \x);
    }
    \foreach \x in {0,...,4} {
        \draw[gray] (12, \x) -- (12.3, \x);
    }

    \node at (-5, 3.5) {$\ldots$};
    \node at (13, 3.5) {$\ldots$};

    \end{tikzpicture}
    \caption{(Left) A shape $\lambda/\mu$ with $L(\lambda/\mu) \cup R(\lambda/\mu)$; (Right) $\theta_3 \# \theta_3$ with $L(\theta_3 \# \theta_3) \cup R(\theta_3 \# \theta_3)$.}
    \label{pic:ribbon_with_p_cells}
\end{figure}
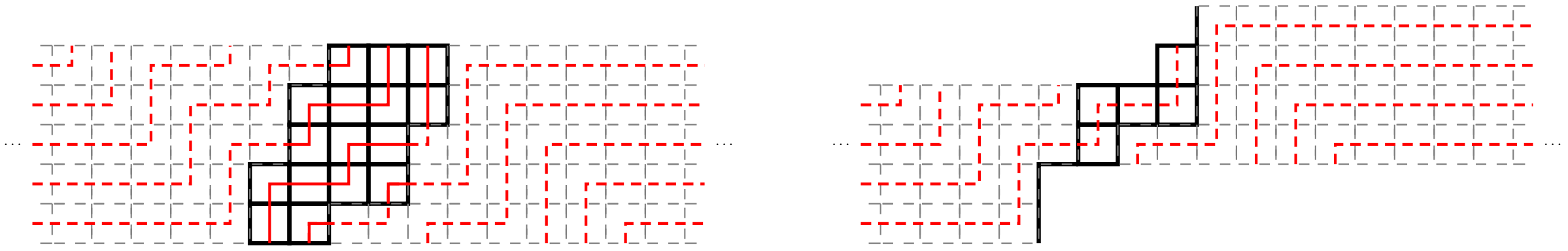

\section{Ribbon decomposition determinants for flagged Grothendieck enumerators}\label{sec:det}

\subsection{Induced ribbon flags}

\begin{definition}
    Let $\lambda/\mu$ be skew shape of signatures from $\mathsf{Sign}_n$ and $\mathbf{a} = (a_i), \mathbf{b} = (b_i)$ be $\mathsf{G}$-flags. 
    Define the following sets:  
    \begin{align*}
    L(\lambda/\mu) &:= \{(i, j)  \ |\ i\in [n],\, a_i \ne -\infty, \, j {\le \mu_i} \} \subset \mathbb{Z}^2, \\ 
    R(\lambda/\mu) &:= \{(i, j) \ |\ i\in [n],\, b_i \ne -\infty,\, j {\ge \lambda_i + 1} \} \subset \mathbb{Z}^2.
    \end{align*}

For $c \in \mathbb{Z}$  
we define $\max(\lambda/\mu, c)$ as the rightmost cell in $L(\lambda/\mu)$ with content $c$, i.e. $\max(\lambda/\mu, c) := (p, q) \in L(\lambda/\mu)$ with the largest $q$ such that $q - p = c$, and if there is no cell in $L(\lambda/\mu)$ with content $c$ we set $\max(\lambda/\mu, c) =$ undefined. Similarly, 
we define $\min(\lambda/\mu, c)$ as the leftmost cell in $R(\lambda/\mu)$ with content $c$, i.e. $\min(\lambda/\mu, c) := (p, q) \in R(\lambda/\mu)$ with the smallest $q$ such that $q - p = c$, and  if there is no cell $R(\lambda/\mu)$ with content $c$ we set $\min(\lambda/\mu, c) = $ undefined.
\end{definition}

\begin{definition}[Induced ribbon $\mathsf{G}$-flags]
Given outer ribbon decomposition $(\theta_1, \ldots, \theta_k)$. 
Let $n^{ij}$ be number of rows in the shape $\theta_i \# \theta_j$. We define {\it induced ribbon $\mathsf{G}$-flags} $\mathbf{a}^{ij} = (a^{ij}_r), \mathbf{b}^{ij} = (b^{ij}_r)$ both of length $n^{ij}$ for the shape $\theta_i \# \theta_j$ as follows: 
Let $\zeta_r$ be the rightmost cell in $r$-th row of $\mathcal{L}(\theta_i \# \theta_j)$, and let $\eta_r$ be the leftmost cell in $r$-th row of $\mathcal{R}(\theta_i \# \theta_j)$. Denote $M_r = \max(\lambda/\mu, c(\zeta_r)),\ m_r = \min(\lambda/\mu, c(\eta))$. Then we set 
\begin{align*}
    a^{ij}_r = 
    \begin{cases}
        -\infty, & \text{\text{if $M_r$ is undefined}}, \\
        a_{\mathrm{row}(M_r)}, & \text{otherwise},
    \end{cases}
    \qquad 
    b^{ij}_r = 
    \begin{cases}
        -\infty & \text{if $m_r$ is undefined}, \\
        b_{\mathrm{row}(m_r)} & \text{otherwise}.
    \end{cases}
\end{align*}
Where by $\mathrm{row}(\gamma)$ we denote the row number of the cell $\gamma$. 

In the same way we define induced ribbon $\mathsf{G}$-flags $\mathbf{a}^{\varepsilon}, \mathbf{b}^{\varepsilon}$ for an empty ribbon $\varepsilon$. 
\end{definition}

It is easy to check that induced ribbon $\mathsf{G}$-flags are indeed $\mathsf{G}$-flags. 
\begin{proposition}
The induced ribbon $\mathsf{G}$-flags $\mathbf{a}^{ij}, \mathbf{b}^{ij}$ are $\mathsf{G}$-flags for the shape $\theta_i \# \theta_j$. 
\end{proposition}

\begin{example}
    \label{ex:p-cells}
    Let $\lambda = (5,5,4,4,2),\ \mu = (2,1,1,0,0)$ with $\mathsf{G}$-flags $\mathbf{a} = (a_1,a_2,a_3,a_4,a_5),$ 
    $\mathbf{b} = (b_1, b_2, b_3, b_4, b_5)$. 
    Define pipe vector $\pi$ as follows: 
    \[
    \begin{array}{rclrclrclrcl}
        -5 & \mapsto & (0, 1) \quad & -2 & \mapsto & (0, 0) \quad & 1 & \mapsto & (0, 0) \quad & 4 & \mapsto & (1, 1) \\ 
        -4 & \mapsto & (1, 1) \quad & -1 & \mapsto & (0, 1) \quad & 2 & \mapsto & (0, 1) \quad & 5 & \mapsto & (1, 0) \\ 
        -3 & \mapsto & (1, 0) \quad &  0 & \mapsto & (1, 0) \quad & 3 & \mapsto & (1, 1) \quad & c & \mapsto & (0, 0), 
    \end{array}
    \]
    for $c \in \mathbb{Z} \setminus [-5, 5]$. 
     In Fig.~\ref{pic:ribbon_with_p_cells} (left) we show $\lambda/\mu$ with the sets $L(\lambda/\mu), R(\lambda/\mu)$.
    Next we find the induced ribbon $\mathsf{G}$-flags $\mathbf{a}^{3,3}, \mathbf{b}^{3,3}$ for the shape $\theta_3 \# \theta_3$, see Fig.~\ref{pic:ribbon_with_p_cells} (right). 
    The number of rows in $\theta_3 \# \theta_3$ is $n_{3,3} = 6$ and we have 
    \begin{align*}
        \mathbf{a}^{3,3} = (-\infty, -\infty, a_2, a_3, a_4, a_5), \quad \mathbf{b}^{3,3} = (b_1, b_2, b_3,  b_5, -\infty, -\infty). 
    \end{align*}
\end{example}


\subsection{Main formula}

We are now ready to present our main Hamel--Goulden-type formulas for flagged Grothendieck numerators. 

\begin{theorem}[Hamel--Goulden-type formulas for flagged Grothendieck numerators]
    \label{thm:HG-sfG}
    Let $\Theta = (\theta_1, \ldots, \theta_k)$ be an outer ribbon decomposition of skew signature $\lambda/\mu$  with associated pipe function $\pi$. Let $\mathbf{a}, \mathbf{b}$ be $\mathsf{G}$-flags and $\mathbf{a}^{ij}, \mathbf{b}^{ij}$ be induced ribbon $\mathsf{G}$-flags for the ribbons $\theta_i \# \theta_j$. Then the following formula holds:  
    \begin{align*}
        \mathsf{G}^{\mathbf{a}, \mathbf{b}}_{\lambda/\mu}(\mathbf{x} ; \mathbf{U} , \mathbf{V}) 
        = (-1)^{\mathsf{hp}^{\pi}_{\lambda/\mu}} \cdot 
        \mathsf{G}^{\mathbf{a}^{\varepsilon}, \mathbf{b}^{\varepsilon}}_{\varepsilon}(\mathbf{x} ;  \mathbf{U} , \mathbf{V})^{-k+1} 
        \cdot \det \left[ (-1)^{\mathsf{hp}^{\pi}_{\theta_i \# \theta_j}} \cdot \mathsf{G}^{\mathbf{a}^{ij}, \mathbf{b}^{ij}}_{\theta_i \# \theta_j}(\mathbf{x} ; \mathbf{U} , \mathbf{V})\right]_{1 \le i,j \le k}. 
    \end{align*}
\end{theorem}

\medskip

By taking $G$-specialization and using Proposition~\ref{prop:gspec} we immediately obtain the following Hamel--Goulden-type formulas for canonical refined stable Grothendieck polynomials. 

\begin{corollary}[Hamel--Goulden-type formulas for canonical refined stable Grothendieck polynomials]
    Let $(\theta_1, \ldots, \theta_k)$ be an outer ribbon decomposition of skew signature $\lambda/\mu$ for $\lambda, \mu \in \mathsf{Par}_n$ with associated pipe function $\pi$. Let $\mathbf{a} = (i), \mathbf{b} = (\lambda_i)$ be $\mathsf{G}$-flags and $\mathbf{a}^{ij}, \mathbf{b}^{ij}$ be corresponding induced ribbon $\mathsf{G}$-flags for the ribbons $\theta_i \# \theta_j$. Then the following formula holds:  
    \begin{align*}
        {G}^{}_{\lambda/\mu}(\mathbf{x}_m ; \boldsymbol{\alpha}, \boldsymbol{\beta}) 
        =  \prod_{i \in [n]} \prod_{j \in [m]} (1 - \beta_i x_j) \cdot 
        \mathsf{G}^{\mathbf{a}^{\varepsilon}, \mathbf{b}^{\varepsilon}}_{\varepsilon}(\mathbf{x}_m ; \boldsymbol{\alpha}, \boldsymbol{\beta})^{-k+1} \cdot 
        \det \left[(-1)^{\mathsf{hp}^{\pi}_{\theta_i \# \theta_j}} \cdot \mathsf{G}^{\mathbf{a}^{ij}, \mathbf{b}^{ij}}_{\theta_i \# \theta_j}(\mathbf{x}_m ; \boldsymbol{\alpha}, \boldsymbol{\beta}) \right]_{1 \le i,j \le k}. 
    \end{align*}
\end{corollary}

\begin{figure}
    \centering
    \begin{subfigure}[c]{0.45\textwidth}
    \centering
    \begin{tikzpicture}[y=-1cm, scale=0.4, font=\small]
	\foreach \x in {3,...,4} {
        \draw[thick] (\x, 1) rectangle ++(1, 1);
    }
    \foreach \x in {2,...,4} {
        \draw[thick] (\x, 2) rectangle ++(1, 1);
    }
    \foreach \x in {2,...,3} {
        \draw[thick] (\x, 3) rectangle ++(1, 1);
    }

    \foreach \x in {-1,...,2} {
        \draw[gray, dashed] (\x, 1) rectangle ++(1, 1);
    }
    \foreach \x in {-1,...,1} {
        \draw[gray, dashed] (\x, 2) rectangle ++(1, 1);
    }
    \foreach \x in {-1,...,1} {
        \draw[gray, dashed] (\x, 3) rectangle ++(1, 1);
    }


    \foreach \x in {5,...,7} {
        \draw[gray, dashed] (\x, 1) rectangle ++(1, 1);
    }
    \foreach \x in {5,...,7} {
        \draw[gray, dashed] (\x, 2) rectangle ++(1, 1);
    }
    \foreach \x in {4,...,7} {
        \draw[gray, dashed] (\x, 3) rectangle ++(1, 1);
    }

    \foreach \x in {1,...,4} {
        \draw[gray] (-1, \x) -- (-1.3, \x);
    }
    \foreach \x in {1,...,4} {
        \draw[gray] (8, \x) -- (8.3, \x);
    }

    \draw [red, thick] (2.5, 4) -- (2.5, 3.5) -- (4, 3.5);
    \draw [red, thick] (4.5, 3) -- (4.5, 1.5) -- (5, 1.5);
    \draw [red, thick] (2, 2.5) -- (3.5, 2.5) -- (3.5, 1);
    \draw [red, thick, densely dashed] (-1.5, 3.5) -- (1.5, 3.5) -- (1.5, 2.5) -- (2, 2.5);
    \draw [red, thick, densely dashed] (-1.5, 2.5) -- (0.5, 2.5) -- (0.5, 1.5) -- (2.5, 1.5) -- (2.5, 1);
    \draw [red, thick, densely dashed] (-1.5, 1.5) -- (-0.5, 1.5) -- (-0.5, 1);
    \draw [red, thick, densely dashed] (4, 3.5) -- (4.5, 3.5) -- (4.5, 3);
    \draw [red, thick, densely dashed] (5, 1.5) -- (8.5, 1.5);
    \draw [red, thick, densely dashed] (5.5, 4) -- (5.5, 2.5) -- (8.5, 2.5);
    \draw [red, thick, densely dashed] (6.5, 4) -- (6.5, 3.5) -- (8.5, 3.5);


    \node at (2.5, 4.75) {$\theta_1$};
    \node at (3.5, 0.25) {$\theta_2$};
    \node at (4.5, 0.25) {$\theta_3$};
    \end{tikzpicture}
    \end{subfigure}
    \hspace{0.5cm}
    \begin{subfigure}[c]{0.45\textwidth}
        \centering
        \begin{tikzpicture}[y=-1cm, scale=0.4]
        \draw[black, thick] (2, 4) -- (2, 3) -- (4, 3) -- (4, 1);

        \foreach \x in {0,...,3} {
            \draw[gray, dashed] (\x, 2) rectangle ++(1, 1);
        }
        \foreach \x in {0,...,1} {
            \draw[gray, dashed] (\x, 3) rectangle ++(1, 1);
        }


        \foreach \x in {4,...,5} {
            \draw[gray, dashed] (\x, 1) rectangle ++(1, 1);
        }
        \foreach \x in {4,...,5} {
            \draw[gray, dashed] (\x, 2) rectangle ++(1, 1);
        }

        \foreach \x in {2,...,4} {
            \draw[gray] (-0, \x) -- (-0.3, \x);
        }
        \foreach \x in {1,...,3} {
            \draw[gray] (6, \x) -- (6.3, \x);
        }

        \draw [red, thick, densely dashed] (-0.5, 3.5) -- (1.5, 3.5) -- (1.5, 2.5) -- (3.5, 2.5) -- (3.5, 2);
        \draw [red, thick, densely dashed] (-0.5, 2.5) -- (0.5, 2.5) -- (0.5, 2);
        \draw [red, thick, densely dashed] (4.5, 3) -- (4.5, 1.5) -- (5, 1.5) -- (6.5, 1.5);
        \draw [red, thick, densely dashed] (5.5, 3) -- (5.5, 2.5) -- (6.5, 2.5);
        \end{tikzpicture}
    \end{subfigure}
    \caption{A skew shape $\lambda/\mu$ with outer ribbon decomposition (on the left) and corresponding empty ribbon $\varepsilon$ (on the right).}
    \label{fig:example-HG-1}
\end{figure}
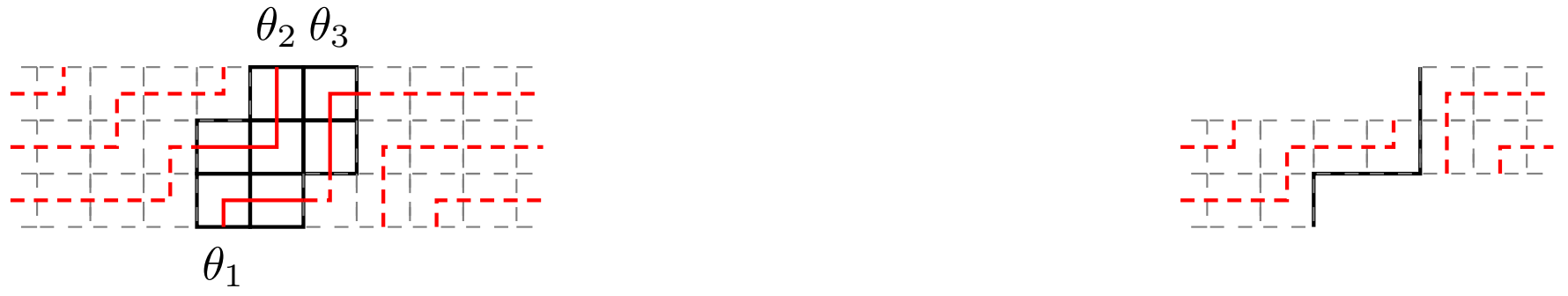
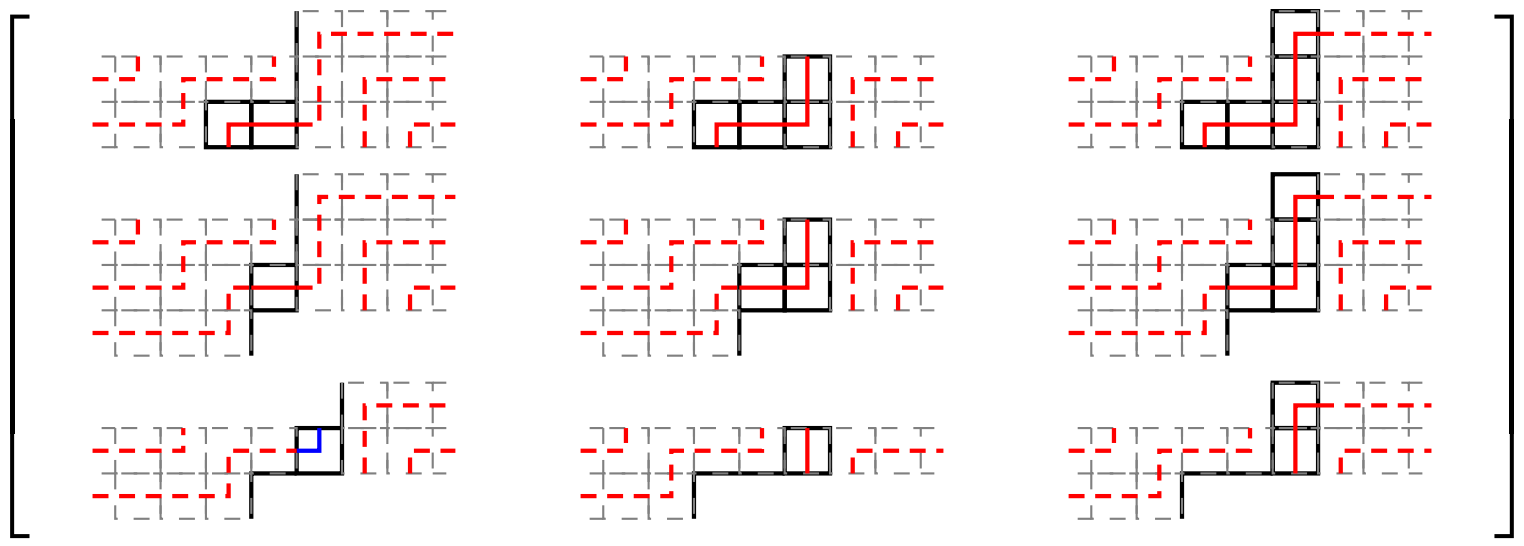
\begin{figure}
\[\arraycolsep=0.4cm\def\arraystretch{1}
\left[
\begin{array}{ccc}
    \begin{tikzpicture}[y=-1cm, scale=0.3]
    \draw[black, thick] (4, 3) -- (4, 1);

    \foreach \x in {2,...,3} {
        \draw[thick] (\x, 3) rectangle ++(1, 1);
    }


    \foreach \x in {0,...,3} {
        \draw[gray, thin, dashed] (\x, 2) rectangle ++(1, 1);
    }
    \foreach \x in {0,...,1} {
        \draw[gray, thin, dashed] (\x, 3) rectangle ++(1, 1);
    }


    \foreach \x in {4,...,6} {
        \draw[gray, thin, dashed] (\x, 1) rectangle ++(1, 1);
    }
    \foreach \x in {4,...,6} {
        \draw[gray, thin, dashed] (\x, 2) rectangle ++(1, 1);
    }
    \foreach \x in {4,...,6} {
        \draw[gray, thin, dashed] (\x, 3) rectangle ++(1, 1);
    }

    \foreach \x in {2,...,4} {
        \draw[gray] (0, \x) -- (-0.3, \x);
    }
    \foreach \x in {1,...,4} {
        \draw[gray] (7, \x) -- (7.3, \x);
    }

    \draw [red, thick] (2.5, 4) -- (2.5, 3.5) -- (4, 3.5);
    \draw [red, thick, densely dashed] (-0.5, 3.5) -- (1.5, 3.5) -- (1.5, 2.5) -- (2, 2.5) -- (3.5, 2.5) -- (3.5, 2);
    \draw [red, thick, densely dashed] (-0.5, 2.5) -- (0.5, 2.5) -- (0.5, 2);
    \draw [red, thick, densely dashed] (4, 3.5) -- (4.5, 3.5) -- (4.5, 3);
    \draw [red, thick, densely dashed] (4.5, 3) -- (4.5, 1.5) -- (5, 1.5) -- (7.5, 1.5);
    \draw [red, thick, densely dashed] (5.5, 4) -- (5.5, 2.5) -- (7.5, 2.5);
    \draw [red, thick, densely dashed] (6.5, 4) -- (6.5, 3.5) -- (7.5, 3.5);
    \end{tikzpicture} & 
    \begin{tikzpicture}[y=-1cm, scale=0.3]

    \foreach \x in {4,...,4} {
        \draw[thick] (\x, 2) rectangle ++(1, 1);
    }
    \foreach \x in {2,...,4} {
        \draw[thick] (\x, 3) rectangle ++(1, 1);
    }


    \foreach \x in {0,...,3} {
        \draw[gray, thin, dashed] (\x, 2) rectangle ++(1, 1);
    }
    \foreach \x in {0,...,1} {
        \draw[gray, thin, dashed] (\x, 3) rectangle ++(1, 1);
    }


    \foreach \x in {4,...,6} {
        \draw[gray, thin, dashed] (\x, 2) rectangle ++(1, 1);
    }
    \foreach \x in {4,...,6} {
        \draw[gray, thin, dashed] (\x, 3) rectangle ++(1, 1);
    }
    \foreach \x in {2,...,4} {
        \draw[gray] (0, \x) -- (-0.3, \x);
    }
    \foreach \x in {2,...,4} {
        \draw[gray] (7, \x) -- (7.3, \x);
    }

    \draw [red, thick] (2.5, 4) -- (2.5, 3.5) -- (4.5, 3.5) -- (4.5, 2);
    \draw [red, thick, densely dashed] (-0.5, 3.5) -- (1.5, 3.5) -- (1.5, 2.5) -- (2, 2.5) -- (3.5, 2.5) -- (3.5, 2);
    \draw [red, thick, densely dashed] (-0.5, 2.5) -- (0.5, 2.5) -- (0.5, 2);
    \draw [red, thick, densely dashed] (5.5, 4) -- (5.5, 2.5) -- (7.5, 2.5);
    \draw [red, thick, densely dashed] (6.5, 4) -- (6.5, 3.5) -- (7.5, 3.5);
    \end{tikzpicture} &
    \begin{tikzpicture}[y=-1cm, scale=0.3]
    \draw[black, thick] (4, 3) -- (4, 1);

    \foreach \x in {4,...,4} {
        \draw[thick] (\x, 1) rectangle ++(1, 1);
    }
    \foreach \x in {4,...,4} {
        \draw[thick] (\x, 2) rectangle ++(1, 1);
    }
    \foreach \x in {2,...,4} {
        \draw[thick] (\x, 3) rectangle ++(1, 1);
    }


    \foreach \x in {0,...,3} {
        \draw[gray, thin, dashed] (\x, 2) rectangle ++(1, 1);
    }
    \foreach \x in {0,...,1} {
        \draw[gray, thin, dashed] (\x, 3) rectangle ++(1, 1);
    }


    \foreach \x in {4,...,6} {
        \draw[gray, thin, dashed] (\x, 1) rectangle ++(1, 1);
    }
    \foreach \x in {4,...,6} {
        \draw[gray, thin, dashed] (\x, 2) rectangle ++(1, 1);
    }
    \foreach \x in {4,...,6} {
        \draw[gray, thin, dashed] (\x, 3) rectangle ++(1, 1);
    }

    \foreach \x in {2,...,4} {
        \draw[gray] (0, \x) -- (-0.3, \x);
    }
    \foreach \x in {1,...,4} {
        \draw[gray] (7, \x) -- (7.3, \x);
    }

    \draw [red, thick] (2.5, 4) -- (2.5, 3.5) -- (4.5, 3.5) -- (4.5, 1.5) -- (5, 1.5);
    \draw [red, thick, densely dashed] (-0.5, 3.5) -- (1.5, 3.5) -- (1.5, 2.5) -- (2, 2.5) -- (3.5, 2.5) -- (3.5, 2);
    \draw [red, thick, densely dashed] (-0.5, 2.5) -- (0.5, 2.5) -- (0.5, 2);
    \draw [red, thick, densely dashed] (5, 1.5) -- (7.5, 1.5);
    \draw [red, thick, densely dashed] (5.5, 4) -- (5.5, 2.5) -- (7.5, 2.5);
    \draw [red, thick, densely dashed] (6.5, 4) -- (6.5, 3.5) -- (7.5, 3.5);
    \end{tikzpicture} \\

    \begin{tikzpicture}[y=-1cm, scale=0.3]
    \draw[black, thick] (3, 3) -- (3, 1);
    \foreach \x in {2,...,2} {
        \draw[thick] (\x, 3) rectangle ++(1, 1);
    }
    \draw[black, thick] (2, 4) -- (2, 5);


    \foreach \x in {-1,...,3} {
        \draw[gray, thin, dashed] (\x, 2) rectangle ++(1, 1);
    }
    \foreach \x in {-1,...,1} {
        \draw[gray, thin, dashed] (\x, 3) rectangle ++(1, 1);
    }
    \foreach \x in {-1,...,1} {
        \draw[gray, thin, dashed] (\x, 4) rectangle ++(1, 1);
    }


    \foreach \x in {3,...,5} {
        \draw[gray, thin, dashed] (\x, 1) rectangle ++(1, 1);
    }
    \foreach \x in {3,...,5} {
        \draw[gray, thin, dashed] (\x, 2) rectangle ++(1, 1);
    }
    \foreach \x in {3,...,5} {
        \draw[gray, thin, dashed] (\x, 3) rectangle ++(1, 1);
    }

    \foreach \x in {2,...,4} {
        \draw[gray] (-1, \x) -- (-1.3, \x);
    }
    \foreach \x in {1,...,4} {
        \draw[gray] (6, \x) -- (6.3, \x);
    }

    \draw [red, thick] (2, 3.5) -- (3, 3.5);
    \draw [red, thick, densely dashed] (3, 3.5) -- (3.5, 3.5) -- (3.5, 1.5) -- (6.5, 1.5);
    \draw [red, thick, densely dashed] (-1.5, 4.5) -- (1.5, 4.5) -- (1.5, 3.5) -- (2, 3.5);
    \draw [red, thick, densely dashed] (-1.5, 3.5) -- (0.5, 3.5) -- (0.5, 2.5) -- (2.5, 2.5) -- (2.5, 2);
    \draw [red, thick, densely dashed] (-1.5, 2.5) -- (-0.5, 2.5) -- (-0.5, 2);
    \draw [red, thick, densely dashed] (4.5, 4) -- (4.5, 2.5) -- (6.5, 2.5);
    \draw [red, thick, densely dashed] (5.5, 4) -- (5.5, 3.5) -- (6.5, 3.5);
    \end{tikzpicture} &
    \begin{tikzpicture}[y=-1cm, scale=0.3]
    \foreach \x in {3,...,3} {
        \draw[thick] (\x, 2) rectangle ++(1, 1);
    }
    \foreach \x in {2,...,3} {
        \draw[thick] (\x, 3) rectangle ++(1, 1);
    }
    \draw[black, thick] (2, 4) -- (2, 5);


    \foreach \x in {-1,...,3} {
        \draw[gray, thin, dashed] (\x, 2) rectangle ++(1, 1);
    }
    \foreach \x in {-1,...,1} {
        \draw[gray, thin, dashed] (\x, 3) rectangle ++(1, 1);
    }
    \foreach \x in {-1,...,1} {
        \draw[gray, thin, dashed] (\x, 4) rectangle ++(1, 1);
    }


    \foreach \x in {4,...,5} {
        \draw[gray, thin, dashed] (\x, 2) rectangle ++(1, 1);
    }
    \foreach \x in {4,...,5} {
        \draw[gray, thin, dashed] (\x, 3) rectangle ++(1, 1);
    }

    \foreach \x in {2,...,4} {
        \draw[gray] (-1, \x) -- (-1.3, \x);
    }
    \foreach \x in {2,...,4} {
        \draw[gray] (6, \x) -- (6.3, \x);
    }

    \draw [red, thick] (2, 3.5) -- (3.5, 3.5) -- (3.5, 2);
    \draw [red, thick, densely dashed] (-1.5, 4.5) -- (1.5, 4.5) -- (1.5, 3.5) -- (2, 3.5);
    \draw [red, thick, densely dashed] (-1.5, 3.5) -- (0.5, 3.5) -- (0.5, 2.5) -- (2.5, 2.5) -- (2.5, 2);
    \draw [red, thick, densely dashed] (-1.5, 2.5) -- (-0.5, 2.5) -- (-0.5, 2);
    \draw [red, thick, densely dashed] (4.5, 4) -- (4.5, 2.5) -- (6.5, 2.5);
    \draw [red, thick, densely dashed] (5.5, 4) -- (5.5, 3.5) -- (6.5, 3.5);
    \end{tikzpicture} &
    \begin{tikzpicture}[y=-1cm, scale=0.3]
    \foreach \x in {3,...,3} {
        \draw[thick] (\x, 1) rectangle ++(1, 1);
    }
    \foreach \x in {3,...,3} {
        \draw[thick] (\x, 2) rectangle ++(1, 1);
    }
    \foreach \x in {2,...,3} {
        \draw[thick] (\x, 3) rectangle ++(1, 1);
    }
    \draw[black, thick] (2, 4) -- (2, 5);


    \foreach \x in {-1,...,3} {
        \draw[gray, thin, dashed] (\x, 2) rectangle ++(1, 1);
    }
    \foreach \x in {-1,...,1} {
        \draw[gray, thin, dashed] (\x, 3) rectangle ++(1, 1);
    }
    \foreach \x in {-1,...,1} {
        \draw[gray, thin, dashed] (\x, 4) rectangle ++(1, 1);
    }

    \foreach \x in {4,...,5} {
        \draw[gray, thin, dashed] (\x, 1) rectangle ++(1, 1);
    }
    \foreach \x in {4,...,5} {
        \draw[gray, thin, dashed] (\x, 2) rectangle ++(1, 1);
    }
    \foreach \x in {4,...,5} {
        \draw[gray, thin, dashed] (\x, 3) rectangle ++(1, 1);
    }

    \foreach \x in {2,...,4} {
        \draw[gray] (-1, \x) -- (-1.3, \x);
    }
    \foreach \x in {1,...,4} {
        \draw[gray] (6, \x) -- (6.3, \x);
    }

    \draw [red, thick] (2, 3.5) -- (3.5, 3.5) -- (3.5, 1.5) -- (4, 1.5);
    \draw [red, thick, densely dashed] (4, 1.5) -- (6.5, 1.5);
    \draw [red, thick, densely dashed] (-1.5, 4.5) -- (1.5, 4.5) -- (1.5, 3.5) -- (2, 3.5);
    \draw [red, thick, densely dashed] (-1.5, 3.5) -- (0.5, 3.5) -- (0.5, 2.5) -- (2.5, 2.5) -- (2.5, 2);
    \draw [red, thick, densely dashed] (-1.5, 2.5) -- (-0.5, 2.5) -- (-0.5, 2);
    \draw [red, thick, densely dashed] (4.5, 4) -- (4.5, 2.5) -- (6.5, 2.5);
    \draw [red, thick, densely dashed] (5.5, 4) -- (5.5, 3.5) -- (6.5, 3.5);
    \end{tikzpicture} \\

    \begin{tikzpicture}[y=-1cm, scale=0.3]
    \draw[black, thick] (4, 3) -- (4, 2);
    \foreach \x in {3,...,3} {
        \draw[thick] (\x, 3) rectangle ++(1, 1);
    }
    \draw[black, thick] (2, 5) -- (2, 4) -- (3, 4);


    \foreach \x in {-1,...,2} {
        \draw[gray, thin, dashed] (\x, 3) rectangle ++(1, 1);
    }
    \foreach \x in {-1,...,1} {
        \draw[gray, thin, dashed] (\x, 4) rectangle ++(1, 1);
    }


    \foreach \x in {4,...,5} {
        \draw[gray, thin, dashed] (\x, 2) rectangle ++(1, 1);
    }
    \foreach \x in {4,...,5} {
        \draw[gray, thin, dashed] (\x, 3) rectangle ++(1, 1);
    }

    \foreach \x in {3,...,4} {
        \draw[gray] (-1, \x) -- (-1.3, \x);
    }
    \foreach \x in {2,...,4} {
        \draw[gray] (6, \x) -- (6.3, \x);
    }

    \draw [blue, thick] (3, 3.5) -- (3.5, 3.5) -- (3.5, 3);
    \draw [red, thick, densely dashed] (-1.5, 4.5) -- (1.5, 4.5) -- (1.5, 3.5) -- (3, 3.5);
    \draw [red, thick, densely dashed] (-1.5, 3.5) -- (0.5, 3.5) -- (0.5, 3);
    \draw [red, thick, densely dashed] (4.5, 4) -- (4.5, 2.5) -- (6.5, 2.5);
    \draw [red, thick, densely dashed] (5.5, 4) -- (5.5, 3.5) -- (6.5, 3.5);
    \end{tikzpicture} &
    \begin{tikzpicture}[y=-1cm, scale=0.3]
    \foreach \x in {4,...,4} {
        \draw[thick] (\x, 2) rectangle ++(1, 1);
    }
    \draw[black, thick] (4, 3) -- (2, 3) -- (2, 4);


    \foreach \x in {0,...,3} {
        \draw[gray, thin, dashed] (\x, 2) rectangle ++(1, 1);
    }
    \foreach \x in {0,...,1} {
        \draw[gray, thin, dashed] (\x, 3) rectangle ++(1, 1);
    }


    \foreach \x in {4,...,6} {
        \draw[gray, thin, dashed] (\x, 2) rectangle ++(1, 1);
    }

    \foreach \x in {2,...,4} {
        \draw[gray] (0, \x) -- (-0.3, \x);
    }
    \foreach \x in {2,...,3} {
        \draw[gray] (7, \x) -- (7.3, \x);
    }

    \draw [red, thick] (4.5, 3) -- (4.5, 2);
    \draw [red, thick, densely dashed] (-0.5, 3.5) -- (1.5, 3.5) -- (1.5, 2.5) -- (2, 2.5) -- (3.5, 2.5) -- (3.5, 2);
    \draw [red, thick, densely dashed] (-0.5, 2.5) -- (0.5, 2.5) -- (0.5, 2);
    \draw [red, thick, densely dashed] (5.5, 3) -- (5.5, 2.5) -- (7.5, 2.5);
    \end{tikzpicture} &
    \begin{tikzpicture}[y=-1cm, scale=0.3]
    \draw[black, thick] (4, 3) -- (4, 1);

    \foreach \x in {4,...,4} {
        \draw[thick] (\x, 1) rectangle ++(1, 1);
    }
    \foreach \x in {4,...,4} {
        \draw[thick] (\x, 2) rectangle ++(1, 1);
    }
    \draw[black, thick] (4, 3) -- (2, 3) -- (2, 4);

    \foreach \x in {0,...,3} {
        \draw[gray, thin, dashed] (\x, 2) rectangle ++(1, 1);
    }
    \foreach \x in {0,...,1} {
        \draw[gray, thin, dashed] (\x, 3) rectangle ++(1, 1);
    }


    \foreach \x in {4,...,6} {
        \draw[gray, thin, dashed] (\x, 1) rectangle ++(1, 1);
    }
    \foreach \x in {4,...,6} {
        \draw[gray, thin, dashed] (\x, 2) rectangle ++(1, 1);
    }

    \foreach \x in {2,...,4} {
        \draw[gray] (0, \x) -- (-0.3, \x);
    }
    \foreach \x in {1,...,3} {
        \draw[gray] (7, \x) -- (7.3, \x);
    }

    \draw [red, thick] (4.5, 3) -- (4.5, 1.5) -- (5, 1.5);
    \draw [red, thick, densely dashed] (-0.5, 3.5) -- (1.5, 3.5) -- (1.5, 2.5) -- (2, 2.5) -- (3.5, 2.5) -- (3.5, 2);
    \draw [red, thick, densely dashed] (-0.5, 2.5) -- (0.5, 2.5) -- (0.5, 2);
    \draw [red, thick, densely dashed] (5, 1.5) -- (7.5, 1.5);
    \draw [red, thick, densely dashed] (5.5, 3) -- (5.5, 2.5) -- (7.5, 2.5);
    \end{tikzpicture}
\end{array}
\right]
\]
    \caption{An example of a ribbon matrix $(\theta_i\#\theta_j)_{i,j}^3$}
    \label{fig:example-HG-2}
\end{figure}

\begin{example}\label{ex:hg-ex-1}
    We give example of application of Theorem~\ref{thm:HG-sfG}. Let $\lambda = (3,3,2),\ \mu = (1,0,0),\ \mathbf{a} = (a_i) = (1,2,3), \mathbf{b} = (b_i) = (3,3,2)$. Let $\Theta = (\theta_1, \theta_2, \theta_3)$ be the outer decomposition shown in Fig.~\ref{fig:example-HG-1} where $k = 3$. 
    The ribbon matrix of shapes $\theta_i \# \theta_j$ with $i,j \in [3]$ is the following:
    \[
    \left[
    \begin{array}{ccc}
        (2,2,2) / (2,2,0) & (2,2) / (1,-1) & (3,3,3) / (2,2,0) \\[0.3cm] 
        (2,2,2,1) / (2,2,1,1) & (2,2,0) / (1,0,0) & (3,3,3,1) / (2,2,1,1) \\[0.3cm] 
        (2,1,0) / (2,2,0) & (2, -1) / (1, -1) & (3,3,0) / (2,2,0)
    \end{array}
    \right],
    \]
    with diagrams of shapes $\theta_i \# \theta_j$ shown in Fig.~\ref{fig:example-HG-2}.
    The $\varepsilon = \kappa/\kappa$ where $\kappa = (2,2,0)$, the ribbon $\mathsf{G}$-flags are $\mathsf{a}^{\varepsilon} = (-\infty, a_1, a_3), \mathsf{b}^{\varepsilon} = (b_2, b_3, -\infty)$ (see Fig.~\ref{fig:example-HG-1}). All other induced ribbon $\mathsf{G}$-flags are given by 
    {\scriptsize 
    \begin{alignat*}{4}
    \mathbf{a}^{11} &= (-\infty, a_1, a_3) 
    &\qquad \mathbf{b}^{11} &= (b_2, b_3, b_3)
    \\
    \mathbf{a}^{12} &= (a_1, a_3) 
    &\qquad \mathbf{b}^{12} &= (b_2, b_3) 
    \\
    \mathbf{a}^{13} &= (-\infty, a_1, a_3)  
    &\qquad \mathbf{b}^{13} &= (b_1, b_2, b_3)
    \\
    \mathbf{a}^{21} &= (-\infty, a_1, a_2, a_3)  
    &\qquad \mathbf{b}^{21} &= (b_2, b_3, b_3)
    \\ 
    \mathbf{a}^{22} &= (a_1, a_2, a_3)
    &\qquad \mathbf{b}^{22} &= (b_2, b_3)
    \\
    \mathbf{a}^{23} &= (a_1, a_2, a_3)
    &\qquad \mathbf{b}^{23} &= (b_1, b_2, b_3)
    \\ 
    \mathbf{a}^{31} &= (-\infty, a_1, a_3) 
    &\qquad \mathbf{b}^{31} &= (b_2, b_3)
    \\
    \mathbf{a}^{32} &= (a_1, a_3)
    &\qquad \mathbf{b}^{32} &= (b_2)
    \\
    \mathbf{a}^{33} &= (-\infty, a_1, a_3) 
    &\qquad \mathbf{b}^{33} &= (b_1, b_2).
    \end{alignat*}
    }
    
    Applying Theorem~\ref{thm:HG-sfG} we get:
    \begin{align*}
        \mathsf{G}^{\mathbf{a}, \mathbf{b}}_{\lambda/\mu}(\mathbf{x}; \mathbf{U}, \mathbf{V}) = 
        - \mathsf{G}^{\mathbf{a}^{\varepsilon}, \mathbf{b}^{\varepsilon}}_{\varepsilon}(\mathbf{x}; \mathbf{U}, \mathbf{V})^{-2} \cdot 
        \det \left[(-1)^{\mathsf{hp}^{\pi}_{\theta_i \# \theta_j}}\ \mathsf{G}^{\mathbf{a}^{ij}, \mathbf{b}^{ij}}_{\theta_i \# \theta_j}(\mathbf{x}; \mathbf{U}, \mathbf{V})\right]_{1 \le i,j \le 3}.
    \end{align*}
    Note that $(-1)^{\mathsf{hp}^{\pi}_{\theta_i \# \theta_j}} = -1$ if $i = 3, j = 1$, and $(-1)^{\mathsf{hp}^{\pi}_{\theta_i \# \theta_j}} = 1$ otherwise, as $\theta_3 \# \theta_1$ is the only negative shape among $(\theta_i \# \theta_j)$ with an odd numbers of hook-pipes. 
    
    In particular, applying Theorem~\ref{thm:HG-sfG} and applying the $G$-specialization we also get:
    \begin{align*}
        G_{\lambda/\mu}(\mathbf{x}_m; \boldsymbol{\alpha}, \boldsymbol{\beta}) = C \cdot \left[\mathsf{G}^{\mathbf{a}^{\varepsilon}, \mathbf{b}^{\varepsilon}}_{\varepsilon}(\mathbf{x}_m; \boldsymbol{\alpha}, \boldsymbol{\beta})\right]^{-2} \cdot \det \left[(-1)^{\mathsf{hp}^{\pi}_{\theta_i \# \theta_j}} \cdot \mathsf{G}^{\mathbf{a}^{ij}, \mathbf{b}^{ij}}_{\theta_i \# \theta_j}(\mathbf{x}_m; \boldsymbol{\alpha}, \boldsymbol{\beta})\right]_{1 \le i,j \le 3},
    \end{align*}
    where $C = \prod_{i = 1}^n \prod_{j = 1}^m (1 - \beta_i x_j)$. 
\end{example}

\section{Talaska's formula}\label{sec:tal}
To prove the main formulas stated in previous section, we will use a generalization of the Lindstr\"om--Gessel--Viennot lemma due to Talaska \cite{tal}, which allows cycles in directed graphs. 
To describe this formula we need several definitions. 

\begin{definition}
(Underlying graph) Let $G$ be a directed weighted graph. 
Each directed edge $e$ of $G$ has weight $w_e$ which is a formal variable of a commutative ring. 


(Walks matrix) 
For a walk $P$ in $G$, its weight is defined as $\mathrm{wt}(P) := \prod_{e \in P} w_e$ the product of edges (we also have {trivial walks} from $v$ to $v$ with no edges and weight $1$).
Define the {\it weighted walk matrix} $W = [w_{a,b}]_{a,b \in V}$, where 
$$
    w_{a,b} := \sum_{\text{walks } P \text{ from $a$ to $b$}} \mathrm{wt}(P),
$$
which is a formal power series in variables $\{ w_e\}$.

(Cycle systems) Let $\mathcal{C}$ be the set of all systems $\mathbf{C} = \{C_1,\ldots, C_k \}$ of non-intersecting (i.e. vertex-disjoint) simple cycles $C_1,\ldots, C_k$ in $G$. The weight of cycle system is defined as $\mathrm{wt}(\mathbf{C}) := \mathrm{wt}(C_1) \cdots \mathrm{wt}(C_k)$ (where the weight of each cycle is the product of edge weights). We also define an empty system of weight $1$ to be the element of $\mathcal{C}$. For nonempty cycle system, its cycles must be nonempty.  

(Path systems) Let ${A} = \{a_1, \dots, a_m\},\ {B} = \{b_1, \ldots, b_m\}$ be subsets of vertices of $G$. Let $\mathcal{{P}}({A}, {B}, \sigma)$ be the set of nonintersecting (i.e. vertex-disjoint) path systems $\mathbf{P} = (P_1, \ldots, P_m)$, where $P_i$ is a path (i.e. a self-avoiding walk) from $A_i$ to $B_{\sigma(i)}$ and $\sigma \in S_m$ is a permutation. The weight of path system is defined as $\mathrm{wt}(\mathbf{P}) := \mathrm{wt}(P_1) \cdots \mathrm{wt}(P_m)$. Let $\mathcal{{P}}({A}, {B}) = \bigcup\limits_{\sigma \in S_m} \mathcal{{P}}({A}, {B}, \sigma)$. 

(Flows) Let $\mathcal{F}({A}, {B})$ be the set of {\it flows}, where each flow $\mathbf{F} \in \mathcal{F}({A}, {B})$ is a pair $(\mathbf{P},  \mathbf{C})$ of a path system $\mathbf{P}  \in \mathcal{P}({A}, {B})$ and a cycle system $\mathbf{C} \in \mathcal{C}$ such that 
$\mathbf{P}$ and $\mathbf{C}$ are vertex-disjoint.
The weight of a flow is defined as 
$\mathrm{wt}(\mathbf{F}) := \mathrm{wt}(\mathbf{P}) \cdot \mathrm{wt}(\mathbf{C})$.

(Signs) For a path system $\mathbf{P} \in \mathcal{{P}}({A}, {B}, \sigma)$, the sign is defined by  $\mathrm{sgn}(\mathbf{P}) := \mathrm{sgn}(\sigma)$. The sign of cycle system is defined as $\mathrm{sgn}(\mathbf{C}) := (-1)^{\left| \mathbf{C} \right|}$, where $\left| \mathbf{C} \right|$ is the number of cycles in $\mathbf{C}$. For a flow $\mathbf{F} = (\mathbf{P}, \mathbf{C}) \in \mathcal{F}(A, B)$, the sign is defined as $\mathrm{sgn}(\mathbf{F}) := \mathrm{sgn}(\mathbf{P}) \cdot \mathrm{sgn}(\mathbf{C})$.
\end{definition}

\begin{theorem}[\cite{tal}]
    \label{thm:Talaska}
    Let $W_{A, B}$ be the submatrix of $W$ indexed by rows $A$ and columns $B$. 
    Then the following formula holds: 
    \begin{align*}
        \det  W_{A, B}  = \frac{\sum\limits_{\mathbf{F} \in \mathcal{F}({A}, {B})} \mathrm{sgn(\mathbf{F}) \cdot wt(\mathbf{F})}}{\sum\limits_{\mathbf{C} \in \mathcal{C}} \mathrm{sgn(\mathbf{C}) \cdot wt(\mathbf{C})}}.
    \end{align*}
\end{theorem}

More specifically, we will use the following corollary of this formula.  
\begin{corollary}\label{cor:tal}
Let 
\begin{align*}
    &F_{i,j} := \sum\limits_{\mathbf{F} \in \mathcal{F}(\{a_i\}, \{b_j\})} \mathrm{sgn(\mathbf{F}) \cdot wt(\mathbf{F})}. 
\end{align*}
Then we have: 
    \begin{align*}
        \sum\limits_{\mathbf{F} \in \mathcal{F}({A}, {B})} \mathrm{sgn(\mathbf{F}) \cdot {wt}(\mathbf{F})} = \frac{\det\left[{F}_{i, j}\right]_{1 \le i,j \le m}}{\left(\sum\limits_{\mathbf{C} \in \mathcal{C}} \mathrm{sgn(\mathbf{C}) \cdot wt(\mathbf{C})}\right)^{m - 1}}.
    \end{align*}
\end{corollary}
\begin{proof}
    Since we have 
    $$
    w_{a_i, b_j} = \frac{F_{i, j}}{\sum\limits_{\mathbf{C} \in \mathcal{C}} \mathrm{sgn(\mathbf{C}) \cdot wt(\mathbf{C})}},
    $$
    then from Theorem~\ref{thm:Talaska} we obtain that 
    \begin{align*}
        \frac{\sum\limits_{\mathbf{F} \in \mathcal{F}(\mathbf{A}, \mathbf{B})} \mathrm{sgn(\mathbf{F}) \cdot wt(\mathbf{F})}}{\sum\limits_{\mathbf{C} \in \mathcal{C}} \mathrm{sgn(\mathbf{C}) \cdot wt(\mathbf{C})}} = \det W_{A,B} =  \frac{1}{\left(\sum\limits_{\mathbf{C} \in \mathcal{C}} \mathrm{sgn(\mathbf{C}) \cdot wt(\mathbf{C})}\right)^m} \det\left[ F_{i, j} \right]_{1 \le i, j \le m}, 
    \end{align*}
    from which the result follows.
\end{proof}

\section{Proof}\label{sec:proof}


In this section we give a proof of the main result. 

\subsection{Summary of the proof.} \label{sec:proofplan}
\subsubsection{T-lattice}
In \textsection\,\ref{subsec:t-lattice} we introduce the {\it T-lattice} denoted by $TL_{\lambda/\mu}(\pi, m, \mathbf{a}, \mathbf{b})$ which is a lattice graph defined by the shape $\lambda/\mu$, pipe function $\pi$ (which induces an outer ribbon decomposition $\Theta = (\theta_1, \ldots, \theta_k)$), and row $\mathsf{G}$-flags $\mathbf{a}, \mathbf{b}$. We place endpoints $\mathbf{A} = (A_1, \ldots, A_k), \mathbf{B} = (B_1, \ldots, B_k)$ on the T-lattice and consider flows $\mathcal{F}(\mathbf{A}, \mathbf{B})$. It is worth noting that the weights of some edges (more precisely, the edges of content $c$ such that $\pi(c)$ is a hook pipe) are negative. 

\subsubsection{Bijection}
In \textsection\,\ref{subsec:bijection} we define a map $\phi : \mathcal{F}(\mathbf{A}, \mathbf{B}) \rightarrow TT_{\lambda/\mu}(\mathbf{a}, \mathbf{b}, m)$ between the flows on the T-lattice and trio tableaux with flags $\mathbf{a}, \mathbf{b}$ and prove that it is a bijection by constructing the inverse map $\psi = \phi^{-1}$ of $\phi$. 
The classical technique for such type of problems is to establish a bijection between a single path $P$ from $A_i$ to $B_j$ and a single filling of some ribbon $\theta_i \# \theta_j, i,j \in [k]$, which then extends the bijection to the whole shape and system of paths by filling each ribbon $\theta_i$ with single paths from $A_i$ to $B_i, i \in [k]$, and collecting all the ribbons $\theta_i$ to form the shape. However this technique fails here for multiple reasons such as: (1) a single path $P$ does not fill the ribbon from an outer decomposition $\Theta$; instead it fills {\it multiple} ribbons from outer decompositions of shapes $\lambda/\tau,\ \nu/\mu,\ \tau/\nu$ where $\nu,\tau \in \mathsf{Sign}_n$ and $\nu \subseteq \mu,\ \lambda \subseteq \tau$; (2) the flows consist of not only paths from $A_i$ to $B_j$ but also cycles that fill some ribbons, contributing to the {\it sign} of the flow. 
Thus, we define the map $\phi$ differently: we fix some content $c$ and fill the cells of $L(\lambda/\mu), R(\lambda/\mu)$ of content $c$ with edges of content $c$ from top, bottom parts of the flow $\mathbf{F}$; then we fill the cells in $\tau/\nu$ (which is well defined, since the cells in $L(\lambda/\mu), R(\lambda/\mu)$ have already been filled) in the same way. 
The map $\phi$ constructs the whole trio tableau from the flow straightaway. We prove that the map $\phi$ indeed constructs some trio tableau (Lemma~\ref{lem:flow-is-trio-tableau}) using Lemmas~\ref{lem:two-cells-two-edges}, \ref{lem:left-tableau-is-of-some-shape}, \ref{lem:right-tableau-is-of-some-shape}, \ref{lem:central-tableau-is-of-some-shape}. 
\subsubsection{Inverse map}
Then we define the map $\psi$ that constructs flow from trio tableaux as follows: we fix $\nu, \tau$ and using them together with $\lambda, \mu$ define {\it minor endpoints}; after that we construct non-intersecting paths between minor endpoints, and using Lemmas~\ref{lem:minor-is-either-major-or-another-minor}, \ref{lem:minor-endpoint-at-most-two-paths}, \ref{lem:systems-of-paths-form-flow} prove that the constructed paths combined give some flow and the minor endpoints become either the endpoints $\mathbf{A}, \mathbf{B}$ or some intermediate nodes. Then we prove that $\phi^{-1} = \psi$ in Lemma~\ref{lem:psi-is-inverse-of-phi}.

\subsubsection{Weight-preserving properties} 
In \textsection~\ref{subsec:bijection-respects-weights} we prove that the defined bijection is weight-preserving up to some sign, where sign is determined by the shape and pipe function $\pi$. Surpisingly, this part is the most challenging; mainly because flows may consist of paths 
that start at $A_i$ and terminate at $B_{\sigma(i)}$ where $\sigma$ is a {\it non-identity} permutation, whose sign has to be accounted as well. In \textsection~\ref{subsec:prelude} we prove auxiliary Lemmas~\ref{lem:cycle-corresponds-to-special-ribbons} (using Lemma~\ref{lem:cycle-contains-only-two-types-of-edges}), \ref{lem:restrictions-for-ribbons}, \ref{lem:ribbons-filled-by-same-path-in-same-copy}, \ref{lem:th-i-th-si-are-in-the-same-ribbon} (using  Lemma~\ref{lem:ribbons-filled-by-same-path-in-same-copy}), \ref{lem:cousin-paths-avoid-each-others-endpoints} (using Lemma~\ref{lem:ribbons-filled-by-same-path-in-same-copy}) which we use in further sections. Especially Lemma~\ref{lem:cousin-paths-avoid-each-others-endpoints} is heavily used in arguments by contradiction. 

In \textsection~\ref{subsec:weigth-preserving-for-ribbons} we prove weight-preserving property for single path $P \in \mathcal{P}(A_i, B_{\sigma(i)})$ and the ribbons filled by $P$ in Lemmas~\ref{lem:weight-of-ribbons-filled-by-same-path-positive} (for positive $\theta_i \# \theta_j$), \ref{lem:weight-of-ribbons-filled-by-same-path-negative} (for negative $\theta_i \# \theta_j$). 

In \textsection~\ref{subsec:weigth-preserving-for-positive-shapes} we prove weight-preserving property for positive shapes. First 
we show that the weight of paths $(P_1, \ldots, P_{\ell})$ (where $P_i \in \mathcal{P}(A_{a_i}, B_{\sigma(a_i)})$ and $\sigma_0 = (a_1, \ldots, a_{\ell})$ is a cycle contained in $\sigma$) coincides with the weight of all ribbons filled by the paths $P_i (i \in [\ell])$ up to some determined sign (Lemma~\ref{lem:sgn-rule-for-one-cycle-positive}) using Lemma~\ref{lem:intersection-of-pos-ribbon-not-exist}. Then we derive the result for general positive shape in Lemma~\ref{lem:weight-preserving-positive}. 

In \textsection~\ref{subsec:weigth-preserving-for-negative-shapes} we handle the negative shapes  (Lemmas~\ref{lem:sgn-rule-for-one-cycle}, \ref{lem:weight-preserving-negative}) using auxiliary Lemmas~\ref{lem:intersection-of-neg-ribbon-is-one-special-ribbon}, \ref{lem:ribbon-between-is-always-odd}. It is worth noting that the negative case is more challenging than the positive one; the reason is that paths from $A_i$ to $B_j$ is more `tricky' if $B_j$ is placed to the left of $A_i$ (i.e. if $\theta_i \# \theta_j$ is negative) rather than if $B_j$ is placed to the right of $A_i$ (i.e. if $\theta_i \# \theta_j$ is positive), to see this compare Lemmas~\ref{lem:weight-of-ribbons-filled-by-same-path-positive} and \ref{lem:weight-of-ribbons-filled-by-same-path-negative}. Finally, we prove that the defined map is weight-preserving in Lemma~\ref{lem:main-weight-preserving-bijection}. 

\subsubsection{Completing the proof}
In \textsection~\ref{subsec:hg-formulas-proof} we give the proof of Hamel--Goulden formulas (Theorem~\ref{thm:HG-sfG}) using the weight-preserving bijection.

\subsection{T-lattice}
\label{subsec:t-lattice}

We construct a lattice graph using the given skew shape $\lambda/\mu$, outer decomposition $\Theta = (\theta_1, \ldots, \theta_k)$ induced by a corresponding pipe function $\pi$ and $\mathsf{G}$-flags $\mathbf{a}, \mathbf{b}$, which are used throughout the section. The lattice graph is fully defined by presenting all of its directed edges. For a content $c \in \mathbb{Z}$, define the following value 
\[
f_c := 
\begin{cases}
    0,\quad &\text{if}\quad c = 0, \\
    |\{i : 0 < i \le c,\ \pi(i) = (1, \cdot)\}|, \quad &\text{if}\quad c > 0, \\
    -|\{i : c\le i < 0,\ \pi(i) = (\cdot, 1)\}|, \quad &\text{if}\quad c < 0.
\end{cases}
\]



Let $c \in \mathbb{Z}$ be some content 
and $m_c = \min(\lambda/\mu, c), M_c = \max(\lambda/\mu, c)$. Define 
\begin{align*}
	\overline{a}_i := a_{\mathrm{row}(M_i)}, \qquad \overline{b}_i := b_{\mathrm{row}(m_i)}, 
\end{align*}
and set $\overline{a}_i = -\infty$ if $M_i$ is undefined, and $\overline{b}_i = -\infty$ if $m_i$ is undefined.

There are four types of edges in the graph: {\it right, left, up, down} (they are defined below); we call right and left edges as {\it horizontal}; up and down edges as {\it vertical}.

\begin{definition}[T-lattice]
    Let $(r', s'),\ (r, s)$ be top-right and bottom-left cells (resp.) in $\lambda/\mu$. 
    Define positive integers $m, n^{\circ}, n^{\bullet}$ such that $n^{\bullet} > m + f_{c(\zeta)},\ n^{\circ} < 1 + f_{c(\eta)}$ where $\zeta = (r', s' + 1), \eta = (r, s - 1)$.
    
    

    We define the \textit{T-lattice} $TL_{\lambda/\mu}(\pi, m, \mathbf{a}, \mathbf{b})$ as a directed graph with nodes from $\mathbb{Z}^2$ and the set $E$ of directed edges. Denote $i^{\circ} := n^{\circ} - i$ where $i \in \mathbb{Z}$ and 
    $i < 1 + f_{c(\eta)}$ so that $i^{\circ} \ge n^{\circ} - f_{c(\eta)}$; 
    similarly, let $i^{\bullet} := i - n^{\bullet}$ where $i \in \mathbb{Z}$ and 
    $i > m + f_{c(\zeta)}$ so that $i^{\bullet} > m + f_{c(\zeta)} - n^{\bullet}$. Let $c_M, c_m$ be the maximal, minimal contents in $L(\lambda/\mu), R(\lambda/\mu)$, respectively.

    Now we describe all edges in $E$ togehter with their types:

    \begin{enumerate}
    \item $(i,j)\rightarrow (i+1,j),\ i \in \mathbb{Z},\ j \in [1 + f_i, m + f_i]$, right 
    \item $(i + 1,\ j^{\circ}) \rightarrow (i,\ j^{\circ}),\ i \le c_M,\ j \in [n^{\circ} - f_{c(\eta)},\  \overline{a}_{i} + f_i]$, left 
    \item $(i + 1,\ j^{\bullet}) \rightarrow (i,\ j^{\bullet}),\ i \ge c_m,\ j \in [m + f_{c(\zeta)} + 1 - n^{\bullet},\  \overline{b}_{i} - f_i]$, left 
    \item $(i, j) \rightarrow (i, j + 1),\ i, j \in \mathbb{Z}$  
    if $\pi(i) = (0, \cdot)$, up 
    \item $(i, j) \rightarrow (i, j - 1),\ i, j \in \mathbb{Z}$ 
    if $\pi(i) = (1, \cdot)$, down 
    \end{enumerate}
    We let $n^{\bullet} \rightarrow \infty,\ n^{\circ} \rightarrow -\infty$. 
    
    Let $e = (c + 1, x) \rightarrow (c, x)$ or $e = (c, x) \rightarrow (c + 1, x)$ then we define the content of the edge $e$ as $c(e) := c$.

    Define the weight $wt(e)$ for horizontal edges $e$ as follows:
    \begin{align*}
        wt(e) := 
        \begin{cases}
            x_a, & e = (c, a + f_c) \rightarrow (c + 1, a + f_c) \\
            u_{a, a - c}, & e = (c + 1, (a - f_c)^{\bullet}) \rightarrow (c, (a - f_c)^{\bullet}),\ \pi(c) \not\in \mathcal{P}_H \\
            -u_{a, a - c}, & e = (c + 1, (a - f_c)^{\bullet}) \rightarrow (c, (a - f_c)^{\bullet}),\ \pi(c) \in \mathcal{P}_H \\
            v_{a,a + c}, & e = (c + 1, (a - f_c)^{\circ}) \rightarrow (c, (a - f_c)^{\circ}),\ \pi(c) \not\in \mathcal{P}_H \\
            -v_{a,a + c}, & e = (c + 1, (a - f_c)^{\circ}) \rightarrow (c, (a - f_c)^{\circ}),\ \pi(c) \in \mathcal{P}_H \\
        \end{cases}
    \end{align*}
\end{definition}

Out lattice consists of three layers: {\it top, middle, bottom}. Top and bottom layers consist of nodes with coordinates $(i, j^{\bullet}),\ (i, j^{\circ})$ (resp.) for $j \in \mathbb{Z}$; the middle layer contains all other nodes. We say that the edge $e = (x, y) \rightarrow (x', y')$ belongs to bottom (resp. top) layer if both nodes $(x, y), (x', y')$ belong to bottom (resp. top) layer; and $e$ belongs to middle layer if it belongs to neither top nor bottom layer. 

To consider paths on the lattice above we define positions of starting points 
$\mathbf{A}=(A_1,...,A_k)$ and ending points $\mathbf{B}=(B_1,...,B_k)$, where $k$
is the number of ribbons in the given outer outer decomposition $\Theta=(\theta_1,...,\theta_k)$.
Let $\delta_i$ be the head and $\gamma_i$ be the tail of the ribbon $\theta_i$. 
The endpoints $A_i$ and $B_i$ are defined by the ribbon $\theta_i$ as follows:

(i) If $\lambda/\mu$ is positive shape, then 
\begin{align*}
    &A_i = 
\begin{cases}
	(c(\delta_i),\ \infty^{\circ}) &\text{ if } \pi(\delta_i) = (0, \cdot), \\
	(c(\delta_i),\ \infty^{\bullet}) &\text{ if } \pi(\delta_i) = (1, \cdot).
\end{cases}\\
    &B_i = 
\begin{cases}
    (c(\gamma_i) + 1,\ \infty^{\bullet}) &\text{ if } \pi(\gamma_i) = (\cdot, 0), \\
	(c(\gamma_i) + 1,\ \infty^{\circ}) &\text{ if } \pi(\gamma_i) = (\cdot, 1).
\end{cases}
\end{align*}

(ii) If $\lambda/\mu$ is negative shape, then  
\begin{align*}
    &A_i = 
\begin{cases}
	(c(\delta_i) + 1,\ \infty^{\circ}) &\text{ if } \pi(\delta_i) = (\cdot, 0) \\
	(c(\delta_i) + 1,\ \infty^{\bullet}) &\text{ if } \pi(\delta_i) = (\cdot, 1)
\end{cases}\\
    &B_i = 
\begin{cases}
    (c(\gamma_i),\ \infty^{\bullet}) &\text{ if } \pi(\gamma_i) = (0, \cdot) \\
	(c(\gamma_i),\ \infty^{\circ}) &\text{ if } \pi(\gamma_i) = (1, \cdot)
\end{cases}
\end{align*}

Let $P$ be a path on the lattice. Let $P' = (e_1, \ldots, e_s)$ be a sub-path of $P$ such that 

(1) each $e_i, i \in [s]$ belongs to bottom (middle, top resp.) layer; 

(2) if there is edge right before $e_1$ in $P$ then it does not belong to bottom (middle, top resp.) layer; 

(3) if there is edge right after $e_s$ in $P$ then it does not belong to bottom (middle, top resp.) layer.

Then $P'$ is called a {\it minor sub-path} of $P$ (or simply {\it minor path} if $P$ is not specified), and the endpoints of $P'$ are called {\it minor endpoints}. Let $(P'_1, \ldots, P'_s)$ be all minor sub-paths of paths in the flow $\mathbf{F}$ that belong to bottom (middle, top resp.) layer, then we call the set paths $(P'_i)$ a {\it bottom (middle, top resp.) part} of the flow $\mathbf{F}$. 

\begin{figure}
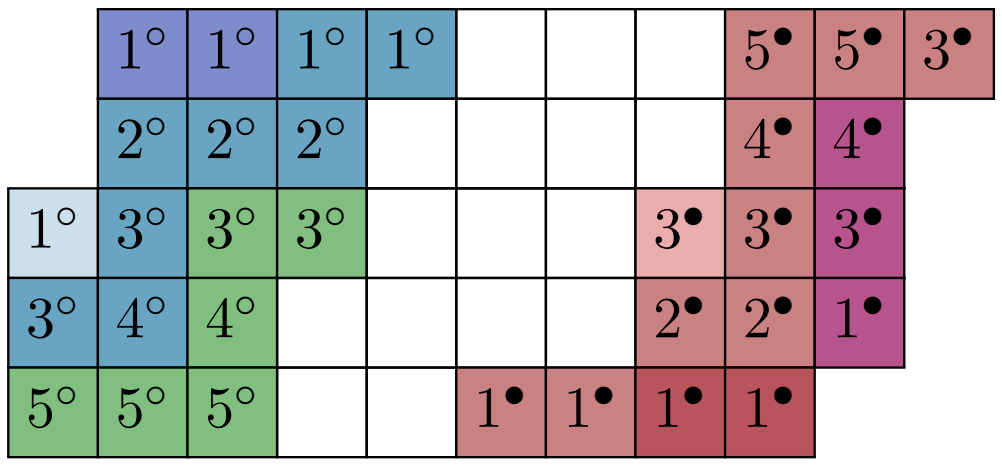
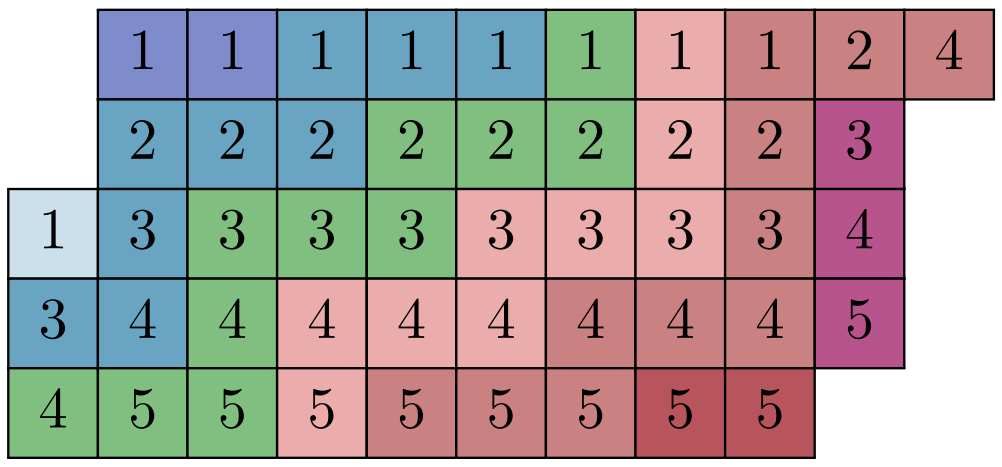

    \centering
    \begin{subfigure}[b]{1\textwidth}
    \centering
    \begin{ytableau}
        \none & *(red 7) 1^{\circ} & *(red 7) 1^{\circ} & *(red 6) 1^{\circ} & *(red 6) 1^{\circ} & & & & *(red 4) 5^{\bullet} &  *(red 4) 5^{\bullet} &  *(red 4) 3^{\bullet} \\
        \none & *(red 6) 2^{\circ} & *(red 6) 2^{\circ} & *(red 6) 2^{\circ} & & & & & *(red 4) 4^{\bullet} &  *(red 1) 4^{\bullet}\\
        *(red 2) 1^{\circ} & *(red 6) 3^{\circ} & *(red 8) 3^{\circ} & *(red 8) 3^{\circ} & & & & *(red 3) 3^{\bullet} & *(red 4) 3^{\bullet} &  *(red 1) 3^{\bullet} \\
        *(red 6) 3^{\circ} & *(red 6) 4^{\circ} & *(red 8) 4^{\circ} & & & & & *(red 4) 2^{\bullet} & *(red 4) 2^{\bullet} &  *(red 1) 1^{\bullet}\\
        *(red 8) 5^{\circ} & *(red 8) 5^{\circ} & *(red 8) 5^{\circ} & & & *(red 4) 1^{\bullet} & *(red 4) 1^{\bullet} & *(red 5) 1^{\bullet} & *(red 5) 1^{\bullet} 
    \end{ytableau}
    \subcaption{$T^{\circ}$ and $T^{\bullet}$}
    \label{fig:trio-tableau-example-b}
    \end{subfigure}  
    \begin{subfigure}[b]{1\textwidth}
    \centering
    \begin{ytableau}
        \none & *(red 7) 1 & *(red 7) 1 & *(red 6) 1 & *(red 6) 1 & *(red 6) 1 & *(red 8) 1 & *(red 3) 1 & *(red 4) 1 & *(red 4) 2 & *(red 4) 4 \\
        \none & *(red 6) 2 & *(red 6) 2 & *(red 6) 2 & *(red 8) 2 & *(red 8) 2 & *(red 8) 2 & *(red 3) 2 & *(red 4) 2 & *(red 1) 3  \\
        *(red 2) 1 & *(red 6) 3 & *(red 8) 3 & *(red 8) 3 & *(red 8) 3 & *(red 3) 3 & *(red 3) 3 & *(red 3) 3 & *(red 4) 3 & *(red 1) 4 \\
        *(red 6) 3 & *(red 6) 4 & *(red 8) 4 & *(red 3) 4 & *(red 3) 4 & *(red 3) 4 & *(red 4) 4 & *(red 4) 4 & *(red 4) 4 & *(red 1) 5 \\
        *(red 8) 4 & *(red 8) 5 & *(red 8) 5 & *(red 3) 5 & *(red 4) 5 & *(red 4) 5 & *(red 4) 5 & *(red 5) 5 & *(red 5) 5
    \end{ytableau}
    \subcaption{$\widetilde{T}$}
    \label{fig:trio-tableau-example-a}
\end{subfigure}
\caption{An example of trio tableau.}
\end{figure}

\begin{figure}[!h]
    \centering
\begin{tikzpicture}[scale=0.6, transform shape]
        \foreach \x in {-8,...,8} {
        \foreach \y in {-8,...,15} {
            \draw[dotted, thin] (\x,\y) rectangle ++(1,1);
        }
        }
        \foreach \x in {-2, ..., 8} {
            \fill[black] (-8.5,\x) circle (0pt) node[left] {$\x$};
        }
        \fill[blue, opacity=0.05] (-4,-8) rectangle (-2, 16);
        \fill[blue, opacity=0.05] (0,-8) rectangle (1, 16);
        \fill[blue, opacity=0.05] (3, -8) rectangle (6, 16);

        \draw[line width=0.2mm, black, dashed] (-8, -8) -- (-1, -8) -- (-1, -7) -- (2, -7) -- (2, 1) -- (3, 1) -- (3, 2) -- (4, 2) -- (4, 3) -- (5, 3) -- (5, 4) -- (9, 4);
        \draw[line width=0.2mm, black, dashed] (-8, 2) -- (-4, 2) -- (-4, 3) -- (-3, 3) -- (-3, 4) -- (-2, 4) -- (-2, 15) -- (0, 15) -- (0, 14) -- (1, 14) -- (1, 16) -- (3, 16) -- (3, 15) -- (5, 15) -- (5, 14) -- (9, 14);

        \foreach \x in {-8,...,9} {
            \fill[black] (\x, -9) circle (0pt) node[above] {$\x$};
        }

        \fill[black] (-9, 9.6) circle (0pt) node {$\vdots$};
        \fill[black] (-9, -3.4) circle (0pt) node {$\vdots$};

        \fill[black] (-4, 16) circle (2pt) node[below right] {$A_1$};
        \fill[black] (-3, 16) circle (2pt) node[below right] {$A_2$};
        \fill[black] (-1, -8) circle (2pt) node[above right] {$A_3$};
        \fill[black] (0, 16) circle (2pt) node[below right] {$A_4$};
        \fill[black] (3, 16) circle (2pt) node[below right] {$A_5$};
        \fill[black] (2, -8) circle (2pt) node[below right] {$A_6$};
        \fill[black] (2, 16) circle (2pt) node[below right] {$B_1$};
        \fill[black] (-2, 16) circle (2pt) node[below right] {$B_2$};
        \fill[black] (4, -8) circle (2pt) node[below right] {$B_3$};
        \fill[black] (1, 16) circle (2pt) node[below right] {$B_4$};
        \fill[black] (5, -8) circle (2pt) node[below right] {$B_5$};
        \fill[black] (3, -8) circle (2pt) node[below right] {$B_6$};

        \foreach \x in {11, ..., 16} {
            \pgfmathsetmacro{\y}{int(\x - 12)}
            \fill[black] (-8.5,\x) circle (0pt) node[left] {$\y^{\bullet}$};
        }
        \foreach \x in {-8, ..., -5} {
            \pgfmathsetmacro{\y}{int(-\x - 6)}
            \fill[black] (-8.5,\x) circle (0pt) node[left] {$\y^{\circ}$};
        }

        \draw[red 3, very thick, -{Stealth}] (-4, 16) -- (-4, 3);
        \draw[red 3, very thick, -{Stealth}] (-4, 3) -- (5, 3);
        \draw[red 3, very thick, -{Stealth}] (5, 3) -- (5, -8);
        \fill[black] (-3.5, 3) node[above] {$5$};
        \fill[black] (-2.5, 3) node[above] {$4$};
        \fill[black] (-1.5, 3) node[above] {$4$};
        \fill[black] (-0.5, 3) node[above] {$4$};
        \fill[black] (0.5, 3) node[above] {$3$};
        \fill[black] (1.5, 3) node[above] {$3$};
        \fill[black] (2.5, 3) node[above] {$3$};
        \fill[black] (3.5, 3) node[above] {$2$};
        \fill[black] (4.5, 3) node[above] {$1$};

        \draw[red 4!80!black, very thick, -{Stealth}] (-3, 16) -- (-3, 4);
        \draw[red 4!80!black, very thick, -{Stealth}] (-3, 4) -- (6, 4);
        \draw[red 4!80!black, very thick, -{Stealth}] (6, 4) -- (6, 5) -- (7, 5) -- (7, 7) -- (8, 7) -- (8, 12) -- (7, 12) -- (7, 14);
        \draw[red 4!80!black, very thick, -{Stealth}] (7, 14) -- (1, 14) -- (1, 16);

        \fill[black] (-2.5, 4) node[above] {$5$};
        \fill[black] (-1.5, 4) node[above] {$5$};
        \fill[black] (-0.5, 4) node[above] {$5$};
        \fill[black] (0.5, 4) node[above] {$4$};
        \fill[black] (1.5, 4) node[above] {$4$};
        \fill[black] (2.5, 4) node[above] {$4$};
        \fill[black] (3.5, 4) node[above] {$3$};
        \fill[black] (4.5, 4) node[above] {$2$};
        \fill[black] (5.5, 4) node[above] {$1$};
        \fill[black] (6.5, 5) node[above] {$2$};
        \fill[black] (7.5, 7) node[above] {$4$};
        
        \fill[black] (7.5, 12) node[above] {$3^{\bullet}$};
        \fill[black] (6.5, 14) node[above] {$5^{\bullet}$};
        \fill[black] (5.5, 14) node[above] {$5^{\bullet}$};
        \fill[black] (4.5, 14) node[above] {$4^{\bullet}$};
        \fill[black] (3.5, 14) node[above] {$3^{\bullet}$};
        \fill[black] (2.5, 14) node[above] {$2^{\bullet}$};
        \fill[black] (1.5, 14) node[above] {$2^{\bullet}$};

        \draw[red 4!80!black, very thick, -{Stealth}] (0, 16) -- (0, 14) -- (-2, 14) -- (-2, 16);

        \fill[black] (-0.5, 14) node[above] {$1^{\bullet}$};
        \fill[black] (-1.5, 14) node[above] {$1^{\bullet}$};

        \draw[red 8, very thick, -{Stealth}] (-1, -8) -- (-7, -8);
        \draw[red 8, very thick, -{Stealth}] (-7, -8) -- (-7, 1);
        \draw[red 8, very thick, -{Stealth}] (-7, 1) -- (-6, 1) -- (-6, 2) -- (4, 2);
        \draw[red 8, very thick, -{Stealth}] (4, 2) -- (4, -8);
        \fill[black] (-1.5, -8) node[above] {$3^{\circ}$};
        \fill[black] (-2.5, -8) node[above] {$3^{\circ}$};
        \fill[black] (-3.5, -8) node[above] {$4^{\circ}$};
        \fill[black] (-4.5, -8) node[above] {$5^{\circ}$};
        \fill[black] (-5.5, -8) node[above] {$5^{\circ}$};
        \fill[black] (-6.5, -8) node[above] {$5^{\circ}$};

        \fill[black] (-6.5, 1) node[above] {$4$};
        \fill[black] (-5.5, 2) node[above] {$5$};
        \fill[black] (-4.5, 2) node[above] {$5$};
        \fill[black] (-3.5, 2) node[above] {$4$};
        \fill[black] (-2.5, 2) node[above] {$3$};
        \fill[black] (-1.5, 2) node[above] {$3$};
        \fill[black] (-0.5, 2) node[above] {$3$};
        \fill[black] (0.5, 2) node[above] {$2$};
        \fill[black] (1.5, 2) node[above] {$2$};
        \fill[black] (2.5, 2) node[above] {$2$};
        \fill[black] (3.5, 2) node[above] {$1$};

        \draw[red 6, very thick, -{Stealth}] (2, -8) -- (2, -7) -- (-5, -7) -- (-5, -6) -- (-6, -6);
        \draw[red 6, very thick, -{Stealth}] (-6, -6) -- (-6, 0);
        \draw[red 6, very thick, -{Stealth}] (-6, 0) -- (-5, 0) -- (-5, 1) -- (3, 1);
        \draw[red 6, very thick, -{Stealth}] (3, 1) -- (3, -8);
        \fill[black] (1.5, -7) node[above] {$1^{\circ}$};
        \fill[black] (0.5, -7) node[above] {$1^{\circ}$};
        \fill[black] (-0.5, -7) node[above] {$2^{\circ}$};
        \fill[black] (-1.5, -7) node[above] {$2^{\circ}$};
        \fill[black] (-2.5, -7) node[above] {$2^{\circ}$};
        \fill[black] (-3.5, -7) node[above] {$3^{\circ}$};
        \fill[black] (-4.5, -7) node[above] {$4^{\circ}$};
        \fill[black] (-5.5, -6) node[above] {$3^{\circ}$};
        
        \fill[black] (-5.5, 0) node[above] {$3$};
        \fill[black] (-4.5, 1) node[above] {$4$};
        \fill[black] (-3.5, 1) node[above] {$3$};
        \fill[black] (-2.5, 1) node[above] {$2$};
        \fill[black] (-1.5, 1) node[above] {$2$};
        \fill[black] (-0.5, 1) node[above] {$2$};
        \fill[black] (0.5, 1) node[above] {$1$};
        \fill[black] (1.5, 1) node[above] {$1$};
        \fill[black] (2.5, 1) node[above] {$1$};

        \draw[red 3, very thick, -{Stealth}] (3, 16) -- (3, 15) -- (2, 15) -- (2, 16);
        \fill[black] (2.5, 15) node[above] {$3^{\bullet}$};

        \draw[red 7, thick, -{Stealth}] (-2, -6) -- (-2, 0) -- (0, 0) -- (0, -6) -- (-2, -6);
        \fill[black] (-0.5, -6) node[above] {$1^{\circ}$};
        \fill[black] (-1.5, -6) node[above] {$1^{\circ}$};
        
        \fill[black] (-1.5, 0) node[above] {$1$};
        \fill[black] (-0.5, 0) node[above] {$1$};

        \draw[red 2!70!black, thick, -{Stealth}] (-5, -4) -- (-5, -2) -- (-4, -2) -- (-4, -4) -- (-5, -4);
        \fill[black] (-4.5, -4) node[above] {$1^{\circ}$};
        \fill[black] (-4.5, -2) node[above] {$1$};

        \draw[red 5, thick, -{Stealth}] (0, 5) -- (2, 5) -- (2, 13) -- (0, 13) -- (0, 5);
        \fill[black] (0.5, 5) node[above] {$5$};
        \fill[black] (1.5, 5) node[above] {$5$};

        \fill[black] (1.5, 13) node[above] {$1^{\bullet}$};
        \fill[black] (0.5, 13) node[above] {$1^{\bullet}$};

        \draw[red 1!80!black, thick, -{Stealth}] (3, 6) -- (6, 6) -- (6, 13) -- (4, 13) -- (4, 12) -- (3, 12) -- (3, 6);
        \fill[black] (3.5, 6) node[above] {$5$};
        \fill[black] (4.5, 6) node[above] {$4$};
        \fill[black] (5.5, 6) node[above] {$3$};

        \fill[black] (5.5, 13) node[above] {$4^{\bullet}$};
        \fill[black] (4.5, 13) node[above] {$3^{\bullet}$};
        \fill[black] (3.5, 12) node[above] {$1^{\bullet}$};

        \foreach \x in {-2,...,0} {
            \foreach \y in {10,...,15} {
                \fill[black] (\x,\y) circle (1.5pt) node[left] {};
            }
        }
        \foreach \x in {1,...,3} {
            \foreach \y in {10,...,16} {
                \fill[black] (\x,\y) circle (1.5pt) node[left] {};
            }
        }
        \foreach \x in {4,...,5} {
            \foreach \y in {10,...,15} {
                \fill[black] (\x,\y) circle (1.5pt) node[left] {};
            }
        }
        \foreach \x in {6,...,9} {
            \foreach \y in {10,...,14} {
                \fill[black] (\x,\y) circle (1.5pt) node[left] {};
            }
        }

        \foreach \x in {-8,...,-4} {
            \foreach \y in {-2,...,2} {
                \fill[black] (\x,\y) circle (1.5pt) node[left] {};
            }
        }
        \foreach \y in {-1,...,3} {
            \fill[black] (-3, \y) circle (1.5pt) node[left] {};
        }
        \foreach \x in {-2,...,-1} {
            \foreach \y in {0,...,4} {
                \fill[black] (\x, \y) circle (1.5pt) node[left] {};
            }
        }
        \foreach \y in {0,...,5} {
            \fill[black] (0, \y) circle (1.5pt) node[left] {};
        }
        \foreach \x in {1,...,2} {
            \foreach \y in {1,...,5} {
                \fill[black] (\x, \y) circle (1.5pt) node[left] {};
            }
        }
        \foreach \y in {1,...,6} {
                \fill[black] (3, \y) circle (1.5pt) node[left] {};
            }
        \foreach \y in {2,...,7} {
            \fill[black] (4, \y) circle (1.5pt) node[left] {};
        }
        \foreach \y in {3,...,8} {
            \fill[black] (5, \y) circle (1.5pt) node[left] {};
        }
        \foreach \x in {6,...,9} {
            \foreach \y in {4,...,8} {
                \fill[black] (\x, \y) circle (1.5pt) node[left] {};
            }
        }

        \foreach \x in {-8,...,-1} {
            \foreach \y in {-4,...,-8} {
                \fill[black] (\x,\y) circle (1.5pt) node[left] {};
            }
        }
        \foreach \x in {0,...,2} {
            \foreach \y in {-4,...,-7} {
                \fill[black] (\x,\y) circle (1.5pt) node[left] {};
            }
        }

    \end{tikzpicture}
    \caption{A flow on the $T$-lattice.}
    \label{fig:lattice-flow-main}
\end{figure}
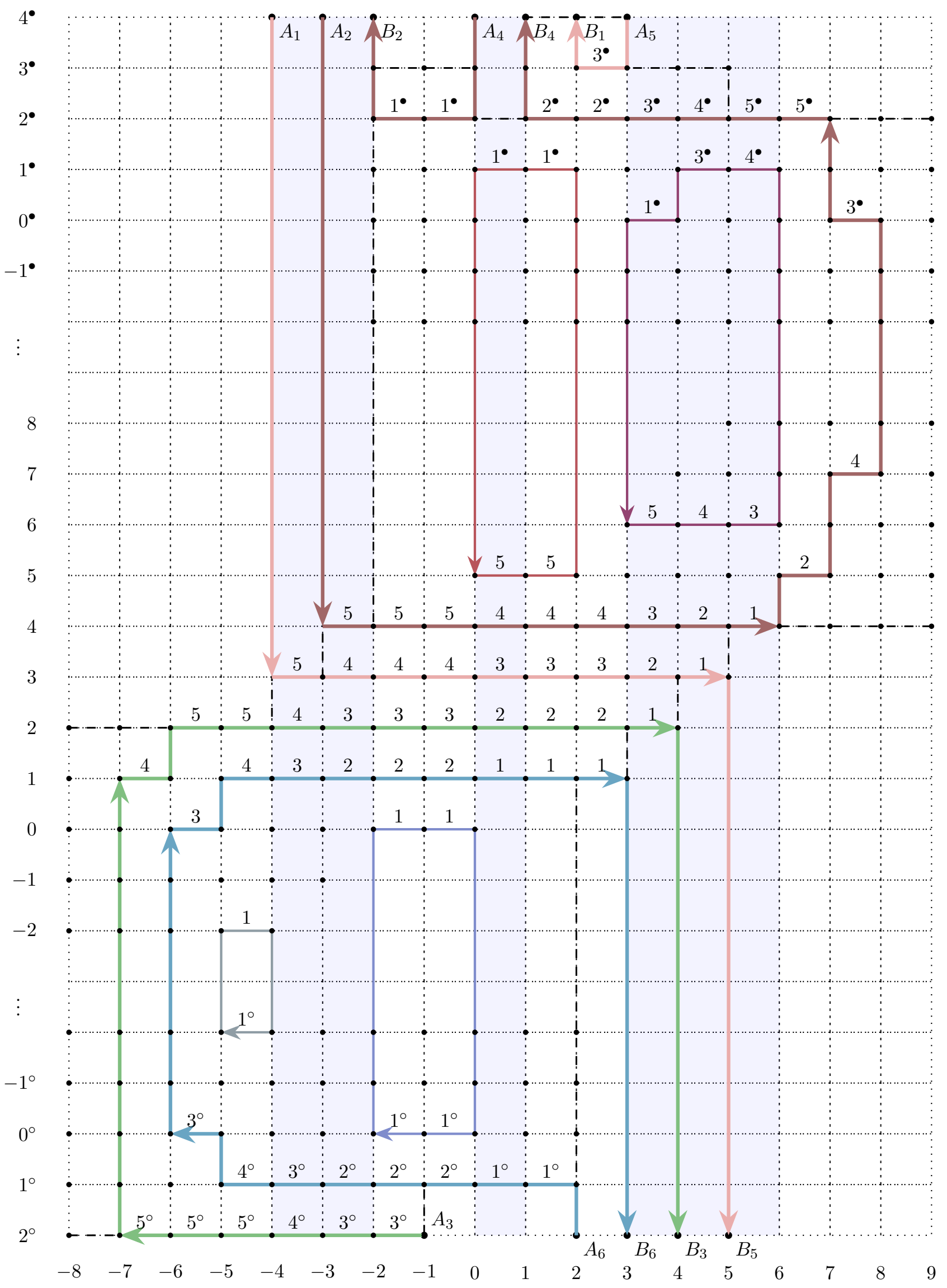

\begin{example}
    We refer to Example~\ref{ex:p-cells}. A trio tableau $(T^{\circ}, \widetilde{T}, T^{\bullet}) \in TT_{\lambda/\mu}(\mathbf{a}, \mathbf{b}, 5)$ is shown in Fig.~\ref{fig:trio-tableau-example-a}, \ref{fig:trio-tableau-example-b} and the corresponding flow on the T-lattice in Fig.~\ref{fig:lattice-flow-main}.
\end{example}

Let $X = (a, b)$ be some endpoint on the lattice; we use $X^1, X^2$ to denote the first and second coordinates of $X$, i.e. $X^1 = a, X^2 = b$. 



\subsection{Bijection}
\label{subsec:bijection}


In this subsection we establish weight-preserving bijection between set of trio-tablueax $TT_{\lambda/\mu}(\mathbf{a}, \mathbf{b}, m)$ and the set of flows $\mathcal{F}(\mathbf{A}, \mathbf{B})$ defined on the T-lattice $TL_{\lambda/\mu}(\pi, m, \mathbf{a}, \mathbf{b})$. 



We define the map $\phi : \mathcal{F}(\mathbf{A}, \mathbf{B}) \rightarrow TT_{\lambda/\mu}(\mathbf{a}, \mathbf{b}, m)$ by introducing an algorithm which converts the flow $\mathbf{F} = (\mathbf{P}, \mathbf{C}) \in \mathcal{F}(\mathbf{A}, \mathbf{B})$ into a trio tableau $(T^{\circ}, \widetilde{T}, T^{\bullet})$. 
The algorithm is the following:

Construct the left tableau $T^{\circ}$ of shape $\nu/\mu$ where  $\nu \in \mathsf{Sign}_n$ such that $\nu \subseteq \mu$. Let $\eta_1, \dots, \eta_r$ be all cells with content $c$ in $L(\lambda/\mu)$ sorted by row number in descending order, and let $e_1, \dots, e_s$ be edges in $\mathbf{F}$ such that $e_i = (c + 1, (p_i + f_c)^{\circ}) \rightarrow (c, (p_i + f_c)^{\circ})$ for all $i$, sorted by $p_i$ in descending order. Then the cell $\eta_i$ is filled with $p_i^{\circ}$ for all $i \in [s]$. 

Construct the right tableau $T^{\bullet}$ of shape $\lambda/\tau$ where $\tau \in \mathsf{Sign}_n$ such that $\lambda \subseteq \tau$. Let $\eta_1, \dots, \eta_r$ be all cells with content $c$ in $R(\lambda/\mu)$ sorted by row number in ascending order, and let $e_1, \dots, e_s$ be edges in $\mathbf{F}$ such that $e_i = (c + 1, (p_i - f_c)^{\bullet}) \rightarrow (c, (p_i - f_c)^{\bullet})$ for all $i$, sorted by $p_i$ in descending order. Then the cell $\eta_i$ is filled with $p_i^{\bullet}$ for all $i \in [s]$. 

Finally, construct the main tableau $\widetilde{T}$ of shape $\tau/\nu$. Let $\eta_1, \dots, \eta_r$ be all cells of content $c$ in $\tau/\nu$ sorted by row number in descending order, and let $e_1, \dots, e_s$ be edges in $\mathbf{F}$ such that $e_i = (c + 1, p_i + f_c) \rightarrow (c, p_i + f_c)$ for all $i$, sorted by $p_i$ in descending order. Then the cell $\eta_i$ is filled with $p_i$ for all $i \in [s]$. 

We then define $\phi(\mathbf{F}) := \mathbf{T} = (T^{\circ}, \widetilde{T}, T^{\bullet})$. It will be shown that $\nu \subseteq \tau$. Also $s \le r$ must hold for left, right tableaux; $s = r$ must hold for central tableau; and that $\mathbf{T}$ is, indeed, a trio-tableau.

\begin{definition}
    Define the relation $\mathrel{R}$ on the set of cells $\mathbb{Z}^2$ as follows: 
    $\eta {R} \kappa \iff \eta, \kappa$ are contained in the same tile $\hat{\theta}_k, k \in \mathbb{Z}$. 
    
\end{definition}

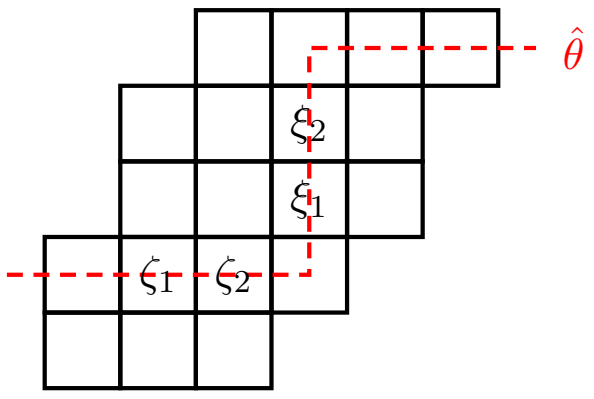
\begin{figure}
    \centering
        \begin{tikzpicture}[y=-1cm, scale=0.5, transform shape, >={Stealth[scale=0.5]}, every node/.append style={scale=1.5}]
            \foreach \x in {6,...,9} {
                \draw[thick] (\x, 1) rectangle ++(1, 1);
            }
            \foreach \x in {5,...,8} {
                \draw[thick] (\x, 2) rectangle ++(1, 1);
            }
            \foreach \x in {5,...,8} {
                \draw[thick] (\x, 3) rectangle ++(1, 1);
            }
            \foreach \x in {4,...,7} {
                \draw[thick] (\x, 4) rectangle ++(1, 1);
            }
            \foreach \x in {4,...,6} {
                \draw[thick] (\x, 5) rectangle ++(1, 1);
            }

            \draw [red, thick, densely dashed] (3.5, 4.5) -- (7.5, 4.5) -- (7.5, 2.5) -- (7.5, 1.5) -- (10.5, 1.5);

            \fill (5.5, 4.5) node {$\zeta_1$};
            \fill (6.5, 4.5) node {$\zeta_2$};
            \fill (7.5, 3.5) node {$\xi_1$};
            \fill (7.5, 2.5) node {$\xi_2$};
            \fill[red] (11, 1.5) node {$\hat{\theta}$};
        \end{tikzpicture}
    \caption{An example for cells $\eta$.}
        \label{fig:neigh-cells-neigh-arrows}
\end{figure}

\begin{figure}
    \centering
    \begin{subfigure}[b]{0.45\textwidth}
        \centering
        \begin{tikzpicture}[y=-1cm, scale=0.4, transform shape, >={Stealth[scale=0.5]}, every node/.append style={scale=1.5}]
            \draw[densely dashed, gray] (3, -1) -- (8, -1) -- (8, 4) -- (6, 4) -- (6, 6) -- (-1, 6) -- (-1, 1) -- (3, 1) -- (3, -1);
            \foreach \x in {3,...,4} {
                \draw[thick] (\x, 5) rectangle ++(1, 1);
            }
            \foreach \x in {2,...,3} {
                \draw[thin, dashed, gray] (\x, 4) rectangle ++(1, 1);
            }
            \foreach \x in {1,...,2} {
                \draw[thin, dashed, gray] (\x, 3) rectangle ++(1, 1);
            }
            \fill (3.5, 5.5) node {$\eta$};
            \fill (4.5, 5.5) node {$\kappa$};
            \fill (1.25, 2) node {$\ddots$};  
        \end{tikzpicture}
        \subcaption{$\eta, \kappa$ are the bottommost cells with contents $c-1, c$. }
        \label{fig:neigh-cells-neigh-arrows-bottommost-cells}
    \end{subfigure} \quad 
    \begin{subfigure}[b]{0.5\textwidth}
        \centering
        \begin{tikzpicture}[y=-1cm, scale=0.4, transform shape, >={Stealth[scale=0.5]}, every node/.append style={scale=1.5}]
            \draw[dashed, gray] (3, -1) -- (8, -1) -- (8, 4) -- (6, 4) -- (6, 6) -- (-1, 6) -- (-1, 1) -- (3, 1) -- (3, -1);
            \foreach \x in {5,...,5} {
                \draw[thin, dashed, gray] (\x, 5) rectangle ++(1, 1);
            }
            \foreach \x in {4,...,5} {
                \draw[thick] (\x, 4) rectangle ++(1, 1);
            }
            \foreach \x in {3,...,4} {
                \draw[thin, dashed, gray] (\x, 3) rectangle ++(1, 1);
            }

            \fill (4.5, 4.5) node {$\eta$};
            \fill (5.5, 4.5) node {$\kappa$};
            \fill (5.5, 5.5) node {$\zeta$};
            \fill (3.25, 2) node {$\ddots$};
            
        \end{tikzpicture}
        \subcaption{$\kappa$ is the bottommost cell with content $c$, but $\eta$ is row above the bottomost cell with content $c - 1$}
        \label{fig:neigh-cells-neigh-arrows-bottommost-cells-2}
    \end{subfigure}
    \caption{An example for cells $\eta, \kappa$. }
\end{figure}
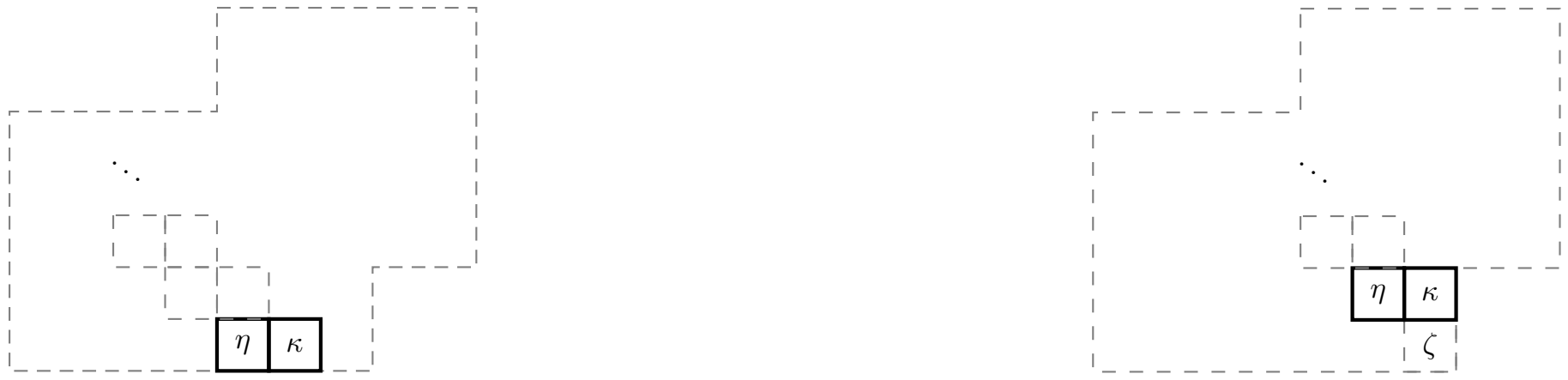

\begin{lemma}
    \label{lem:two-cells-two-edges}
    Let $\mathbf{F} \in \mathcal{F}(\mathbf{A}, \mathbf{B})$ and $(T^{\circ}, \widetilde{T}, T^{\bullet}) = \phi(\mathbf{F})$.
    Let $\eta, \kappa$ be cells in $T' \in (T^{\circ}, T^{\bullet})$ filled by edges $e_{\eta}, e_{\kappa}$ (resp.), and $\eta$ be left or bottom neighbor of $\kappa$. The edges $e_{\eta}, e_{\kappa}$ are contained in the same path $P$ if and only if $\eta R \kappa$.  
\end{lemma}
\begin{proof}
    We prove the statement for $T' = T^{\circ}$; for $T' = T^{\bullet}$ it can be done similarly.

    Let $c(\kappa) = c$, then $c(\eta) = c - 1$. Also let $e_{\eta} = (c, (p + f_c)^{\circ}) \rightarrow (c - 1, (p + f_c)^{\circ})$ and $e_{\kappa} = (c, (q + f_c)^{\circ}) \rightarrow (c - 1, (q + f_c)^{\circ})$.
    
    First we prove the necessity using induction on the row number. Let $\eta R \kappa$, we have either $\pi(c - 1) = (\cdot, 0),\ \pi(c) = (0, \cdot)$ (see Fig.~\ref{fig:neigh-cells-neigh-arrows}, $\eta = \zeta_1,\ \kappa = \zeta_2$) or $\pi(c - 1) = (\cdot, 1),\ \pi(c) = (1, \cdot)$ (see Fig.~\ref{fig:neigh-cells-neigh-arrows}, $\eta = \xi_1,\ \kappa = \xi_2$). We give proof for the first case using induction. There are two base cases two consider: (1) $\eta, \kappa$ are both bottommost cells of contents $c - 1, c$ (resp.) (Fig.~\ref{fig:neigh-cells-neigh-arrows-bottommost-cells}); (2) $\kappa$ is the bottommost cell of content $c$, but $\eta$ is one row above the bottomost cell of content $c - 1$ (Fig.~\ref{fig:neigh-cells-neigh-arrows-bottommost-cells-2}). 

    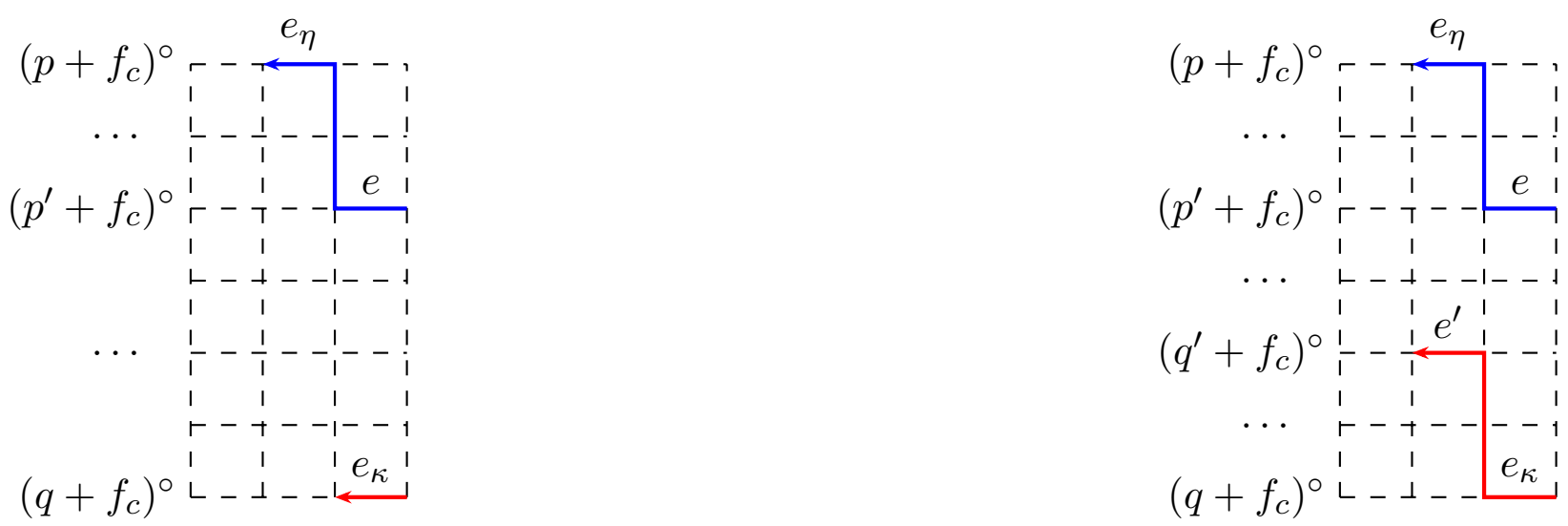
\begin{figure}
        \centering
        \begin{subfigure}[b]{0.5\textwidth}
        \centering
        \begin{tikzpicture}[scale=0.5, transform shape, >={Stealth[scale=0.5]}, every node/.append style={scale=1.5}]
            \draw[dashed, thin] (-1, 2) grid (2, 8);
            \draw[->, blue, thick] (2, 6) -- (1, 6) -- (1, 8) -- (0, 8);
            \draw[->, red, thick] (2, 2) -- (1, 2);
            \fill (1.5, 2) node[above] {$e_{\kappa}$};
            \fill (1.5, 6) node[above] {$e$};
            \fill (0.5, 8) node[above] {$e_{\eta}$};
            \fill (-1, 8) node[left] {$(p + f_c)^{\circ}$};
            \fill (-1, 6) node[left] {$(p' + f_c)^{\circ}$};
            \fill (-1, 2) node[left] {$(q + f_c)^{\circ}$};
            \fill (-1.5, 7) node[left] {$\ldots$};
            \fill (-1.5, 4) node[left] {$\ldots$};
        \end{tikzpicture}
        \subcaption{}
        \label{fig:one-path-one-ribbon-a}
        \end{subfigure}
        \begin{subfigure}[b]{0.45\textwidth}
        \centering
        \begin{tikzpicture}[scale=0.5, transform shape, >={Stealth[scale=0.5]}, every node/.append style={scale=1.5}]
            \draw[dashed, thin] (-1, 2) grid (2, 8);
            \draw[->, blue, thick] (2, 6) -- (1, 6) -- (1, 8) -- (0, 8);
            \draw[->, red, thick] (2, 2) -- (1, 2) -- (1, 4) -- (0, 4);
            \fill (1.5, 2) node[above] {$e_{\kappa}$};
            \fill (1.5, 6) node[above] {$e$};
            \fill (0.5, 4) node[above] {$e'$};
            \fill (0.5, 8) node[above] {$e_{\eta}$};
            \fill (-1, 8) node[left] {$(p + f_c)^{\circ}$};
            \fill (-1, 6) node[left] {$(p' + f_c)^{\circ}$};
            \fill (-1, 2) node[left] {$(q + f_c)^{\circ}$};
            \fill (-1, 4) node[left] {$(q' + f_c)^{\circ}$};
            \fill (-1.5, 7) node[left] {$\ldots$};
            \fill (-1.5, 5) node[left] {$\ldots$};
            \fill (-1.5, 3) node[left] {$\ldots$};
        \end{tikzpicture}
        \subcaption{}
        \label{fig:one-path-one-ribbon-b}
        \end{subfigure}
        \caption{An example of edge locations.}
    \end{figure}
    
    In the first case it follows that $p$ is maximal among all edges of content $c - 1$, and $q$ is maximal among all edges of content $c$. By contradiction argument assume that $e_{\eta}, e_{\kappa}$ are not contained in the same path. As $\eta$ is the bottommost cell which is not rightmost in $L(\lambda/\mu)$, there is no starting node at $(c, \cdot)$ in bottom layer, 
    thus there must be edge $e = (c + 1, (p' + f_c)^{\circ}) \rightarrow (c, (p' + f_c)^{\circ})$ which is contained in the same path as $e_{\eta}$. Due to definition of the map $\phi$, $p' < q$ must hold (see Fig.~\ref{fig:one-path-one-ribbon-a}). Then to avoid intersection of paths, there must be edge $e' = (c, (q' + f_c)^{\circ}) \rightarrow (c - 1, (q' + f_c)^{\circ})$ such that $q' > p$, which contradicts with the maximality of $p$ (see Fig.~\ref{fig:one-path-one-ribbon-b}). 

    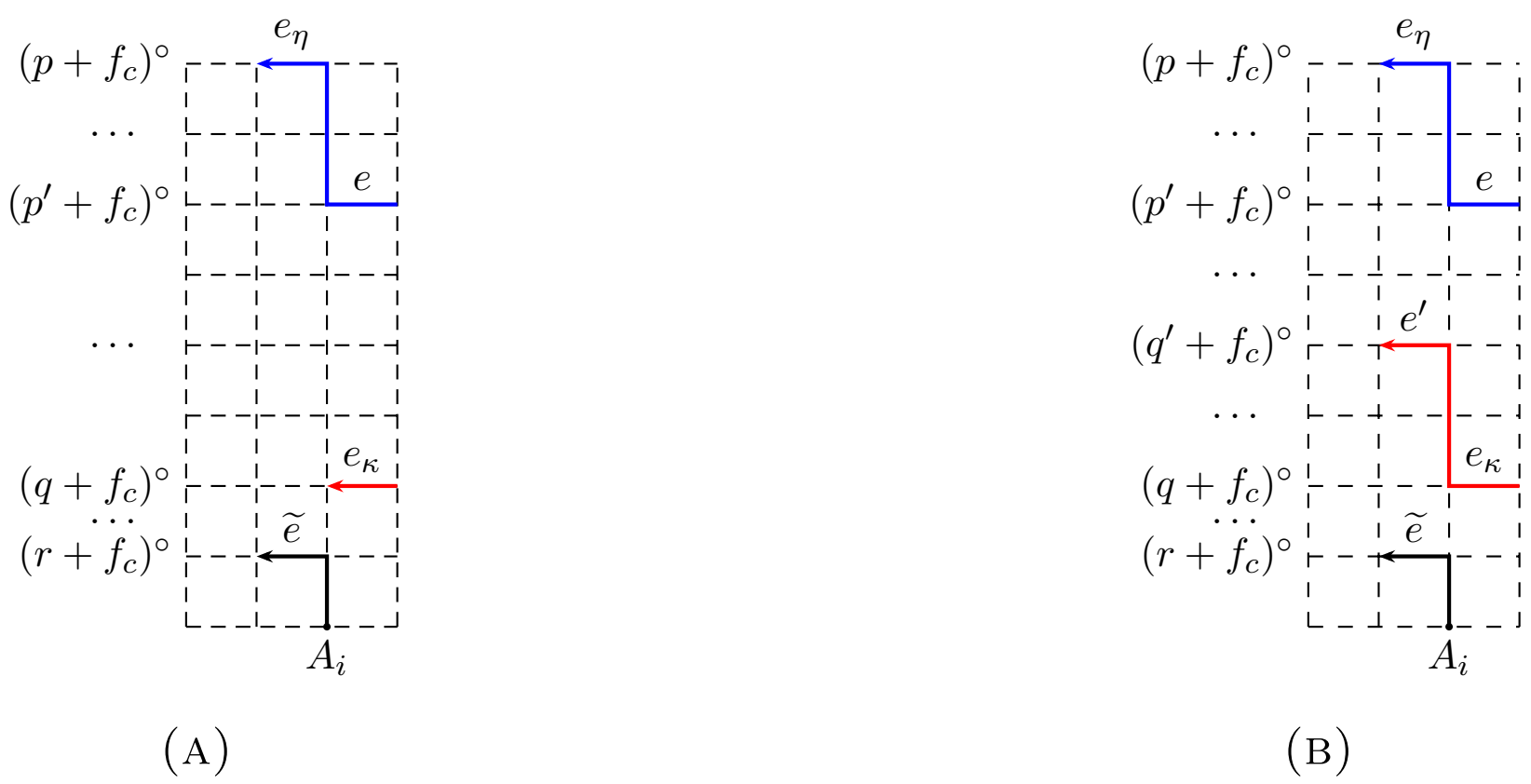
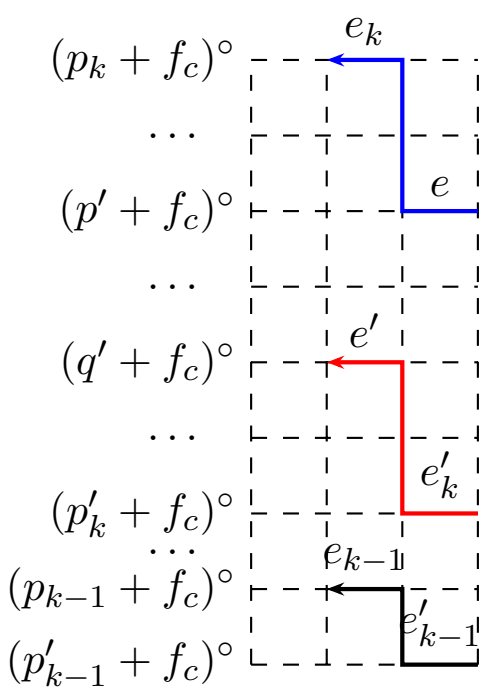
\begin{figure}
    \centering
    \begin{subfigure}[b]{0.5\textwidth}
    \centering
    \begin{tikzpicture}[scale=0.5, transform shape, >={Stealth[scale=0.5]}, every node/.append style={scale=1.5}]
        \draw[densely dashed, thin] (-1, 0) grid (2, 8);
        \draw[->, blue, thick] (2, 6) -- (1, 6) -- (1, 8) -- (0, 8);
        \draw[->, red, thick] (2, 2) -- (1, 2);
        \draw[->, black, thick] (1, 0) -- (1, 1) -- (0, 1);
        \fill (1.5, 2) node[above] {$e_{\kappa}$};
        \fill (1.5, 6) node[above] {$e$};
        \fill (0.5, 8) node[above] {$e_{\eta}$};
        \fill (0.5, 1) node[above] {$\widetilde{e}$};
        \fill (-1, 8) node[left] {$(p + f_c)^{\circ}$};
        \fill (-1, 6) node[left] {$(p' + f_c)^{\circ}$};
        \fill (-1, 2) node[left] {$(q + f_c)^{\circ}$};
        \fill (-1, 1) node[left] {$(r + f_c)^{\circ}$};
        \fill (-1.5, 7) node[left] {$\ldots$};
        \fill (-1.5, 4) node[left] {$\ldots$};
        \fill (-1.5, 1.5) node[left] {$\ldots$};
        \fill (1, 0) circle (1.5pt) node[left, below] {$A_i$};
    \end{tikzpicture}
    \subcaption{}
    \label{fig:one-path-one-ribbon-c}
    \end{subfigure}
    \begin{subfigure}[b]{0.45\textwidth}
    \centering
    \begin{tikzpicture}[scale=0.5, transform shape, >={Stealth[scale=0.5]}, every node/.append style={scale=1.5}]
        \draw[dashed, thin] (-1, 0) grid (2, 8);
        \draw[->, blue, thick] (2, 6) -- (1, 6) -- (1, 8) -- (0, 8);
        \draw[->, red, thick] (2, 2) -- (1, 2) -- (1, 4) -- (0, 4);
        \draw[->, black, thick] (1, 0) -- (1, 1) -- (0, 1);
        \fill (1.5, 2) node[above] {$e_{\kappa}$};
        \fill (1.5, 6) node[above] {$e$};
        \fill (0.5, 4) node[above] {$e'$};
        \fill (0.5, 8) node[above] {$e_{\eta}$};
        \fill (0.5, 1) node[above] {$\widetilde{e}$};
        \fill (-1, 8) node[left] {$(p + f_c)^{\circ}$};
        \fill (-1, 6) node[left] {$(p' + f_c)^{\circ}$};
        \fill (-1, 2) node[left] {$(q + f_c)^{\circ}$};
        \fill (-1, 4) node[left] {$(q' + f_c)^{\circ}$};
        \fill (-1, 1) node[left] {$(r + f_c)^{\circ}$};
        \fill (-1.5, 7) node[left] {$\ldots$};
        \fill (-1.5, 5) node[left] {$\ldots$};
        \fill (-1.5, 3) node[left] {$\ldots$};
        \fill (-1.5, 1.5) node[left] {$\ldots$};
        \fill (1, 0) circle (1.5pt) node[left, below] {$A_i$};
    \end{tikzpicture}
    \subcaption{}
    \label{fig:one-path-one-ribbon-d}
    \end{subfigure}
    \caption{An example of edge locations.}
    \end{figure}

    In the second case it follows that the bottommost cell $\zeta$ of content $c - 1$ is also located on the right perimeter of $\nu/\mu$. As $\pi(c - 1) = (\cdot, 0)$, there is a starting node $A_i \in \mathbf{A}$ at $(c, \cdot)$ in the bottom layer, and a path that starts at $A$ contains edge $\widetilde{e} = (c, (r + f_{c - 1})^{\circ}) \rightarrow (c - 1, (r + f_{c - 1})^{\circ})$ which fills the cell $\zeta$, thus to avoid intersection $p < r,\ q < r$ must hold. By contradiction argument, assume that $e_{\eta}, e_{\kappa}$ are not contained in the same path. As $e_{\eta}$ is not contained in the path which starts from node $A$, there must be edge $e = (c + 1, (p' + f_c)^{\circ}) \rightarrow (c, (p' + f_c)^{\circ})$ which is contained in the same path with $e_{\eta}$ (see Fig.~\ref{fig:one-path-one-ribbon-c}). Due to definition of the map $\phi$, $p' < q$ must hold. Then to avoid intersection of paths with edge $e$, there must be edge $e' = (c, (q' + f_c)^{\circ}) \rightarrow (c - 1, (q' + f_c)^{\circ})$ in the path that contains $e_{\kappa}$ such that $r > q' > p$, but it means that $\eta$ is should have been filled by the edge $e'$ instead of $e_{\eta}$ (see Fig.~\ref{fig:one-path-one-ribbon-d}), a contradiction.

    \begin{figure}
    \centering
    \begin{subfigure}[b]{0.45\textwidth}
        \centering
        \begin{tikzpicture}[y=-1cm, scale=0.5, transform shape, >={Stealth[scale=0.5]}, every node/.append style={scale=1.5}]
            \draw[dashed, gray] (3, -1) -- (8, -1) -- (8, 4) -- (6, 4) -- (6, 6) -- (-1, 6) -- (-1, 1) -- (3, 1) -- (3, -1);
            \foreach \x in {3,...,4} {
                \draw[thick] (\x, 5) rectangle ++(1, 1);
            }
            \foreach \x in {2,...,3} {
                \draw[thick] (\x, 4) rectangle ++(1, 1);
            }
            \foreach \x in {1,...,2} {
                \draw[thick] (\x, 3) rectangle ++(1, 1);
            }

            \fill (3.5, 5.5) node {$\eta_1$};
            \fill (4.5, 5.5) node {$\kappa_1$};
            \fill (2.5, 4.5) node {$\eta_2$};
            \fill (3.5, 4.5) node {$\kappa_2$};
            \fill (1.5, 3.5) node {$\eta_3$};
            \fill (2.5, 3.5) node {$\kappa_3$};
            \fill (1.25, 2) node {$\ddots$};
            
        \end{tikzpicture}
        \subcaption{}
        \label{fig:neigh-cells-neigh-arrows-all-cells}
    \end{subfigure}
    \begin{subfigure}[b]{0.5\textwidth}
        \centering
        \begin{tikzpicture}[y=-1cm, scale=0.5, transform shape, >={Stealth[scale=0.5]}, every node/.append style={scale=1.5}]
        \draw[dashed, gray] (3, -1) -- (8, -1) -- (8, 4) -- (6, 4) -- (6, 6) -- (-1, 6) -- (-1, 1) -- (3, 1) -- (3, -1);
        \foreach \x in {3,...,4} {
            \draw[thin, dashed, gray] (\x, 5) rectangle ++(1, 1);
        }
        \foreach \x in {1,...,2} {
            \draw[thick] (\x, 3) rectangle ++(1, 1);
        }
        \foreach \x in {2,...,3} {
            \draw[thin, dashed, gray] (\x, 4) rectangle ++(1, 1);
        }
        \foreach \x in {0,...,1} {
            \draw[thin, dashed, gray] (\x, 2) rectangle ++(1, 1);
        }
        \foreach \x in {-1,...,0} {
            \draw[thin, dashed, gray] (\x, 1) rectangle ++(1, 1);
        }
        \draw[thick] (2, 4) rectangle ++(1, 1);
        
        \draw [red, thick, densely dashed] (0.5, 3.5) -- (1.5, 3.5) -- (1.5, 2.5);
        \draw [red, thick, densely dashed] (2.5, 4.5) -- (2.5, 3.5) -- (3.5, 3.5);
        
        \fill (1.5, 3.5) node {$\eta$};
        \fill (2.5, 3.5) node {$\kappa$};
        \fill (2.5, 4.5) node {$\zeta$};
    \end{tikzpicture}
    \subcaption{} 
    \label{fig:one-path-one-ribbon-h}
    \end{subfigure}
\caption{An example of cell locations.}
    \end{figure}
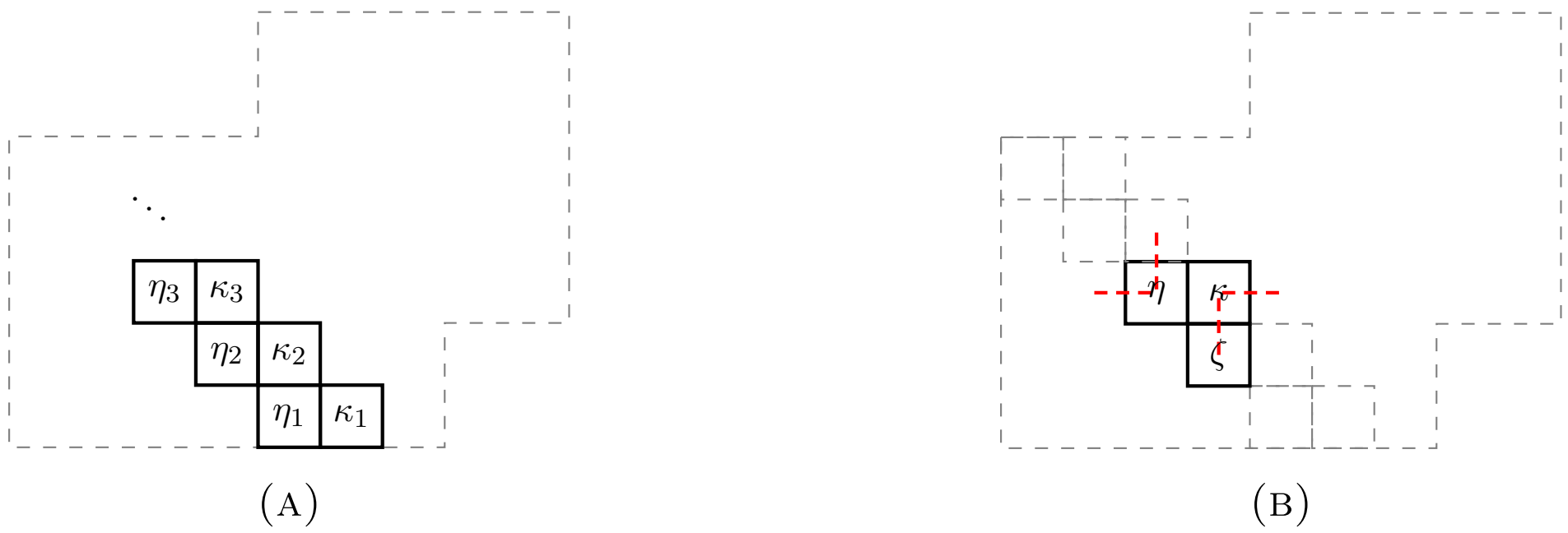

    Now, define pairs of cells $\eta_i, \kappa_i$ such that $c(\eta_i) = c - 1,\ c(\kappa_i) = c$; $\eta_i$ is left neighbor of $\kappa_i$; $\eta_i R \kappa_i$ and they are filled by edges $e_i, e'_i$ of contents $c - 1, c$ (resp.). The pairs $\eta_i, \kappa_i$ are sorted by the row number in descending order, i.e. $\mathrm{row}(\eta_i) > \mathrm{row}(\eta_{i + 1})$ for all $i$ (see Fig.~\ref{fig:neigh-cells-neigh-arrows-all-cells}).
    
    \begin{figure}
        \centering
        \centering
        \begin{tikzpicture}[scale=0.5, transform shape, >={Stealth[scale=0.5]}, every node/.append style={scale=1.5}]
            \draw[dashed, thin] (-1, 0) grid (2, 8);
            \draw[->, blue, thick] (2, 6) -- (1, 6) -- (1, 8) -- (0, 8);
            \draw[->, red, thick] (2, 2) -- (1, 2) -- (1, 4) -- (0, 4);
            \draw[->, black, thick] (2, 0) -- (1, 0) -- (1, 1) -- (0, 1);
            \fill (1.5, 2) node[above] {$e'_k$};
            \fill (1.5, 6) node[above] {$e$};
            \fill (0.5, 8) node[above] {$e_k$};
            \fill (0.5, 4) node[above] {$e'$};
            \fill (0.5, 1) node[above] {$e_{k-1}$};
            \fill (1.5, 0) node[above] {$e'_{k-1}$};
            \fill (-1, 8) node[left] {$(p_k + f_c)^{\circ}$};
            \fill (-1, 6) node[left] {$(p' + f_c)^{\circ}$};
            \fill (-1, 2) node[left] {$(p'_k + f_c)^{\circ}$};
            \fill (-1, 1) node[left] {$(p_{k-1} + f_c)^{\circ}$};
            \fill (-1, 0) node[left] {$(p'_{k-1} + f_c)^{\circ}$};
            \fill (-1, 4) node[left] {$(q' + f_c)^{\circ}$};
            \fill (-1.5, 7) node[left] {$\ldots$};
            \fill (-1.5, 5) node[left] {$\ldots$};
            \fill (-1.5, 3) node[left] {$\ldots$};
            \fill (-1.5, 1.5) node[left] {$\ldots$};
        \end{tikzpicture}
        \caption{An example of edge locations.}
                \label{fig:one-path-one-ribbon-f}

    \end{figure}

    Let $k > 1$, by induction suppose that for all $i < k$, $e_i, e'_i$ are contained in the same paths $P_i$. We show that $e_k, e'_k$ are contained in the same path $P_k$. Let $e_k = (c, (p_k + f_{c - 1})^{\circ}) \rightarrow (c - 1, (p_k + f_{c - 1})^{\circ})$ and $e'_k = (c + 1, (p'_k + f_{c - 1})^{\circ}) \rightarrow (c, (p'_k + f_{c - 1})^{\circ})$. By contradiction argument assume that $e_k, e'_k$ are not contained in the same path. There must be edge $e = (c + 1, (p' + f_c)^{\circ}) \rightarrow (c, (p' + f_c)^{\circ})$ which is contained in the same path with $e_k$. Due to definition of the map $\phi$, $p' < p'_k$ must hold. Then to avoid intersection of paths, there must be edge $e' = (c, (q' + f_c)^{\circ}) \rightarrow (c - 1, (q' + f_c)^{\circ})$ in the path that contains $e'_k$ such that $q' > p_k$, but it means that $\eta_k$ should have been filled by the edge $e'$ instead of $e_k$ (see Fig.~\ref{fig:one-path-one-ribbon-f}), a contradiction.
    
    The case $\pi(c - 1) = (\cdot, 1),\ \pi(c) = (1, \cdot)$ can be proved similarly.

    Now let us prove the sufficiency part. We argue by contradiction. Assume that the edges $e_{\eta}, e_{\kappa}$ are contained in the same path and $\eta \not\mathrel{R} \kappa$.  Similarly we have two cases: (1) $\pi(c - 1) = (\cdot, 1),\ \pi(c) = (1, \cdot)$; (2) $\pi(c - 1) = (\cdot, 0),\ \pi(c) = (0, \cdot)$. We consider the first case, thus $\eta$ must be left neighbor of $\kappa$ (as $\eta \not\mathrel{R} \kappa$). If there is a cell $\zeta \in T^{\circ}$ with $\mathrm{row}(\zeta) = \mathrm{row}(\eta) + 1$ and $\mathrm{col}(\zeta) = \mathrm{col}(\eta) + 1$ (see Fig.~\ref{fig:one-path-one-ribbon-h}), then $\zeta \mathrel{R} \kappa$ and by the nesessity part of current lemma, it follows that $\zeta$ must be filled with edge of content $c - 1$ which is contained in the path $P$ together with $e_{\kappa}$, which is a contradiction, as we have $e_{\eta}$ of content $c - 1$ and no two edges in bottom layer of content $c - 1$ can be contained in the same path $P$.

    Now if $\zeta \notin T^{\circ}$ then there must be ending node $B$ in the bottom layer with coordinates $(c, \cdot)$, such that the path $P$, that fills the cell $\kappa$, must terminate; but $P$ contains the edge $e_{\eta}$ of content $c - 1$, thus $P$ does not terminate at node $B$, contradiction. 

    The case $\pi(c - 1) = (\cdot, 0),\ \pi(c) = (0, \cdot)$ can be proved similarly.
\end{proof}

\begin{lemma}
    \label{lem:left-tableau-is-of-some-shape}
    Let $\eta_1, \dots, \eta_r$ be all cells of content $c$ in $L(\lambda/\mu)$ sorted by row number in descending order, and let $e_1, \dots, e_s$ be all edges of the bottom part of flow $\mathbf{F}$ such that $e_i = (c + 1, (p_i + f_c)^{\circ}) \rightarrow (c, (p_i + f_c)^{\circ})$ for all $i$, sorted by $p_i$ in descending order. 
    
    (1) We must have $r \ge s$, i.e. the number of cells of content $c$ in $L(\lambda/\mu)$ is at least the number of edges of content $c$ in the bottom part of the flow $\mathbf{F}$.
    
    (2) The tableau $T^{\circ}$ has shape $\nu/\mu$, for some shape $\nu \subset L(\lambda/\mu)$. 
\end{lemma}

\begin{figure}
    \centering
    \begin{tikzpicture}[y=-1cm, scale=0.6]

    \draw[thick] (0, 5) -- (0, 4) -- (1, 4) -- (1, 3) -- (2, 3) -- (2, 1);

    \foreach \x in {-4,...,1} {
        \draw[gray, dashed] (\x, 1) rectangle ++(1, 1);
    }
    \foreach \x in {-4,...,1} {
        \draw[gray, dashed] (\x, 2) rectangle ++(1, 1);
    }
    \foreach \x in {-4,...,0} {
        \draw[gray, dashed] (\x, 3) rectangle ++(1, 1);
    }
    \foreach \x in {-4,...,-1} {
        \draw[gray, dashed] (\x, 4) rectangle ++(1, 1);
    }

    \fill (1.5, 3.5) node {$\mu$};
    \fill (1.5, 1.5) node {$\eta_1$};

    \end{tikzpicture}
    \caption{An example of $\mu$ and $\eta_1$.}
    \label{fig:filling-is-shape-c}
\end{figure}
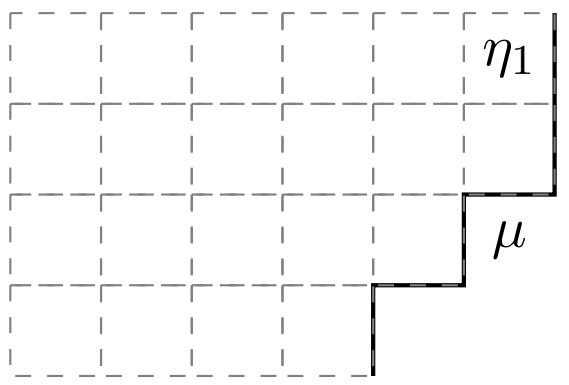

\begin{figure}
    \centering
    \begin{subfigure}[b]{0.5\textwidth}
    \centering
    \begin{tikzpicture}[scale=0.5, transform shape, >={Stealth[scale=0.5]}, every node/.append style={scale=1.5}]
        \draw[dashed, thin] (-1, 0) grid (2, 3);
        \draw[->, red, thick] (1, 0) -- (1, 2) -- (0, 2);
        \fill (0.5, 2) node[above] {$e$};
        \fill (1.5, 0) node[below] {$c_M + 1$};
        \fill (0, 0) node[below] {$c_M$};
        \fill (1, 0) circle (1.5pt) node[above right] {$A$};
    \end{tikzpicture}
    \subcaption{}
    \label{fig:filling-is-shape-d}
    \end{subfigure}
    \begin{subfigure}[b]{0.45\textwidth}
    \centering
    \begin{tikzpicture}[scale=0.5, transform shape, >={Stealth[scale=0.5]}, every node/.append style={scale=1.5}]
        \draw[dashed, thin] (-1, 0) grid (2, 3);
        \draw[->, red, thick] (1, 3.5) -- (1, 1) -- (0, 1);
        \fill (0.5, 1) node[above] {$e$};
        \fill (1.5, 0) node[below] {$c_M + 1$};
        \fill (0, 0) node[below] {$c_M$};
    \end{tikzpicture}
    \subcaption{}
    \label{fig:filling-is-shape-e}
    \end{subfigure}
    \caption{An example of edge locations.}
\end{figure}
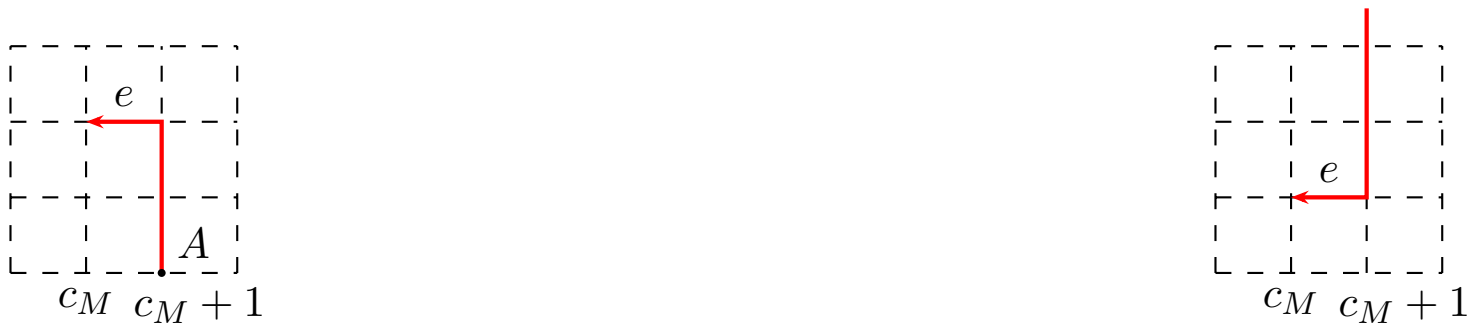

\begin{figure}
    \centering
    \begin{subfigure}[b]{0.5\textwidth}
    \centering
    \begin{tikzpicture}[scale=0.5, transform shape, >={Stealth[scale=0.5]}, every node/.append style={scale=1.5}]
        \draw[dashed, thin] (-1, 0) grid (2, 5);
        \draw[->, blue, thick] (0.95, 5.5) -- (0.95, 1) -- (0, 1);
        \draw[->, red, thick] (1.05, 5.5) -- (1.05, 2) -- (0, 2);
        \fill (0.5, 2) node[above] {$e$};
        \fill (0.5, 1) node[below] {$e''$};
        \fill (0.85, 1.5) node[left] {$\dots$};
    \end{tikzpicture}
    \subcaption{}
    \label{fig:filling-is-shape-a}
    \end{subfigure}
    \begin{subfigure}[b]{0.45\textwidth}
    \centering
    \begin{tikzpicture}[scale=0.5, transform shape, >={Stealth[scale=0.5]}, every node/.append style={scale=1.5}]
        \draw[dashed, thin] (-1, 0) grid (2, 5);
        \draw[->, blue, thick] (0.95, 5.5) -- (0.95, 0);
        \draw[->, red, thick] (1.05, 5.5) -- (1.05, 2) -- (0, 2);
        \fill (0.5, 2) node[above] {$e$};
        \fill (1, 0) node[below] {$B$};
    \end{tikzpicture}
    \subcaption{}
    \label{fig:filling-is-shape-b}
    \end{subfigure}
    \caption{An example of edge locations.}
\end{figure}
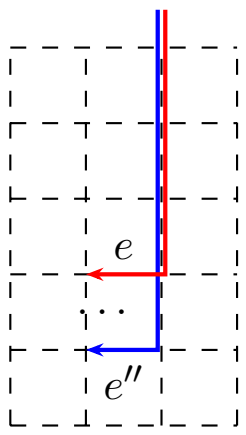
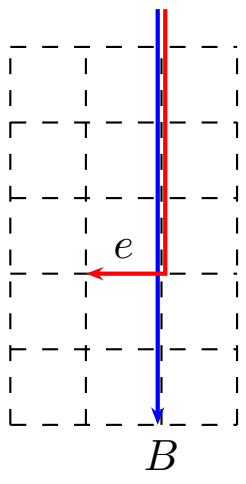

\begin{proof}
    We prove part (1) by induction on the content number starting from the maximal content in $L(\lambda/\mu)$. For the base case, let $c = c_M$ be the maximum content in $L(\lambda/\mu)$ (see e.g. Fig~\ref{fig:filling-is-shape-c}), there is only one cell with content $c_M$ in $L(\lambda/\mu)$, i.e. $r = 1$. Assume that $\pi(c_M) = (\cdot, 0)$, by the construction of the lattice there are no horizontal edges of the form $(c + 1, x^{\circ}) \rightarrow (c, x^{\circ})$ for all $c > c_M$ and $x \in \mathbb{Z}$, thus $\eta_1$ is filled by the path $P$ that starts at the node $A = (c_M + 1, \infty^{\circ}) \in \mathbf{A}$, and there is only one possible edge of content $c_M$, i.e. $r \ge s$ (see Fig.~\ref{fig:filling-is-shape-d}). Now assume that $\pi(c_M) = (\cdot, 1)$. By the same reason as above and as there is no starting node of the form $(c_M + 1, \infty^{\circ})$, there can be only one edge $e$ of content $c_M$ (i.e. $r \ge s$), otherwise a path containing $e$ will inevitably intersect any other edge of content $c_M$ (see Fig.~\ref{fig:filling-is-shape-d}).

    Define content $c' < c_M$. Now suppose that for all contents $c$ with $c_M \ge c > c'$, the condition $r \ge s$ holds. We prove that it also must hold for $c = c'$. Let $r'$ be the number of cells of content $c + 1$ in $L(\lambda/\mu)$, then either $r = r'$ or $r = r' + 1$. Let $s'$ be the number of edges $(e'_i)$ of content $c + 1$ in the bottom part of $\mathbf{F}$, and also let $P_1, \ldots, P_{s'}$ be paths that contain the edges $(e'_i)$. Assume that $r = r'$,  we consider two cases: (1) $\pi(c) = (\cdot, 0)$; and (2) $\pi(c) = (\cdot, 1)$. In the first case, there is no starting node at $(c, \infty^{\circ})$ (as $r = r'$), and every edge $e_i$ is contained in $P_i$ for some $i \in [s']$, thus $s' \ge s$, and $r \ge s$. In the second case, there is ending node for $B \in \mathbf{B}$ at $(c + 1, \infty^{\circ})$, but there can be new edge $e' = (c + 1, (p + f_c)^{\circ}) \rightarrow (c, (p + f_c)^{\circ})$ such that $p < p_i$ for all $p_i$ and $e'$ is not contained in any of paths $P_i$, thus $s' \ge s$ and $r \ge s$. 
    For $r = r' + 1$ the statement can be checked similarly.
    
    Now we prove part (2). Let $\eta = (i, j)$ be some cell in $T^{\circ}$, it suffices to show that if $\eta$ is not rightmost cell in $L(\lambda/\mu)$, the cell $\kappa = (i, j + 1)$ must be in $T^{\circ}$; and, if $\eta$ is not bottommost cell in $L(\lambda/\mu)$, the cell $(i + 1, j)$ must be in $T^{\circ}$. We prove the fist part; the second part can be proved similarly. Let $c(\eta) = c - 1,\ c(\kappa) = c$. 
    
    
    If $\eta \mathrel{R} \kappa$, then $\pi(c - 1) = (\cdot, 0),\ \pi(c) = (0, \cdot)$. There can be starting node at $A = (c, \cdot)$ in bottom layer, but as $\eta$ is not rightmost cell in $L(\lambda/\mu)$,
    it cannot be filled by a path that starting at $A$, 
    thus there must be some edge $e'$ of content $c$ contained in some path $P$ together with $e$. By Lemma~\ref{lem:two-cells-two-edges} it follows that $\kappa$ must be filled by $e'$.

    If $\eta \not\mathrel{R} \kappa$, then $\pi(c - 1) = (\cdot, 1),\ \pi(c) = (1, \cdot)$. Assume, by contradiction, that $\kappa$ is not filled by any edge. There are two cases to consider: (1) the cell $\zeta := (i + 1, j + 1) \in L(\lambda/\mu)$; and (2) $\zeta \not\in L(\lambda/\mu)$. In the first case, $\zeta$ must be filled by some edge $e''$, and the path that contains $e''$ intersects with the edge $e$, which contradicts the definition of a flow (Fig.~\ref{fig:filling-is-shape-a}). In the second case, it follows that there is no edge of content $c$ in bottom part of $\mathbf{F}$, and there is ending node $B = (c, \cdot) \in \mathbf{B}$ in the bottom layer,
    thus the path that ends on $B$ intersects the edge $e$, a contradiction (Fig.~\ref{fig:filling-is-shape-b}). 
\end{proof}

\begin{lemma}
    \label{lem:right-tableau-is-of-some-shape}
    Let $\eta_1, \dots, \eta_r$ be all cells of content $c$ in $R(\lambda/\mu)$ sorted by row number in ascending order, and let $e_1, \dots, e_s$ be all edges of the top part of $\mathbf{F}$ such that $e_i = (c + 1, (p_i + f_c)^{\bullet}) \rightarrow (c, (p_i + f_c)^{\bullet})$ for all $i$, sorted by $p_i$ in ascending order. 
    
    (1) We must have $r \ge s$, i.e. the number of cells of content $c$ in $R(\lambda/\mu)$ is at least the number of edges of content $c$ in the top part of the flow $\mathbf{F}$.
    
    (2) The tableau $T^{\bullet}$ is of shape $\lambda/\tau$, for some shape $\tau \subset R(\lambda/\mu)$. 
\end{lemma}
\begin{proof}
    Similar to the proof of Lemma~\ref{lem:left-tableau-is-of-some-shape}. 
\end{proof}

Next our aim is to prove that $\nu \subseteq \tau$. 
\begin{lemma}
    \label{lem:central-tableau-is-of-some-shape}
    Let $\eta_1, \dots, \eta_r$ be all cells of content $c$ in $\tau/\nu$ sorted by row number in descending order, and let $e_1, \dots, e_s$ be all edges of the middle part of $\mathbf{F}$ such that $e_i = (c + 1, p_i + f_c) \rightarrow (c, p_i + f_c)$ for all $i$, sorted by $p_i$ in descending order. 
    
    (1) We must have $r = s$, i.e. the number of cells of content $c$ in $\tau/\nu$ is equal to the number of edges of content $c$ in the middle part of the flow $\mathbf{F}$.

    (2) For adjacent cells $\alpha, \beta$ in $\tau/\mu$ filled by edges $e, e'$ we have that $\alpha \mathrel{R} \beta$ if and only if $e, e'$ are in the same path. 
    
    (3) The tableau $\widetilde{T}$ is of shape $\tau/\nu$, where $\nu \subseteq \tau$. 
\end{lemma}

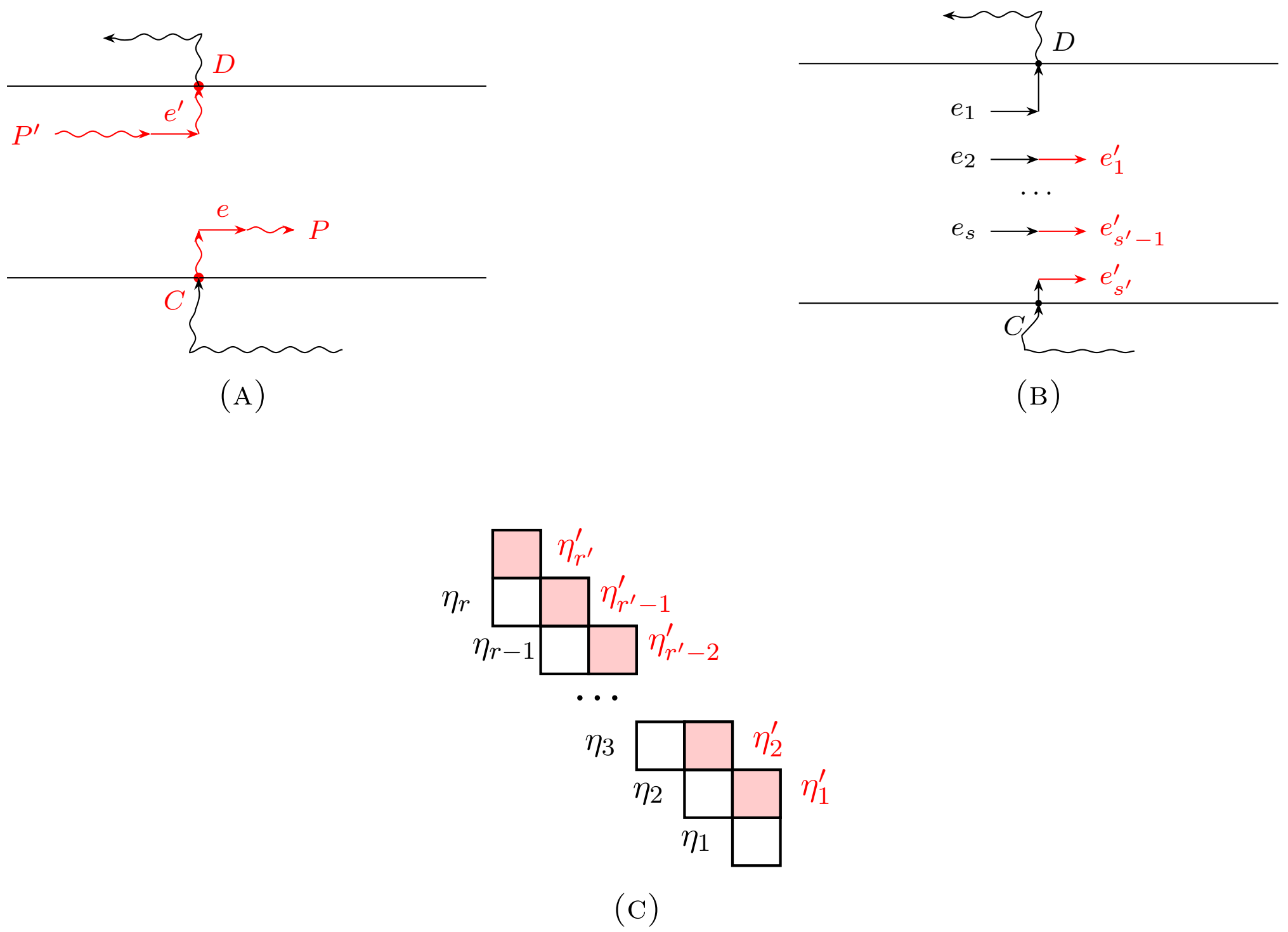
\begin{figure}
    \centering
    \begin{subfigure}[b]{0.5\textwidth}
        \centering
        \begin{tikzpicture}[scale=0.5, font=\scriptsize, >={Stealth[scale=0.75]}]


            \fill[red] (2, 1) circle (3pt) node[below left] {$C$};
            \fill[red] (2, 5) circle (3pt) node[above right] {$D$};
            \fill[red] (2.5, 2) circle (0pt) node[above] {$e$}; 
            \fill[red] (4, 2) circle (0pt) node[right] {$P$}; 
            \fill[red] (1.5, 4) circle (0pt) node[above] {$e'$}; 
            \fill[red] (-1, 4) circle (0pt) node[left] {$P'$}; 

            \draw (-2, 5) -- (8, 5);
            \draw (-2, 1) -- (8, 1);

            \draw[decorate, decoration={snake, amplitude=0.3mm, segment length=3mm}, ->, black] (5, -0.5) -- (2, -0.5) -- (2, 1);
            \draw[decorate, decoration={snake, amplitude=0.3mm, segment length=3mm}, ->, red] (2, 1) -- (2, 2);
            \draw[->, red] (2,2) -- (3, 2); 
            \draw[decorate, decoration={snake, amplitude=0.3mm, segment length=3mm}, ->, red] (3, 2) -- (4, 2);
            \draw[->, red] (1, 4) -- (2, 4);
            \draw[decorate, decoration={snake, amplitude=0.3mm, segment length=3mm}, ->, red] (-1, 4) -- (1, 4);
            \draw[decorate, decoration={snake, amplitude=0.3mm, segment length=3mm}, ->, red] (2, 4) -- (2, 5);
            \draw[decorate, decoration={snake, amplitude=0.3mm, segment length=3mm}, ->, black] (2, 5) -- (2, 6) -- (0, 6);
        \end{tikzpicture}
        \subcaption{}
        \label{fig:central-tableau-is-of-some-shape-a}
    \end{subfigure}
    \begin{subfigure}[b]{0.49\textwidth}
        \centering
        \begin{tikzpicture}[scale=0.5, font=\scriptsize, >={Stealth[scale=0.75]}]

            \fill (3, 6) circle (2pt) node[above right] {$D$};
            \fill (3, 1) circle (2pt) node[below left] {$C$};

            
            \fill (2, 2.5) circle (0pt) node[left] {$e_s$};
            \fill (2, 4) circle (0pt) node[left] {$e_2$};
            \fill (2, 5) circle (0pt) node[left] {$e_1$};
            \fill[red] (4, 1.5) circle (0pt) node[right] {$e'_{s'}$};
            \fill[red] (4, 2.5) circle (0pt) node[right] {$e'_{s' - 1}$};
            \fill[red] (4, 4) circle (0pt) node[right] {$e'_1$};
            \fill[black] (3, 3) circle (0pt) node[above] {$\dots$};

            \draw (-2, 6) -- (8, 6);
            \draw (-2, 1) -- (8, 1);

            \draw[decorate, decoration={snake, amplitude=0.15mm, segment length=3mm}, ->, black] (5, 0) -- (3, 0) -- (3, 1);
            \draw[decorate, decoration={snake, amplitude=0.3mm, segment length=3mm}, ->, black] (3, 6) -- (3, 7) -- (1, 7);
            \draw[->, black] (3, 1) -- (3, 1.5); 
            \draw[->, black] (2, 2.5) -- (3, 2.5); 
            \draw[->, black] (2, 4) -- (3, 4); 
            \draw[->, black] (2, 5) -- (3, 5); 
            \draw[->, black] (3, 5) -- (3, 6); 
            \draw[->, red] (3, 1.5) -- (4, 1.5); 
            \draw[->, red] (3, 2.5) -- (4, 2.5); 
            \draw[->, red] (3, 4) -- (4, 4); 
 
        \end{tikzpicture}
        \subcaption{}
        \label{fig:central-tableau-is-of-some-shape-b}
    \end{subfigure}\\[1cm]
    \begin{subfigure}[b]{0.45\textwidth}
        \centering
        \begin{tikzpicture}[y=-1cm, transform shape, scale=0.5, every node/.append style={scale=1.85}]
            \foreach \x in {1,...,1} {
                \draw[thick] (\x, 0) rectangle ++(1, 1);
            }
            \foreach \x in {1,...,2} {
                \draw[thick] (\x, 1) rectangle ++(1, 1);
            }
            \foreach \x in {2,...,3} {
                \draw[thick] (\x, 2) rectangle ++(1, 1);
            }
            \foreach \x in {4,...,5} {
                \draw[thick] (\x, 4) rectangle ++(1, 1);
            }
            \foreach \x in {5,...,6} {
                \draw[thick] (\x, 5) rectangle ++(1, 1);
            }
            \foreach \x in {6,...,6} {
                \draw[thick] (\x, 6) rectangle ++(1, 1);
            }

            \foreach \x in {1,...,1} {
                \draw[thick, fill=red!20] (\x, 0) rectangle ++(1, 1);
            }
            \foreach \x in {2,...,2} {
                \draw[thick, fill=red!20] (\x, 1) rectangle ++(1, 1);
            }
            \foreach \x in {3,...,3} {
                \draw[thick, fill=red!20] (\x, 2) rectangle ++(1, 1);
            }
            \foreach \x in {5,...,5} {
                \draw[thick, fill=red!20] (\x, 4) rectangle ++(1, 1);
            }
            \foreach \x in {6,...,6} {
                \draw[thick, fill=red!20] (\x, 5) rectangle ++(1, 1);
            }



            \fill (5.25, 6.5) node {$\eta_1$};
            \fill (4.25, 5.5) node {$\eta_2$};
            \fill (3.25, 4.5) node {$\eta_3$};
            \fill (1.25, 2.5) node {$\eta_{r - 1}$};
            \fill (0.25, 1.5) node {$\eta_r$};
            \fill (3.25, 3.5) node {$\boldsymbol{\dots}$};

            \fill[red] (7.75, 5.35) node {$\eta'_1$};
            \fill[red] (6.75, 4.35) node {$\eta'_2$};
            \fill[red] (5, 2.35) node {$\eta'_{r' - 2}$};
            \fill[red] (4, 1.35) node {$\eta'_{r' - 1}$};
            \fill[red] (2.75, 0.35) node {$\eta'_{r'}$};


        \end{tikzpicture}
        \subcaption{}
        \label{fig:central-tableau-is-of-some-shape-c}
    \end{subfigure}
    \caption{An example of paths and cells $\eta_i, \eta'_i$.}
    \label{fig:central-tableau-is-of-some-shape}
\end{figure}

\begin{proof}
    First assume that ${D}(\tau, \nu) = \varnothing$, it follows that $\tau/\mu$ consists of only null-strips. 
    If there is no middle part of $\mathbf{F}$ 
    there is nothing to prove. Suppose that there are some minor sub-paths of some paths in $\mathbf{F}$ in the middle layer.
    We claim that the minor paths in the middle layer do not contain right-edges. By contradiction, 
    let $P$ be the minor path in the middle layer that contains the leftmost right-edge $e = (c, x) \rightarrow (c + 1, x)$ for some $c, x \in \mathbb{Z}$ (i.e. content of $e$ is minimal among all right-edges in the middle layer). Let $C = (c, \cdot)$ be the minor starting node from which the path $P$ starts. As $\tau/\nu$ consists of only null-strips, it follows that there must be ending node $D = (c, \cdot)$ on the same vertical as $C$. As the path $P$ contains right-edge $e$, it does not terminate at $D$, thus there must be another path $P'$ in the middle layer which terminates at $D$. As $P'$ does not start at $C$ it follows that $P'$ must contain the right-edge $e' = (c - 1, y) \rightarrow (c, y)$ (Fig.~\ref{fig:central-tableau-is-of-some-shape-a}), which contradicts with the minimality of $c$, thus there are no paths with right-edges in the middle layer; it follows that $r = s = 0$ and $\nu \subseteq \tau$. 

    Now suppose that $D(\tau, \nu) \neq \varnothing$ and it consists of connected components $D_1, \ldots, D_k, k \ge 1$. Let $D_1$ be the bottommost connected component, and $\alpha$ be the bottom-left cell in $D_1$ with content $c = c(\alpha)$. It can be shown using a similar argument by contradiction as above, that there are no minor paths in the middle layer that contain right-edge of content strictly less than $c$. 
    
    Suppose that row numbers of cells in the diagram $D_1$ form the set $[a, b]$;
    we claim that $\nu_i < \tau_i, i \in [a, b]$. Assume the contrary, then particularly we have that $\nu_b > \tau_b$. 
    There are two cases: (1) $\pi(c) = (0, \cdot)$; (2) $\pi(c) = (1, \cdot)$. 
    In the first case, 
    then there must be minor path $P$ in the middle layer that terminates at $D = (c, \cdot)$; as there are no paths that contain right edge of content less than $c$ it follows that $P$ starts at minor starting node $C = (c, \cdot)$; i.e. $P$ does not contain right-edge. Thus, there must be $j \in [n]$ such that $\nu_j - j = \tau_b - b$. If $j > b$, then as there are no cells lower than $\alpha$ in $D(\tau, \nu)$
    it follows that $\tau_j = \nu_j = \tau_b - b + j$, thus we have that $\tau_j > \tau_b$ where $j > b$, consequently $\tau \notin \mathsf{Sign}_n$, a contradiction. If $j \le b$, then $\nu_j = \tau_b - b + j < \nu_b - b + j$, and $\nu_j < \nu_b$ where $j \le b$; a contradiction. 
    Thus, if $\pi(c) = (0, \cdot)$ we have that $\nu_b < \tau_b$, consequently (as $D_1$ is a connected diagram), $\nu_i < \tau_i, i \in [a, b]$. It can be shown in a similar way that if $\pi(c) = (1, \cdot)$, then $\nu_i < \tau_i, i \in [a, b]$ must hold. The path $P$ starts at minor starting node $C = (c, \cdot)$, and it must contain at least one right-edge $e_1$ of content $c$ (otherwise it leads to similar contradiction as above). 
    Thus, the edge $e_1$ fills the cell $\alpha$. Let $\beta$ be the top-right cell of $D_1$, then contents of cells in $D_1$ lie in the interval $[c(\alpha), c(\beta)]$. Now we show that the number of cells of content $c \in [c(\alpha), c(\beta)]$ is equal to the number of right-edges of content $c$ 
    using induction on $c \in [c(\alpha), c(\beta)]$. For $c = c(\alpha)$ we have already shown it. Suppose that the statement holds for $c \in [c(\alpha), c(\beta) - 1]$, we show that it also must hold for $c + 1$. Let $\eta_1, \ldots, \eta_r, r \ge 1$ be cells of content $c$ in $D_1$ sorted by row number in descending order, and $\eta'_1, \ldots, \eta'_{r'}$ be cells of content $c + 1$ in $D_1$ sorted in the same way (Fig.~\ref{fig:central-tableau-is-of-some-shape-c}). 
    Similarly, let $e_1, \ldots, e_s, s \ge 1$ be right-edges of content $c$, where $e_i = (c, p_i + f_c) \rightarrow (c + 1, p_i + f_c), i \in [s]$, sorted by $p_i$ in descending order, 
    and $e'_1, \ldots, e'_{s'}$ be right-edges of content $c + 1$ sorted in the same way (Fig.~\ref{fig:central-tableau-is-of-some-shape-b}). 
    By induction hypothesis, we have $r = s$. 
    Suppose that $\pi(c + 1) = (0, \cdot)$, there are four cases induced by presence or absence of minor ending node at $(c + 1, \cdot)$, and minor starting node at $(c + 1, \cdot)$. Suppose that both are present (Fig~\ref{fig:central-tableau-is-of-some-shape-b}), other cases can be handled similarly. Then, it follows that $\eta_1$ is on the right perimeter of $D_1$, and there must be cell $\eta'_{r'}$ on top of $\eta_r$ of content $c + 1$, and $\eta'_{r'}$ is on left perimeter of $D_1$ (Fig~\ref{fig:central-tableau-is-of-some-shape-c}); thus, $r' = r$, and $s' = s$ (as every edge of content $c$ except $e_1$ continues, and there is additional right-edge $e'_{s'}$ in the path that starts from the minor endpoint at $(c + 1, \cdot)$), thus $r' = s'$ which concludes the induction. The case $\pi(c + 1) = (1, \cdot)$ is handled by similar argument, which completes the induction. 

    Now we prove (2) for cells in $D_1$. Let $c \in [c(\alpha), c(\beta) - 1]$. First suppose that $\pi(c + 1) = (0, \cdot)$; again, there are four cases induced by presence or absence of minor ending node at $(c + 1, \cdot)$, and minor starting node at $(c + 1, \cdot)$. Suppose that both are present, other cases can be handled similarly. Define the cells $(\eta_1, \ldots, \eta_r), (\eta'_1, \ldots, \eta'_{r'})$ and edges $(e_1, \ldots, e_s), (e'_1, \ldots, e'_{s'})$ as above, then 
    by definition of the map $\phi$, the cell $\eta'_i$ is filled by the edge $e'_i, i \in [r']$. From the induction above it follows that the cell $\eta'_i$ is right neighbor of $\eta_{i + 1}, i < r$ (i.e. $\eta'_i \mathrel{R} \eta_{i + 1}, i < r$), and $e'_i$ is in the same path as $e_{i + 1}$ for $i < s$. The case $\pi(c + 1) = (1, \cdot)$ can be shown similarly. 

    Now we claim that $\nu \subseteq \tau$ and the number of cells of content $c \in \mathbb{Z}$ in diagrams $D_i, i \in [k]$ is equal to the number of right-edges of content $c$ on the lattice. 

    The statement is true for $D_1$. We continue by induction, assume that the statement is true for the diagram $D_i, i \in [k - 1]$, then we prove that it must hold for the diagram $D_{i + 1}$. Let $\beta'$ be the top-right cell of $D_i$. From induction hypothesis it follows that there is single path $P$ that contains right-edge of content $c(\beta')$ (that fills $\beta'$ as well), and $P$ must terminate at minor ending node at $(c(\beta') + 1, \cdot)$. Consequently, there is no path that contains right-edge of content $c(\beta') + 1$. Thus, 
    applying similar arguments used to prove the statement for $D_1$, it can be shown that the statement holds for $D_{i + 1}$, completing the induction. Then (2) can be proved for cells in diagrams $D_i, 1 < i \le k$ in the same way as it was done for diagram $D_1$, which completes the proof. 
\end{proof}

\begin{figure}
    \centering
    \begin{subfigure}[b]{0.33\textwidth}
        \centering
        \begin{tikzpicture}[scale=0.7, transform shape]
            \draw[dashed, thin] (-1, 1) grid (2, 5);
            \draw[->, red, thick] (2, 2) -- (1, 2) -- (1, 4) -- (0, 4);
            \fill (1.5, 2) node[above] {$e_{\kappa}$};
            \fill (0.5, 4) node[above] {$e_{\eta}$};
            \fill (-1, 2) node[left] {$(q + f_c)^{\circ}$};
            \fill (-1, 4) node[left] {$(p + f_{c - 1})^{\circ}$};
        \end{tikzpicture}
        \subcaption{}
        \label{fig:tableau-is-semistandard-a}
    \end{subfigure}
    \begin{subfigure}[b]{0.33\textwidth}
        \centering
        \begin{tikzpicture}[scale=0.7, transform shape]
            \draw[dashed, thin] (-1, 1) grid (2, 6);
            \draw[->, blue, thick] (2, 4) -- (1, 4) -- (1, 2) -- (0, 2);
            \draw[->, red, thick] (1, 6) -- (1, 5) -- (0, 5);
            \fill (1.5, 4) node[above] {$e_{\kappa}$};
            \fill (0.5, 5) node[above] {$e_{\eta}$};
            \fill (0.5, 2) node[above] {$e_{\zeta}$};
            \fill (-1, 2) node[left] {$(p' + f_{c - 1})^{\circ}$};
            \fill (-1, 4) node[left] {$(q + f_c)^{\circ}$};
            \fill (-1, 5) node[left] {$(p + f_{c - 1})^{\circ}$};
        \end{tikzpicture}
        \subcaption{}
        \label{fig:tableau-is-semistandard-b}
    \end{subfigure}
    \begin{subfigure}[b]{0.3\textwidth}
        \centering
        \begin{tikzpicture}[scale=0.7, transform shape]
            \draw[dashed, thin] (-1, 1) grid (2, 6);
            \draw[->, blue, thick] (2, 4) -- (1, 4) -- (1, 1);
            \draw[->, red, thick] (1, 6) -- (1, 5) -- (0, 5);
            \fill (1.5, 4) node[above] {$e_{\kappa}$};
            \fill (0.5, 5) node[above] {$e_{\eta}$};
            \fill (1, 1) circle (1.5pt) node[above right] {$B$};
            \fill (-1, 4) node[left] {$(q + f_c)^{\circ}$};
            \fill (-1, 5) node[left] {$(p + f_{c - 1})^{\circ}$};
        \end{tikzpicture}
        \subcaption{}
        \label{fig:tableau-is-semistandard-c}
    \end{subfigure}
    \caption{An example of edge locations.}
    \label{fig:}
\end{figure}
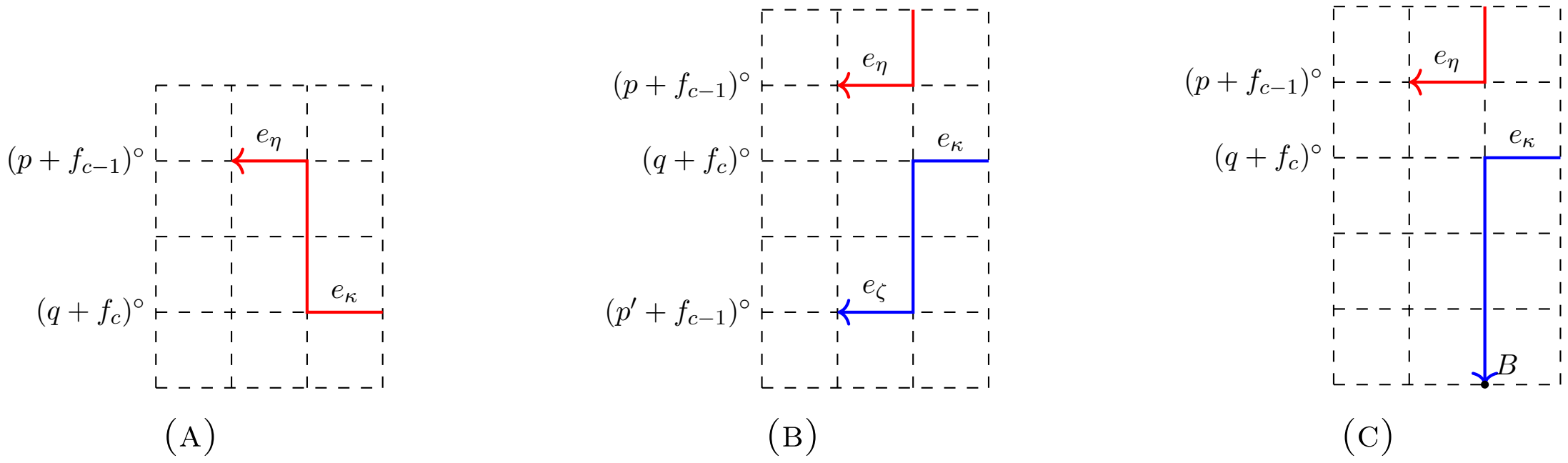

\begin{lemma}
    \label{lem:flow-is-trio-tableau}
    We have 
    
    (1) $T^{\circ}$ is left tableau; i.e. it is semistandard; and respects the row $\mathsf{G}$-flag $\mathbf{a} = (a_i)$.

    (2) $T^{\bullet}$ is right tableau; i.e. it is semistandard; and respects the row $\mathsf{G}$-flag $\mathbf{b} = (b_i)$.

    (3) $\widetilde{T}$ is central tableau; i.e. $\widetilde{T}$ is semistandard.

    Thus, $(T^{\circ}, \widetilde{T}, T^{\bullet}) \in TT_{\lambda/\mu}(\mathbf{a}, \mathbf{b}, m)$.
\end{lemma}
\begin{proof}
    We already know that $T^{\circ}, \widetilde{T}, T^{\bullet}$ are tableaux of shapes $\nu/\mu,\ \tau/\nu,\ \lambda/\tau$ (resp.) where $\nu \subseteq \mu,\ \lambda \subseteq \tau,\ \nu \subseteq \tau$.

    First we prove that $T^{\circ}$ is semistandard. Let $\eta, \kappa \in T^{\circ}$ such that $c(\eta) = c - 1,\ c(\kappa) = c$; $\eta = (i, j)$ is left neighbor of $\kappa = (i, j + 1)$; and edges which they filled by are $e_{\eta} = (c, (p + f_{c - 1})^{\circ}) \rightarrow (c - 1, (p + f_{c - 1})^{\circ})$ and $e_{\kappa} = (c + 1, (q + f_{c})^{\circ}) \rightarrow (c, (q + f_{c})^{\circ})$. If $\eta \mathrel{R} \kappa$, we have $\pi(c - 1) = (\cdot, 0),\ \pi(c) = (0, \cdot)$, thus by construction of the lattice $p + f_{c - 1} \le q + f_c$, and as $f_c = f_{c - 1}$ we have $p \le q$ (see Fig.~\ref{fig:tableau-is-semistandard-a}), thus $T^{\circ}_{\eta}\le T^{\circ}_{\kappa}$. 

    Now if $\eta \not\mathrel{R} \kappa$, then we must have $\pi(c - 1) = (\cdot, 1),\ \pi(c) = (1, \cdot)$. If there is a cell $\zeta = (i + 1, j + 1) \in T^{\circ}$ filled by some edge $e_{\zeta} = (c,  (p' + f_{c - 1})^{\circ}) \rightarrow (c - 1, (p' + f_{c - 1})^{\circ})$ then $\zeta \mathrel{R} \kappa$ must hold, and thus $p < p'$ by definition of the map $\phi$, moreover, to avoid intersection $p + f_{c - 1} < q + f_c$ must hold (see Fig.~\ref{fig:tableau-is-semistandard-b}). If there is no cell $\zeta \in T^{\circ}$, then there must be ending node $B = (c, \cdot)$ in the lower layer and the path $P$ which contains edge $e_{\kappa}$ must end on the node $B$, thus to avoid intersection, again, $p + f_{c - 1} < q + f_c$ must hold (see Fig.~\ref{fig:tableau-is-semistandard-c}), as $f_{c - 1} = f_c - 1$ we have $p \le q$, thus $T^{\circ}_{\eta}\le T^{\circ}_{\kappa}$.

    In the same way we can prove that $T^{\circ}_{\eta} < T^{\circ}_{\kappa}$ where $\eta = (i, j)$ is upper neighbor of $\kappa = (i + 1, j)$.

    Now we prove that $T^{\circ}$ respects the flag. Let $\eta$ be the rightmost cell on the $i$-row of $T^{\circ}$, and let $c(\eta) = c$. Suppose that the cell $\eta$ is filled by edge $e = (c + 1, (p + f_c)^{\circ}) \rightarrow (c, (p + f_c)^{\circ})$. As $\eta$ is the rightmost cell $\mathrm{max}(\lambda/\mu, c) = \eta$. By construction of the T-lattice the maximum possible value of $p$ is $a_{\mathrm{row}(\eta)}$, thus maximum value of $T^{\circ}_{\eta}$ is $a_{\mathrm{row}(\eta)}^{\circ}$.  

    The parts (2) and (3) can be proved similarly. 
\end{proof}

\subsubsection{Inverse map}

Let $\mathbf{T} = (T^{\circ}, \widetilde{T}, T^{\bullet}) \in TT_{\lambda/\mu}(\mathbf{a}, \mathbf{b}, m)$ of shapes $\nu/\mu,\ \tau/\nu,\ \lambda/\tau$ (resp.) where $\nu \subseteq \mu,\ \lambda \subseteq \tau,\ \nu \subseteq \tau$. Next we define the map 
$$\psi : TT_{\lambda/\mu}(\mathbf{a}, \mathbf{b}, m) \rightarrow \mathcal{F}(\mathbf{A}, \mathbf{B})$$ 
by constructing a flow $\mathbf{F} \in \mathcal{F}(\mathbf{A}, \mathbf{B})$ that corresponds to $\mathbf{T}$. We consider outer decompositions $\Theta^{\circ} = (\theta^{\circ}_i),\ \widetilde{\Theta} = (\widetilde{\theta}_i),\ \Theta^{\bullet} = (\theta^{\bullet}_i)$ of shapes $\nu/\mu,\ \tau/\nu,\ \lambda/\tau$ (resp.) induced by pipe function $\pi$. Let $\delta^{\circ}_i, \widetilde{\delta}_i, \delta_i^{\bullet}$ denote the heads of $\theta^{\circ}_i, \widetilde{\theta}_i, \theta^{\bullet}_i$ (resp.); similarly, let $\gamma^{\circ}_i, \widetilde{\gamma}_i, \gamma_i^{\bullet}$ denote their tails. Note that $\theta^{\circ}, \theta^{\bullet}$ are negative ribbons. 

In addition to starting and ending nodes $\mathbf{A}, \mathbf{B}$ for flows, we consider starting and ending nodes for three systems of paths on the T-lattice such that, each system correpond to each tableau $T^{\circ}, \widetilde{T}, T^{\bullet}$. Then we prove that systems of paths combined, yield a flow $\mathbf{F} \in \mathcal{F}(\mathbf{A}, \mathbf{B})$, and we set $\psi(\mathbf{T}) = \mathbf{F}$. 

The endpoints $\mathbf{C}^{\circ} = (C^{\circ}_i),\ \mathbf{D}^{\circ} = (D^{\circ}_i)$ for the left tableau $T^{\circ}$ are defined as follows:
\begin{align*}
	&C^{\circ}_i = 
    \begin{cases}
        (c(\delta^{\circ}_i)+1,\ \infty^{\circ}) &\text{ if } \pi(\delta^{\circ}_i) = (\cdot, 0) \\
        (c(\delta^{\circ}_i)+1,\ 1 + f_{c(\delta^{\circ}_i)} - 1) &\text{ if } \pi(\delta^{\circ}_i) = (\cdot, 1) \\
    \end{cases}\\
    &D^{\circ}_i = 
    \begin{cases}
        (c(\gamma^{\circ}_i),\ 1 + f_{c(\gamma^{\circ}_i)} - 1) &\text{ if } \pi(\gamma^{\circ}_i) = (0, \cdot) \\
        (c(\gamma^{\circ}_i),\ \infty^{\circ}) &\text{ if } \pi(\gamma^{\circ}_i) = (1, \cdot)
    \end{cases}
\end{align*}

The endpoints $\widetilde{\mathbf{C}} = (\widetilde{C}_i),\ \mathbf{\widetilde{D}} = (\widetilde{D}_i)$ for the central tableau $\widetilde{T}$ are defined as follows: 
\begin{align*}
    &\widetilde{C}_i = 
    \begin{cases}
        (c(\widetilde{\delta}_i),\ 1 + f_{c(\widetilde{\delta}_i)} - 1) &\text{ if } \pi(\widetilde{\delta}_i) = (0, \cdot) \\
        (c(\widetilde{\delta}_i),\ m + f_{c(\widetilde{\delta}_i)} + 1) &\text{ if } \pi(\widetilde{\delta}_i) = (1, \cdot) \\
    \end{cases} \\
	&\widetilde{D}_i = 
    \begin{cases}
        (c(\widetilde{\gamma}_i) + 1,\ m + f_{c(\widetilde{\gamma}_i)} + 1) &\text{ if } \pi(\widetilde{\gamma}_i) = (\cdot, 0) \\
        (c(\widetilde{\gamma}_i) + 1,\ 1 + f_{c(\widetilde{\gamma}_i)} - 1) &\text{ if } \pi(\widetilde{\gamma}_i) = (\cdot, 1)
    \end{cases}
\end{align*}

The endpoints $\mathbf{C}^{\bullet} = (C^{\bullet}_i),\ \mathbf{D}^{\bullet} = (D^{\bullet}_i)$ for the right tableau $T^{\bullet}$ are defined as follows:
\begin{align*}
    &C^{\bullet}_i = 
    \begin{cases}
        (c(\delta^{\bullet}_i)+1,\ m + f_{c(\delta^{\bullet}_i)} + 1) &\text{ if } \pi(\delta^{\bullet}_i) = (\cdot, 0) \\
        (c(\delta^{\bullet}_i)+1,\ \infty^{\bullet}) &\text{ if } \pi(\delta^{\bullet}_i) = (\cdot, 1) \\
    \end{cases} \\
    &D^{\bullet}_i = 
    \begin{cases}
        (c(\gamma^{\bullet}_i),\ \infty^{\bullet}) &\text{ if } \pi(\gamma^{\bullet}_i) = (0, \cdot) \\ 
        (c(\gamma^{\bullet}_i),\ m + f_{c(\gamma^{\bullet}_i)} + 1) &\text{ if } \pi(\gamma^{\bullet}_i) = (1, \cdot) 
    \end{cases}
\end{align*}

We call $\mathbf{A}, \mathbf{B}$ {\it major endpoints}, while the newly defined endpoints are called {\it minor}.

Define paths system $\mathbf{P} = (P_i) \in \mathcal{P}(\mathbf{C}^{\circ}, \mathbf{D}^{\circ})$ that corresponds to the left tableau $T^{\circ}$, where each $P_i \in \mathcal{P}(\{C^{\circ}_i\}, \{D^{\circ}_i\})$ is defined as follows:  $P_i$ is a path from $C^{\circ}_i$ to $D^{\circ}_i$ that consists of $r = c(\delta^{\circ}_i) - c(\gamma^{\circ}_i) + 1$ horizontal edges $e_{c(\gamma^{\circ}_i)}, \ldots, e_{c(\delta^{\circ}_i)}$ such that $e_c = (c + 1, (a + f_c)^{\circ}) \rightarrow (c, (a + f_c)^{\circ}),\ c \in [c(\gamma^{\circ}_i), c(\delta^{\circ}_i)]$, where $a$ is value in the cell of content $c$ in $T^{\circ}_{\theta^{\circ}_i}$. 

Define paths system $\mathbf{P} = (P_i) \in \mathcal{P}(\mathbf{C}^{\bullet}, \mathbf{D}^{\bullet})$ that corresponds to the right tableau $T^{\bullet}$, where each $P_i \in \mathcal{P}(\{C^{\bullet}_i\}, \{D^{\bullet}_i\})$ is defined as follows:  $P_i$ is a path from $C^{\bullet}_i$ to $D^{\bullet}_i$ that consists of $r = c(\delta^{\bullet}_i) - c(\gamma^{\bullet}_i) + 1$ horizontal edges $e_{c(\gamma^{\bullet}_i)}, \ldots, e_{c(\delta^{\bullet}_i)}$ such that $e_c = (c + 1, (a + f_c)^{\bullet}) \rightarrow (c, (a + f_c)^{\bullet}),\ c \in [c(\gamma^{\bullet}_i), c(\delta^{\bullet}_i)]$, where $a$ is value in the cell of content $c$ in $T^{\bullet}_{\theta^{\bullet}_i}$. 

Finally, define paths system $\mathbf{P} = (P_i) \in \mathcal{P}(\mathbf{\widetilde{C}}, \mathbf{\widetilde{D}})$ that corresponds to the central tableau $\widetilde{T}$, where each $P_i \in \mathcal{P}(\{\widetilde{C}_i\}, \{\widetilde{D}_i\})$ is defined as follows:  $P_i$ is a path from $\widetilde{C}_i$ to $\widetilde{D}_i$ that consists of $r = c(\widetilde{\gamma}_i) - c(\widetilde{\delta}_i) + 1$ horizontal edges $e_{c(\widetilde{\delta}_i)}, \ldots, e_{c(\widetilde{\gamma}_i)}$ such that $e_c = (c, a + f_c) \rightarrow (c + 1, a + f_c),\ c \in [c(\widetilde{\delta}_i), c(\widetilde{\gamma}_i)]$, where $a$ is value in the cell of content $c$ in $\widetilde{T}_{\widetilde{\theta}_i}$.

The paths in $\mathcal{P}(C^{\circ}_i, D^{\circ}_i), \mathcal{P}(\widetilde{C}_i, \widetilde{D}_i), \mathcal{P}(C^{\bullet}_i, D^{\bullet}_i)$ are indeed minor paths defined above. 

\begin{proposition}
    Let $\mathbf{P}^{\circ} \in \mathcal{P}(\mathbf{C}^{\circ}, \mathbf{D}^{\circ}), \mathbf{\widetilde{P}} \in \mathcal{P}(\mathbf{\widetilde{C}}, \mathbf{\widetilde{D}}), \mathbf{P}^{\bullet} \in \mathcal{P}(\mathbf{C}^{\bullet}, \mathbf{D}^{\bullet})$ be systems of paths constructed from trio-tableau $(T^{\circ}, \widetilde{T}, T^{\bullet})$ (resp.). Then the paths systems $\mathbf{P}^{\circ}, \mathbf{\widetilde{P}}, \mathbf{P}^{\bullet}$ exist on the T-lattice and they are non-intersecting. 
\end{proposition}
\begin{proof}
    We prove for $T^{\circ}$, for other tableaux it can be done similarly. Let $\theta^{\circ}_i \in \Theta^{\circ}$ and $P_i \in \mathbf{P}^{\circ}$ be the path filled by $T^{\circ}_{\theta^{\circ}_i}$ with horizontal edges $e_{c(\gamma^{\circ}_i)}, \ldots, e_{c(\delta^{\circ}_i)}$. 
    Let $\alpha$ be some cell in $\theta^{\circ}_i$ and $M = \mathrm{max}(\lambda/\mu, c(\alpha))$ then $T^{\circ}_{\alpha} \le T^{\circ}_{M} \le a_{\mathrm{row}(M)}$, thus the edge $e = (c(\alpha) + 1, (T^{\circ}_{\alpha} + f_{c(\alpha)})^{\circ}) \rightarrow (c(\alpha), (T^{\circ}_{\alpha} + f_{c(\alpha)})^{\circ})$ exists in the lattice. Then from construction of the lattice it follows that there is a path consisting of only vertical ribbons between head of $e_i$ and tail of $e_{i + 1}, c(\delta) \le i < c(\gamma)$. Thus, the path $P_i$ exists on the T-lattice. The paths in $\mathbf{P}$ are non-intersecting as $T^{\circ}$ is semistandard. 
\end{proof}




\begin{lemma}
    \label{lem:minor-is-either-major-or-another-minor}
    (1) If $C^{\circ}_i \in \mathbf{C}^{\circ}$, then either $C^{\circ}_i \in \mathbf{A}$ or $C^{\circ}_i \in \mathbf{\widetilde{D}}$; and $C^{\circ}_i$ does not belong to other sets of endpoints.  

    (2) If $D^{\circ}_i \in \mathbf{D}^{\circ}$, then either $D^{\circ}_i \in \mathbf{B}$ or $D^{\circ}_i \in \mathbf{\widetilde{C}}$; and $D^{\circ}_i$ does not belong to other sets of endpoints.  

    (3) If $C^{\bullet}_i \in \mathbf{C}^{\bullet}$, then either $C^{\bullet}_i \in \mathbf{A}$ or $C^{\bullet}_i \in \mathbf{\widetilde{D}}$; and $C^{\bullet}_i$ does not belong to other sets of endpoints.  

    (4) If $D^{\bullet}_i \in \mathbf{D}^{\bullet}$, then either $D^{\bullet}_i \in \mathbf{B}$ or $D^{\bullet}_i \in \mathbf{\widetilde{C}}$; and $D^{\bullet}_i$ does not belong to other sets of endpoints.  

    (5) If $\widetilde{C}_i \in \mathbf{\widetilde{C}}$, then either $\widetilde{C}_i \in \mathbf{D}^{\circ}$ or $\widetilde{C}_i \in \mathbf{D}^{\bullet}$; and $\widetilde{C}_i$ does not belong to other sets of endpoints.  

    (6) If $\widetilde{D}_i \in \mathbf{\widetilde{D}}$, then either $\widetilde{D}_i \in \mathbf{C}^{\circ}$ or $\widetilde{D}_i \in \mathbf{C}^{\bullet}$; and $\widetilde{D}_i$ does not belong to other sets of endpoints.  

    (7) If $A_i \in \mathbf{A}$, then either $A_i \in \mathbf{C}^{\circ}$ or $A_i \in \mathbf{C}^{\bullet}$; and $A_i$ does not belong to other sets of endpoints.  

    (8) If $B_i \in \mathbf{B}$, then either $B_i \in \mathbf{D}^{\circ}$ or $B_i \in \mathbf{D}^{\bullet}$; and $B_i$ does not belong to other sets of endpoints.  
\end{lemma}
\begin{proof}
    We prove the first statement, others can be proved similarly. Let $\delta^{\circ}_i = (r, s)$. If $\pi(\delta^{\circ}_i) = (\cdot, 0)$ then the cell $(r, s + 1)$ (or the cell $(r, s)$, if $\lambda/\mu$ is negative) is a head of some ribbon $\theta_j \in \Theta$; it is not hard to see that $a_{\mathrm{row}(\delta^{\circ}_i)} = \overline{a}_{c(\delta_j) - 1}$ (or that $a_{\mathrm{row}(\delta^{\circ})} = \overline{a}_{c(\delta_j)}$ if $\lambda/\mu$ is negative), where $\delta_j$ is head of $\theta_j$, 
    thus $C^{\circ}_i = A_j$. If $\pi(\delta^{\circ}_i) = (\cdot, 1)$, then the cell $\delta^{\circ}_i$ is also tail of some ribbon $\widetilde{\theta}_j \in \widetilde{\Theta}$, thus $C^{\circ}_i = \widetilde{D}_j$. By definition of endpoints, $C^{\circ}_i$ can not be in other sets of endpoints.
\end{proof}

\begin{lemma}
    \label{lem:minor-endpoint-at-most-two-paths}
    Any minor node can be an endpoint for at most two minor paths. 
\end{lemma}
\begin{proof}
    Follows from Lemma~\ref{lem:minor-is-either-major-or-another-minor}. 
\end{proof}

\begin{lemma}
    \label{lem:systems-of-paths-form-flow}
    Let $\mathbf{\widetilde{P}} \in \mathcal{P}(\mathbf{\widetilde{C}}, \mathbf{\widetilde{D}}),\ \mathbf{P}^{\circ} \in \mathcal{P}(\mathbf{C}^{\circ}, \mathbf{D}^{\circ}),\ \mathbf{P}^{\bullet} \in \mathcal{P}(\mathbf{C}^{\bullet}, \mathbf{D}^{\bullet})$. Combining all non-intersecting paths $\mathbf{\widetilde{P}},\ \mathbf{P}^{\circ},\ \mathbf{P}^{\bullet}$ results in a flow $\mathbf{F} \in \mathcal{F}(\mathbf{A}, \mathbf{B})$. 
\end{lemma}
\begin{proof}
    Let $C_1 \in \mathbf{\widetilde{C}} \cup \mathbf{C}^{\circ} \cup \mathbf{C}^{\bullet}$ be some minor starting node. By Lemma~\ref{lem:minor-is-either-major-or-another-minor}, any minor starting node is either a major starting node, or minor ending node. 

    First suppose that $C_1$ is a major starting node $A \in \mathbf{A}$. By Lemma~\ref{lem:minor-is-either-major-or-another-minor}, any minor ending node of some path is either a major ending node or minor starting node of another path. Let us consider paths $(P_1, P_2, \ldots)$ where $P_i \in \mathbf{\widetilde{P}} \cup \mathbf{P}^{\circ} \cup \mathbf{P}^{\bullet}$ for all $i \ge 1$; and each path $P_i$ starts at $C_i$ and ends at minor ending node $D_i$, where $D_i$ is also a minor starting node $C_{i + 1}$ of $P_{i + 1}$, for all $i \ge 1$. Thus the paths $(P_i)$ combined together results in some `bigger' path $P$ that starts at $A$. As the systems of paths are non-intersecting and by Lemma~\ref{lem:minor-endpoint-at-most-two-paths}, the path $P$ is self-avoiding. As the number of minor nodes is finite, and $P$ is not a cycle (as it starts at major endpoint $A$), the path $P$ must terminate at some ending node $B$ which is not starting; thus $B \in \mathbf{B}$. 
    
    Now suppose that $C_1$ is some minor ending node. Consider paths $(P_1, P_2, \ldots)$ where $P_i \in \mathbf{\widetilde{P}} \cup \mathbf{P}^{\circ} \cup \mathbf{P}^{\bullet}$ for all $i \ge 1$; and each $P_i$ starts at $C_i$ and ends at minor ending node $D_i$, where $D_i$ is also a minor starting node $C_{i + 1}$ of $P_{i + 1}$ for all $i \ge 1$. Thus the paths $(P_i)$ combined together results in some path $P$ that starts at $C_1$, which is also self-avoiding by Lemma~\ref{lem:minor-endpoint-at-most-two-paths} with the fact that systems of paths are non-intersecting. As the number of minor nodes is finite, we must have finite number of paths $(P_1, \ldots, P_k)$ contained in $P$. If $P$ terminates at some $D_k \in \mathbf{B}$ then we have the first case where $P$ is path from $A$ to $B$ for some $A \in \mathbf{A}$; and $C_1$ is some intermediate node in $P$. If all nodes $C_i, D_i, i \in [k]$ are not major, then $C_1 = D_k$, and $P$ must be a self-avoiding cycle. 

    Finally, let $\mathbf{P} = (P_i)$ be a system of paths (combined from paths in $\mathbf{\widetilde{P}} \cup \mathbf{P}^{\circ} \cup \mathbf{P}^{\bullet}$) such that $P_i$ is a path that starts at $A_i$ and ends at $B_{j_i}$, and $\mathbf{C}$ be a cycles collection that contain all minor nodes which are not in paths $(P_i)$. The $\mathbf{P}$ and $\mathbf{C}$ are non-intersecting and all paths in $\mathbf{P}$ and $\mathbf{C}$ are self-avoiding, thus $(\mathbf{P}, \mathbf{C})$ form a flow $\mathbf{F} \in \mathcal{F}(\mathbf{A}, \mathbf{B})$.
\end{proof}
\begin{lemma}
    \label{lem:psi-is-inverse-of-phi}
    The map $\phi$ is inverse of $\psi$, i.e. $\psi^{-1} = \phi$. Thus, $\psi$ is a bijection. 
\end{lemma}
\begin{proof}
    Let $\mathbf{F} = \psi((T^{\circ}, \widetilde{T}, T^{\bullet}))$ where $T^{\circ} \in TT^{\circ}_{\nu/\mu}(\mathbf{a})$, $T^{\bullet} \in TT^{\bullet}_{\lambda/\tau}(\mathbf{b})$, $\widetilde{T} \in TT_{\tau/\nu}(m)$. Also let $\phi(\mathbf{F}) = (Q^{\circ}, Q, Q^{\bullet})$. Let $\zeta_1, \ldots, \zeta_r \in \nu/\mu$ be cells of content $c$ sorted by the row number in descending order, i.e. $\mathrm{row}(\zeta_i) > \mathrm{row}(\zeta_{i + 1})$, let $e_1, \ldots, e_r$ be edges in $\mathbf{F}$ such that $e_i = (c + 1, (a_i + f_c)^{\circ}) \rightarrow (c, (a_i + f_c)^{\circ})$ and $e_i$ is constructed by $T^{\circ}_{\zeta_i}$, i.e. $a_i = T^{\circ}_{\zeta_i}$ for all $i$. It follows from the definition of left tableau that $a_i > a_{i + 1}$ for all $i < r$, thus by definition of the map $\phi$ we have that $Q^{\circ}_{\zeta_i} = a_i$ for all $i < r$, consequently $Q^{\circ} = T^{\circ}$.
    In the same way it can be checked that $\widetilde{T} = Q, T^{\bullet} = Q^{\bullet}$, which concludes the proof.
\end{proof}

\subsection{The bijection respects weights.}
\label{subsec:bijection-respects-weights}

\subsubsection{Preliminary facts.}
\label{subsec:prelude}

\begin{lemma}
    \label{lem:cycle-contains-only-two-types-of-edges}
    Let $C$ be a (simple) cycle on the T-lattice. Let $e_1, \ldots, e_r$ be the left edges in $C$ i.e.  $e_i = (c_i + 1, x_i) \rightarrow (c_i, x_i)$ for all $i \in [r]$. Then either $x_i \in \mathbb{Z}^{\circ}$ or $x_i \in \mathbb{Z}^{\bullet}$ for all $i \in [r]$. In other words, all edges $e_i$ are located either on the top or bottom layer. 
\end{lemma}
\begin{proof}
    As $C$ is a simple cycle, it should contain at least one left edge and one right edge. All right edges of $C$ lie in the middle layer. Let $v_1 = (p_1, q_1)$ be the node in the middle layer such that it is contained in the cycle $C$ and $p_1$ is minimal (i.e. $v_1$ is leftmost).
    
    By contradiction, suppose that some left edges are on top layer and other left edges are on bottom layer, then there are nodes $v = (p, q), v' = (p', q')$ contained in $C$ such that $q \in \mathbb{Z}^{\circ}$ and $q' \in \mathbb{Z}^{\bullet}$, there is a sub-path $P$ of $C$ from $v$ to $v_1$ that contains vertical edge $e = (p_1, q_1 - 1) \rightarrow (p_1, q_1)$, and there is path $P'$ from $v'$ to $v_1$ that contains vertical edge $e' = (p_1, q_1 + 1) \rightarrow (p_1, q_1)$, but this is not possible as directions of $e, e'$ must be the same by the construction of the lattice.
\end{proof}

\begin{lemma}
    \label{lem:cycle-corresponds-to-special-ribbons}
    Let $\mathbf{F} = (\mathbf{P}, \mathbf{C}) = \psi(\mathbf{T})$. Let $C \in \mathbf{C}$, then either of the following two cases must hold:

    (1) the cycle $C$ is constructed by $\widetilde{T}_{\widetilde{\theta}}$ and $T^{\bullet}_{\theta^{\bullet}}$ where $\widetilde{\theta} \in \widetilde{\Theta}$ and $\theta^{\bullet} \in \Theta^{\bullet}$ such that $\widetilde{\theta} = \theta^{\bullet}$ and $\pi(\widetilde{\delta}) = (1, \cdot),\ \pi(\widetilde{\gamma}) = (\cdot, 0)$.

    (2) the cycle $C$ is constructed by $\widetilde{T}_{\widetilde{\theta}}$ and $T^{\circ}_{\theta^{\circ}}$ where $\widetilde{\theta} \in \widetilde{\Theta}$ and $\theta^{\circ} \in \Theta^{\circ}$ such that $\widetilde{\theta} = \theta^{\circ}$ and $\pi(\widetilde{\delta}) = (0, \cdot),\ \pi(\widetilde{\gamma}) = (\cdot, 1)$.

\end{lemma}
\begin{proof}
    Suppose that $C$ contains right edges $e_1, \ldots, e_r$ and left edges $e'_1, \ldots, e'_r$. Let $v = (v_1, v_2), u = (u_1, u_2)$ be the leftmost ($v_1$ is minimal) and rightmost ($u_1$ is maximal) nodes contained in the right edges. Let $v' = (v'_1, v'_2), u' = (u'_1, u'_2)$ be the leftmost ($v'_1$ is minimal) and rightmost ($u'_1$ is maximal) nodes contained in left edges. It follows that $v_1 = v'_1,\ u_1 = u'_1$. As the nodes $v, u, v', u'$ are contained in the cycle, there must be sub-path $P$ of $C$ from $v$ to $u$ that does not contain left edges and sub-path $P'$ of $C$ from $u'$ to $v'$ that does not contain right edges. 

    By Lemma~\ref{lem:cycle-contains-only-two-types-of-edges} the edges $(e'_i)$ contained either in the top layer or in bottom layer. 
    
    In the first case, let $\widetilde{\theta} \in \widetilde{\Theta}$ and $\theta^{\bullet} \in \Theta^{\bullet}$ be ribbons from fillings using which the paths $P, P'$ are constructed. As contents of horizontal edges in both $P, P'$ are in the same range, it follows that $\widetilde{\theta} = \theta^{\bullet}$.
    Moreover, there must be sub-path $Q_1$ of $C$ from $v'$ to $v$ that consists of only down edges; and there must be sub-path $Q_2$ of $C$ from $u$ to $u'$ that consists of only up edges; thus $\pi(\widetilde{\delta}) = (1, \cdot),\ \pi(\widetilde{\gamma}) = (\cdot, 0)$.  
    

    In the second case, let $\widetilde{\theta} \in \widetilde{\Theta}$ and $\theta^{\circ} \in \Theta^{\circ}$ be ribbons from fillings using which the paths $P, P'$ are constructed. As contents of horizontal edges in both $P, P'$ are in the same range, it follows that $\widetilde{\theta} = \theta^{\circ}$.
    Moreover, there must be sub-path $Q_1$ of $C$ from $v'$ to $v$ that consists of only up edges; and there must be sub-path $Q_2$ of $C$ from $u$ to $u'$ that consists of only down edges, thus $\pi(\widetilde{\delta}) = (0, \cdot),\ \pi(\widetilde{\gamma}) = (\cdot, 1)$. 
\end{proof}

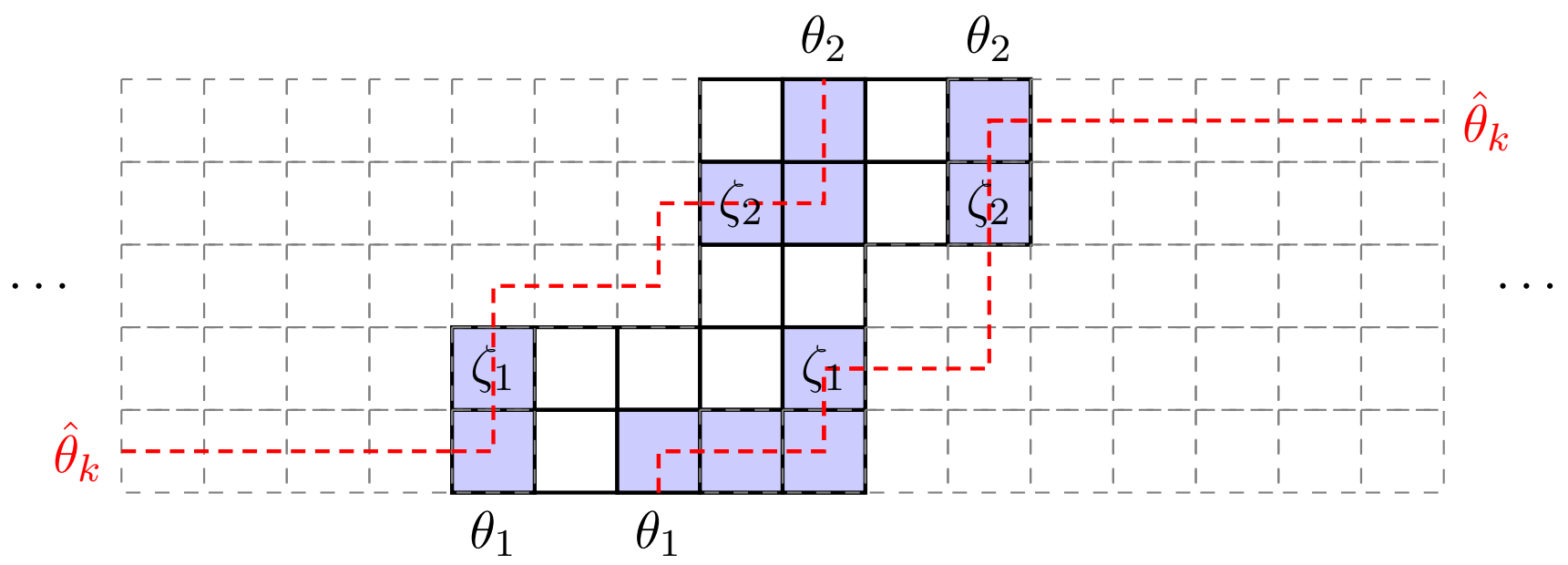
\begin{figure}
    \centering
    \begin{tikzpicture}[y=-1cm, scale=0.6]
	\foreach \x in {3,...,6} {
        \draw[thick] (\x, 1) rectangle ++(1, 1);
    }
    \foreach \x in {3,...,6} {
        \draw[thick] (\x, 2) rectangle ++(1, 1);
    }
    \foreach \x in {3,...,4} {
        \draw[thick] (\x, 3) rectangle ++(1, 1);
    }
    \foreach \x in {0,...,4} {
        \draw[thick] (\x, 4) rectangle ++(1, 1);
    }
    \foreach \x in {0,...,4} {
        \draw[thick] (\x, 5) rectangle ++(1, 1);
    }

    \draw[thick, fill=blue!20] (0, 5) rectangle ++(1, 1);
    \draw[thick, fill=blue!20] (0, 4) rectangle ++(1, 1);
    \draw[thick, fill=blue!20] (3, 2) rectangle ++(1, 1);
    \draw[thick, fill=blue!20] (4, 2) rectangle ++(1, 1);
    \draw[thick, fill=blue!20] (4, 1) rectangle ++(1, 1);

    \draw[thick, fill=blue!20] (2, 5) rectangle ++(1, 1);
    \draw[thick, fill=blue!20] (3, 5) rectangle ++(1, 1);
    \draw[thick, fill=blue!20] (4, 5) rectangle ++(1, 1);
    \draw[thick, fill=blue!20] (4, 4) rectangle ++(1, 1);
    \draw[thick, fill=blue!20] (6, 2) rectangle ++(1, 1);
    \draw[thick, fill=blue!20] (6, 1) rectangle ++(1, 1);

    \foreach \x in {-4,...,2} {
        \draw[gray, dashed] (\x, 1) rectangle ++(1, 1);
    }
    \foreach \x in {-4,...,2} {
        \draw[gray, dashed] (\x, 2) rectangle ++(1, 1);
    }
    \foreach \x in {-4,...,2} {
        \draw[gray, dashed] (\x, 3) rectangle ++(1, 1);
    }
    \foreach \x in {-4,...,0} {
        \draw[gray, dashed] (\x, 4) rectangle ++(1, 1);
    }
    \foreach \x in {-4,...,0} {
        \draw[gray, dashed] (\x, 5) rectangle ++(1, 1);
    }


    \foreach \x in {6,...,11} {
        \draw[gray, dashed] (\x, 1) rectangle ++(1, 1);
    }
    \foreach \x in {6,...,11} {
        \draw[gray, dashed] (\x, 2) rectangle ++(1, 1);
    }
    \foreach \x in {5,...,11} {
        \draw[gray, dashed] (\x, 3) rectangle ++(1, 1);
    }
    \foreach \x in {5,...,11} {
        \draw[gray, dashed] (\x, 4) rectangle ++(1, 1);
    }
    \foreach \x in {3,...,11} {
        \draw[gray, dashed] (\x, 5) rectangle ++(1, 1);
    }

    \draw [red, thick, densely dashed] (-4, 5.5) -- (0, 5.5);
    \draw [red, thick, densely dashed] (0, 5.5) -- (0.5, 5.5) -- (0.5, 4); 
    \draw [red, thick, densely dashed] (0.5, 4) -- (0.5, 3.5) -- (2.5, 3.5) -- (2.5, 2.5) -- (3, 2.5);
    \draw [red, thick, densely dashed] (3, 2.5) -- (4.5, 2.5) -- (4.5, 1);
  
    \draw [red, thick, densely dashed] (2.5, 6) -- (2.5, 5.5) -- (4.5, 5.5) -- (4.5, 4.5) -- (5, 4.5);
    \draw [red, thick, densely dashed] (6.5, 3) -- (6.5, 1.5) -- (7, 1.5);
    \draw [red, thick, densely dashed] (2.5, 6) -- (2.5, 5.5) -- (4.5, 5.5) -- (4.5, 4.5) -- (6.5, 4.5) -- (6.5, 1.5) -- (12, 1.5);

    \node at (-5, 3.5) {$\ldots$};
    \node at (13, 3.5) {$\ldots$};

    \fill (0.5, 6) node[below] {$\theta_1$};
    \fill (0.5, 4.5) node {$\zeta_1$};
    \fill (4.5, 1) node[above] {$\theta_2$};
    \fill (3.5, 2.5) node {$\zeta_2$};
    \fill[red] (-4, 5.5) node[left] {$\hat{\theta}_k$};

    \fill (2.5, 6) node[below] {$\theta_1$};
    \fill (4.5, 4.5) node {$\zeta_1$};
    \fill (6.5, 1) node[above] {$\theta_2$};
    \fill (6.5, 2.5) node {$\zeta_2$};
    \fill[red] (12, 1.5) node[right] {$\hat{\theta}_k$};

    \end{tikzpicture}
    \caption{Two examples of $\theta_1$ and $\theta_2$.}
    \label{fig:ribbons-connections-types}
\end{figure}

\begin{figure}
    \centering
    \begin{tikzpicture}[y=-1cm, scale=0.5]
    \foreach \x in {10,...,12} {
        \draw[thick] (\x, -3) rectangle ++(1, 1);
    }
    \foreach \x in {9,...,11} {
        \draw[thick] (\x, -2) rectangle ++(1, 1);
    }
    \foreach \x in {9,...,10} {
        \draw[thick] (\x, -1) rectangle ++(1, 1);
    }

	\foreach \x in {3,...,5} {
        \draw[thick] (\x, 1) rectangle ++(1, 1);
    }
    \foreach \x in {3,...,5} {
        \draw[thick] (\x, 2) rectangle ++(1, 1);
    }
    \foreach \x in {2,...,4} {
        \draw[thick] (\x, 3) rectangle ++(1, 1);
    }

    \draw[thick, fill=blue!20] (2, 3) rectangle ++(1, 1);
    \draw[thick, fill=blue!20] (3, 3) rectangle ++(1, 1);
    \draw[thick, fill=blue!20] (4, 3) rectangle ++(1, 1);
    \draw[thick, fill=blue!20] (4, 2) rectangle ++(1, 1);
    \draw[thick, fill=blue!20] (5, 2) rectangle ++(1, 1);
    \draw[thick, fill=blue!20] (9, -2) rectangle ++(1, 1);
    \draw[thick, fill=blue!20] (10, -2) rectangle ++(1, 1);
    \draw[thick, fill=blue!20] (11, -2) rectangle ++(1, 1);

    \draw [red, thick, densely dashed] (1, 3.5) -- (4.5, 3.5) -- (4.5, 2.5) -- (7.5, 2.5) -- (7.5, 0.5) -- (8.5, 0.5) -- (8.5, -1.5) -- (13, -1.5);
    \draw (6, 1) -- (9, 1) -- (9, -1);
    \draw [ultra thick] (7, 1) -- (8, 1);

    \fill (8, 1.5) node {$\theta_3$};

    \fill (3.5, 4) node[below] {$\theta_1$};
    \fill (5.5, 2.5) node {$\zeta_1$};
    \fill (11.5, -0.5) node {$\theta_2$};
    \fill (9.5, -1.5) node {$\zeta_2$};
    \fill[red] (0.5, 3.5) node {$\hat{\theta}_k$};

    \end{tikzpicture}
    \caption{Contradiction example.}
    \label{fig:contra-with-null-strip}
\end{figure}
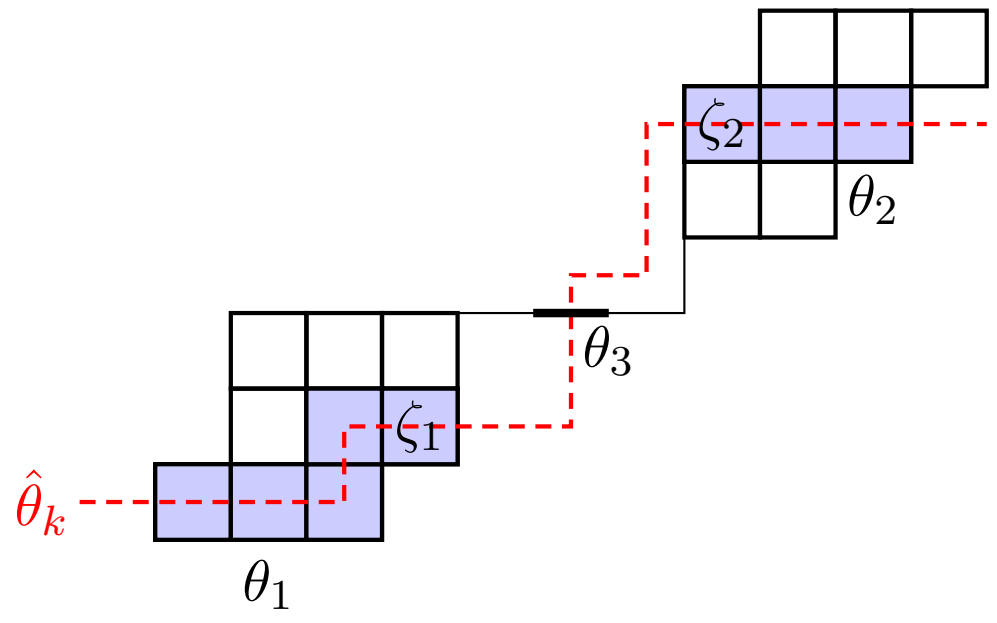

\begin{lemma}
    \label{lem:restrictions-for-ribbons}
    Let $\theta_1, \theta_2 \in \Theta$ (from given outer decomposition of $\lambda/\mu$). 
    Let $\zeta_1$ be the top-right cell of $\theta_1$ and $\zeta_2$ be the bottom-left cell of $\theta_2$. Assume that $c(\zeta_1) < c(\zeta_2)$ and they are both contained in the same tile $\hat{\theta}_k, k \in \mathbb{Z}$ such that there is no other ribbon $\theta \in \Theta$ between $\theta_1$ and $\theta_2$ in $\hat{\theta}_k$ (see Fig.~\ref{fig:ribbons-connections-types}). 
    Then we have either $\pi(\zeta_1) = (\cdot, 0),\ \pi(\zeta_2) = (1, \cdot)$ or $\pi(\zeta_1) = (\cdot, 1),\ \pi(\zeta_2) = (0, \cdot)$.
\end{lemma}
\begin{proof}
    By contradiction, suppose the statement does not hold, then we must have either $\pi(\zeta_1) = (\cdot, 0),\ \pi(\zeta_2) = (0, \cdot)$ or $\pi(\zeta_1) = (\cdot, 1),\ \pi(\zeta_2) = (1, \cdot)$. 
    In both cases the shape must be disconnected, and then it leads to contradiction, as there must be null-strip $\theta_3$ (i.e. ribbon) between $\theta_1$ and $\theta_2$ which contradicts with definition of chosen ribbons (see Fig.~\ref{fig:contra-with-null-strip}).
\end{proof}

\begin{definition}
    Let $\theta \in \Theta \cup \Theta^{\circ} \cup \widetilde{\Theta} \cup \Theta^{\bullet}$ be a ribbon with the bottom-left cell $\delta$ and the top-right cell $\gamma$. 
    We say that $\theta$ is {\it even} if $\pi(\delta) = (x, \cdot),\ \pi(\gamma) = (\cdot, x),\ x \in \{0, 1\}$; otherwise we say that $\theta$ is {\it odd}. 
\end{definition}

\begin{figure}
    \centering
    \begin{subfigure}[b]{0.5\textwidth}
        \centering
        \begin{tikzpicture}[scale=0.5]

            \fill[red] (3, 4) node {$P_{i + 1}$};
            \fill[blue] (3, -0.4) node {$P_i$};

            \draw (-1, 5) -- (7, 5);
            \draw (-1, 1) -- (7, 1);
            \draw[decorate, decoration={snake, amplitude=0.3mm, segment length=3mm}, ->, blue] (5, -1.2) -- (0, -1.2) -- (0, 1);
            \draw[decorate, decoration={snake, amplitude=0.3mm, segment length=3mm}, ->, red] (0, 1) -- (0, 3) -- (4.6, 3);
        \end{tikzpicture}
        \subcaption{}
    \end{subfigure}
    \begin{subfigure}[b]{0.45\textwidth}
        \centering
            \begin{tikzpicture}[scale=0.5]

            \fill[blue] (3, 3.85) node {$P_i$};
            \fill[red] (3, -0.5) node {$P_{i + 1}$};

            \draw (-1, 5) -- (7, 5);
            \draw (-1, 1) -- (7, 1);
            \draw[decorate, decoration={snake, amplitude=0.3mm, segment length=3mm}, ->, blue] (0, 3) -- (5, 3) -- (5, 1);
            \draw[decorate, decoration={snake, amplitude=0.3mm, segment length=3mm}, ->, red] (5, 1) -- (5, -1.2) -- (0, -1.2);
        \end{tikzpicture}
        \subcaption{}
    \end{subfigure}\\[1cm]
    \begin{subfigure}[b]{0.5\textwidth}
        \centering
        \begin{tikzpicture}[scale=0.5]

            \fill[red] (3, 3.5) node {$P_{i + 1}$};
            \fill[blue] (3, 7.5) node {$P_i$};

            \draw (-1, 5) -- (7, 5);
            \draw (-1, 1) -- (7, 1);
            \draw[decorate, decoration={snake, amplitude=0.3mm, segment length=3mm}, ->, blue] (5, 7) -- (0, 7) -- (0, 5);
            \draw[decorate, decoration={snake, amplitude=0.3mm, segment length=3mm}, ->, red] (0, 5) -- (0, 3) -- (5, 3);
        \end{tikzpicture}
        \subcaption{}
    \end{subfigure}
    \begin{subfigure}[b]{0.45\textwidth}
        \centering
            \begin{tikzpicture}[scale=0.5]

            \fill[blue] (3, 3.85) node {$P_i$};
            \fill[red] (3, 7.5) node {$P_{i + 1}$};

            \draw (-1, 5) -- (7, 5);
            \draw (-1, 1) -- (7, 1);
            \draw[decorate, decoration={snake, amplitude=0.3mm, segment length=3mm}, ->, blue] (0, 3) -- (5, 3) -- (5, 5);
            \draw[decorate, decoration={snake, amplitude=0.3mm, segment length=3mm}, ->, red] (5, 5) -- (5, 7) -- (0, 7);
        \end{tikzpicture}
        \subcaption{}
    \end{subfigure}
    \caption{Four different cases for paths.}
    \label{fig:parts-of-one-path}
\end{figure}
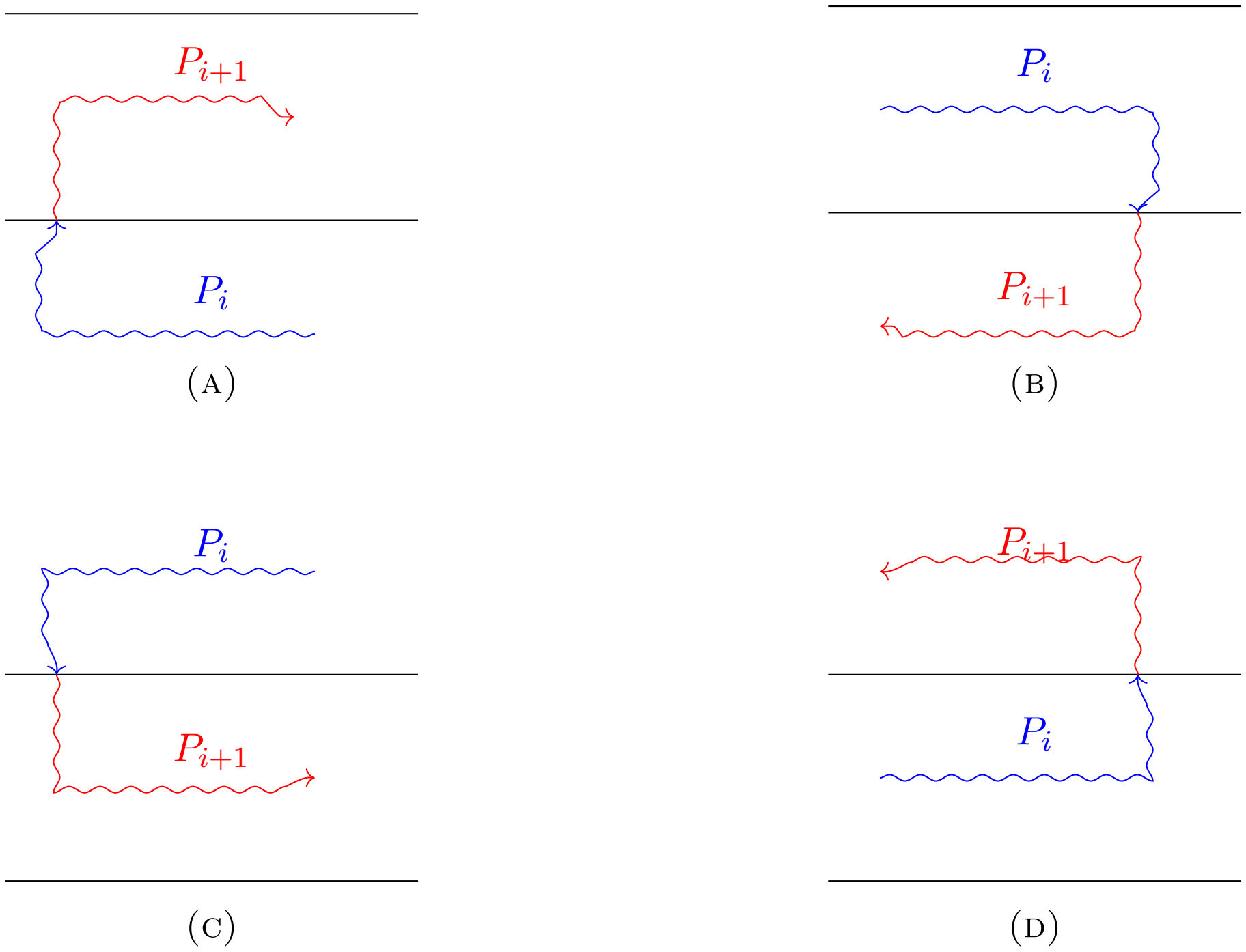

\begin{lemma}
    \label{lem:ribbons-filled-by-same-path-in-same-copy}
    Let $\mathbf{F} = (\mathbf{P}, \mathbf{C}) = \psi(\mathbf{T}) = \psi((T^{\circ}, \widetilde{T}, T^{\bullet}))$ such that $\mathbf{P} \in \mathcal{P}(\mathbf{A}, \mathbf{B}, \sigma)$ where $\sigma$ is some permutation. Let $P \in \mathcal{P}(A_i, B_{\sigma(i)})$ be some path and $\rho_1, \ldots, \rho_s \in \Theta^{\circ} \cup \Theta \cup \Theta^{\bullet}$ be ribbons which are filled by the path $P$. Then all the ribbons $\rho_1, \ldots, \rho_s$ are contained the same tile $\hat{\theta}_k, k \in \mathbb{Z}$. 
\end{lemma}
\begin{proof}
    The path $P$ consists of minor paths $P_1, \ldots, P_s$ such that each $\rho_i$ is filled by the minor path $P_i$ for $i \in [s]$; and we sort the ribbons $(\rho_i)$ so that, the ending node of $P_i$ coincides with starting node of $P_{i + 1}$. It follows that $P_i, P_{i + 1}$ must be in different layers of the lattice.
    
    Fix integer $i$ with $i < s$. We show that the ribbons $\rho_i, \rho_{i + 1}$ are both in $\hat{\theta}_k$. Let $(p_1, \ldots, p_a)$ and $(q_1, \ldots, q_b)$ are horizontal edges of $P_i$ and $P_{i + 1}$ (resp.) where $a, b \ge 0$. 
    There are four cases to consider, each case for the paths $P_i, P_{i + 1}$ being in: (a) bottom, middle; (b) middle, bottom; (c) top, middle; (d) middle, top layers (see Fig.~\ref{fig:parts-of-one-path}). We work out the first case, and other cases can be proved in similarly. In the first case, $\rho_i \in \Theta^{\circ},\ \rho_{i + 1} \in \widetilde{\Theta}$. Let $\gamma^{\circ}$ be tail of $\rho_i$, and $\widetilde{\delta}$ be head of $\rho_{i + 1}$. It follows that $c(\gamma^{\circ}) = c(\widetilde{\delta})$.
    
    Note that $\gamma^{\circ}$ (or left border of $\gamma^{\circ}$ if $a = 0$) is located on the left perimeter of the $\nu/\mu$, thus $\gamma^{\circ}$ (or left border of $\gamma^{\circ}$ if $b = 0$) is located on the left perimeter of $\tau/\nu$; and $\gamma^{\circ}$ is head of some ribbon in $\widetilde{\Theta}$. As there are no two ribbons with the same head content in $\widetilde{\Theta}$, it follows that $\gamma^{\circ} = \widetilde{\delta}$ which implies $\rho_i, \rho_{i + 1}$ are in the same $\hat{\theta}_k$; and the statement of the lemma follows.
\end{proof}

\begin{lemma}
    \label{lem:th-i-th-si-are-in-the-same-ribbon}
    Let $\mathbf{F} = (\mathbf{P}, \mathbf{C}) = \psi(\mathbf{T}) = \psi((T^{\circ}, \widetilde{T}, T^{\bullet}))$ such that $\mathbf{P} \in \mathcal{P}(\mathbf{A}, \mathbf{B}, \sigma)$ where $\sigma$ is some permutation. Let $\rho_1, \ldots, \rho_r \in \Theta^{\circ} \cup \widetilde{\Theta} \cup \Theta^{\bullet}$ be ribbons filled by path $P \in \mathbf{P}$ where $P \in \mathcal{P}(A_i, B_{\sigma(i)})$.
    Then the ribbons $\theta_i, \theta_{\sigma(i)}$ and all the ribbons $(\rho_i)$ are contained in the same tile $\hat{\theta}_k, k \in \mathbb{Z}$. 
    
\end{lemma}

\begin{proof}
    From Lemma~\ref{lem:ribbons-filled-by-same-path-in-same-copy} it follows that ribbons $(\rho_j)$ are sub-ribbons of $\hat{\theta}_k$. Let us specify the coordinates of $A_i, B_j$: $A_i = (c_a, y_a), B_j = (c_b, y_b)$. 
    
    There are two cases to consider: (1) $\lambda/\mu$ is positive; and (2) $\lambda/\mu$ is negative. 
    We show the first case, and the second case can be handled similarly. 

    Let $\rho_s, s \in [r]$ be negative ribbon (i.e. $\rho_s \in \Theta^{\circ} \cup \Theta^{\bullet}$)  filled with the minor path $P'$ contained in $P$ which starts at $A_i$; thus its head content is $c_a - 1$; As head content of $\theta_i$ is $c_a$, it follows that $\rho_s$
    is contained in the same tile 
    with $\theta_i$.

    
    Also there must be negative ribbon $\rho_t, t \in [r]$ with tail content $c_b + 1$ that is contained in the same tile with the ribbon $\theta_{\sigma(i)}$ (as tail content of $\theta_{\sigma(i)}$ is $c_b$). As $\rho_s, \rho_t$ both contained in $\hat{\theta}_k$, it follows that $\theta_i, \theta_{\sigma(i)}$ are contained in $\hat{\theta}_k$ too. 
\end{proof}

\begin{definition}
    Let $P_1, P_2$ be minor paths on the lattice. We say that they \textit{overlap} if there are nodes $u_1 \in P_1, u_2 \in P_2$ in the same layer such that their first coordinates coincide, i.e. if $u_1 = (x, y_1),\ u_2 = (x, y_2),\ y_1, y_2 \in \mathbb{Z}^{\circ} \cup \mathbb{Z} \cup \mathbb{Z}^{\bullet}$.
\end{definition}

\begin{lemma}
    \label{lem:cousin-paths-avoid-each-others-endpoints}
    Let $\mathbf{F} = (\mathbf{P}, \mathbf{C}) = \psi(\mathbf{T}) = \psi((T^{\circ}, \widetilde{T}, T^{\bullet}))$ such that $\mathbf{P} \in \mathcal{P}(\mathbf{A}, \mathbf{B}, \sigma)$ where $\sigma$ is some permutation. Assume that $\theta_r, \theta_s \in \Theta$ are in the same tile. Let $P'_r, P'_s$ be the minor paths contained in the paths $P_r, P_s \in \mathbf{P}$ (resp.) where $P_i \in \mathcal{P}(A_i, B_{\sigma(i)}),\ i \in \{r, s\}$, such that $P'_r, P'_s$ are in the same layer. Then, $P'_r, P'_s$ do not overlap.
\end{lemma}
\begin{proof}
    Assume on the contrary that paths $P'_r, P'_s$ overlap, i.e. $P'_r$ contains node $u_1 = (x, y_1)$, $P'_s$ contains $u_2 = (x, y_2)$ for some $y_1, y_2 \in \mathbb{Z}^{\circ} \cup \mathbb{Z} \cup \mathbb{Z}^{\bullet}$. Thus there are ribbons $\rho_1, \rho_2$ such that both of them belong to $X$ where $X$ is either $\Theta^{\circ}$ or $\widetilde{\Theta}$ or $\Theta^{\bullet}$, and $\rho_1, \rho_2$ are filled by the paths $P'_r, P'_s$. Note that by Lemma~\ref{lem:ribbons-filled-by-same-path-in-same-copy} $\rho_1, \rho_2$ are contained in the same tile $\hat{\theta}_k, k \in \mathbb{Z}$. As $P'_r, P'_s$ overlap, it follows that $\rho_1, \rho_2$ share some common cells (see Fig.~\ref{fig:cousin-paths-avoid-each-others-endpoints-1}), or $\rho_1, \rho_2$ are placed consecutively in $\hat{\theta}_k$ such that that there are no cells between $\rho_1, \rho_2$ (see Fig.~\ref{fig:cousin-paths-avoid-each-others-endpoints-2}); this contradicts the fact that $X$ is an outer decomposition. 
\end{proof}

\begin{figure}
    \centering
    \begin{subfigure}[b]{0.5\textwidth}
        \centering
        \begin{tikzpicture}[scale=0.5]
            \fill (3, 4) circle (2pt) node[above] {$u_2$};
            \fill (3, 2) circle (2pt) node[below] {$u_1$};

            \fill[blue] (1.5, 2.5) node {$P_r$};
            \fill[red] (4.5, 4.5) node {$P_s$};

            \draw (-1, 5) -- (7, 5);
            \draw (-1, 1) -- (7, 1);

            \draw[decorate, decoration={snake, amplitude=0.3mm, segment length=3mm}, ->, red] (2, 4) -- (6, 4);

            \draw[decorate, decoration={snake, amplitude=0.3mm, segment length=3mm}, ->, blue] (0, 2) -- (4, 2);
        \end{tikzpicture}
        \subcaption{$P'_r, P'_s$ overlap at more than one node.}
    \end{subfigure}
    \begin{subfigure}[b]{0.45\textwidth}
        \centering
        \begin{tikzpicture}[y=-1cm, scale=0.3]
            \foreach \x in {7,...,9} {
                \draw[thick] (\x, 1) rectangle ++(1, 1);
            }
            \foreach \x in {7,...,7} {
                \draw[thick, fill=red!20] (\x, 1) rectangle ++(1, 1);
            }
            \foreach \x in {6,...,7} {
                \draw[thick, fill=red!20] (\x, 2) rectangle ++(1, 1);
            }
            \foreach \x in {6,...,6} {
                \draw[thick, fill=blue!50!red] (\x, 3) rectangle ++(1, 1);
            }
            \foreach \x in {6,...,6} {
                \draw[thick, fill=blue!50!red] (\x, 4) rectangle ++(1, 1);
            }
            \foreach \x in {3,...,6} {
                \draw[thick, fill=blue!20] (\x, 5) rectangle ++(1, 1);
            }
            \foreach \x in {5,...,6} {
                \draw[thick, fill=blue!50!red] (\x, 5) rectangle ++(1, 1);
            }
            \foreach \x in {2,...,3} {
                \draw[thick, fill=blue!20] (\x, 6) rectangle ++(1, 1);
            }
            \foreach \x in {1,...,2} {
                \draw[thick] (\x, 7) rectangle ++(1, 1);
            }

            \foreach \x in {7,...,8} {
                \draw[thick] (1, \x) -- (0.7, \x);
            }
            \foreach \x in {1,...,2} {
                \draw[thick] (10, \x) -- (10.3, \x);
            }


            \fill (0, 7.5) node {$\hat{\theta}_k$};

        \end{tikzpicture}
        \subcaption{Intersection of $\rho_1, \rho_2$. The blue cells belong only to $\rho_1$, red cells belong only to $\rho_2$, purple cells belong to both.}
    \end{subfigure}
    \caption{Paths and ribbon example.}
    \label{fig:cousin-paths-avoid-each-others-endpoints-1}
\end{figure}
\begin{figure}
    \centering
    \begin{subfigure}[b]{0.5\textwidth}
        \centering
        \begin{tikzpicture}[scale=0.5]
            \fill (3, 4) circle (2pt) node[left] {$u_2$};
            \fill (3, 2) circle (2pt) node[right] {$u_1$};

            \fill[blue] (1.5, 2.5) node {$P_r$};
            \fill[red] (4.5, 4.5) node {$P_s$};

            \draw (-1, 5) -- (7, 5);
            \draw (-1, 1) -- (7, 1);

            \draw[decorate, decoration={snake, amplitude=0.3mm, segment length=3mm}, ->, red] (3, 5) -- (3, 4) -- (6, 4);

            \draw[decorate, decoration={snake, amplitude=0.3mm, segment length=3mm}, ->, blue] (0, 2) -- (3, 2) -- (3, 1);
        \end{tikzpicture}
        \subcaption{$P'_r, P'_s$ overlap at exactly one node.}
    \end{subfigure}
    \begin{subfigure}[b]{0.45\textwidth}
        \centering
        \begin{tikzpicture}[y=-1cm, scale=0.3]
            \foreach \x in {7,...,9} {
                \draw[thick] (\x, 1) rectangle ++(1, 1);
            }
            \foreach \x in {7,...,7} {
                \draw[thick, fill=red!20] (\x, 1) rectangle ++(1, 1);
            }
            \foreach \x in {6,...,7} {
                \draw[thick, fill=red!20] (\x, 2) rectangle ++(1, 1);
            }
            \foreach \x in {6,...,6} {
                \draw[thick, fill=red!20] (\x, 3) rectangle ++(1, 1);
            }
            \foreach \x in {6,...,6} {
                \draw[thick, fill=red!20] (\x, 4) rectangle ++(1, 1);
            }
            \foreach \x in {3,...,6} {
                \draw[thick, fill=blue!20] (\x, 5) rectangle ++(1, 1);
            }
            \foreach \x in {2,...,3} {
                \draw[thick, fill=blue!20] (\x, 6) rectangle ++(1, 1);
            }
            \foreach \x in {1,...,2} {
                \draw[thick] (\x, 7) rectangle ++(1, 1);
            }

            \foreach \x in {7,...,8} {
                \draw[thick] (1, \x) -- (0.7, \x);
            }
            \foreach \x in {1,...,2} {
                \draw[thick] (10, \x) -- (10.3, \x);
            }


            \fill (0, 7.5) node {$\hat{\theta}_k$};

        \end{tikzpicture}
        \subcaption{Ribbons $\rho_1, \rho_2$ placed consecutively. The blue cells belong to $\rho_1$, red cells belong to $\rho_2$.}
    \end{subfigure}
    \caption{Paths and ribbon example.}
    \label{fig:cousin-paths-avoid-each-others-endpoints-2}
\end{figure}

\begin{figure}
    \centering
    \begin{subfigure}[b]{0.5\textwidth}
        \centering
        \begin{tikzpicture}[scale=0.5, font=\scriptsize]

            \fill[blue] (-1.75, -1.25) node {$P_1$};

            \fill (0, -0.75) circle (2pt) node[below right] {$A_i$};
            \fill (0, -0.75) circle (2pt) node[above right] {$C^{\circ}_1$};
            \fill (-1.5, 1) circle (2pt) node[above left] {$\widetilde{C}_1$};
            \fill (-1.5, 1) circle (2pt) node[above right] {$D^{\circ}_1$};
            \fill (5.5, -0.5) circle (2pt) node[below left] {$B_{\sigma(i)}$};

            \draw (-2, 5) -- (8, 5);
            \draw (-2, 1) -- (8, 1);

            \draw[decorate, decoration={snake, amplitude=0.3mm, segment length=3mm}, ->, blue] (0, -0.75) -- (-1.5, -0.75) -- (-1.5, 1);
 
        \end{tikzpicture}
        \subcaption{}
        \label{fig:weight-of-ribbons-filled-by-same-path-positive-1a}
    \end{subfigure}
    \begin{subfigure}[b]{0.45\textwidth}
        \centering
        \begin{tikzpicture}[scale=0.5, font=\scriptsize]

            \fill[blue] (-1.75, -1.25) node {$P_1$};
            \fill[red] (5, 1.75) node {$P_2$};
            \fill[cyan] (5, 7) node {$P_3$};
            \fill[brown] (1.5, 3) node {$P'$};

            \fill (0, -0.75) circle (2pt) node[below right] {$A_i$};
            \fill (0, -0.75) circle (2pt) node[above right] {$C^{\circ}_1$};

            \fill (-1.5, 1) circle (2pt) node[above left] {$\widetilde{C}_1$};
            \fill (-1.5, 1) circle (2pt) node[above right] {$D^{\circ}_1$};

            \fill (5, 5) circle (2pt) node[above left] {$\widetilde{D}_1$};
            \fill (5, 5) circle (2pt) node[above right] {$C^{\bullet}_1$};

            \fill (2, 5) circle (2pt) node[below left] {$D^{\bullet}_1$};
            \fill (2, 5) circle (2pt) node[below right] {$\widetilde{C}_2$};

            \fill (5.5, -0.5) circle (2pt) node[below left] {$B_{\sigma(i)}$};

            \draw (-2, 5) -- (8, 5);
            \draw (-2, 1) -- (8, 1);

            \draw[decorate, decoration={snake, amplitude=0.3mm, segment length=3mm}, ->, blue] (0, -0.75) -- (-1.5, -0.75) -- (-1.5, 1);

            \draw[decorate, decoration={snake, amplitude=0.3mm, segment length=3mm}, ->, red] (-1.5, 1) -- (-1.5, 2.25) -- (5, 2.25) -- (5, 5);

            \draw[decorate, decoration={snake, amplitude=0.3mm, segment length=3mm}, ->, cyan] (5, 5) -- (5, 6.25) -- (2, 6.25) -- (2, 5);

            \draw[decorate, decoration={snake, amplitude=0.3mm, segment length=3mm}, ->, brown] (2, 5) -- (2, 3) -- (4, 3);
        \end{tikzpicture}
        \subcaption{}
        \label{fig:weight-of-ribbons-filled-by-same-path-positive-1b}
    \end{subfigure} \\[1cm]
    \begin{subfigure}[b]{0.5\textwidth}
        \centering
        \begin{tikzpicture}[scale=0.5, font=\scriptsize]

            \fill[blue] (-1.75, -1.25) node {$P_1$};
            \fill[red] (8.25, 3) node {$P_2$};
            \fill[cyan] (8.25, -1) node {$P_3$};

            \fill (0, -0.75) circle (2pt) node[below right] {$A_i$};
            \fill (0, -0.75) circle (2pt) node[above right] {$C^{\circ}_1$};

            \fill (-1.5, 1) circle (2pt) node[above left] {$\widetilde{C}_1$};
            \fill (-1.5, 1) circle (2pt) node[above right] {$D^{\circ}_1$};

            \fill (7.5, 1) circle (2pt) node[above left] {$\widetilde{D}_1$};
            \fill (7.5, 1) circle (2pt) node[above right] {$C^{\circ}_2$};

            \fill (5.5, -0.5) circle (2pt) node[above left] {$D^{\circ}_2$};
            \fill (5.5, -0.5) circle (2pt) node[below left] {$B_{\sigma(i)}$};

            \draw (-2, 5) -- (8, 5);
            \draw (-2, 1) -- (8, 1);

            \draw[decorate, decoration={snake, amplitude=0.3mm, segment length=3mm}, ->, blue] (0, -0.75) -- (-1.5, -0.75) -- (-1.5, 1);

            \draw[decorate, decoration={snake, amplitude=0.3mm, segment length=3mm}, ->, red] (-1.5, 1) -- (-1.5, 3) -- (5, 3) -- (7.5, 3) -- (7.5, 1);

            \draw[decorate, decoration={snake, amplitude=0.3mm, segment length=3mm}, ->, cyan] (7.5, 1) -- (7.5, -0.5) -- (5.5, -0.5);
        \end{tikzpicture}
        \subcaption{}
        \label{fig:weight-of-ribbons-filled-by-same-path-positive-1c}
    \end{subfigure}
    \caption{Paths examples.}
    \label{fig:weight-of-ribbons-filled-by-same-path-positive-1}
\end{figure}

\subsubsection{Weight-preserving property for ribbons.}
\label{subsec:weigth-preserving-for-ribbons}

For convenience we introduce the following definition.  
\begin{definition}
    Let $D$ be some diagram, denote the sign 
    \[
    \mathrm{sgn}_{\pi}(D) := (-1)^{\mathsf{hp}^{\pi}_{D}}.
    \]
    Similarly for shape $\lambda/\mu$, define the sign 
    \[
    \mathrm{sgn}_{\pi}(\lambda/\mu) := (-1)^{\mathsf{hp}^{\pi}_{\lambda/\mu}}. 
    \]
\end{definition}

\begin{lemma}
    \label{lem:weight-of-ribbons-filled-by-same-path-positive}
    Let $\mathbf{F} = (\mathbf{P}, \mathbf{C}) = \psi(\mathbf{T}) = \psi((T^{\circ}, \widetilde{T}, T^{\bullet}))$ such that $\mathbf{P} \in \mathcal{P}(\mathbf{A}, \mathbf{B}, \sigma)$ where $\sigma$ is some permutation, and recall that $\Theta = (\theta_i)$ is the given outer decomposition. 
    Assume that $\theta_i \# \theta_{\sigma(i)}$ is positive. Let $P \in \mathcal{P}(A_i, B_{\sigma(i)})$ be 
    some path and $\rho^{\circ}_1, \ldots, \rho^{\circ}_{s_1} \in \Theta^{\circ}$, $\rho^{\bullet}_1, \ldots, \rho^{\bullet}_{s_2} \in \Theta^{\bullet}$, $\widetilde{\rho}_1, \ldots, \widetilde{\rho}_{s_3} \in \widetilde{\Theta},\ s_1, s_2, s_3 \ge 0$ be ribbons filled by the path $P$. Then we have only one ribbon $\widetilde{\rho}_1$ filled by $P$ that belongs to $\widetilde{\Theta}$ (i.e. $s_3 = 1$); and $\theta_i \# \theta_{\sigma(i)}, \rho^{\circ}_r, \rho^{\bullet}_s$ are sub-ribbons of $\widetilde{\rho}_1$ for all $r \in [s_1], s \in [s_2]$. Moreover, the following identity holds: 
    \begin{equation*}
        \mathrm{wt}(P) = \mathrm{sgn}_{\pi}(\theta_i \# \theta_{\sigma(i)}) \prod^{s_1}_{r = 1} \mathrm{wt}(T^{\circ}_{\rho^{\circ}_r}) \prod^{s_2}_{r = 1} \mathrm{wt}(T^{\bullet}_{\rho^{\bullet}_r}) \prod^{s_3}_{r = 1} \mathrm{wt}(\widetilde{T}_{\widetilde{\rho}_{\sigma(i)}}).
    \end{equation*}
\end{lemma}
\begin{proof}
    It follows from the construction of the lattice that 
    \begin{equation*}
        |\mathrm{wt}(P)| = \prod^{s_1}_r \mathrm{wt}(T^{\circ}_{\rho^{\circ}_r}) \prod^{s_2}_r \mathrm{wt}(T^{\bullet}_{\rho^{\bullet}_r}) \prod^{s_3}_r \mathrm{wt}(\widetilde{T}_{\widetilde{\rho}_{\sigma(i)}}).
    \end{equation*}
    Define the sign $\alpha(P)$ as follows: 
    \begin{equation*}
        \mathrm{wt}(P) = \alpha(P) \cdot |\mathrm{wt}(P)|.
    \end{equation*}
    Note that $\mathrm{sgn}_{\pi}(\theta_i \# \theta_{\sigma(i)}) = 1$, and $\alpha(P) = (-1)^N$ where $N$ is the number of cells with hook pipes in all the ribbons $(\rho^{\circ}_j), (\rho^{\bullet}_j)$. We prove that $\alpha(P) = 1$.
    
    There are four cases induced by the placement of endpoints $A_i, B_{\sigma(i)}$ to the top and bottom layers. We consider the case in which both $A_i, B_{\sigma(i)}$ are in bottom layer; and other cases can be proved similarly. 

    From Lemma~\ref{lem:ribbons-filled-by-same-path-in-same-copy} it follows that all the ribbons $(\rho^{\circ}_i), (\widetilde{\rho}_i), (\rho^{\bullet}_i)$ are in the same tile $\hat{\theta}_k, k \in \mathbb{Z}$. 
    
    Recall that the ribbons in $\Theta^{\circ} \cup \widetilde{\Theta} \cup \Theta^{\bullet}$ are filled by the minor paths. Denote the minor starting and ending node of the path that fills the ribbon $\rho^{\circ}_r$ as $C^{\circ}_r, D^{\circ}_r,\ r \in [s_1]$. In the same way define the endpoints $C^{\bullet}_r, D^{\bullet}_r$ for the ribbon $\rho^{\bullet}_r,\ r \in [s_2]$; and $\widetilde{C}_r, \widetilde{D}_r$ for the ribbon $\widetilde{\rho}_r,\ r \in [s_3]$.
    
    As $A_i, B_{\sigma(i)}$ are in the bottom layer, there must be minor endpoint $D^{\circ}_1 = (x_1, \cdot)$ such that the minor path $P_1$ contained in $P$ which starts at $A_i = C^{\circ}_1$ terminates at $D^{\circ}_1$; thus $x_1 < A^1_i$. The path $P_1$ fills some ribbon $\rho^{\circ}_1$ from $\Theta^{\circ}$ with $\pi(\delta^{\circ}_1) = (\cdot, 0),\ \pi(\gamma^{\circ}_1) = (0, \cdot)$ (Fig.~\ref{fig:weight-of-ribbons-filled-by-same-path-positive-1a}). Then the next minor path $P_2$ contained in $P$ which stars at $\widetilde{C}_1 = D^{\circ}_1$ terminates at $\widetilde{D}_1 = (x_2, \cdot)$, and it fills the ribbon $\widetilde{\rho}_1$. Suppose that $\widetilde{D}_1 = C^{\bullet}_1$, then we will have minor path $P_3$ which starts at $C^{\bullet}_1$ and terminates at $D^{\bullet}_1$. Then the ending minor node also must be some minor starting node $\widetilde{C}_2 = (x_3, \cdot)$ with $x_3 < x_2$; it follows that the sub-path $P'$ of $P$ which starts at $\widetilde{C}_2$ and terminates at $B_{\sigma(i)}$ contains minor path in the middle layer that overlaps with $P_2$, which is a contradiction by Lemma~\ref{lem:cousin-paths-avoid-each-others-endpoints}. 
    Thus, $\widetilde{D}_1 \neq C^{\bullet}_1$. We have the same contradiction if $\widetilde{D}_1 = C^{\circ}_2$ and $x_2 < B^1_{\sigma(i)}$, thus $x_2 \ge B^1_{\sigma(i)}$ must hold. Now suppose that the part $P_3$ that starts at $C^{\circ}_2$ and terminates at $D^{\circ}_2$. If $D^{\circ}_2 = \widetilde{C}_2$ then we have the above contradiction, thus the only choice for $D^{\circ}_2$ is $B_{\sigma(i)}$; hence the part $P_3$ fills some ribbon $\rho^{\circ}_2$ with $\pi(\delta^{\circ}_2) = (\cdot, 1),\ \pi(\gamma^{\circ}_2) = (1, \cdot)$. The only ribbons filled by path $P$ are $\rho^{\circ}_1, \rho^{\circ}_2, \widetilde{\rho}_1$. 
    
    Note that the contents of cells in $\widetilde{\rho}_1$ form the range $[x_1, x_2 - 1]$. Also the contents of cells in $\rho^{\circ}_1, \rho^{\circ}_2$ are in $[x_1, A^1_i - 1],\ [B^1_{\sigma(i)}, x_2 - 1]$; thus $\rho^{\circ}_1, \rho^{\circ}_2$ are sub-ribbons of $\widetilde{\rho}_1$. Similarly, $\theta_i \# \theta_{\sigma(i)}$ is sub-ribbon of $\widetilde{\rho}_1$ as the contents of cells of $\theta_i \# \theta_{\sigma(i)}$ are in range $[A^1_i, B^1_{\sigma(i)}]$.
    
    As the ribbons $\rho^{\circ}_1, \rho^{\circ}_2$ are even, the overall number of hook pipes in them is even, and the result follows.
\end{proof}

\begin{lemma}
    \label{lem:weight-of-ribbons-filled-by-same-path-negative}
    Let $\mathbf{F} = (\mathbf{P}, \mathbf{C}) = \psi(\mathbf{T}) = \psi((T^{\circ}, \widetilde{T}, T^{\bullet}))$ such that $\mathbf{P} \in \mathcal{P}(\mathbf{A}, \mathbf{B}, \sigma)$ where $\sigma$ is some permutation.
    Assume that $\theta_i \# \theta_{\sigma(i)}$ is negative. Let $P \in \mathcal{P}(A_i, B_{\sigma(i)})$ be some path and $\rho^{\circ}_1, \ldots, \rho^{\circ}_{s_1} \in \Theta^{\circ}$, $\rho^{\bullet}_1, \ldots, \rho^{\bullet}_{s_2} \in \Theta^{\bullet}$, $\widetilde{\rho}_1, \ldots, \widetilde{\rho}_{s_3} \in \widetilde{\Theta},\ s_1, s_2, s_3 \ge 0$ be ribbons filled by the path $P$. Then the following identity holds: 
    \begin{equation*}
        \mathrm{wt}(P) = \mathrm{sgn}_{\pi}(\theta_i \# \theta_{\sigma(i)}) \prod^{s_1}_{r = 1} \mathrm{wt}(T^{\circ}_{\rho^{\circ}_r}) \prod^{s_2}_{r = 1} \mathrm{wt}(T^{\bullet}_{\rho^{\bullet}_r}) \prod^{s_3}_{r =1} \mathrm{wt}(\widetilde{T}_{\widetilde{\rho}_r}).
    \end{equation*}
\end{lemma}

\begin{figure}
    \centering
    \begin{subfigure}[b]{0.5\textwidth}
        \centering
        \begin{tikzpicture}[scale=0.5, font=\scriptsize]

            \fill[red] (3, -1.5) node {$P^{\circ}_1$};

            \fill (7, -0.75) circle (2pt) node[below right] {$A_i$};
            \fill (7, -0.75) circle (2pt) node[above right] {$C^{\circ}_1$};
            \fill (-1.5, -0.75) circle (2pt) node[above left] {$D^{\circ}_1$};
            \fill (-1.5, -0.75) circle (2pt) node[below left] {$B_{\sigma(i)}$};

            \draw (-2, 5) -- (8, 5);
            \draw (-2, 1) -- (8, 1);

            \draw[decorate, decoration={snake, amplitude=0.3mm, segment length=3mm}, ->, red] (7, -0.75) -- (-1.5, -0.75);
 
        \end{tikzpicture}
        \subcaption{}
        \label{fig:weight-of-ribbons-filled-by-same-path-negative-a}
    \end{subfigure}
    \begin{subfigure}[b]{0.49\textwidth}
        \centering
        \begin{tikzpicture}[scale=0.5, font=\scriptsize]

            \fill[red] (4, -1.25) node {$P^{\circ}_1$};

            \fill (7, -0.75) circle (2pt) node[below right] {$A_i$};
            \fill (7, -0.75) circle (2pt) node[above right] {$C^{\circ}_1$};
            \fill (4, 1) circle (2pt) node[above left] {$\widetilde{C}_1$};
            \fill (4, 1) circle (2pt) node[above right] {$D^{\circ}_1$};
            \fill (-1.5, -0.5) circle (2pt) node[below left] {$B_{\sigma(i)}$};

            \draw (-2, 5) -- (8, 5);
            \draw (-2, 1) -- (8, 1);

            \draw[decorate, decoration={snake, amplitude=0.3mm, segment length=3mm}, ->, red] (7, -0.75) -- (4, -0.75) -- (4, 1);
 
        \end{tikzpicture}
        \subcaption{}
        \label{fig:weight-of-ribbons-filled-by-same-path-negative-b}
    \end{subfigure}\\[1cm]
    \begin{subfigure}[b]{0.5\textwidth}
        \centering
        \begin{tikzpicture}[scale=0.5, font=\scriptsize]

            \fill[black] (4, -1.25) node {$P^{\circ}_1$};

            \fill (7, -0.75) circle (2pt) node[below right] {$A_i$};
            \fill (4, 1) circle (2pt) node[above left] {$\widetilde{C}_1$};
            \fill (7.5, 1) circle (2pt) node[above right] {$\widetilde{D}_1$};
            \fill (6, 1) circle (2pt) node[above] {$\widetilde{C}_2$};
            \fill (-1.5, -0.5) circle (2pt) node[below left] {$B_{\sigma(i)}$};

            \draw (-2, 5) -- (8, 5);
            \draw (-2, 1) -- (8, 1);

            \draw[decorate, decoration={snake, amplitude=0.3mm, segment length=3mm}, ->, black] (7, -0.75) -- (4, -0.75) -- (4, 1);
            \draw[decorate, decoration={snake, amplitude=0.3mm, segment length=3mm}, densely dotted, thick, ->, red] (4, 1) -- (4, 3) -- (7.5, 3) -- (7.5, 1); 
            \draw[decorate, decoration={snake, amplitude=0.3mm, segment length=3mm}, densely dotted, thick, ->, red] (7.5, 1) -- (7.5, 0.25) -- (6, 0.25) -- (6, 1);
 
        \end{tikzpicture}
        \subcaption{}
        \label{fig:weight-of-ribbons-filled-by-same-path-negative-c}
    \end{subfigure}
    \begin{subfigure}[b]{0.49\textwidth}
    \centering
    \begin{tikzpicture}[scale=0.5, font=\scriptsize]

            \fill[black] (4, -1.25) node {$P^{\circ}_1$};
            \fill[red] (6.5, 3.5) node {$\widetilde{P}_1$};
            \fill[red] (3.75, 6.5) node {$P^{\bullet}_1$};
            \fill[red] (1.5, 3.5) node {$\widetilde{P}_2$};

            \fill (7, -0.75) circle (2pt) node[below right] {$A_i$};
            \fill (4, 1) circle (2pt) node[above left] {$\widetilde{C}_1$};
            \fill (6, 5) circle (2pt) node[above right] {$C^{\bullet}_1$};
            \fill (1, 5) circle (2pt) node[above right] {$D^{\bullet}_1$};
            \fill (1, 5) circle (2pt) node[below left] {$\widetilde{C}_2$};
            \fill (6, 5) circle (2pt) node[below left] {$\widetilde{D}_1$};
            \fill (2, 1) circle (2pt) node[above left] {$\widetilde{D}_2$};
            \fill (2, 1) circle (2pt) node[below] {$C^{\circ}_2$};
            \fill (-1.5, -0.5) circle (2pt) node[below left] {$B_{\sigma(i)}$};

            \draw (-2, 5) -- (8, 5);
            \draw (-2, 1) -- (8, 1);

            \draw[decorate, decoration={snake, amplitude=0.3mm, segment length=3mm}, ->, black] (7, -0.75) -- (4, -0.75) -- (4, 1);
            \draw[decorate, decoration={snake, amplitude=0.3mm, segment length=3mm}, ->, red] (4, 1) -- (4, 3) -- (6, 3) -- (6, 5);
            \draw[decorate, decoration={snake, amplitude=0.3mm, segment length=3mm}, ->, red] (6, 5) -- (6, 6) -- (1, 6) -- (1, 5);
            \draw[decorate, decoration={snake, amplitude=0.3mm, segment length=3mm}, ->, red] (1, 5) -- (1, 3) -- (2, 3) -- (2, 1);
 
        \end{tikzpicture}
        \subcaption{}
        \label{fig:weight-of-ribbons-filled-by-same-path-negative-d}
    \end{subfigure}\\[1cm]
    \caption{Path examples.}
    \label{fig:weight-of-ribbons-filled-by-same-path-negative}
\end{figure}

\begin{proof}
    It follows from the construction of the lattice that 
    \begin{equation*}
        |\mathrm{wt}(P)| = \prod^{s_1}_r \mathrm{wt}(T^{\circ}_{\rho^{\circ}_r}) \prod^{s_2}_r \mathrm{wt}(T^{\bullet}_{\rho^{\bullet}_r}) \prod^{s_3}_r \mathrm{wt}(\widetilde{T}_{\widetilde{\rho}_{\sigma(i)}}).
    \end{equation*}
    The non-trivial part is the sign of $\mathrm{wt}(P)$. Define the sign $\alpha(P)$ as follows: 
    \begin{equation*}
        \mathrm{wt}(P) = \alpha(P) \cdot |\mathrm{wt}(P)|.
    \end{equation*}
    We aim to prove that $\alpha(P) = \mathrm{sgn}_{\pi}(\theta_i \# \theta_{\sigma(i)})$. We prove it first by establishing some properties of the path $P$, and consequently properties of ribbons in $\Theta^{\circ}, \widetilde{\Theta}, \Theta^{\bullet}$. 
    
    There are four cases induced by the placement of endpoints $A_i, B_{\sigma(i)}$ to top and bottom layers. We consider the case in which both $A_i, B_{\sigma(i)}$ are in bottom layer; and other cases can be proved similarly. Note that, if $A_i, B_{\sigma(i)}$ are in the lower layer, then $\theta_i \# \theta_{\sigma(i)}$ is an odd ribbon, and hence $\mathrm{sgn}_{\pi}(\theta_i \# \theta_{\sigma(i)}) = -1$. 
    
    Recall that the ribbons in $\Theta^{\circ} \cup \widetilde{\Theta} \cup \Theta^{\bullet}$ are filled by  minor paths. Denote the minor starting and ending node of the minor path that fills the ribbon $\rho^{\circ}_r$ as $C^{\circ}_r, D^{\circ}_r,\ r \in [s_1]$. In the same way define the endpoints $C^{\bullet}_r, D^{\bullet}_r$ for the ribbon $\rho^{\bullet}_r,\ r \in [s_2]$; and $\widetilde{C}_r, \widetilde{D}_r$ for the ribbon $\widetilde{\rho}_r,\ r \in [s_3]$. Define the minor path contained in $P$ that starts at $C^{\circ}_i$ and terminates at $D^{\circ}_i$ as $P^{\circ}_i$. In the same way define the minor paths $P^{\bullet}_i,i \in [s_2],\ \widetilde{P}_i, i \in [s_3]$. We also denote the coordinates of minor endpoints as follows: $C^{\bullet}_r = (x^{\bullet}_r, \cdot), D^{\bullet}_r = (y^{\bullet}_r, \cdot), r \in [s_2]\ C^{\circ}_r = (x^{\circ}_r, \cdot), D^{\circ}_r = (y^{\circ}_r, \cdot), r \in [s_1],\ \widetilde{C}_r = (\widetilde{x}_r, \cdot), \widetilde{D}_r = (\widetilde{y}_r, \cdot), r \in [s_3]$
    

    We index the minor endpoints along the path $P$ according to their order of occurrence, starting from index $1$ and node $A_i$. As $A_i$ is located in lower layer, we set $C^{\circ}_1 = A_i$. The path $P$ terminates at $B_{\sigma(i)}$ which is located in the lower layer, thus $B_{\sigma(i)} = D^{\circ}_{s_1}$. Let $P^{\circ}_1$ be part of $P$ that starts at $C^{\circ}_1$ and terminates at $D^{\circ}_1$. If $D^{\circ}_1 = B_{\sigma(i)}$ (i.e. $s_1 = 1$), then we are done (Fig.~\ref{fig:weight-of-ribbons-filled-by-same-path-negative-a}); otherwise from Lemma~\ref{lem:minor-is-either-major-or-another-minor} it follows that $D^{\circ}_1 = \widetilde{C}_1$ (Fig.~\ref{fig:weight-of-ribbons-filled-by-same-path-negative-b}). Let $\widetilde{P}_1$ be minor path which starts at $\widetilde{C}_1$ and terminates at $\widetilde{D}_1$. There are two choices for $\widetilde{D}_1$: either $\widetilde{D}_1 = C^{\circ}_2$ or $\widetilde{D}_1 = C^{\bullet}_1$. In the first case, we have that $D^{\circ}_2 = \widetilde{C}_2$ and this endpoint is located between $\widetilde{C}_1, \widetilde{D}_1$; i.e. $\widetilde{x}_1 < \widetilde{x}_2 < \widetilde{y}_2$ and it follows that the paths $\widetilde{P}_1, \widetilde{P}_2$ overlap, which is a contradiction (Fig.~\ref{fig:weight-of-ribbons-filled-by-same-path-negative-c}); thus, the only choice for $\widetilde{D}_1$ is $C^{\bullet}_1$. Next, it follows that $D^{\bullet}_1 = \widetilde{C}_2$, as there are no other choice for $D^{\bullet}_1$ (Lemma~\ref{lem:minor-is-either-major-or-another-minor}, Fig.~\ref{fig:weight-of-ribbons-filled-by-same-path-negative-d}). The inequality $\widetilde{y}_2 < \widetilde{x}_1$ must hold, since otherwise the paths $\widetilde{P}_1, \widetilde{P}_2$ overlap. Now, there are two choices for $\widetilde{D}_2$: either $\widetilde{D}_2 = C^{\bullet}_2$ or $\widetilde{D}_2 = C^{\circ}_2$. In the first case we have a contradiction similar to the case when we had choice for $\widetilde{D}_1$; thus $\widetilde{D}_2 = C^{\circ}_2$ (Fig.~\ref{fig:weight-of-ribbons-filled-by-same-path-negative-d}). 

\begin{figure}
    \centering
    \begin{subfigure}[b]{0.5\textwidth}
        \centering
        \begin{tikzpicture}[scale=0.5, font=\scriptsize]


            \fill (7, -0.75) circle (2pt) node[below right] {$A_i$};
            \fill (7, -0.75) circle (2pt) node[above right] {$C^{\circ}_1$};
            \fill (5.5, 1) circle (2pt) node[below left] {$C^{\circ}_2$};
            \fill (-1.5, -0.75) circle (2pt) node[below left] {$B_{\sigma(i)}$};

            \draw (-2, 5) -- (8, 5);
            \draw (-2, 1) -- (8, 1);

            \draw[decorate, decoration={snake, amplitude=0.25mm, segment length=3.5mm}, ->, black] (7, -0.75) -- (6.5, -0.75) -- (6.5, 1) -- (6.5, 3) -- (7.5, 3) -- (7.5, 5) -- (7.5, 5.5) -- (5, 5.5) -- (5, 3) -- (5.5, 3) -- (5.5, 1);
            \draw[decorate, decoration={snake, amplitude=0.25mm, segment length=3.5mm}, ->, red] (5.5, 1) -- (5.5, -0.75) -- (-1.5, -0.75);
 
        \end{tikzpicture}
        \subcaption{}
        \label{fig:weight-of-ribbons-filled-by-same-path-negative-2-a}
    \end{subfigure}
    \begin{subfigure}[b]{0.49\textwidth}
        \centering
        \begin{tikzpicture}[scale=0.5, font=\scriptsize]


            \fill (7, -0.75) circle (2pt) node[below right] {$A_i$};
            \fill (7, -0.75) circle (2pt) node[above right] {$C^{\circ}_1$};
            \fill (5.5, 1) circle (2pt) node[below left] {$C^{\circ}_2$};
            \fill (3, 1) circle (2pt) node[below left] {$D^{\circ}_2$};
            \fill (-1.5, -0.75) circle (2pt) node[below left] {$B_{\sigma(i)}$};

            \draw (-2, 5) -- (8, 5);
            \draw (-2, 1) -- (8, 1);

            \draw[decorate, decoration={snake, amplitude=0.25mm, segment length=3.5mm}, ->, black] (7, -0.75) -- (6.5, -0.75) -- (6.5, 1) -- (6.5, 3) -- (7.5, 3) -- (7.5, 5) -- (7.5, 5.5) -- (5, 5.5) -- (5, 3) -- (5.5, 3) -- (5.5, 1);
            \draw[decorate, decoration={snake, amplitude=0.25mm, segment length=3mm}, ->, red] (5.5, 1) -- (5.5, 0) -- (3, 0) -- (3, 1);
 
        \end{tikzpicture}
        \subcaption{}
        \label{fig:weight-of-ribbons-filled-by-same-path-negative-2-b}
    \end{subfigure}\\[1cm]
    \begin{subfigure}[b]{0.5\textwidth}
        \centering
        \begin{tikzpicture}[scale=0.5, font=\scriptsize]

            \fill[black] (3.75, -1.5) node {$P^{\circ}_2$};
            \fill[red] (2.55, 0) node {$P^{\circ}_3$};

            \fill (7, -0.75) circle (2pt) node[below right] {$A_i$};
            \fill (7, -0.75) circle (2pt) node[above right] {$C^{\circ}_1$};
            \fill (5.5, 1) circle (2pt) node[below left] {$C^{\circ}_2$};
            \fill (2, 1) circle (2pt) node[below left] {$D^{\circ}_2$};
            \fill (2, 1) circle (2pt) node[above left] {$\widetilde{C}_3$};
            \fill (4, 1) circle (2pt) node[above left] {$\widetilde{D}_3$};
            \fill (-1.5, -0.75) circle (2pt) node[below left] {$B_{\sigma(i)}$};

            \draw (-2, 5) -- (8, 5);
            \draw (-2, 1) -- (8, 1);

            \draw[decorate, decoration={snake, amplitude=0.25mm, segment length=3.5mm}, ->, black] (7, -0.75) -- (6.5, -0.75) -- (6.5, 1) -- (6.5, 3) -- (7.5, 3) -- (7.5, 5) -- (7.5, 5.5) -- (5, 5.5) -- (5, 3) -- (5.5, 3) -- (5.5, 1);
            \draw[decorate, decoration={snake, amplitude=0.25mm, segment length=3mm}, ->, black] (5.5, 1) -- (5.5, -0.5) -- (2, -0.5) -- (2, 1);
            \draw[decorate, decoration={snake, amplitude=0.25mm, segment length=3mm}, ->, red, densely dotted, thick] (2, 1) -- (2, 3) -- (4, 3) -- (4, 1);
            \draw[decorate, decoration={snake, amplitude=0.25mm, segment length=3mm}, ->, red, densely dotted, thick] (4, 1) -- (4, 0) -- (3, 0);
 
        \end{tikzpicture}
        \subcaption{}
        \label{fig:weight-of-ribbons-filled-by-same-path-negative-2-c}
    \end{subfigure}
    \begin{subfigure}[b]{0.49\textwidth}
    \centering
        \begin{tikzpicture}[scale=0.5, font=\scriptsize]

            \fill[red] (2.5, 6.25) node[above] {$P^{\bullet}_2$};
            \fill[red] (3.5, 2.5) node {$\widetilde{P}_3$};
            \fill[red] (4, -0.25) node[below] {$P^{\circ}_2$};

            \fill (7, -0.75) circle (2pt) node[below right] {$A_i$};
            \fill (7, -0.75) circle (2pt) node[above right] {$C^{\circ}_1$};
            \fill (5.5, 1) circle (2pt) node[below left] {$C^{\circ}_2$};
            \fill (3, 1) circle (2pt) node[below left] {$D^{\circ}_2$};
            \fill (3, 1) circle (2pt) node[above left] {$\widetilde{C}_3$};
            \fill (4, 5) circle (2pt) node[below left] {$\widetilde{D}_3$};
            \fill (4, 5) circle (2pt) node[above left] {$C^{\bullet}_2$};
            \fill (1, 5) circle (2pt) node[above left] {$D^{\bullet}_2$};
            \fill (1, 5) circle (2pt) node[below left] {$\widetilde{C}_4$};
            \fill (-1.5, -0.75) circle (2pt) node[below left] {$B_{\sigma(i)}$};

            \draw (-2, 5) -- (8, 5);
            \draw (-2, 1) -- (8, 1);

            \draw[decorate, decoration={snake, amplitude=0.25mm, segment length=3.5mm}, ->, black] (7, -0.75) -- (6.5, -0.75) -- (6.5, 1) -- (6.5, 3) -- (7.5, 3) -- (7.5, 5) -- (7.5, 5.5) -- (5, 5.5) -- (5, 3) -- (5.5, 3) -- (5.5, 1);
            \draw[decorate, decoration={snake, amplitude=0.25mm, segment length=2.5mm}, ->, red] (5.5, 1) -- (5.5, 0) -- (3, 0) -- (3, 1);
            \draw[decorate, decoration={snake, amplitude=0.25mm, segment length=2.5mm}, ->, red] (3, 1) -- (3, 3) -- (4, 3) -- (4, 5);
            \draw[decorate, decoration={snake, amplitude=0.25mm, segment length=2.5mm}, ->, red] (4, 5) -- (4, 6) -- (1, 6) -- (1, 5);
 
        \end{tikzpicture}
        \subcaption{}
        \label{fig:weight-of-ribbons-filled-by-same-path-negative-2-d}
    \end{subfigure}\\[1cm]
    \caption{Path examples.}
    \label{fig:weight-of-ribbons-filled-by-same-path-negative-2}
\end{figure}

    Visually the sub-path of $P$ connecting $C^{\circ}_1$ and $C^{\circ}_2$ is `$\Omega$-shaped'. We claim that the sub-paths connecting $C^{\circ}_{r}, C^{\circ}_{r + 1}, 2 \le r < s_1$ are `$\Omega$-shaped' as well. More precisely, we claim that  
    \begin{align*}
        &C^{\circ}_r = \widetilde{D}_{2r - 2},\ 1 < r \le s_1 \qquad D^{\circ}_r = \widetilde{C}_{2r - 1},\ 1 \le r < s_1 \\ 
        &C^{\bullet}_r = \widetilde{D}_{2r - 1},\ r \in [s_2] \qquad D^{\bullet}_r = \widetilde{C}_{2r},\ r \in [s_2].
    \end{align*}
    Hence we have $s_3 = 2s_1 - 2, s_3 = 2s_2$, and consequently $s_1 - 1 = s_2$. 
    
    Let us check statement for $r = 2$. We already know that $C^{\circ}_2 = \widetilde{D}_2$. Consider the path $P^{\circ}_2$; if $r = s_1$ then $D^{\circ}_2 = B_{\sigma(i)}$ and there is nothing to check as path $P$ terminates at $B_{\sigma(i)}$ (Fig.~\ref{fig:weight-of-ribbons-filled-by-same-path-negative-2-a}). If $r < s_1$ it follows from Lemma~\ref{lem:minor-is-either-major-or-another-minor} that $D^{\circ}_2 = \widetilde{C}_3$ from which the part $\widetilde{P}_3$ of path $P$ starts (Fig.~\ref{fig:weight-of-ribbons-filled-by-same-path-negative-2-b}). Then there are two choices for $\widetilde{D}_3$: either $\widetilde{D}_3 = C^{\circ}_3$ or $\widetilde{D}_3 = C^{\bullet}_3$; In the first case, the inequality $y^{\circ}_2 < \widetilde{y}_3 < x^{\circ}_2$ must hold; but then $P^{\circ}_3$ overlaps with $P^{\circ}_2$ (Fig.~\ref{fig:weight-of-ribbons-filled-by-same-path-negative-2-c}), contradiction. 
    Thus, $\widetilde{D}_3 = C^{\bullet}_2$; then, as $D^{\bullet}_2$ is not major endpoint, we have $D^{\bullet}_2 = \widetilde{C}_4$ (Fig.~\ref{fig:weight-of-ribbons-filled-by-same-path-negative-2-d}), proving the statement for $r = 2$. 
    
    We continue with induction on $r$. Suppose the statement holds for all $r \le N - 1$ where $N \le s_1$. Let $r = N$. By induction hypothesis we have $\widetilde{C}_{2N - 2} = D^{\bullet}_{N - 1}$. We have two choices for $\widetilde{D}_{2N - 2}$: either $\widetilde{D}_{2N - 2} = C^{\bullet}_N$ or $\widetilde{D}_{2N - 2} = C^{\circ}_N$. In the first case, we have  $y^{\bullet}_{N - 1} < \widetilde{y}_{2N - 2} < x^{\bullet}_{N - 1}$, and consequently the path $P^{\bullet}_N$ overlaps with $P^{\bullet}_{N - 1}$, a contradiction; thus $\widetilde{D}_{2N - 2} = C^{\circ}_N$. If $N = s_1$, then $D^{\circ}_N = B_{\sigma(i)}$ and there is nothing to check further. If $N < s_1$, then $D^{\circ}_N = \widetilde{C}_{2N - 1}$. There are two choices for $\widetilde{D}_{2N - 1}$: either $\widetilde{D}_{2N - 1} = C^{\circ}_{N + 1}$ or $\widetilde{D}_{2N - 1} = C^{\bullet}_N$. In the first case we have that $y^{\circ}_N < \widetilde{y}_{2N - 1} < z^{\circ}_N$, thus the path $P^{\circ}_{N + 1}$ overlaps with $P^{\circ}_N$, a contradiction; thus $\widetilde{D}_{2N - 1} = C^{\bullet}_N$. Since $D^{\bullet}_N \neq B_{\sigma(i)}$ it follows that $D^{\bullet}_N = \widetilde{C}_{2N}$; and by induction the statement follows. 
    
    Let $\delta^{\circ}_r, \gamma^{\circ}_r$ are head and tail of the ribbon $\rho^{\circ}_r,\ r \in [s_1]$. 
    From the proved statement it follows that $\pi(\delta^{\circ}_1) = (\cdot, 0),\ \pi(\gamma^{\circ}_{s_1}) = (1, \cdot)$; additionally we have $\pi(\gamma^{\circ}_1) = (0, \cdot),\ \pi(\delta^{\circ}_{s_1}) = (\cdot, 1)$ if $s_1 > 1$. For all other ribbons $\rho^{\circ}_r,\ 1 < r < s_1$ we have that $\pi(\delta^{\circ}_r) = (\cdot, 1)$ and $\pi(\gamma^{\circ}_r) = (0, \cdot)$. For the ribbons in $\rho^{\bullet}_r, r \in [s_2]$ we have $\pi(\delta^{\bullet}_r) = (\cdot, 0),\ \pi(\gamma^{\bullet}_r) = (1, \cdot)$, where $\delta^{\bullet}_r, \gamma^{\bullet}_r$ are head and tail of the ribbon $\rho^{\bullet}_r,\ r \in [s_2]$. We also have $\pi(\widetilde{\delta}_r) = (x, \cdot),\ \pi(\widetilde{\gamma}_r) = (\cdot, x),\ r \in [s_3],\ x \in \{0, 1\}$, where $\widetilde{\delta}_r, \widetilde{\gamma}_r$ are head and tail of the ribbon $\widetilde{\rho}_r,\ r \in [s_3]$; consequently, the number of hook pipes in every ribbon $\widetilde{\rho}_r$ is even (Fig.~\ref{fig:weight-of-ribbons-filled-by-same-path-negative-3}).

    Now we are ready to count the number of hook pipes in the ribbons $(\rho^{\circ}_r), (\rho^{\bullet}_r)$; the sign of the weight $\mathrm{wt}(P)$ is affected only by hook pipes in these ribbons. 
    Note that $\alpha(P)$ is $(-1)^N$ where $N$ is the number of hook pipes in all ribbons $(\rho^{\circ}_r), (\rho^{\bullet}_r)$; in other words 
    \begin{equation*}
        \alpha(P) = \prod_r \mathrm{sgn}_{\pi}(\rho^{\circ}_r) \prod_r \mathrm{sgn}_{\pi}(\rho^{\bullet}_r).
    \end{equation*}
    If $s_1 = 1$ we have that $s_2 = s_3 = 0$ and $\alpha(P) = \mathrm{sgn}_{\pi}(\rho^{\circ}_1) = -1 = \mathrm{sgn}_{\pi}(\theta_i \# \theta_{\sigma(i)})$ and we are done. If $s_1 > 1$ we have   
    \begin{equation*}
        \mathrm{sgn}_{\pi}(\rho^{\circ}_r) = 
        \begin{cases}
            -1, \quad 1 < r < s_1, \\
            +1, \quad \text{otherwise},
        \end{cases}
        \qquad 
        \mathrm{sgn}_{\pi}(\rho^{\bullet}_r) = -1, \quad r \in [s_2].
    \end{equation*}
    Hence recalling that $s_2 = s_1 - 1$ we have 
    \begin{equation*}
        \alpha(P) = (-1)^{s_1 - 2 + s_2} = (-1)^{2s_1 - 3} = (-1) = \mathrm{sgn}_{\pi}(\theta_i \# \theta_{\sigma(i)}), 
    \end{equation*}
    which concludes proof. 
\end{proof}

\begin{figure}
    \centering
    \begin{subfigure}[b]{0.5\textwidth}
        \centering
        \begin{tikzpicture}[scale=0.5, font=\scriptsize]

            \fill[blue] (4, -0.25) node[below] {$P^{\circ}_r$};
            \fill[black] (4, 3) node[below] {$\widetilde{P}_{2r - 1}$};
            \fill[red] (4, 6.5) node[above] {$P^{\bullet}_r$};

            \fill (2.5, 1) circle (2pt) node[above left] {$\widetilde{C}_{2r - 1}$};
            \fill (5, 5) circle (2pt) node[below right] {$\widetilde{D}_{2r - 1}$};

            \draw (-2, 5) -- (8, 5);
            \draw (-2, 1) -- (8, 1);

            \draw[decorate, decoration={snake, amplitude=0.25mm, segment length=3mm}, ->, blue] (7, 0) -- (2.5, 0) -- (2.5, 1);
            \draw[decorate, decoration={snake, amplitude=0.25mm, segment length=3.5mm}, ->, purple] (2.5, 1) -- (2.5, 3) -- (5, 3) -- (5, 5);
            \draw[decorate, decoration={snake, amplitude=0.25mm, segment length=3.5mm}, ->, red] (5, 5) -- (5, 6) -- (0, 6);
            

        \end{tikzpicture}
        \subcaption{}
        \label{fig:weight-of-ribbons-filled-by-same-path-negative-3-a}
    \end{subfigure}
    \begin{subfigure}[b]{0.49\textwidth}
        \centering
        \begin{tikzpicture}[y=-1cm, scale=0.25]
            \foreach \x in {7,...,9} {
                \draw[thick] (\x, 1) rectangle ++(1, 1);
            }
            \foreach \x in {7,...,8} {
                \draw[thick, fill=blue!20] (\x, 1) rectangle ++(1, 1);
            }
            \foreach \x in {5,...,7} {
                \draw[thick] (\x, 2) rectangle ++(1, 1);
            }
            \foreach \x in {7,...,7} {
                \draw[thick, fill=blue!20] (\x, 2) rectangle ++(1, 1);
            }
            \foreach \x in {5,...,6} {
                \draw[thick, fill=sredpurple] (\x, 2) rectangle ++(1, 1);
            }

            \foreach \x in {5,...,5} {
                \draw[thick, fill=sredpurple] (\x, 3) rectangle ++(1, 1);
            }

            \foreach \x in {4,...,5} {
                \draw[thick, fill=sredpurple] (\x, 4) rectangle ++(1, 1);
            }
            \foreach \x in {3,...,3} {
                \draw[thick, fill=red!20] (\x, 4) rectangle ++(1, 1);
            }

            \foreach \x in {2,...,3} {
                \draw[thick, fill=red!20] (\x, 5) rectangle ++(1, 1);
            }
            \foreach \x in {1,...,2} {
                \draw[thick] (\x, 5) rectangle ++(1, 1);
            }

            \foreach \x in {7,...,8} {
                \draw[thick, fill=blue!20] (\x, -5) rectangle ++(1, 1);
            }
            \foreach \x in {5,...,7} {
                \draw[thick, fill=blue!20] (\x, -4) rectangle ++(1, 1);
            }

            \foreach \x in {5,...,5} {
                \draw[thick, fill=blue!20] (\x, -3) rectangle ++(1, 1);
            }

            \foreach \x in {4,...,5} {
                \draw[thick, fill=blue!20] (\x, -2) rectangle ++(1, 1);
            }

            \foreach \x in {13,...,14} {
                \draw[thick, fill=sredpurple] (\x, 2) rectangle ++(1, 1);
            }

            \foreach \x in {13,...,13} {
                \draw[thick, fill=sredpurple] (\x, 3) rectangle ++(1, 1);
            }

            \foreach \x in {12,...,13} {
                \draw[thick, fill=sredpurple] (\x, 4) rectangle ++(1, 1);
            }
            
            \foreach \x in {7,...,9} {
                \draw[thick] (\x, 1) rectangle ++(1, 1);
            }
            \foreach \x in {7,...,8} {
                \draw[thick, fill=blue!20] (\x, 1) rectangle ++(1, 1);
            }
            
            \foreach \x in {5,...,6} {
                \draw[thick, fill=red!20] (\x, 8) rectangle ++(1, 1);
            }

            \foreach \x in {5,...,5} {
                \draw[thick, fill=red!20] (\x, 9) rectangle ++(1, 1);
            }

            \foreach \x in {4,...,5} {
                \draw[thick, fill=red!20] (\x, 10) rectangle ++(1, 1);
            }
            \foreach \x in {3,...,3} {
                \draw[thick, fill=red!20] (\x, 11) rectangle ++(1, 1);
            }

            \foreach \x in {2,...,3} {
                \draw[thick, fill=red!20] (\x, 12) rectangle ++(1, 1);
            }

            \foreach \x in {5,...,6} {
                \draw[thick] (1, \x) -- (0.7, \x);
            }
            \foreach \x in {1,...,2} {
                \draw[thick] (10, \x) -- (10.3, \x);
            }

            \draw[dashed] (4, 4) -- (4, -2);
            \draw[dashed] (9, 1) -- (9, -4);
            \draw[dashed] (2, 6) -- (2, 12);
            \draw[dashed] (7, 3) -- (7, 8);


            \fill (0, 5.5) node {$\hat{\theta}_k$};
            \fill (3, -1.5) node {$\rho^{\circ}_r$};
            \fill (1, 12.5) node {$\rho^{\bullet}_r$};
            \fill (17, 2.5) node {$\widetilde{\rho}_{2r - 1}$};

        \end{tikzpicture}
        \subcaption{}
        \label{fig:weight-of-ribbons-filled-by-same-path-negative-3-b}
    \end{subfigure}
    \caption{Path and ribbons example.}
    \label{fig:weight-of-ribbons-filled-by-same-path-negative-3}
\end{figure}

\subsubsection{Weight-preserving property for positive shapes.}
\label{subsec:weigth-preserving-for-positive-shapes}
\begin{lemma}
    \label{lem:intersection-of-pos-ribbon-not-exist}
    Let $\mathbf{F} = (\mathbf{P}, \mathbf{C}) = \psi(\mathbf{T}) = \psi((T^{\circ}, \widetilde{T}, T^{\bullet}))$ such that $\mathbf{P} \in \mathcal{P}(\mathbf{A}, \mathbf{B}, \sigma)$ where $\sigma$ is some permutation. Let $\sigma_0 = (a_1, \ldots, a_{\ell}), \ell \ge 1$ be some cycle in cycle decomposition of $\sigma$. Then we have: 
    
    (1) There is a positive ribbon $\theta_{a_j} \# \theta_{\sigma(a_j)},\ j \in [\ell]$ among ribbons $(\theta_{a_i} \# \theta_{\sigma(a_i)})$ such that every other ribbon $\theta_{a_i} \# \theta_{\sigma(a_i)}, i \in [\ell] \setminus \{j\}$ is a sub-ribbon of $\theta_{a_j} \# \theta_{\sigma(a_j)}$.
    
    (2) All ribbons $\theta_{a_i} \# \theta_{\sigma(a_i)}, i \in [\ell] \setminus \{j\}$ are negative and odd.
    
    
\end{lemma}

\begin{figure}
    \centering
    \begin{subfigure}[b]{0.5\textwidth}
        \centering
        \begin{tikzpicture}[scale=0.5, font=\scriptsize]
            \fill[blue] (4, 2) node[above] {$P'_j$};
            \fill[red] (3, 4) node[below] {$P'_s$};

            \fill (3, 0) circle (2pt) node[below right] {$A_{a_j}$};
            \fill (0, 6) circle (2pt) node[below right] {$A_{a_s}$};
            \fill (6, 6) circle (2pt) node[below left] {$B_{\sigma(a_s)}$};

            \draw (-2, 5) -- (8, 5);
            \draw (-2, 1) -- (8, 1);

            \draw[decorate, decoration={snake, amplitude=0.25mm, segment length=3mm}, ->, blue] (3, 0) -- (2, 0) -- (2, 1);
            
            \draw[decorate, decoration={snake, amplitude=0.25mm, segment length=3mm}, ->, blue] (2, 1) -- (2, 2) -- (6, 2);


            \draw[decorate, decoration={snake, amplitude=0.25mm, segment length=3mm}, ->, red] (0, 6) -- (-1, 6) -- (-1, 5);

            \draw[decorate, decoration={snake, amplitude=0.25mm, segment length=3mm}, ->, red] (-1, 5) -- (-1, 4) -- (7, 4) -- (7, 5);

            \draw[decorate, decoration={snake, amplitude=0.25mm, segment length=3mm}, ->, red] (7, 5) -- (7, 6) -- (6, 6);
            

        \end{tikzpicture}
        \subcaption{}
        \label{fig:intersection-of-pos-ribbon-not-exist-1-a}
    \end{subfigure}
    \begin{subfigure}[b]{0.49\textwidth}
        \centering
        \begin{tikzpicture}[scale=0.5, font=\scriptsize]
            \fill[blue] (4, 2) node[above] {$P'_j$};
            \fill[red] (3, 4) node[below] {$P'_s$};

            \fill (3, 0) circle (2pt) node[below left] {$B_{a_j}$};
            \fill (0, 6) circle (2pt) node[below right] {$A_{a_s}$};
            \fill (6, 6) circle (2pt) node[below left] {$B_{\sigma(a_s)}$};

            \draw (-2, 5) -- (8, 5);
            \draw (-2, 1) -- (8, 1);

            
            \draw[decorate, decoration={snake, amplitude=0.25mm, segment length=3mm}, ->, blue] (1, 2) -- (5, 2) -- (5, 1);

            \draw[decorate, decoration={snake, amplitude=0.25mm, segment length=3mm}, ->, blue] (5, 1) -- (5, 0) -- (3, 0);


            \draw[decorate, decoration={snake, amplitude=0.25mm, segment length=3mm}, ->, red] (0, 6) -- (-1, 6) -- (-1, 5);

            \draw[decorate, decoration={snake, amplitude=0.25mm, segment length=3mm}, ->, red] (-1, 5) -- (-1, 4) -- (7, 4) -- (7, 5);

            \draw[decorate, decoration={snake, amplitude=0.25mm, segment length=3mm}, ->, red] (7, 5) -- (7, 6) -- (6, 6);
            

        \end{tikzpicture}
        \subcaption{}
        \label{fig:intersection-of-pos-ribbon-not-exist-1-b}
    \end{subfigure}
    \caption{Paths examples.}
    \label{fig:intersection-of-pos-ribbon-not-exist-1}
\end{figure}
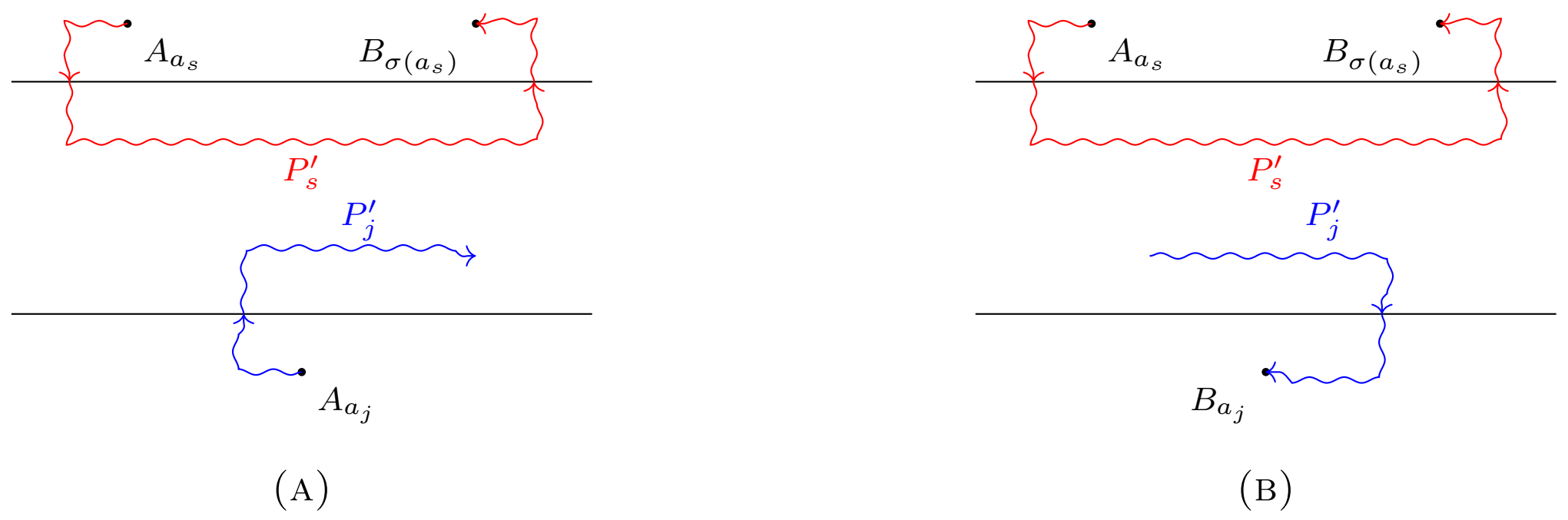

\begin{proof}
    Let $\delta_{a_i}, \gamma_{a_i}$ be the head and tail of $\theta_{a_i}, i \in [\ell]$. 
    From Lemma~\ref{lem:th-i-th-si-are-in-the-same-ribbon} it follows that all the ribbons $(\theta_{a_i})$ are in the same tile $\hat{\theta}_k, k \in \mathbb{Z}$. Hence they do not share cells with common contents. Note also that $(\theta_{a_i})$ are positive, since the shape $\lambda/\mu$ is positive; thus it follows that there must be some positive ribbon $\theta_{a_j} \# \theta_{\sigma(a_j)}$ for some $j \in [\ell]$. We claim that every ribbon $\theta_{a_i} \# \theta_{\sigma(a_i)}$ is a sub-ribbon of $\theta_{a_j} \# \theta_{\sigma(a_j)}$; i.e. as the ribbons $(\theta_{a_i})$ do not intersect, it is enough to show that every ribbon $\theta_{a_i}, i \in [\ell] \setminus \{j\}$ is between the ribbons $\theta_{a_j}$ and $\theta_{\sigma(a_j)}$; more precisely we need to show that $c(\gamma_{a_j}) < c(\delta_{a_i}) < c(\delta_{\sigma(a_j)}), i \in [\ell] \setminus \{j\}$. By contradiction, suppose that there is a ribbon $\theta_{a_r}, r \in [\ell]$ such that $\theta_{a_r}$ is not between $\theta_{a_j}$ and $\theta_{\sigma(a_j)}$. We have $a_r = \sigma^{t}(a_j)$ for some $t \in [\ell - 1]$. There are two cases to consider: (1) $c(\gamma_{a_r}) < c(\delta_{a_j})$; and (2) $c(\gamma_{\sigma(a_j)}) < c(\delta_{a_r})$. In the first case, there must be $s$ such that $t \le s < \ell$ and $c(\delta_{\sigma^{s}(a_j)}) < c(\delta_{a_j}) \le c(\gamma_{\sigma^{s + 1}(a_j)}) + 1$; thus the paths $P_j \in \mathcal{P}(A_{a_j}, B_{\sigma(a_j)}),\ P_s \in \mathcal{P}(A_{\sigma^s(a_j)}, B_{\sigma^{s + 1}(a_j)})$ contain minor paths $P'_j, P'_s$ that belong to central layer and overlap, which is a contradiction by Lemma~\ref{lem:cousin-paths-avoid-each-others-endpoints} (Fig.~\ref{fig:intersection-of-pos-ribbon-not-exist-1-a}). In the second case there must be $s$ such that $1 \le s \le t$ and $c(\delta_{\sigma^s(a_j)}) \le c(\gamma_{\sigma(a_j)}) + 1 < c(\gamma_{\sigma^{s + 1}(a_j)}) + 1$; thus the paths $P_j \in \mathcal{P}(A_{a_j}, B_{\sigma(a_j)}),\ P_s \in \mathcal{P}(A_{\sigma^s(a_j)}, B_{\sigma^{s + 1}(a_j)})$ contain the minor paths $P'_j, P'_s$ that belong to central layer and overlap, which is a contradiction by Lemma~\ref{lem:cousin-paths-avoid-each-others-endpoints} (Fig.~\ref{fig:intersection-of-pos-ribbon-not-exist-1-b}).

    Now we prove to the second part of the lemma. First we prove that the ribbons $\theta_{a_i} \# \theta_{\sigma(a_i)}, i \in [\ell] \setminus \{j\}$ are negative. By contrary, suppose there is a positive ribbon $\theta_{a_r} \# \theta_{\sigma(a_r)}, r \in [\ell] \setminus \{j\}$, then it must be contained in $\theta_{a_j} \# \theta_{\sigma(a_j)}$, but it follows that there are minor paths $P'_r, P'_j$ contained in the paths $P_r \in \mathcal{P}(A_{a_r}, B_{\sigma(a_r)}),\ P_j \in \mathcal{P}(A_{a_j}, B_{\sigma(a_j)})$ (resp.) which belong to the middle layer, and they overlap, 
    a contradiction. 


    Now, suppose that there is even ribbon $\theta_{a_i} \# \theta_{\sigma(a_i)}, i \in [\ell] \setminus \{j\}$; it follows that the endpoints $A_{a_i}, B_{\sigma(a_i)}$ are in different layers, thus the path $P_i \in \mathcal{P}(A_{a_i}, B_{\sigma(a_i)})$ contains minor path that belongs to the middle layer, and it overlaps with minor path contained in the path $P_j \in \mathcal{P}(A_{a_j}, B_{\sigma(a_j)})$ and belongs to middle layer, a contradiction. 
\end{proof}

\begin{lemma}
    \label{lem:sgn-rule-for-one-cycle-positive}
    Assume that $\lambda/\mu$ is a positive shape. Let $\mathbf{F} = (\mathbf{P}, \mathbf{C}) = \psi(\mathbf{T}) = \psi((T^{\circ}, \widetilde{T}, T^{\bullet}))$ such that $\mathbf{P} \in \mathcal{P}(\mathbf{A}, \mathbf{B}, \sigma)$ where $\sigma$ is some permutation. Let $\sigma_0 = (a_1, \ldots, a_{\ell}), \ell \ge 1$ be some cycle in cycle decomposition of $\sigma$. Let $\rho^{\circ}_1, \ldots, \rho^{\circ}_{s_1} \in \Theta^{\circ}$, $\rho^{\bullet}_1, \ldots, \rho^{\bullet}_{s_2} \in \Theta^{\bullet}$, $\widetilde{\rho}_1, \ldots, \widetilde{\rho}_{s_3} \in \widetilde{\Theta},\ s_1, s_2, s_3 \ge 0$ be ribbons filled by paths $P_1, \ldots, P_{\ell}$ where $P_i \in \mathcal{P}(A_{a_i}, B_{\sigma(a_i)}),\ i \in [\ell]$. Then the following identity holds:
    \begin{equation*}
        \prod_{i = 1}^{\ell} \mathrm{wt}(P_i) = \mathrm{sgn}(\sigma_0) \prod_{i = 1}^{s_1} \mathrm{wt}(T^{\circ}_{\rho^{\circ}_i})\prod_{i = 1}^{s_2} \mathrm{wt}(\widetilde{T}_{\widetilde{\rho}_i})\prod_{i = 1}^{s_3} \mathrm{wt}(T^{\bullet}_{\rho^{\bullet}_i}).
    \end{equation*}
\end{lemma}
\begin{proof}
    It follows from the construction of the lattice that 
    \begin{equation*}
        \prod_{i = 1}^{\ell} |\mathrm{wt}(P_i)| = \prod_{i = 1}^{s_1} \mathrm{wt}(T^{\circ}_{\rho^{\circ}_i})\prod_{i = 1}^{s_2}\mathrm{wt}(\widetilde{T}_{\widetilde{\rho}_i})\prod_{i = 1}^{s_3} \mathrm{wt}(T^{\bullet}_{\rho^{\bullet}_i}).
    \end{equation*}
    Define $\alpha(P_i)$ as the sign of $\mathrm{wt}(P_i), i \in [\ell]$; i.e. 
    \begin{equation*}
        \mathrm{wt}(P_i) = \alpha(P_i) \cdot |\mathrm{wt}(P_i)|.
    \end{equation*}
    From Lemmas~\ref{lem:weight-of-ribbons-filled-by-same-path-positive} and \ref{lem:weight-of-ribbons-filled-by-same-path-negative} it follows that
    \begin{equation*}
        \alpha(P_i) = \mathrm{sgn}_{\pi}(\theta_{a_i} \# \theta_{\sigma(a_i)}).
    \end{equation*}
    Now we aim to prove that
    \begin{equation*}
        \prod_{i = 1}^{\ell}\mathrm{sgn}_{\pi}(\theta_{a_i} \# \theta_{\sigma(a_i)}) = \mathrm{sgn}(\sigma_0).
    \end{equation*}
    From Lemma~\ref{lem:intersection-of-pos-ribbon-not-exist} it follows that there is only one positive ribbon $\theta_{a_j} \# \theta_{\sigma(a_j)}, j \in [\ell]$ among all ribbon $(\theta_{a_i} \# \theta_{\sigma(a_i)})$, and others are negative and odd. Hence,     
    \begin{equation*}
         \prod_{i = 1}^{\ell}\mathrm{sgn}_{\pi}(\theta_{a_i} \# \theta_{\sigma(a_i)}) = (-1)^{\ell - 1} = \mathrm{sgn}(\sigma_0),
    \end{equation*}
    which concludes proof.
\end{proof}

\begin{lemma}
    \label{lem:weight-preserving-positive}
    Assume that $\lambda/\mu$ is a positive shape. Let $\mathbf{F} = (\mathbf{P}, \mathbf{C}) = \psi(\mathbf{T}) = \psi((T^{\circ}, \widetilde{T}, T^{\bullet}))$. Then the following identity holds: 
    \begin{equation*}
        \mathrm{sgn}(\mathbf{F}) \mathrm{wt}(\mathbf{F}) = \mathrm{sgn}_{\pi}(\lambda/\mu) \cdot \mathrm{wt}(\mathbf{T}).
    \end{equation*}
\end{lemma}
\begin{proof}
    From the construction of the lattice it follows that 
    \begin{equation*}
        |\mathrm{wt}(\mathbf{P})| \cdot |\mathrm{wt}(\mathbf{C}| = \mathrm{wt}(\mathbf{T}).
    \end{equation*}
    Define $\alpha(\mathbf{P}), \alpha(\mathbf{C})$ as follows:
    \begin{equation*}
        |\mathrm{wt}(\mathbf{P})| = \alpha(\mathbf{P}) \cdot \mathrm{wt}(\mathbf{P}), \qquad |\mathrm{wt}(\mathbf{C})| = \alpha(\mathbf{C}) \cdot \mathrm{wt}(\mathbf{C}).
    \end{equation*}
    Assume that $\mathbf{P} \in \mathcal{P}(\mathbf{A}, \mathbf{B}, \sigma)$ where $\sigma$ is some permutation. Let $\sigma = \sigma_1 \ldots \sigma_{\ell}$ be cycle decomposition of $\sigma$. From Lemma~\ref{lem:sgn-rule-for-one-cycle-positive} it follows that
    \begin{align*}
        \alpha(\mathbf{P}) = \prod_{i = 1}^{\ell} \mathrm{sgn}(\sigma_i) = \mathrm{sgn}(\sigma).
    \end{align*}
    Next we discuss signs generated by cycles, i.e. $\alpha(\mathbf{C})$. From Lemma~\ref{lem:cycle-corresponds-to-special-ribbons} it follows that the ribbons from $\Theta^{\circ} \cup \Theta^{\bullet}$ filled by cycles are odd and hence  
    \begin{equation*}
        \alpha(\mathbf{C}) = (-1)^{|\mathbf{C}|} = \mathrm{sgn}(\mathbf{C}).
    \end{equation*}
    Recall that $\mathrm{sgn}(\mathbf{F}) = \mathrm{sgn}(\sigma)\, \mathrm{sgn}(\mathbf{C})$. Combining the preceding observations we have:
    \begin{align*}
        \mathrm{sgn}(\mathbf{F}) \cdot \alpha(\mathbf{P}) \cdot \alpha(\mathbf{C}) &= \mathrm{sgn}(\sigma)\, \mathrm{sgn}(\mathbf{C}) \cdot \mathrm{sgn}(\sigma) \cdot \mathrm{sgn}(\mathbf{C}) = 1 = \mathrm{sgn}_{\pi}(\lambda/\mu), 
    \end{align*}
    as needed. 
\end{proof}

\subsubsection{Weight-preserving property for negative shapes.}
\label{subsec:weigth-preserving-for-negative-shapes}

\begin{lemma}
    \label{lem:intersection-of-neg-ribbon-is-one-special-ribbon}
    Assume that  $\lambda/\mu$ is a negative shape. 
    Let $\mathbf{F} = (\mathbf{P}, \mathbf{C}) = \psi(\mathbf{T}) = \psi((T^{\circ}, \widetilde{T}, T^{\bullet}))$ such that $\mathbf{P} \in \mathcal{P}(\mathbf{A}, \mathbf{B}, \sigma)$ where $\sigma$ is some permutation. Let $\sigma_0 = (a_1, \ldots, a_{\ell}), \ell \ge 1$ be some cycle in cycle decomposition of  $\sigma$. Suppose that the distinct ribbons $\theta_{a_r} \# \theta_{\sigma(a_r)}$ and $\theta_{a_s} \# \theta_{\sigma(a_s)}, r, s \in [\ell]$ are negative and they have intersection (i.e. common cells); without loss of generality assume that $c(\delta_{a_r}) > c(\delta_{a_s})$. Then we have $a_s = \sigma(a_r)$, i.e. the intersection is indeed the ribbon $\theta_{\sigma(a_r)}$; and moreover the ribbon $\theta_{\sigma(a_r)}$ is even. 
\end{lemma}

\begin{figure}
    \centering
    \begin{subfigure}[b]{0.5\textwidth}
        \centering
        \begin{tikzpicture}[scale=0.5, font=\scriptsize]
            \fill[black] (6, 3) node {$P_{a_r}$};
            \fill[black] (1, 3) node {$P_{a_s}$};
            \fill[red] (2, -0.5) node[below] {$P$};

            \fill (5, 0) circle (2pt) node[below right] {$A_{a_s}$};
            \fill (7, 0) circle (2pt) node[below right] {$A_{a_r}$};
            \fill (3, -0.5) circle (2pt) node[below right] {$A_{\sigma(a_r)}$};
            \fill (2, 6) circle (2pt) node[above] {$B_{\sigma(a_r)}$};
            \fill (0, 6) circle (2pt) node[left] {$B_{\sigma(a_s)}$};

            \draw (-2, 5) -- (8, 5);
            \draw (-2, 1) -- (8, 1);

            \draw[decorate, decoration={snake, amplitude=0.3mm, segment length=3mm}, ->, black] (5, 0) -- (5, 0.25) -- (0, 0.25) -- (0, 6);
            
            \draw[decorate, decoration={snake, amplitude=0.3mm, segment length=5mm}, ->, black] (7, 0) -- (7, 6) -- (2, 6);

            \draw[decorate, decoration={snake, amplitude=0.25mm, segment length=3mm}, ->, red] (3, -0.5) -- (1, -0.5);
        \end{tikzpicture}
        \subcaption{}
        \label{fig:intersection-of-neg-ribbon-is-one-special-ribbon-1-a}
    \end{subfigure}
    \begin{subfigure}[b]{0.49\textwidth}
        \centering
        \begin{tikzpicture}[scale=0.5, font=\scriptsize]
            \fill[black] (4, 6) node[below] {$P_{a_s}$};
            \fill[black] (1, 3) node {$P_{a_r}$};
            \fill[red] (2, -0.5) node[below] {$P$};

            \fill (6, 6) circle (2pt) node[below right] {$A_{a_s}$};
            \fill (7, 0) circle (2pt) node[below right] {$A_{a_r}$};
            \fill (3, -0.5) circle (2pt) node[below right] {$A_{\sigma(a_s)}$};
            \fill (1.5, 6) circle (2pt) node[above] {$B_{\sigma(a_s)}$};
            \fill (0, 6) circle (2pt) node[left] {$B_{\sigma(a_r)}$};

            \draw (-2, 5) -- (8, 5);
            \draw (-2, 1) -- (8, 1);

            \draw[decorate, decoration={snake, amplitude=0.3mm, segment length=3mm}, ->, black] (6, 6) -- (1.5, 6);
            
            \draw[decorate, decoration={snake, amplitude=0.3mm, segment length=5mm}, ->, black] (7, 0) -- (0, 0) -- (0, 6);

            \draw[decorate, decoration={snake, amplitude=0.25mm, segment length=3mm}, ->, red] (3, -0.5) -- (1, -0.5);
        \end{tikzpicture}
        \subcaption{}
        \label{fig:intersection-of-neg-ribbon-is-one-special-ribbon-1-b}
    \end{subfigure}
    \caption{Paths examples.}
    \label{fig:intersection-of-neg-ribbon-is-one-special-ribbon-1}
\end{figure}
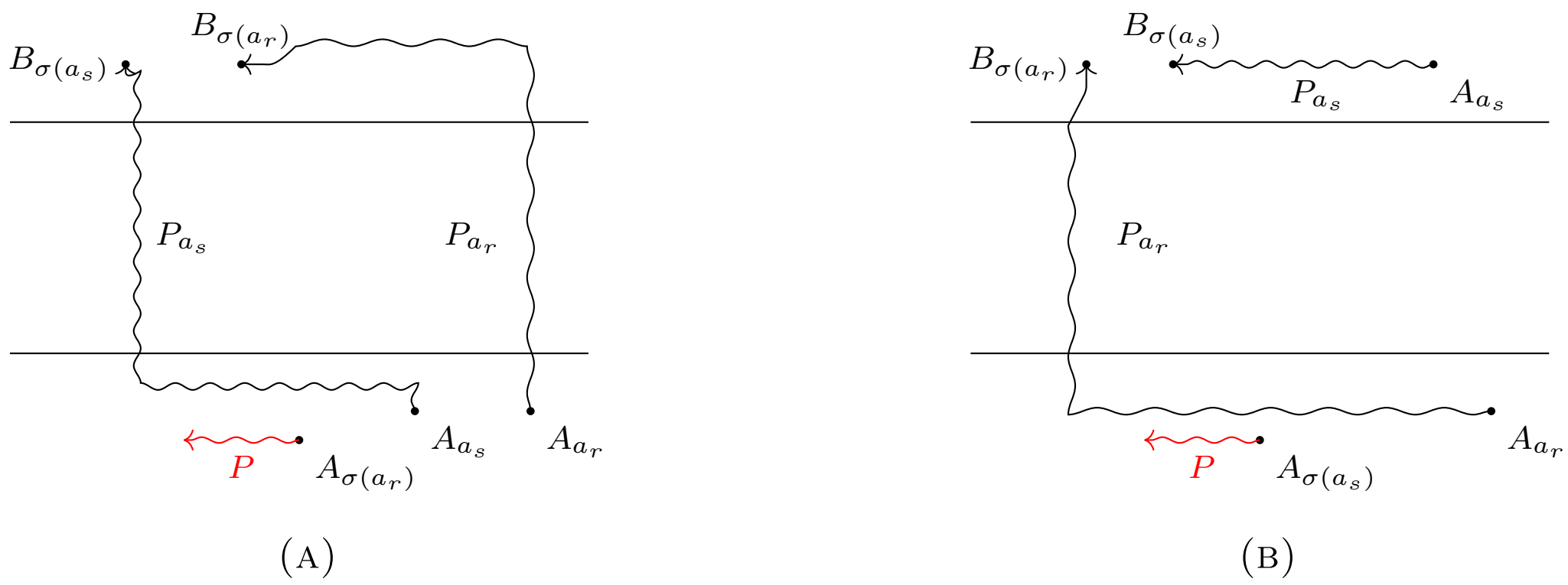

\begin{proof}
    Let $\delta_{a_i}, \gamma_{a_i}$ be the head and tail of $\theta_{a_i}, i \in [\ell]$. From Lemma~\ref{lem:th-i-th-si-are-in-the-same-ribbon} it follows that the ribbons $(\theta_{a_i})$ are in the same tile $\hat{\theta}_k, k \in \mathbb{Z}$.
    Since the shape $\lambda/\mu$ is negative, all ribbons $(\theta_{a_i})$ are negative and hence,  
    \begin{equation}
        \label{eq:intersection-of-neg-ribbon-is-one-special-ribbon-2}
        c(\gamma_{a_i}) + 1 \le c(\delta_{a_i}),\ i \in [\ell]. 
    \end{equation}
    In addition, as $\theta_{a_i} \# \theta_{\sigma(a_i)}$ are negative for $i \in \{r, s\}$ we have 
    \begin{equation}
        \label{eq:intersection-of-neg-ribbon-is-one-special-ribbon-3}
        c(\delta_{\sigma(a_i)}) < c(\delta_{a_i}),\ i \in \{t, s\}. 
    \end{equation}
    
    Assume the contrary, i.e. $a_s \neq \sigma(a_r)$; as $\theta_{a_r} \# \theta_{\sigma(a_r)}$ and $\theta_{a_s} \# \theta_{\sigma(a_s)}$ intersect, then the following inequality must hold: 
    \begin{equation}
        \label{eq:intersection-of-neg-ribbon-is-one-special-ribbon-4}
        c(\delta_{\sigma(a_r)}) < c(\delta_{a_s}).
    \end{equation}
    There are two cases: (1) $c(\gamma_{\sigma(a_r)}) > c(\gamma_{\sigma(a_s)})$; (2) $c(\gamma_{\sigma(a_r)}) < c(\gamma_{\sigma(a_s)})$.
    
    In the first case, using the inequalities~\eqref{eq:intersection-of-neg-ribbon-is-one-special-ribbon-2}, \eqref{eq:intersection-of-neg-ribbon-is-one-special-ribbon-4}, we have 
    \begin{equation*}
        c(\gamma_{\sigma(a_s)}) + 1 < c(\delta_{\sigma(a_r)}) < c(\delta_{a_s}).
    \end{equation*}
    Similarly, using the inequalities~\eqref{eq:intersection-of-neg-ribbon-is-one-special-ribbon-2}, \eqref{eq:intersection-of-neg-ribbon-is-one-special-ribbon-3} we have 
    \begin{equation*}
        c(\gamma_{\sigma(a_r)}) + 1 \le c(\delta_{\sigma(a_r)}) < c(\delta_{a_r}).
    \end{equation*}
    Translate the last two inequalities into the lattice. Specify the coordinates of endpoints as follows: $A_{a_i} = (A^1_{a_i}, A^2_{a_i}),\ B_{a_i} = (B^1_{a_i}, B^2_{a_i}),\ i \in [\ell]$; then
    we have (Fig.~\ref{fig:intersection-of-neg-ribbon-is-one-special-ribbon-1-a}):
    \begin{equation*}
        B^1_{\sigma(a_s)} < A^1_{\sigma(a_r)} < A^1_{a_s} \qquad 
        B^1_{\sigma(a_r)} \le A^1_{\sigma(a_r)} < A^1_{a_r}.
    \end{equation*}
    It follows that some minor subpath of  
    $P \in \mathcal{P}(A_{\sigma(a_r)}, B_{\sigma^2(a_r)})$ that belongs to the top or bottom layer inevitably overlaps with minor subpath of 
    either $P_r \in \mathcal{P}(A_{a_r}, B_{\sigma(a_r)})$ or $P_s \in \mathcal{P}(A_{a_s}, B_{\sigma(a_s)})$, which is a contradiction by Lemma~\ref{lem:cousin-paths-avoid-each-others-endpoints}. 

    In the second case, using the inequalities~\eqref{eq:intersection-of-neg-ribbon-is-one-special-ribbon-2},\eqref{eq:intersection-of-neg-ribbon-is-one-special-ribbon-3} it follows that 
    \begin{equation*}
        c(\gamma_{\sigma(a_s)}) + 1 \le c(\delta_{\sigma(a_s)}) < c(\delta_{a_s}).
    \end{equation*}
    Recall that $c(\delta_{a_s}) < c(\delta_{a_r})$; thus using inequalities~\eqref{eq:intersection-of-neg-ribbon-is-one-special-ribbon-2}, \eqref{eq:intersection-of-neg-ribbon-is-one-special-ribbon-3} we have 
    \begin{equation*}
        c(\gamma_{\sigma(a_r)}) + 1 < c(\delta_{\sigma(a_s)}) < c(\delta_{a_r}).
    \end{equation*}
    Translating the last two inequalities into the lattice, we have (Fig.~\ref{fig:intersection-of-neg-ribbon-is-one-special-ribbon-1-b})
    \begin{equation*}
        B_{\sigma_{a_s}} \le A_{\sigma(a_s)} < A_{a_s} \qquad B_{\sigma(a_r)} < A_{\sigma(a_s)} < A_{a_r}.
    \end{equation*}
    It follows that some minor subpath of 
    $P \in \mathcal{P}(A_{\sigma(a_s)}, B_{\sigma^2(a_s)})$ that belongs to top or bottom layer inevitably  overlaps with minor subpath of either $P_r \in \mathcal{P}(A_{a_r}, B_{\sigma(a_r)})$ or $P_s \in \mathcal{P}(A_{a_s}, B_{\sigma(a_s)})$, which is a contradiction by Lemma~\ref{lem:cousin-paths-avoid-each-others-endpoints}. Thus, $a_s = \sigma(a_r)$. 

    Now we claim that the endpoints $A_{\sigma(a_r)}, B_{\sigma(a_r)}$ must be located in different layers. By contradiction, suppose that $A_{\sigma(a_r)}, B_{\sigma(a_r)}$ are in the same layer, then there are  minor subpaths of $P \in \mathcal{P}(A_{a_r}, B_{\sigma(a_r)}), P' \in \mathcal{P}(A_{\sigma(a_r)}, B_{\sigma^2(a_r)})$ such that they belong to top or bottom layers and they overlap, as the ribbons $\theta_{a_r} \# \theta_{\sigma(a_r)},\ \theta_{\sigma(a_r)} \# \theta_{\sigma^2(a_r)},\ \theta_{\sigma(a_r)}$ are negative and $B^1_{\sigma^2(a_r)} < B^1_{\sigma(a_r)} \le A^1_{\sigma(a_r)} < A^1_{a_r}$; thus the endpoints $A_{\sigma(a_r)}, B_{\sigma(a_r)}$ are in different layers, and the ribbon $\theta_{\sigma(a_r)}$ is even. 
\end{proof}

\begin{definition}
    Let $\theta_i, \theta_j \in \Theta$ be ribbons which belong to the same tile $\hat{\theta}_k, k \in \mathbb{Z}$ and $c(\delta_j) < c(\delta_i)$, and there is no other ribbon $\theta_r \in \Theta$ between $\theta_i, \theta_j$, i.e. there is no ribbon $\theta_r$ contained in $\hat{\theta}_k$ such that $c(\delta_j) < c(\delta_r) < c(\delta_i)$. We then call such ribbons {\it consecutive}.
\end{definition}

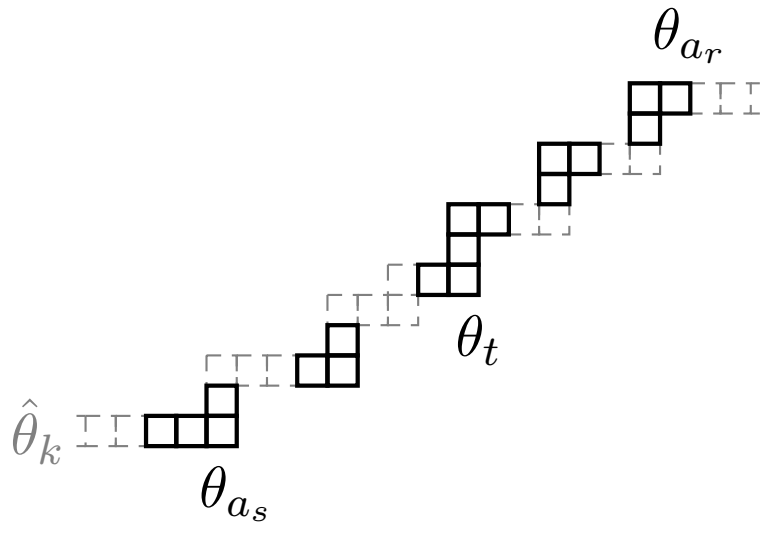
\begin{figure}
    \centering
    \begin{tikzpicture}[y=-1cm, scale=0.2]
        \foreach \x in {7,...,10} {
            \draw[densely dashed, gray] (\x, 1) rectangle ++(1, 1);
        }
        \foreach \x in {7,...,7} {
            \draw[densely dashed, gray] (\x, 2) rectangle ++(1, 1);
        }
        \foreach \x in {4,...,7} {
            \draw[densely dashed, gray] (\x, 3) rectangle ++(1, 1);
        }
        \foreach \x in {4,...,4} {
            \draw[densely dashed, gray] (\x, 4) rectangle ++(1, 1);
        }
        \foreach \x in {1,...,4} {
            \draw[densely dashed, gray] (\x, 5) rectangle ++(1, 1);
        }
        \foreach \x in {1,...,1} {
            \draw[densely dashed, gray] (\x, 6) rectangle ++(1, 1);
        }
        \foreach \x in {-1,...,1} {
            \draw[densely dashed, gray] (\x, 7) rectangle ++(1, 1);
        }
        \foreach \x in {-3,...,-1} {
            \draw[densely dashed, gray] (\x, 8) rectangle ++(1, 1);
        }
        \foreach \x in {-3,...,-3} {
            \draw[densely dashed, gray] (\x, 9) rectangle ++(1, 1);
        }
        \foreach \x in {-7,...,-3} {
            \draw[densely dashed, gray] (\x, 10) rectangle ++(1, 1);
        }
        \foreach \x in {-7,...,-7} {
            \draw[densely dashed, gray] (\x, 11) rectangle ++(1, 1);
        }
        \foreach \x in {-11,...,-7} {
            \draw[densely dashed, gray] (\x, 12) rectangle ++(1, 1);
        }

        \foreach \x in {-7,...,-7} {
            \draw[thick] (\x, 11) rectangle ++(1, 1);
        }
        \foreach \x in {-9,...,-7} {
            \draw[thick] (\x, 12) rectangle ++(1, 1);
        }

        \foreach \x in {-3,...,-3} {
            \draw[thick] (\x, 9) rectangle ++(1, 1);
        }
        \foreach \x in {-4,...,-3} {
            \draw[thick] (\x, 10) rectangle ++(1, 1);
        }

        \foreach \x in {1,...,2} {
            \draw[thick] (\x, 5) rectangle ++(1, 1);
        }
        \foreach \x in {1,...,1} {
            \draw[thick] (\x, 6) rectangle ++(1, 1);
        }
        \foreach \x in {0,...,1} {
            \draw[thick] (\x, 7) rectangle ++(1, 1);
        }

        \foreach \x in {4,...,5} {
            \draw[thick] (\x, 3) rectangle ++(1, 1);
        }
        \foreach \x in {4,...,4} {
            \draw[thick] (\x, 4) rectangle ++(1, 1);
        }

        \foreach \x in {7,...,8} {
            \draw[thick] (\x, 1) rectangle ++(1, 1);
        }
        \foreach \x in {7,...,7} {
            \draw[thick] (\x, 2) rectangle ++(1, 1);
        }

        \foreach \x in {12,...,13} {
            \draw[thin, gray, densely dashed] (-11, \x) -- (-11.3, \x);
        }
        \foreach \x in {1,...,2} {
            \draw[thin, gray, densely dashed] (11, \x) -- (11.3, \x);
        }


        \fill[gray] (-11, 12.5) node[left] {$\hat{\theta}_k$};
        \fill (-6, 13) node[below] {$\theta_{a_s}$};
        \fill (2, 8) node[below] {$\theta_t$};
        \fill (9, 1) node[above] {$\theta_{a_r}$};

    \end{tikzpicture}
    \caption{Ribbon example.}
    \label{fig:ribbon-between-is-always-odd-2}
\end{figure}

\begin{lemma}
    \label{lem:ribbon-between-is-always-odd}
    Assume that $\lambda/\mu$ is a negative shape. 
    Let $\mathbf{F} = (\mathbf{P}, \mathbf{C}) = \psi(\mathbf{T}) = \psi((T^{\circ}, \widetilde{T}, T^{\bullet}))$ such that $\mathbf{P} \in \mathcal{P}(\mathbf{A}, \mathbf{B}, \sigma)$ where $\sigma$ is some permutation. Let $\sigma_0 = (a_1, \ldots, a_{\ell}), \ell \ge 1$ be some cycle in cycle decomposition of  $\sigma$. The ribbons $(\theta_{a_i})$ are in the same tile $\hat{\theta}_k, k \in \mathbb{Z}$ by Lemma~\ref{lem:th-i-th-si-are-in-the-same-ribbon}. 
    Let $\theta_{a_r}, \theta_{a_s}$ be consecutive ribbons in $(\theta_{a_i})$ with $c(\delta_{a_s}) < c(\delta_{a_r})$. 
    Then $\pi(\delta_{a_s}) = (\cdot, x),\ \pi(\gamma_{a_r}) = (\overline{x}, \cdot),\ x \in \{0, 1\},\ \overline{x} \in \{0, 1\} \setminus \{x\},\ i < \ell$.
\end{lemma}

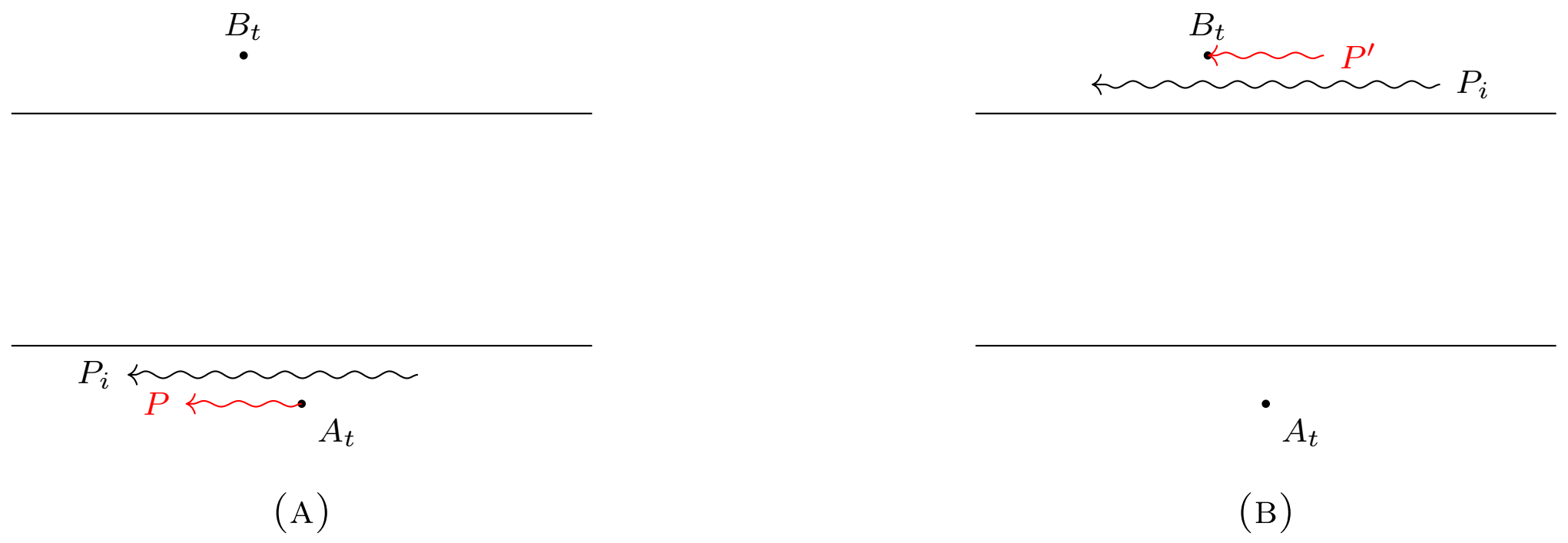
\begin{figure}
    \centering
    \begin{subfigure}[b]{0.5\textwidth}
        \centering
        \begin{tikzpicture}[scale=0.5, font=\scriptsize]
            \fill[black] (0, 0.5) node[left] {$P_i$};
            \fill[red] (1, 0) node[left] {$P$};

            \fill (3, 0) circle (2pt) node[below right] {$A_t$};
            \fill (2, 6) circle (2pt) node[above] {$B_t$};

            \draw (-2, 5) -- (8, 5);
            \draw (-2, 1) -- (8, 1);

            \draw[decorate, decoration={snake, amplitude=0.3mm, segment length=3mm}, ->, black] (5, 0.5) -- (0, 0.5);
            
            \draw[decorate, decoration={snake, amplitude=0.25mm, segment length=3mm}, ->, red] (3, 0) -- (1, 0);
        \end{tikzpicture}
        \subcaption{}
        \label{fig:ribbon-between-is-always-odd-1-a}
    \end{subfigure}
    \begin{subfigure}[b]{0.49\textwidth}
        \centering
        \begin{tikzpicture}[scale=0.5, font=\scriptsize]
            \fill[black] (6, 5.5) node[right] {$P_i$};
            \fill[red] (4, 6) node[right] {$P'$};

            \fill (3, 0) circle (2pt) node[below right] {$A_t$};
            \fill (2, 6) circle (2pt) node[above] {$B_t$};

            \draw (-2, 5) -- (8, 5);
            \draw (-2, 1) -- (8, 1);

            \draw[decorate, decoration={snake, amplitude=0.3mm, segment length=3mm}, ->, black] (6, 5.5) -- (0, 5.5);
            \draw[decorate, decoration={snake, amplitude=0.25mm, segment length=3mm}, ->, red] (4, 6) -- (2, 6);
        \end{tikzpicture}
        \subcaption{}
        \label{fig:ribbon-between-is-always-odd-1-b}
    \end{subfigure}
    \caption{Paths examples.}
    \label{fig:ribbon-between-is-always-odd-1}
\end{figure}

\begin{proof}
    Assume to the contrary, that there are consecutive ribbons $\theta_{a_r}, \theta_{a_s}, r, s \in [\ell]$ such that $c(\delta_{a_s}) < c(\delta_{a_r})$, and $\pi(\delta_{a_s}) = (\cdot, x),\ \pi(\gamma_{a_r}) = (x, \cdot),\ x \in \{0, 1\}$ (Fig.~\ref{fig:ribbon-between-is-always-odd-2}). From Lemma~\ref{lem:restrictions-for-ribbons} it follows that there must be ribbon $\theta_t \in \Theta$ which is contained in the $\hat{\theta}_k$, with the following properties: (1) $\theta_t \notin (\theta_{a_i})$; (2) $\pi(\delta_t) = (\cdot, y),\ \pi(\gamma_t) = (y, \cdot),\ y \in \{0, 1\}$ where $\delta_t, \gamma_t$ are head and tail of $\theta_t$; (3) located between ribbons $\theta_{a_r}, \theta_{a_s}$ (i.e. $c(\delta_{a_s}) < c(\delta_t) < c(\delta_{a_r})$). It follows that $A_t, B_t$ are in different layers. Let $P_i \in \mathcal{P}(A_{a_i}, B_{\sigma(a_i)}), i \in [\ell]$. As $c(\delta_{a_s}) < c(\delta_t) < c(\delta_{a_r})$ and ribbons $\theta_{a_r}, \theta_{a_s}, \theta_t$ do not intersect, there is $P_i, i \in [\ell]$ such that its  minor subpath in top or bottom layer overlaps with minor subpath that starts at $A_t$ (path $P$ in Fig.~\ref{fig:ribbon-between-is-always-odd-1-a}) or terminates at $B_t$ (path $P'$ in Fig.~\ref{fig:ribbon-between-is-always-odd-1-b}), which is a contradiction by Lemma~\ref{lem:cousin-paths-avoid-each-others-endpoints}. 
\end{proof}

\begin{lemma}
    \label{lem:sgn-rule-for-one-cycle}
    Assume that $\lambda/\mu$ is a negative shape. 
    Let $\mathbf{F} = (\mathbf{P}, \mathbf{C}) = \psi(\mathbf{T}) = \psi((T^{\circ}, \widetilde{T}, T^{\bullet}))$ such that $\mathbf{P} \in \mathcal{P}(\mathbf{A}, \mathbf{B}, \sigma)$ where $\sigma$ is some permutation. Let $\sigma_0 = (a_1, \ldots, a_{\ell}), \ell \ge 1$ be some cycle in cycle decomposition of $\sigma$. Let $\rho^{\circ}_1, \ldots, \rho^{\circ}_{s_1} \in \Theta^{\circ}$, $\rho^{\bullet}_1, \ldots, \rho^{\bullet}_{s_2} \in \Theta^{\bullet}$, $\widetilde{\rho}_1, \ldots, \widetilde{\rho}_{s_3} \in \widetilde{\Theta},\ s_1, s_2, s_3 \ge 0$ be ribbons filled by paths $P_1, \ldots, P_{\ell}$ where $P_i \in \mathcal{P}(A_{a_i}, B_{\sigma(a_i)}),\ i \in [\ell]$. Then the following identity holds:
    \begin{equation*}
        \prod_{i = 1}^{\ell} \mathrm{wt}(P_i) = \mathrm{sgn}(\sigma_0)\prod_{i = 1}^{\ell}\mathrm{sgn}_{\pi}(\theta_{a_i}) \prod_{i = 1}^{s_1} \mathrm{wt}(T^{\circ}_{\rho^{\circ}_i})\prod_{i = 1}^{s_2} \mathrm{wt}(\widetilde{T}_{\widetilde{\rho}_i})\prod_{i = 1}^{s_3} \mathrm{wt}(T^{\bullet}_{\rho^{\bullet}_i})
    \end{equation*}
\end{lemma}

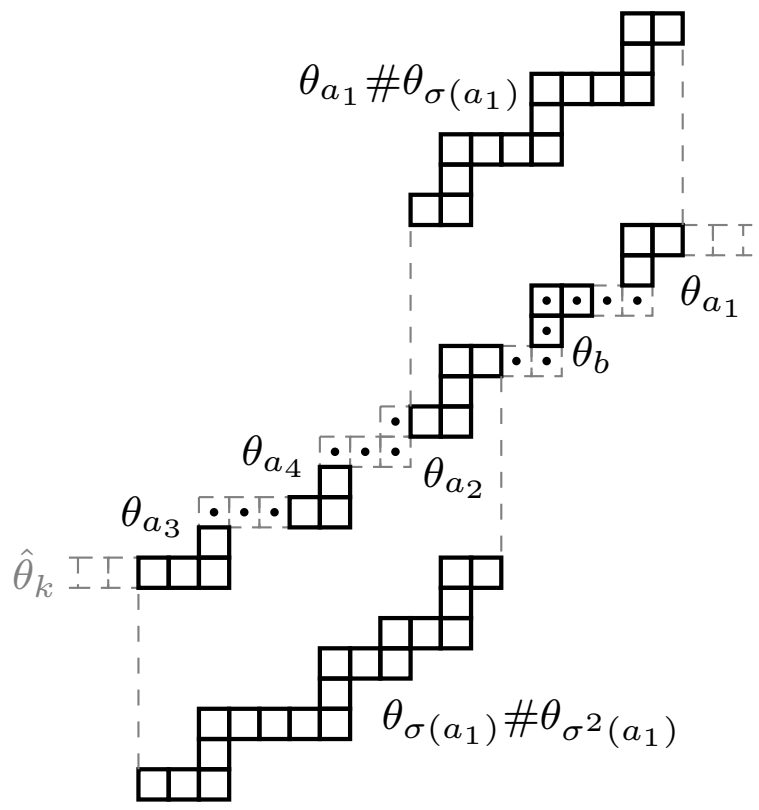
\begin{figure}
    \centering
    \begin{tikzpicture}[y=-1cm, scale=0.2, font=\scriptsize]
        \foreach \x in {7,...,8} {
            \draw[thick] (\x, -6) rectangle ++(1, 1);
        }
        \foreach \x in {7,...,7} {
            \draw[thick] (\x, -5) rectangle ++(1, 1);
        }
        \foreach \x in {4,...,7} {
            \draw[thick] (\x, -4) rectangle ++(1, 1);
        }
        \foreach \x in {4,...,4} {
            \draw[thick] (\x, -3) rectangle ++(1, 1);
        }
        \foreach \x in {1,...,4} {
            \draw[thick] (\x, -2) rectangle ++(1, 1);
        }
        \foreach \x in {1,...,1} {
            \draw[thick] (\x, -1) rectangle ++(1, 1);
        }
        \foreach \x in {0,...,1} {
            \draw[thick] (\x, 0) rectangle ++(1, 1);
        }

        \foreach \x in {7,...,10} {
            \draw[densely dashed, gray] (\x, 1) rectangle ++(1, 1);
        }
        \foreach \x in {7,...,7} {
            \draw[densely dashed, gray] (\x, 2) rectangle ++(1, 1);
        }
        \foreach \x in {4,...,7} {
            \draw[densely dashed, gray] (\x, 3) rectangle ++(1, 1);
        }
        \foreach \x in {4,...,4} {
            \draw[densely dashed, gray] (\x, 4) rectangle ++(1, 1);
        }
        \foreach \x in {1,...,4} {
            \draw[densely dashed, gray] (\x, 5) rectangle ++(1, 1);
        }
        \foreach \x in {1,...,1} {
            \draw[densely dashed, gray] (\x, 6) rectangle ++(1, 1);
        }
        \foreach \x in {-1,...,1} {
            \draw[densely dashed, gray] (\x, 7) rectangle ++(1, 1);
        }
        \foreach \x in {-3,...,-1} {
            \draw[densely dashed, gray] (\x, 8) rectangle ++(1, 1);
        }
        \foreach \x in {-3,...,-3} {
            \draw[densely dashed, gray] (\x, 9) rectangle ++(1, 1);
        }
        \foreach \x in {-7,...,-3} {
            \draw[densely dashed, gray] (\x, 10) rectangle ++(1, 1);
        }
        \foreach \x in {-7,...,-7} {
            \draw[densely dashed, gray] (\x, 11) rectangle ++(1, 1);
        }
        \foreach \x in {-11,...,-7} {
            \draw[densely dashed, gray] (\x, 12) rectangle ++(1, 1);
        }

        \foreach \x in {-6.5,...,-4.5} {
            \fill (\x, 10.5) circle (4pt) node {};
        }
        \foreach \x in {-2.5,...,-0.5} {
            \fill (\x, 8.5) circle (4pt) node {};
        }
        \foreach \x in {-0.5,...,-0.5} {
            \fill (\x, 7.5) circle (4pt) node {};
        }
        \foreach \x in {3.5,...,4.5} {
            \fill (\x, 5.5) circle (4pt) node {};
        }
        \foreach \x in {4.5,...,4.5} {
            \fill (\x, 4.5) circle (4pt) node {};
        }
        \foreach \x in {4.5,...,7.5} {
            \fill (\x, 3.5) circle (4pt) node {};
        }

        \foreach \x in {1,...,2} {
            \draw[thick] (\x, 12) rectangle ++(1, 1);
        }
        \foreach \x in {1,...,1} {
            \draw[thick] (\x, 13) rectangle ++(1, 1);
        }
        \foreach \x in {-1,...,1} {
            \draw[thick] (\x, 14) rectangle ++(1, 1);
        }
        \foreach \x in {-3,...,-1} {
            \draw[thick] (\x, 15) rectangle ++(1, 1);
        }
        \foreach \x in {-3,...,-3} {
            \draw[thick] (\x, 16) rectangle ++(1, 1);
        }
        \foreach \x in {-7,...,-3} {
            \draw[thick] (\x, 17) rectangle ++(1, 1);
        }
        \foreach \x in {-7,...,-7} {
            \draw[thick] (\x, 18) rectangle ++(1, 1);
        }
        \foreach \x in {-9,...,-7} {
            \draw[thick] (\x, 19) rectangle ++(1, 1);
        }

        \foreach \x in {-7,...,-7} {
            \draw[thick] (\x, 11) rectangle ++(1, 1);
        }
        \foreach \x in {-9,...,-7} {
            \draw[thick] (\x, 12) rectangle ++(1, 1);
        }

        \foreach \x in {-3,...,-3} {
            \draw[thick] (\x, 9) rectangle ++(1, 1);
        }
        \foreach \x in {-4,...,-3} {
            \draw[thick] (\x, 10) rectangle ++(1, 1);
        }

        \foreach \x in {1,...,2} {
            \draw[thick] (\x, 5) rectangle ++(1, 1);
        }
        \foreach \x in {1,...,1} {
            \draw[thick] (\x, 6) rectangle ++(1, 1);
        }
        \foreach \x in {0,...,1} {
            \draw[thick] (\x, 7) rectangle ++(1, 1);
        }

        \foreach \x in {4,...,5} {
            \draw[thick] (\x, 3) rectangle ++(1, 1);
        }
        \foreach \x in {4,...,4} {
            \draw[thick] (\x, 4) rectangle ++(1, 1);
        }

        \foreach \x in {7,...,8} {
            \draw[thick] (\x, 1) rectangle ++(1, 1);
        }
        \foreach \x in {7,...,7} {
            \draw[thick] (\x, 2) rectangle ++(1, 1);
        }

        \foreach \x in {12,...,13} {
            \draw[thin, gray, densely dashed] (-11, \x) -- (-11.3, \x);
        }
        \foreach \x in {1,...,2} {
            \draw[thin, gray, densely dashed] (11, \x) -- (11.3, \x);
        }

        \draw[dashed, gray] (9, 1) -- (9, -5);
        \draw[dashed, gray] (0, 7) -- (0, 1);

        \draw[dashed, gray] (3, 6) -- (3, 12);
        \draw[dashed, gray] (-9, 19) -- (-9, 13);


        \fill[gray] (-11, 12.5) node[left] {$\hat{\theta}_k$};
        \fill (-8.5, 12) node[above] {$\theta_{a_3}$};
        \fill (-4.5, 10) node[above] {$\theta_{a_4}$};
        \fill (1.5, 8) node[below] {$\theta_{a_2}$};
        \fill (6, 4) node[below] {$\theta_b$};
        \fill (10, 2) node[below] {$\theta_{a_1}$};
        \fill (0, -5) node[below] {$\theta_{a_1} \# \theta_{\sigma(a_1)}$};
        \fill (4, 16) node[below] {$\theta_{\sigma(a_1)} \# \theta_{\sigma^2(a_1)}$};

    \end{tikzpicture}
    \caption{Here $\sigma_0 = (a_1, a_2, a_3, a_4)$; $\theta_b \in \Theta \setminus (\theta_{a_i})$; $s = 2$; cells of intermediate ribbons are highlighted by dots.}
    \label{fig:sgn-rule-for-one-cycle-1}
\end{figure}

\begin{proof}
    It follows from the construction of the lattice that 
    \begin{equation*}
        \prod_{i = 1}^{\ell} |\mathrm{wt}(P_i)| = \prod_{i = 1}^{s_1} \mathrm{wt}(T^{\circ}_{\rho^{\circ}_i}) \prod_{i = 1}^{s_2} \mathrm{wt}(\widetilde{T}_{\widetilde{\rho}_i}) \prod_{i = 1}^{s_3} \mathrm{wt}(T^{\bullet}_{\rho^{\bullet}_i}).
    \end{equation*}
    Define $\alpha(P_i)$ as the sign of $\mathrm{wt}(P_i), i \in [\ell]$; i.e. 
    \begin{equation*}
        \mathrm{wt}(P_i) = \alpha(P_i)\cdot |\mathrm{wt}(P_i)|, \quad i \in [\ell].
    \end{equation*}
    From Lemmas~\ref{lem:weight-of-ribbons-filled-by-same-path-positive}, \ref{lem:weight-of-ribbons-filled-by-same-path-negative} it follows that 
    \begin{equation*}
        \alpha(P_i) = \mathrm{sgn}_{\pi}(\theta_{a_i} \# \theta_{\sigma(a_i)}), \quad i \in [\ell].
    \end{equation*}
    Now we aim to prove that
    \begin{equation*}
        \prod_{i = 1}^{\ell}\mathrm{sgn}_{\pi}(\theta_{a_i} \# \theta_{\sigma(a_i)}) = \mathrm{sgn}(\sigma_0) \prod_{i = 1}^{\ell} \mathrm{sgn}_{\pi}(\theta_{a_i}).
    \end{equation*}
    We know that $\mathrm{sgn}_{\pi}(\theta_{a_i} \# \theta_{\sigma(a_i)})$ is positive if ribbon $\theta_{a_i} \# \theta_{\sigma(a_i)}$ is positive, and hence the sign is affected only by the ribbons $\theta_{a_i} \# \theta_{\sigma(a_i)}$ which are negative. Without loss of generality assume that $\theta_{a_1}$ is the ribbon with maximum head content among $(\theta_{a_i})$; then it follows from Lemma~\ref{lem:intersection-of-neg-ribbon-is-one-special-ribbon} that there is $s \in [\ell]$ such that all the ribbons $\theta_{a_i} \# \theta_{\sigma(a_i)},\ i \in [s]$ are negative, and all the ribbons $\theta_{a_i} \# \theta_{\sigma(a_i)},\ i > s$ are positive; hence we have
    \begin{equation*}
        \prod_{i = 1}^{\ell}\mathrm{sgn}_{\pi}(\theta_{a_i} \# \theta_{\sigma(a_i)}) = \prod_{i = 1}^{s}\mathrm{sgn}_{\pi}(\theta_{a_i} \# \theta_{\sigma(a_i)}).
    \end{equation*}
    From Lemma~\ref{lem:intersection-of-neg-ribbon-is-one-special-ribbon} it follows that the intersection of negative ribbons $\theta_{a_i} \# \theta_{\sigma(a_i)}$ and $\theta_{\sigma(a_i)} \# \theta_{\sigma^2(a_i)}$ is the ribbon $\theta_{\sigma(a_i)},\ i < s$ and it is even. Fix $j \in [s]$, assume that $\theta_{a_p}, \theta_{a_q},\ p, q \in [\ell]$ are consecutive sub-ribbons of $\theta_{a_j} \# \theta_{\sigma(a_j)}$ with $c(\delta_p) < c(\delta_q)$. 
    Consider the ribbon constructed from cells between $\theta_{a_p}$ and $\theta_{a_q}$; i.e. the ribbon $[c(\gamma_q) - 1, c(\delta_p) + 1]_{\hat{\theta}}$. Let us call ribbons that consist of the cells between consecutive ribbons in $(\theta_{a_i})$ as `intermediate ribbons' (Fig.~\ref{fig:sgn-rule-for-one-cycle-1}). Note that from Lemma~\ref{lem:ribbon-between-is-always-odd} it follows that the intermediate ribbons are odd.
    
    
    
    Suppose that $\theta^j_1, \ldots, \theta^j_{r_j}$ are all the ribbons from $(\theta_{a_i})$ contained in $\theta_{a_j} \# \theta_{\sigma(a_j)}$. 
    In light of the preceding observations, we obtain
    \begin{equation*}
        \mathrm{sgn}_{\pi}(\theta_{a_j} \# \theta_{\sigma(a_j)}) = (-1)^{r_j - 1}\prod_{i = 1}^{r_j}\mathrm{sgn}_{\pi}(\theta^j_i).
    \end{equation*}
    Note that $(-1)^{r_j - 1}$ appears because there are $r_j - 1$ intermediate ribbons in $\theta_{a_j} \# \theta_{\sigma(a_j)}$. Thus, 
    \begin{align*}
        \prod_{j = 1}^{s}\mathrm{sgn}_{\pi}(\theta_{a_j} \# \theta_{\sigma(a_j)}) &= \prod_{j = 1}^{s}(-1)^{r_j - 1}\prod_{i = 1}^{r_j}\mathrm{sgn}_{\pi}(\theta^j_i) \\ 
        &= (-1)^{\ell - 1} \prod_{i = 1}^{\ell} \mathrm{sgn}_{\pi}(\theta_{a_i}) \prod_{i = 1}^{s - 1}\mathrm{sgn}_{\pi}(\theta_{\sigma(a_i)}) \\
        &= \mathrm{sgn}(\sigma_0) \prod_{i = 1}^{\ell} \mathrm{sgn}_{\pi}(\theta_{a_i}).
    \end{align*}
    The last equation follows from the facts that $\theta_{\sigma(a_i)}$ is even for $i < s$; and $\theta_{\sigma(a_i)}, i < s$ are the only possible intersections of negative ribbons in $(\theta_{a_i} \# \theta_{\sigma(a_i)})$; thus their signs are counted twice; but noting that their signs are positive we can ignore them, and the result follows. 
\end{proof}

\begin{lemma}
    \label{lem:weight-preserving-negative}
    Assume that $\lambda/\mu$ is a negative shape. 
    Let $\mathbf{F} = (\mathbf{P}, \mathbf{C}) = \psi(\mathbf{T}) = \psi((T^{\circ}, \widetilde{T}, T^{\bullet}))$. Then the following identity holds: 
    \begin{equation*}
        \mathrm{sgn}(\mathbf{F}) \cdot \mathrm{wt}(\mathbf{F}) = \mathrm{sgn}_{\pi}(\lambda/\mu) \cdot \mathrm{wt}(\mathbf{T}).
    \end{equation*}
\end{lemma}
\begin{proof}
    From the construction of the lattice it follows that 
    \begin{equation*}
        |\mathrm{wt}(\mathbf{P})| \cdot |\mathrm{wt}(\mathbf{C}| = \mathrm{wt}(\mathbf{T}).
    \end{equation*}
    Define $\alpha(\mathbf{P}), \alpha(\mathbf{C})$ as follows:
    \begin{equation*}
        |\mathrm{wt}(\mathbf{P})| = \alpha(\mathbf{P}) \cdot \mathrm{wt}(\mathbf{P}), \qquad |\mathrm{wt}(\mathbf{C})| = \alpha(\mathbf{C}) \cdot \mathrm{wt}(\mathbf{C}).
    \end{equation*}
    Assume that $\mathbf{P} \in \mathcal{P}(\mathbf{A}, \mathbf{B}, \sigma)$ where $\sigma$ is some permutation. Let $\sigma = \sigma_1 \ldots \sigma_{\ell}$ be cycle decomposition of $\sigma$, such that $\sigma_i = (a^i_1, \ldots, a^i_{\ell_i}), \ell_i \ge 1,\ i \in \ell$. From Lemma~\ref{lem:sgn-rule-for-one-cycle} it follows that
    \begin{align*}
        \alpha(\mathbf{P}) = \prod_{i = 1}^{\ell} \mathrm{sgn}(\sigma_i) \prod_{j = 1}^{\ell_i}\mathrm{sgn}_{\pi}(\theta_{a^i_j}) = \mathrm{sgn}(\sigma) \cdot \mathrm{sgn}_{\pi}(\lambda/\mu).
    \end{align*}
    Next we discuss sign generated by cycles, i.e. $\alpha(\mathbf{C})$. From Lemma~\ref{lem:cycle-corresponds-to-special-ribbons} it follows that the ribbons from $\Theta^{\circ} \cup \Theta^{\bullet}$ filled by cycles are odd; thus 
    \begin{equation*}
        \alpha(\mathbf{C}) = (-1)^{|\mathbf{C}|} = \mathrm{sgn}(\mathbf{C}).
    \end{equation*}
    Recall that $\mathrm{sgn}(\mathbf{F}) = \mathrm{sgn}(\sigma)\, \mathrm{sgn}(\mathbf{C})$. Combining the preceding observations we have:
    \begin{align*}
        \mathrm{sgn}(\mathbf{F}) \cdot \alpha(\mathbf{P}) \cdot \alpha(\mathbf{C}) &= \mathrm{sgn}(\sigma) \mathrm{sgn}(\mathbf{C}) \cdot \mathrm{sgn}(\sigma) \cdot \mathrm{sgn}_{\pi}(\lambda/\mu) \cdot \mathrm{sgn}(\mathbf{C}) = \mathrm{sgn}_{\pi}(\lambda/\mu), 
    \end{align*}
    concluding the proof.
\end{proof}

\begin{lemma}
    \label{lem:main-weight-preserving-bijection}
    The map $\psi$ is a weight preserving bijection.
\end{lemma}
\begin{proof}
    Directly follows from Lemmas~\ref{lem:psi-is-inverse-of-phi}, \ref{lem:weight-preserving-positive} and \ref{lem:weight-preserving-negative}.
\end{proof}

\subsection{Concluding the proof of Hamel--Goulden-type formulas.}
\label{subsec:hg-formulas-proof}

The following corollary is a direct consequence of Lemma~\ref{lem:main-weight-preserving-bijection}.
\begin{corollary}
    \label{cor:bijection-major}
    Let $\Theta = (\theta_i)$ be the outer decomposition of $\lambda/\mu$ induced by pipe function $\pi$, and $\mathbf{a}, \mathbf{b}$ be row $\mathsf{G}$-flags. Let $\mathcal{F}(\mathbf{A}, \mathbf{B})$ be set of flows on the lattice $TL_{\lambda/\mu}(\pi, m, \mathbf{a}, \mathbf{b})$. Then we have 
    \begin{equation*}
        \mathrm{sgn}_{\pi}(\lambda/\mu) \cdot \mathsf{G}^{\mathbf{a}, \mathbf{b}}_{\lambda/\mu}(\mathbf{x}_m; \mathbf{y}, \mathbf{z}) = \sum\limits_{\mathbf{F} \in \mathcal{F}(\mathbf{A}, \mathbf{B})}\mathrm{sgn}(\mathbf{F}) \cdot wt(\mathbf{F}).
    \end{equation*}
\end{corollary}

\begin{corollary}
    \label{cor:bijection-minor}
    Let $\theta_i, \theta_j \in \Theta$. There is a weight-preserving bijection between the tableaux in $TT_{\theta_i \# \theta_j}(\mathbf{a}^{i, j}, \mathbf{b}^{i, j}, m)$ and flows in $\mathcal{F}(\{A_i\}, \{B_j\})$ on the lattice $TL_{\lambda/\mu}(\pi, m, \mathbf{a}, \mathbf{b})$, i.e. the following identity holds:
    \begin{equation*}
        \mathrm{sgn}_{\pi}(\theta_i \# \theta_j) \cdot \mathsf{G}^{\mathbf{a}^{ij}, \mathbf{b}^{ij}}_{\theta_i \# \theta_j}(\mathbf{x}_m; \mathbf{y}, \mathbf{z}) = \sum\limits_{\mathbf{F}\in \mathcal{F}(\{A_i\}, \{B_j\})}\mathrm{sgn}(\mathbf{F}) \cdot wt(\mathbf{F}).
    \end{equation*}
\end{corollary}
\begin{proof}
    One can construct the lattice for $TL_{\theta_i \# \theta_j}(\pi, m, \mathbf{a}^{i, j}, \mathbf{b}^{i, j})$, which in fact coincides with the lattice $TL_{\lambda/\mu}(\pi, m, \mathbf{a}, \mathbf{b})$ with all endpoints removed except $A_i, B_j$. Hence, the flows on the lattice $TL_{\theta_i \# \theta_j}(\pi, m, \mathbf{a}^{i, j}, \mathbf{b}^{i, j})$ are indeed the flows in $\mathcal{F}(\{A_i\}, \{B_j\})$. Then the statement follows from Lemma~\ref{lem:main-weight-preserving-bijection}.
\end{proof}

\begin{figure}
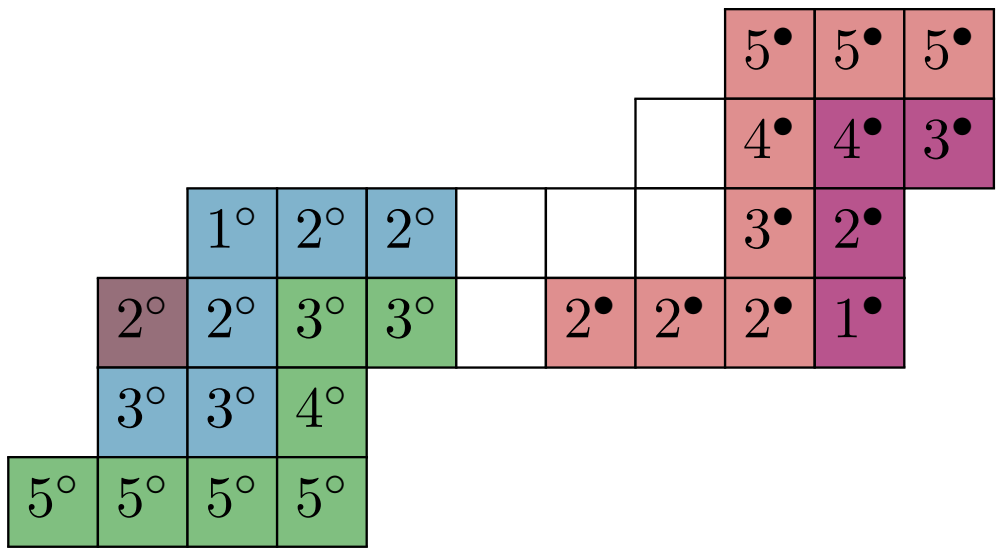

    \centering
    \begin{ytableau}
        \none & \none & \none & \none & \none & \none & \none & \none & *(blue!50!green!25!red!50) 5^{\bullet} & *(blue!50!green!25!red!50) 5^{\bullet} & *(blue!50!green!25!red!50) 5^{\bullet} \\
        \none & \none & \none & \none & \none & \none & \none & *(white) & *(blue!50!green!25!red!50) 4^{\bullet} & *(blue!90!green!40!red!70) 4^{\bullet} & *(blue!90!green!40!red!70) 3^{\bullet} \\
        \none & \none & *(blue!60!green!50) 1^{\circ} & *(blue!60!green!50) 2^{\circ} & *(blue!60!green!50) 2^{\circ} & *(white)  & *(white)  & *(white)  & *(blue!50!green!25!red!50) 3^{\bullet} & *(blue!90!green!40!red!70) 2^{\bullet} \\
        \none & *(blue!55!green!55!red!75) 2^{\circ} & *(blue!60!green!50) 2^{\circ} & *(green!50!black!50) 3^{\circ} & *(green!50!black!50) 3^{\circ} & *(white) & *(blue!50!green!25!red!50) 2^{\bullet} & *(blue!50!green!25!red!50) 2^{\bullet} & *(blue!50!green!25!red!50) 2^{\bullet} & *(blue!90!green!40!red!70) 1^{\bullet} \\
        \none & *(blue!60!green!50) 3^{\circ} & *(blue!60!green!50) 3^{\circ} & *(green!50!black!50) 4^{\circ} \\ 
        *(green!50!black!50) 5^{\circ} & *(green!50!black!50) 5^{\circ} & *(green!50!black!50) 5^{\circ} & *(green!50!black!50) 5^{\circ}
    \end{ytableau}
    \caption{Left and right tableaux $T^{\bullet}$ and $T^{\circ}$.}
    \label{fig:ex-ribbon-filling-vs-path1}
\end{figure}

\begin{figure}
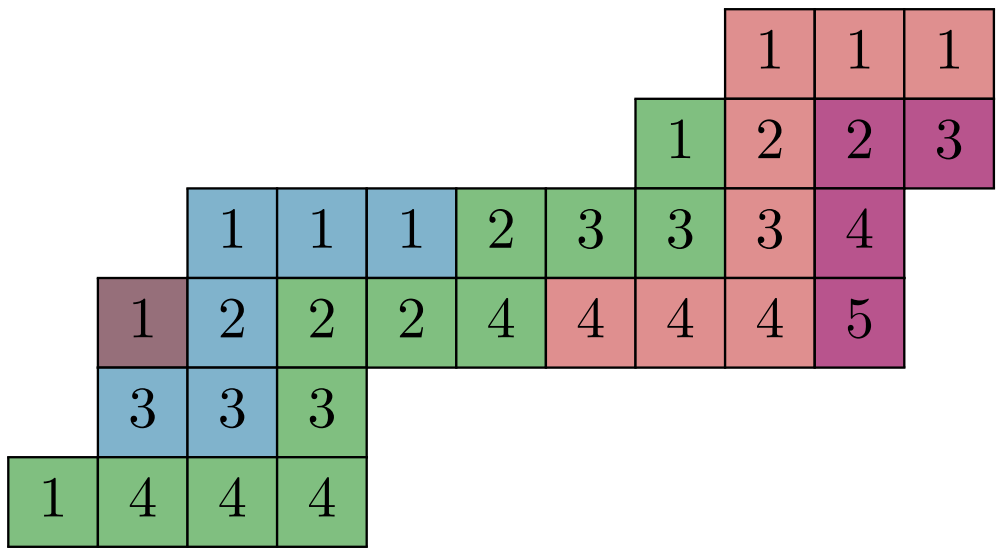

    \centering
    \begin{ytableau}
        \none & \none & \none & \none & \none & \none & \none & \none & *(blue!50!green!25!red!50) 1 & *(blue!50!green!25!red!50) 1 & *(blue!50!green!25!red!50) 1 \\
        \none & \none & \none & \none & \none & \none & \none & *(green!50!black!50) 1 & *(blue!50!green!25!red!50) 2 & *(blue!90!green!40!red!70) 2 & *(blue!90!green!40!red!70) 3 \\
        \none & \none & *(blue!60!green!50) 1 & *(blue!60!green!50) 1 & *(blue!60!green!50) 1 & *(green!50!black!50) 2  & *(green!50!black!50) 3 & *(green!50!black!50) 3  & *(blue!50!green!25!red!50) 3 & *(blue!90!green!40!red!70) 4 \\
        \none & *(blue!55!green!55!red!75) 1 & *(blue!60!green!50) 2 & *(green!50!black!50) 2 & *(green!50!black!50) 2 & *(green!50!black!50) 4 & *(blue!50!green!25!red!50) 4 & *(blue!50!green!25!red!50) 4 & *(blue!50!green!25!red!50) 4 & *(blue!90!green!40!red!70) 5 \\
        \none & *(blue!60!green!50) 3 & *(blue!60!green!50) 3 & *(green!50!black!50) 3 \\ 
        *(green!50!black!50) 1 & *(green!50!black!50) 4 & *(green!50!black!50) 4 & *(green!50!black!50) 4
    \end{ytableau}
    \caption{Central tableau $\widetilde{T}$.}
    \label{fig:ex-ribbon-filling-vs-path}
\end{figure}

\begin{figure}
    \centering
\begin{tikzpicture}[scale=0.6, transform shape]
        \foreach \x in {-8,...,8} {
        \foreach \y in {-8,...,15} {
            \draw[dotted, thin] (\x,\y) rectangle ++(1,1);
        }
        }
        \foreach \x in {-2, ..., 8} {
            \fill[black] (-8.5,\x) circle (0pt) node[left] {$\x$};
        }
        \fill[blue, opacity=0.05] (-4,-8) rectangle (-2, 16);
        \fill[blue, opacity=0.05] (0,-8) rectangle (1, 16);
        \fill[blue, opacity=0.05] (3, -8) rectangle (6, 16);

        \foreach \x in {-2,...,0} {
            \foreach \y in {10,...,15} {
                \fill[black] (\x,\y) circle (1.5pt) node[left] {};
            }
        }
        \foreach \x in {1,...,3} {
            \foreach \y in {10,...,16} {
                \fill[black] (\x,\y) circle (1.5pt) node[left] {};
            }
        }
        \foreach \x in {4,...,5} {
            \foreach \y in {10,...,15} {
                \fill[black] (\x,\y) circle (1.5pt) node[left] {};
            }
        }
        \foreach \x in {6,...,9} {
            \foreach \y in {10,...,14} {
                \fill[black] (\x,\y) circle (1.5pt) node[left] {};
            }
        }

        \foreach \x in {-8,...,-4} {
            \foreach \y in {-2,...,2} {
                \fill[black] (\x,\y) circle (1.5pt) node[left] {};
            }
        }
        \foreach \y in {-1,...,3} {
            \fill[black] (-3, \y) circle (1.5pt) node[left] {};
        }
        \foreach \x in {-2,...,-1} {
            \foreach \y in {0,...,4} {
                \fill[black] (\x, \y) circle (1.5pt) node[left] {};
            }
        }
        \foreach \y in {0,...,5} {
            \fill[black] (0, \y) circle (1.5pt) node[left] {};
        }
        \foreach \x in {1,...,2} {
            \foreach \y in {1,...,5} {
                \fill[black] (\x, \y) circle (1.5pt) node[left] {};
            }
        }
        \foreach \y in {1,...,6} {
                \fill[black] (3, \y) circle (1.5pt) node[left] {};
            }
        \foreach \y in {2,...,7} {
            \fill[black] (4, \y) circle (1.5pt) node[left] {};
        }
        \foreach \y in {3,...,8} {
            \fill[black] (5, \y) circle (1.5pt) node[left] {};
        }
        \foreach \x in {6,...,9} {
            \foreach \y in {4,...,8} {
                \fill[black] (\x, \y) circle (1.5pt) node[left] {};
            }
        }

        \foreach \x in {-8,...,-1} {
            \foreach \y in {-4,...,-8} {
                \fill[black] (\x,\y) circle (1.5pt) node[left] {};
            }
        }
        \foreach \x in {0,...,2} {
            \foreach \y in {-4,...,-7} {
                \fill[black] (\x,\y) circle (1.5pt) node[left] {};
            }
        }
    
        \draw[line width=0.2mm, black, dashed] (-8, -8) -- (-1, -8) -- (-1, -7) -- (2, -7) -- (2, 1) -- (3, 1) -- (3, 2) -- (4, 2) -- (4, 3) -- (5, 3) -- (5, 4) -- (9, 4);
        \draw[line width=0.2mm, black, dashed] (-8, 2) -- (-4, 2) -- (-4, 3) -- (-3, 3) -- (-3, 4) -- (-2, 4) -- (-2, 15) -- (0, 15) -- (0, 14) -- (1, 14) -- (1, 16) -- (3, 16) -- (3, 15) -- (5, 15) -- (5, 14) -- (9, 14);

        \foreach \x in {-8,...,9} {
            \fill[black] (\x, -9) circle (0pt) node[above] {$\x$};
        }

        \fill[black] (-9, 9.6) circle (0pt) node {$\vdots$};
        \fill[black] (-9, -3.4) circle (0pt) node {$\vdots$};

        \fill[red] (-4, 16) circle (2pt) node[below right] {$A_1$};
        \fill[blue] (-3, 16) circle (2pt) node[below right] {$A_2$};
        \fill[green!50!black] (-1, -8) circle (2pt) node[above right] {$A_3$};
        \fill[black] (0, 16) circle (2pt) node[below right] {$A_4$};
        \fill[purple] (3, 16) circle (2pt) node[below right] {$A_5$};
        \fill[red 6] (2, -8) circle (2pt) node[below right] {$A_6$};
        \fill[red] (2, 16) circle (2pt) node[below right] {$B_1$};
        \fill[blue] (-2, 16) circle (2pt) node[below right] {$B_2$};
        \fill[green!50!black] (4, -8) circle (2pt) node[below right] {$B_3$};
        \fill[black] (1, 16) circle (2pt) node[below right] {$B_4$};
        \fill[purple] (5, -8) circle (2pt) node[below right] {$B_5$};
        \fill[red 6] (3, -8) circle (2pt) node[below right] {$B_6$};

        \foreach \x in {11, ..., 16} {
            \pgfmathsetmacro{\y}{int(\x - 12)}
            \fill[black] (-8.5,\x) circle (0pt) node[left] {$\y^{\bullet}$};
        }
        \foreach \x in {-8, ..., -5} {
            \pgfmathsetmacro{\y}{int(-\x - 6)}
            \fill[black] (-8.5,\x) circle (0pt) node[left] {$\y^{\circ}$};
        }


        \draw[green!50!black, thick, -{Stealth}] (-1, -8) --
        node[above, pos=0.0715]{$3^{\circ}$}
        node[above, pos=0.2143]{$3^{\circ}$}
        node[above, pos=0.35715]{$4^{\circ}$}
        node[above, pos=0.5]{$5^{\circ}$}
        node[above, pos=0.643]{$5^{\circ}$}
        node[above, pos=0.786]{$5^{\circ}$}
        node[above, pos=0.92]{$5^{\circ}$}
        (-8, -8);
        \draw[green!50!black, thick, -{Stealth}] (-8, -8) -- (-8, -2);
        \draw[green!50!black, thick, -{Stealth}] (-8, -2) -- 
        node[above, pos=0.5]{$1$}
        (-7, -2);
        \draw[green!50!black, thick, -{Stealth}] (-7, -2) -- (-7, 1);
        \draw[green!50!black, thick] (-7, 1) -- 
        node[above, pos=0.166]{$4$}
        node[above, pos=0.5]{$4$}
        node[above, pos=0.833]{$4$}
        (-4, 1);
        \draw[green!50!black, thick, -{Stealth}] (-4, 1) -- 
        node[above, pos=0.166]{$3$}
        node[above, pos=0.5]{$2$}
        node[above, pos=0.833]{$2$}
        (-1, 1);
        \draw[green!50!black, thick, -{Stealth}] (-1, 1) -- 
        (-1, 3);
        \draw[green!50!black, thick, -{Stealth}] (-1, 3) -- 
        node[above, pos=0.5]{$4$}
        (0, 3);
        \draw[green!50!black, thick, -{Stealth}] (0, 3) -- 
        (0, 2);
        \draw[green!50!black, thick, -{Stealth}] (0, 2) -- 
        node[above, pos=0.5]{$2$}
        (1, 2);
        \draw[green!50!black, thick, -{Stealth}] (1, 2) -- 
        (1, 3);
        \draw[green!50!black, thick, -{Stealth}] (1, 3) -- 
        node[above, pos=0.25]{$3$}
        node[above, pos=0.75]{$3$}
        (3, 3);
        \draw[green!50!black, thick, -{Stealth}] (3, 3) -- 
        (3, 2);
        \draw[green!50!black, thick, -{Stealth}] (3, 2) -- 
        node[above, pos=0.5]{$1$}
        (4, 2);
        \draw[green!50!black, thick, -{Stealth}] (4, 2) -- (4, -8);
        
        \draw[blue!60!green, thick, -{Stealth}] (0, -7) -- 
        node[above, pos=0.25]{$2^{\circ}$}
        node[above, pos=0.73]{$2^{\circ}$}
        (-2, -7);
        \draw[blue!60!green, thick, -{Stealth}] (-2, -7) -- 
        (-2, -6);
        \draw[blue!60!green, thick, -{Stealth}] (-2, -6) -- 
        node[above, pos=0.125]{$1^{\circ}$}
        node[above, pos=0.375]{$2^{\circ}$}
        node[above, pos=0.625]{$3^{\circ}$}
        node[above, pos=0.875]{$3^{\circ}$}
        (-6, -6);
        \draw[blue!60!green, thick, -{Stealth}] (-6, -6) -- 
        (-6, 0);
        \draw[blue!60!green, thick, -{Stealth}] (-6, 0) -- 
        node[above, pos=0.0833]{$3$}
        node[above, pos=0.25]{$3$}
        node[above, pos=0.4167]{$2$}
        node[above, pos=0.5833]{$1$}
        node[above, pos=0.75]{$1$}
        node[above, pos=0.9167]{$1$}
        (0, 0);
        \draw[blue!60!green, thick, -{Stealth}] (0, 0) -- 
        (0, -7);

        \draw[blue!70!green!65!red, thick, -{Stealth}] (-4, -5) -- 
        node[above, pos=0.5]{$2^{\circ}$}
        (-5, -5);
        \draw[blue!70!green!65!red, thick, -{Stealth}] (-5, -5) -- 
        (-5, -2);
        \draw[blue!70!green!65!red, thick, -{Stealth}] (-5, -2) -- 
        node[above, pos=0.5]{$1$}
        (-4, -2);
        \draw[blue!70!green!65!red, thick, -{Stealth}] (-4, -2) -- 
        (-4, -5);

        \draw[blue!50!green!25!red, thick, -{Stealth}] (8, 14) -- 
        node[above, pos=0.0625]{$5^{\bullet}$}
        node[above, pos=0.1875]{$5^{\bullet}$}
        node[above, pos=0.3125]{$5^{\bullet}$}
        node[above, pos=0.4375]{$4^{\bullet}$}
        node[above, pos=0.5625]{$3^{\bullet}$}
        node[above, pos=0.6875]{$2^{\bullet}$}
        node[above, pos=0.8125]{$2^{\bullet}$}
        node[above, pos=0.9375]{$2^{\bullet}$}
        (0, 14);

        \draw[blue!50!green!25!red, thick, -{Stealth}] (0, 14) -- 
        (0, 4);

        \draw[blue!50!green!25!red, thick, -{Stealth}] (0, 4) -- 
        node[above, pos=0.0625]{$4$}
        node[above, pos=0.1875]{$4$}
        node[above, pos=0.3125]{$4$}
        node[above, pos=0.4375]{$3$}
        node[above, pos=0.5625]{$2$}
        node[above, pos=0.6875]{$1$}
        node[above, pos=0.8125]{$1$}
        node[above, pos=0.9375]{$1$}
        (8, 4);

        \draw[blue!50!green!25!red, thick, -{Stealth}] (8, 4) -- 
        (8, 14);

        \draw[blue!90!green!45!red, thick, -{Stealth}] (7, 12) -- 
        node[above, pos=0.5]{$3^{\bullet}$}
        (6, 12);
        \draw[blue!90!green!45!red, thick, -{Stealth}] (6, 12) -- 
        (6, 13);
        \draw[blue!90!green!45!red, thick, -{Stealth}] (6, 13) -- 
        node[above, pos=0.5]{$4^{\bullet}$}
        (5, 13);
        \draw[blue!90!green!45!red, thick, -{Stealth}] (5, 13) -- 
        (5, 12);
        \draw[blue!90!green!45!red, thick, -{Stealth}] (5, 12) -- 
        node[above, pos=0.25]{$2^{\bullet}$}
        node[above, pos=0.75]{$1^{\bullet}$}
        (3, 12);
        \draw[blue!90!green!45!red, thick, -{Stealth}] (3, 12) -- 
        (3, 6);
        \draw[blue!90!green!45!red, thick, -{Stealth}] (3, 6) -- 
        node[above, pos=0.25]{$5$}
        node[above, pos=0.75]{$4$}
        (5, 6);
        \draw[blue!90!green!45!red, thick, -{Stealth}] (5, 6) -- 
        (5, 5);
        \draw[blue!90!green!45!red, thick, -{Stealth}] (5, 5) -- 
        node[above, pos=0.5]{$2$}
        (6, 5);
        \draw[blue!90!green!45!red, thick, -{Stealth}] (6, 5) -- 
        (6, 6);
        \draw[blue!90!green!45!red, thick, -{Stealth}] (6, 6) -- 
        node[above, pos=0.5]{$3$}
        (7, 6);
        \draw[blue!90!green!45!red, thick, -{Stealth}] (7, 6) -- 
        (7, 12);
    \end{tikzpicture}
    \caption{A flow on T-lattice.}
    \label{fig:flowex}
\end{figure}
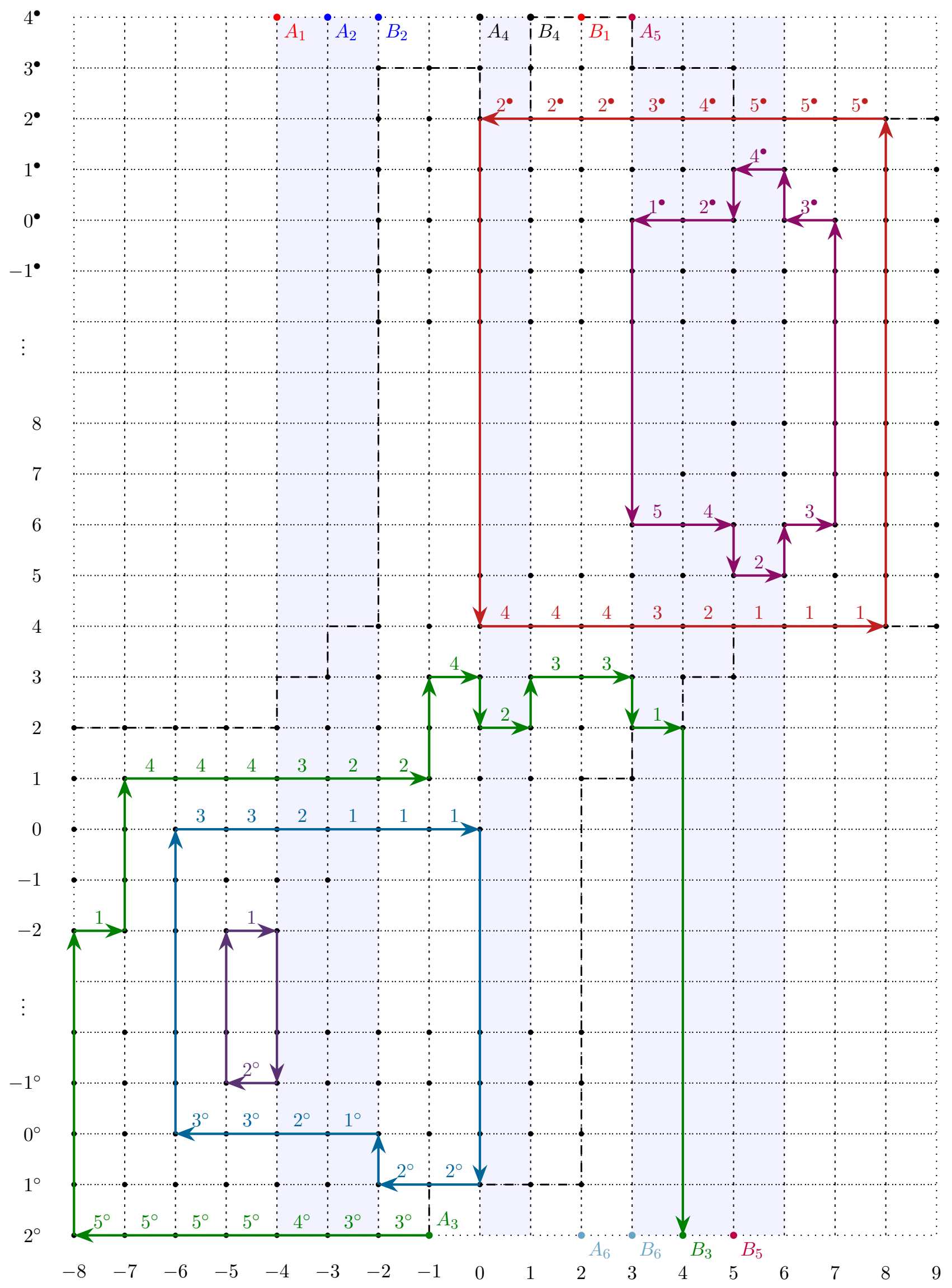

\begin{example}
    Consider Example~\ref{ex:p-cells}, we show an example of a flow from $C_3$ to $D_3$ and associated $\mathbf{T} \in {TT}^*_{\theta_3}(\mathbf{a}^{33}, \mathbf{b}^{33}, m)$, see Fig. \ref{fig:ex-ribbon-filling-vs-path1}, \ref{fig:ex-ribbon-filling-vs-path}, \ref{fig:flowex}. 
\end{example}


\begin{corollary}
    \label{cor:bijection-empty-ribbon}
    Let $\mathcal{C}$ be set of all non-intersecting systems of (simple) cycles on the lattice $TL_{\lambda/\mu}(\pi, m, \mathbf{a}, \mathbf{b})$. Then we have 
    \begin{equation*}
        \mathsf{G}^{\mathbf{a}^{\varepsilon}, \mathbf{b}^{\varepsilon}}_{\varepsilon}(\mathbf{x}_m; \mathbf{y}, \mathbf{z}) = \sum\limits_{\mathbf{C} \in \mathcal{C}} \mathrm{sgn(\mathbf{C}) \cdot wt(\mathbf{C})}.
    \end{equation*}
    In other words, there is a weight-preserving bijection between the sets $TT_{\varepsilon}(\mathbf{a}^{\varepsilon}, \mathbf{b}^{\varepsilon}, m)$ and $\mathcal{C}$.
\end{corollary}
\begin{proof}
    One can construct the lattice for $TL_{\varepsilon}(\pi, m, \mathbf{a}^{\varepsilon}, \mathbf{b}^{\varepsilon})$, which in fact coincides with the lattice $TL_{\lambda/\mu}(\pi, m, \mathbf{a}, \mathbf{b})$ with all endpoints removed. Hence, flows on the lattice $TL_{\varepsilon}(\pi, m, \mathbf{a}^{\varepsilon}, \mathbf{b}^{\varepsilon})$ are indeed collections of cycles, i.e. elements of the set $\mathcal{C}$ (as there are no endpoints, any flow $\mathbf{F}$ on the $TL_{\varepsilon}(\pi, m, \mathbf{a}^{\varepsilon}, \mathbf{b}^{\varepsilon})$ is of the form $(\varnothing, \mathbf{C})$, where $\mathbf{C} \in \mathcal{C}$). Then the proof follows from Lemma~\ref{lem:main-weight-preserving-bijection}.
\end{proof}

\begin{proof}[Proof of Theorem~\ref{thm:HG-sfG}]
    Now follows from Corollaries~\ref{cor:bijection-major}, \ref{cor:bijection-minor}, \ref{cor:bijection-empty-ribbon} together with modified Talaska's formula given in Corollary~\ref{cor:tal}.
\end{proof}

\section{Duality with column flagged functions}\label{sec:dual}
In this section we discuss column flagged Grothendieck enumerators that are dual to flagged Grothendieck enumerators under the standard involution $\omega$ of the ring of symmetric functions extended to completion. 

\subsection{Conjugate shapes and flags}


\begin{definition}[Conjugate shapes]
    Let $\lambda/\mu = ({D}(\lambda, \mu), {N}(\lambda, \mu))$ be some skew signature shape. For a cell $\alpha = (i, j) \in \mathbb{Z}^2$ define its {\it conjugate} as $\alpha' = (j, i)$. Define the conjugates of ${D}(\lambda, \mu), {N}(\lambda, \mu)$ as follows:
    \begin{equation*}
        {D}'(\lambda, \mu) := \{\alpha'\ |\ \alpha \in {D}(\lambda, \mu)\}, \qquad {N}' (\lambda, \mu) := \{(\alpha', \beta')\ |\ (\alpha, \beta) \in {N}(\lambda, \mu)\}.
    \end{equation*}
    Then the conjugate of shape $\lambda/\mu$ is given by:
    \begin{equation*}
        \lambda'/\mu' := ({D}'(\lambda, \mu), {N}'(\lambda, \mu)). 
    \end{equation*}

    

    We shall use primed $(')$ notation only for conjugate shapes. For simplicity, we use the term 'shape' instead of 'conjugate shape'. 
\end{definition}

\begin{definition}[Column $\mathsf{G}$-flags]
    Given skew signature $\lambda'/\mu'$ with $\lambda, \mu \in \mathsf{Sign}_n$. 
    The vectors $\mathbf{a}' = (a'_i) \in \mathbb{Z}_{-}^{n}, \mathbf{b}' = (b'_i) \in \mathbb{Z}_{-}^{n}$ are called {\it column $\mathsf{G}$-flags} if they satisfy the following conditions:
    \begin{align*}
        a'_i - a'_{i + 1} \le \mu_i - \mu_{i + 1} + 1,  \qquad b'_i - b'_{i + 1} \ge \lambda_{i + 1} - \lambda_i - 1, \qquad i \ge 1. 
    \end{align*}
\end{definition}

\begin{definition}[Conjugate flags]
    For column flags $\mathbf{a}', \mathbf{b}'$ define conjugate $\mathsf{G}$-flags $\mathbf{a} = (a_i), \mathbf{b} = (b_i)$ given by:
    \begin{align*}
        a_i = a'_i + c(\gamma_i) \qquad b_i = b'_i - c(\delta_i)
    \end{align*}
    where $\gamma_i = (\mu_i, i)$, and $\delta_i = (\lambda_i + 1, i)$.
    It is not hard to see that conjugate of row flags for $\lambda/\mu$ are column flags for $\lambda'/\mu'$ and vice versa.
\end{definition}


\begin{definition}
    Let $\lambda'/\mu'$ be skew shape of signatures from $\mathsf{Sign}_n$ and $\mathbf{a}' = (a'_i), \mathbf{b}' = (b'_i)$ be column $\mathsf{G}$-flags. 
    Define the following sets:
    \[
        \mathcal{A}(\lambda'/\mu') := \{(i, j)\ |\ j \in [n], i \le \mu_j\}, \qquad \mathcal{B}(\lambda'/\mu') := \{(i, j)\ |\ j \in [n], i \ge \lambda_j + 1\},
    \]  
    and    
    \begin{align*}
    {A}(\lambda'/\mu') &:= \{(i, j)  \ |\ j\in [n],\, a'_j \ne -\infty, \, i {\le \mu_j} \} \subset \mathbb{Z}^2, \\ 
    {B}(\lambda'/\mu') &:= \{(i, j) \ |\ j\in [n],\, b'_j \ne -\infty,\, i {\ge \lambda_j + 1} \} \subset \mathbb{Z}^2.
    \end{align*}
\end{definition}

\subsection{Conjugate flagged trio tableaux}
\begin{definition}[Conjugate flagged trio tableaux]
    Let $\lambda, \mu \in \mathsf{Sign}_n$ and $m$ be nonnegative integer. Let $\mathbf{a}' = (a'_i), \mathbf{b}' = (b'_i)$ be column $\mathsf{G}$-flags for skew shape $\lambda/\mu$. 
    
    A {\it conjugate flagged trio tableau} for shape $\lambda/\mu$ is a triple $\mathbf{T} = (T^{\circ}, T, T^{\bullet})$ of tableaux such that 
    
    (i) (Above tableau) $T^{\circ}$ is a semistandard tableau of some negative shape $\nu'/\mu'$ with $\nu \in \mathsf{Sign}_n$ and $\nu \subseteq \mu$, filled with entries from the set 
    $\mathbb{Z}^{\circ} := \{ i^{\circ} \,|\,  i \in \mathbb{Z} \}$
    such that the elements in $i$-th column are at most $a_i^{'\circ}$.
    
    (ii) (Below tableau) $T^{\bullet}$ is a reverse semistandard tableau of some negative shape $\lambda'/\tau'$ with $\tau \in \mathsf{Sign}_n$ and $\lambda \subseteq \tau$, filled with entries from the set 
    $\mathbb{Z}^{\bullet} := \{ i^{\bullet} \,|\, i \in \mathbb{Z}\}$
    such that the elements in $i$-th column are at least $b_i^{'\bullet}$. 
    
    (iii) (Central tableau) $T$ is a semistandard tableau of shape $\tau'/\nu'$ filled with entries from $[m]$. 
    

    The {\it weight} of $\mathbf{T}$ is defined as 
    $
    \mathrm{wt}(\mathbf{T}) := \mathrm{wt}(T^{\circ}) \cdot \mathrm{wt}(T) \cdot \mathrm{wt}(T^{\bullet}),     
    $
    and 
    \begin{align*}
    \mathrm{wt}(T^{\circ}) := \prod_{(i,j) \in D'(\nu, \mu)}  \mathrm{wt}(T^{\circ}_{i,j}), \quad 
    \mathrm{wt}(T) := \prod_{(i,j) \in D'(\tau, \nu)} \mathrm{wt}(T_{i,j}), \quad 
    \mathrm{wt}(T^{\bullet}) := \prod_{(i,j) \in D'(\lambda, \tau)} \mathrm{wt}(T^{\bullet}_{i,j}),  
    \end{align*}
    where denoting $T^{\circ}_{i,j} = a^{\circ}$, $T^{}_{i,j} = a^{}$, or $T^{\bullet}_{i,j} = a^{\bullet}$ and $c = c(i,j)$ we define the cell weights as
    $$
    \mathrm{wt}(T^{\circ}_{i,j}) := v_{a, a + c}, \quad  
    \mathrm{wt}(T_{i,j}) := x_{a}, \quad 
    \mathrm{wt}(T^{\bullet}_{i,j}) := u_{a, a - c},      
    $$
    for sets of variables $\mathbf{x}_m = (x_1,\ldots, x_m)$, $\mathbf{U} = (u_{i,j})_{i,j \in \mathbb{Z}}$, $\mathbf{V} = (v_{i,j})_{i,j \in \mathbb{Z}}$.
    
    
    We denote by $TT_{\lambda'/\mu'}(\mathbf{a}', \mathbf{b}', m)$ the set of conjugate flagged trio tableaux for shape $\lambda/\mu$ with column $\mathsf{G}$-flags $\mathbf{a}', \mathbf{b}'$.  
    We also denote by $TT^{\circ}_{\nu'/\mu'}(\mathbf{a}')$, $TT^{\bullet}_{\lambda'/\tau'}(\mathbf{b}')$, and $TT_{\tau'/\nu'}(m)$ for the corresponding sets of above, below, and central tableaux, respectively, where $\tau, \mu \in \mathsf{Sign}_n$ with $\nu \subseteq \mu$, $\lambda \subseteq \tau$, $\nu \subseteq \tau$.  
\end{definition}

\subsection{Column flagged Grothendieck enumerators}

\begin{definition}[Column flagged Grothendieck enumerators]
Let $m \in \mathbb{Z}_{\ge 0}$ and $\mathbf{x}_m = (x_1,\ldots, x_m)$, $\mathbf{U} = (u_{i,j})_{i,j \in \mathbb{Z}}$, $\mathbf{V} = (v_{i,j})_{i,j \in \mathbb{Z}}$ be sets of variables. 
    Let $\lambda/\mu$ be skew shape of signatures $\lambda, \mu \in \mathsf{Sign}_n$ with column $\mathsf{G}$-flags $\mathbf{a}', \mathbf{b}'$.
    We define the {\it column flagged Grothendieck enumerators} as follows:  
    \begin{align*}
        \mathsf{G}^{\mathbf{a}', \mathbf{b}'}_{\lambda'/\mu'}(\mathbf{x}_m ; \mathbf{U} , \mathbf{V}) 
        := 
        \sum_{\mathbf{T}\, \in\, TT_{\lambda'/\mu'}(\mathbf{a}', \mathbf{b}', m) } \mathrm{wt}(\mathbf{T}),
    \end{align*}
    i.e. as the weighted sum over conjugate flagged trio tableaux. 
\end{definition}

\subsection{Skew Schur expansions}
    Let us define the following enumerators: 
    \begin{align*}
        L^{\mathbf{a}}_{\nu/\mu}(\mathbf{U}, \mathbf{V}) &:= \sum_{T^{\circ} \in TT^{\circ}_{\nu/\mu}(\mathbf{a})}\ \mathrm{wt}(T^{\circ}) , \qquad
        R^{\mathbf{b}}_{\lambda/\tau}(\mathbf{U}, \mathbf{V}) := \sum_{T^{\bullet} \in TT^{\bullet}_{\lambda/\tau}(\mathbf{b})} \mathrm{wt}(T^{\bullet}), \\
        A^{\mathbf{a}'}_{\nu'/\mu'}(\mathbf{U}, \mathbf{V}) &:= \sum_{T^{\circ} \in TT^{\circ}_{\nu'/\mu'}(\mathbf{a}')} \mathrm{wt}(T^{\circ}), \qquad
        B^{\mathbf{b}'}_{\lambda'/\tau'}(\mathbf{U}, \mathbf{V}) := \sum_{T^{\bullet} \in TT^{\bullet}_{\lambda'/\tau'}(\mathbf{b}')} \mathrm{wt}(T^{\bullet}).
    \end{align*}

Now we will give Schur expansions for the enumerators $\mathsf{G}^{\mathbf{a}, \mathbf{b}}_{\lambda/\mu}(\mathbf{x}; \mathbf{U}, \mathbf{V})$ and $\mathsf{G}^{\mathbf{a}', \mathbf{b}'}_{\lambda'/\mu'}(\mathbf{x}; \mathbf{U}, \mathbf{V})$ that directly follow from the definitions of corresponding enumerators.

\begin{proposition}
    Let $\lambda, \mu \in \mathsf{Sign}_n$. We have the following skew Schur expansions of flagged Grothendieck enumerators:
    \begin{align*}
        \mathsf{G}^{\mathbf{a}, \mathbf{b}}_{\lambda/\mu}(\mathbf{x}; \mathbf{U}, \mathbf{V}) &= \sum_{\nu, \tau} L^{\mathbf{a}}_{\nu/\mu}(\mathbf{U}, \mathbf{V}) \, s_{\tau/\nu}(\mathbf{x}) \, R^{\mathbf{b}}_{\lambda/\tau}(\mathbf{U}, \mathbf{V}), \\
        \mathsf{G}^{\mathbf{a}', \mathbf{b}'}_{\lambda'/\mu'}(\mathbf{x}; \mathbf{U}, \mathbf{V}) &= \sum_{\nu, \tau} A^{\mathbf{a}'}_{\nu'/\mu'}(\mathbf{U}, \mathbf{V}) \, s_{\tau'/\nu'}(\mathbf{x}) \,  B^{\mathbf{b}'}_{\lambda'/\tau'}(\mathbf{U}, \mathbf{V}),
    \end{align*}
    where the summations goes over $\nu, \tau \in \mathsf{Sign}_n$ such that $\nu \subseteq \mu,\ \lambda \subseteq \tau,\ \nu \subseteq \tau$.
\end{proposition}

\subsection{Duality involution}
We consider the standard involution $\omega$ of the ring of symmetric functions, which acts on Schur functions as $\omega: s_{\lambda}(\mathbf{x}) \mapsto s_{\lambda'}(\mathbf{x})$, extended to completion containing infinite formal power series expansions of Schur functions, so that its action on flagged Grothendieck enumerators is well defined.   

\begin{theorem}[Duality]\label{thm:omega}
We have:  
    \begin{equation*}
        \omega_{\mathbf{x}}(\mathsf{G}^{\mathbf{a}, \mathbf{b}}_{\lambda/\mu}(\mathbf{x} ; \mathbf{U} , \mathbf{V})) = \mathsf{G}^{\mathbf{a}', \mathbf{b}'}_{\lambda'/\mu'}(\mathbf{x} ; {\mathbf{U}^{T}} , {\mathbf{V}^T}). 
    \end{equation*}
    Here $\mathbf{U}^T$ denotes the transpose matrix of $\mathbf{U}$. 
\end{theorem}

\begin{lemma}
    \label{lem:bijection-left-above-right-below}
    Let $\lambda, \mu \in \mathsf{Sign}_n$ and $\mathbf{a}, \mathbf{b}$ be row $\mathsf{G}$-flags for $\lambda/\mu$. Then the conjugate flags $\mathbf{a}', \mathbf{b}'$ are column $\mathsf{G}$-flags for $\lambda'/\mu'$. 
    Let $\nu, \tau \in \mathsf{Sign}_n$ such that $\nu \subseteq \mu,\ \lambda \subseteq \tau,\ \nu \subseteq \tau$. There are bijections $\psi^{\circ}, \psi^{\bullet}$:
    \begin{align*}
        \psi^{\circ} : TT^{\circ}_{\nu/\mu}(\mathbf{a}) \rightarrow TT^{\circ}_{\nu'/\mu'}(\mathbf{a}'), \qquad
        \psi^{\bullet} : TT^{\bullet}_{\lambda/\tau}(\mathbf{b}) \rightarrow TT^{\bullet}_{\lambda'/\tau'}(\mathbf{b}'),
    \end{align*}
    so that $Q^{\circ} = \psi^{\circ}(T^{\circ}),\ Q^{\bullet} = \psi^{\bullet}(T^{\bullet})$ are defined as follows:
    \begin{equation*}
        Q^{\circ}_{j, i} = T^{\circ}_{i, j} + c(i, j), \qquad Q^{\bullet}_{j, i} = T^{\bullet}_{i, j} - c(i, j),
    \end{equation*}
    where $c(i, j)$ is the content of cell $(i, j)$ in the shape $\nu/\mu$. 
\end{lemma}
\begin{proof}
    Let us prove the statement for $\psi^{\circ}$, and for $\psi^{\bullet}$ it can be done similarly. It is not hard to see that $Q^{\circ}$ is unique and semistandard tableau of shape $\nu'/\mu'$. Also $Q^{\circ}$ respects the flags $\mathbf{a}'$ as $Q^{\circ}_{j, i} = T^{\circ}_{i, j} + c(i, j) \le a_i + \lambda_i + 1 - i = a'_i$. Thus $Q^{\circ} \in TT^{\circ}_{\nu'/\mu'}(\mathbf{a}')$. 
    
    Next we define the inverse map $\phi^{\circ}$ as follows: $T^{\circ} = \phi^{\circ}(Q^{\circ})$, where $T^{\circ}_{j, i} = Q^{\circ}_{i, j} + c(i, j)$, clearly $T^{\circ}$ is unique and belongs to $TT^{\circ}_{\nu/\mu}(\mathbf{a})$, and also  $Q^{\circ} = \psi^{\circ}(\phi^{\circ}(Q^{\circ}))$ holds, which concludes the proof.
\end{proof}

\begin{lemma}
    \label{lem:left-above-right-below}
    We have: 
    \begin{equation*}
        L^{\mathbf{a}}_{\nu/\mu}(\mathbf{U}, \mathbf{V}) = A^{\mathbf{a}'}_{\nu'/\mu'}({\mathbf{U}^T}, {\mathbf{V}^T}),
        \qquad 
        R^{\mathbf{b}}_{\lambda/\tau}(\mathbf{U}, \mathbf{V}) = B^{\mathbf{b}'}_{\lambda'/\tau'}({\mathbf{U}^T}, {\mathbf{V}^T}).
    \end{equation*}
\end{lemma}
    \begin{proof}
    Let us focus on the first equation. We use the bijection $\psi^{\circ}$ from Lemma~\ref{lem:bijection-left-above-right-below}. Let $T \in TT^{\circ}_{\nu/\mu}(\mathbf{a})$ and $Q = \psi^{\circ}(T) \in TT^{\circ}_{\nu'/\mu'}(\mathbf{a}')$.
    We have 
    \begin{align*}
        \mathrm{wt}(T_{i, j}) = v_{T_{i, j}, T_{i,j} + c(i, j)} = v_{Q_{j,i} + c(j,i), Q_{j,i}} 
        = \mathrm{wt}(Q_{j,i})\ |_{ v_{a,b} \to v_{b,a} }.
    \end{align*}
    Thus, the equation follows. The second equation can be proved similarly.
\end{proof}

\begin{proof}[Proof of Theorem~\ref{thm:omega}]
    Using the above Schur expansions and Lemma~\ref{lem:left-above-right-below} we have 
    \begin{align*}
        \omega(\mathsf{G}^{\mathbf{a}, \mathbf{b}}_{\lambda/\mu}(\mathbf{x}; \mathbf{U}, \mathbf{V})) &= \sum_{\nu, \tau} L^{\mathbf{a}}_{\nu/\mu}(\mathbf{U}, \mathbf{V}) \, \omega(s_{\tau/\nu}(\mathbf{x})) \, R^{\mathbf{b}}_{\lambda/\tau}(\mathbf{U}, \mathbf{V}) \\
        &= \sum_{\nu, \tau} A^{\mathbf{a}'}_{\nu'/\mu'}({\mathbf{U}^T}, {\mathbf{V}^T}) \, s_{\tau'/\nu'}(\mathbf{x}) \, B^{\mathbf{b}'}_{\lambda'/\tau'}({\mathbf{U}^T}, {\mathbf{V}^T}) \\
        &= \mathsf{G}^{\mathbf{a}', \mathbf{b}'}_{\lambda'/\mu'}(\mathbf{x}; {\mathbf{U}^T}, {\mathbf{V}^T}),
    \end{align*}
    which shows the statement. 
\end{proof}

\subsection{Column Hamel--Goulden formulas}
\begin{definition}
    Let $\Theta = A \cup B$ be outer decomposition of shape $\lambda/\mu$, where $A$ is set of ribbons, and $B$ is set of null-strips. Define the following sets:
    \begin{equation*}
        A' = \{\theta'\ |\ \theta \in A\}, \qquad B' = \{(\alpha', \beta')\ |\ (\alpha, \beta) \in B\}.
    \end{equation*}
    We define the outer decomposition $\Theta'$ of the shape $\lambda'/\mu'$ as follows:
    \begin{equation*}
        \Theta' := A' \cup B'.
    \end{equation*}
    Let $\theta'_i, \theta'_j \in \Theta'$. Define the $\#$-operation on $\Theta'$ as follows: 
    \begin{equation*}
        \theta'_i \# \theta'_j := (\theta_j \# \theta_i)'.
    \end{equation*}
    Let $\delta, \gamma$ be head and tail of $\theta_i \in \Theta$, then the head and tail of $\theta'_i$ are defined to be $\gamma', \delta'$. 
\end{definition}

\begin{definition}[Induced column $\mathsf{G}$-flags]
Let $\lambda'/\mu'$ be skew shape with column $\mathsf{G}$-flags $\mathbf{a}' = (a'_i), \mathbf{b}' = (b'_i)$.
Let $c \in \mathbb{Z}$ be some content. Define $\mathrm{max}(\lambda'/\mu', c)$ as the bottommost cell in $A(\lambda'/\mu')$ of content $c$. If $A(\lambda'/\mu')$ does not contain a cell of content $c$, set $\mathrm{max}(\lambda'/\mu', c) = \mathrm{undefined}$. Similarly, define $\mathrm{min}(\lambda'/\mu', c)$ as the topmost cell of content $c$ in $B(\lambda'/\mu')$, if there is no such cell, set $\mathrm{min}(\lambda'/\mu', c) = \mathrm{undefined}$. 

Let $n_{ij}$ be the number of columns in $\theta'_i \# \theta'_j$. Now we define induced column $\mathsf{G}$-flags of length $n_{ij}$, $\mathbf{a}'^{ij} = (a'^{ij}_r), \mathbf{b}'^{ij} = (b'^{ij}_r)$ for the shape $\theta'_i \# \theta'_j$ as follows: 
Let $\zeta_r$ be bottommost cell on $r$-th column from $\mathcal{A}(\theta'_i \# \theta'_j)$, and $\eta_r$ be topmost cell on $r$-th column from $\mathcal{B}(\theta'_i \# \theta_j)$. Let $M'_r = \mathrm{max}(\lambda'/\mu', c(\zeta_r)),\ m'_r = \mathrm{min}(\lambda'/\mu', c(\eta))$, then we have
\begin{align*}
    a'^{ij}_r &= 
    \begin{cases}
        -\infty, & \text{if} \quad \text{$M'_r$ is undefined}, \\
        a'_{\mathrm{col}(M'_r)} + \mathrm{row}(M'_r) - \mu_{\mathrm{col}(M'_r)}, & \text{otherwise},
    \end{cases}
    \\[0.4cm]
    b'^{ij}_r &= 
    \begin{cases}
        -\infty, & \text{if} \quad \text{$m'_r$ is undefined}, \\
        b'_{\mathrm{col}(m'_r)} - \mathrm{row}(m'_r) + \lambda_{\mathrm{col}(m'_r)} + 1, & \text{otherwise}.
    \end{cases}
\end{align*}
In the same way we define induced column $\mathsf{G}$-flags $\mathbf{a}^{\varepsilon}, \mathbf{b}^{\varepsilon}$ for empty ribbon $\varepsilon$. 
\end{definition}


\begin{corollary}
    \label{cor:dual-HG-sfG}
    Let $\Theta' = (\theta'_1, \ldots, \theta'_k)$ be an outer ribbon decomposition of skew signature $\lambda'/\mu'$  with associated pipe function $\pi$. Let $\mathbf{a}', \mathbf{b}'$ be $\mathsf{G}$-flags and $\mathbf{a}'^{ij}, \mathbf{b}'^{ij}$ be induced column ribbon $\mathsf{G}$-flags for the ribbons $\theta'_i \# \theta'_j$. Then the following formula holds:  
    \begin{align*}
        \mathsf{G}^{\mathbf{a}', \mathbf{b}'}_{\lambda'/\mu'}(\mathbf{x} ; \mathbf{U} , \mathbf{V}) 
        = (-1)^{\mathsf{hp}^{\pi}_{\lambda'/\mu'}} \cdot 
        \mathsf{G}^{\mathbf{a}'^{\varepsilon}, \mathbf{b}'^{\varepsilon}}_{\varepsilon}(\mathbf{x} ;  \mathbf{U} , \mathbf{V})^{-k+1} \cdot 
        \det \left[ (-1)^{\mathsf{hp}^{\pi}_{\theta'_i \# \theta'_j}} \cdot \mathsf{G}^{\mathbf{a}'^{ij}, \mathbf{b}'^{ij}}_{\theta'_i \# \theta'_j}(\mathbf{x} ; \mathbf{U} , \mathbf{V})\right]_{1 \le i,j \le k}. 
    \end{align*}
\end{corollary}

\begin{proof}
    Let $\mathbf{a}, \mathbf{b}$ be row $\mathsf{G}$-flags for the shape $\lambda/\mu$ which are conjugate to column flags $\mathbf{a}', \mathbf{b}'$. Let $\mathbf{a}^{ij} = (a^{ij}_r), \mathbf{b}^{ij} = (b^{ij}_r)$ be induced row flags for the shape $\theta_i \# \theta_j$, where $\theta_i, \theta_j \in \Theta$. Then using $\omega$ involution on the Hamel--Goulden-type formula for $\mathsf{G}^{\mathbf{a}, \mathbf{b}}_{\lambda/\mu}$ we have: 
    \begin{equation*}
        \mathsf{G}^{\mathbf{a}', \mathbf{b}'}_{\lambda'/\mu'}(\mathbf{x}; \mathbf{U}, \mathbf{V}) = (-1)^{\mathsf{hp}^{\pi}_{\lambda/\mu}} \cdot 
        \mathsf{G}^{(\mathbf{a}^{\varepsilon})', (\mathbf{b}^{\varepsilon})'}_{\varepsilon'}(\mathbf{x}; \mathbf{U}, \mathbf{V})^{-k+1} \cdot  
        \det \left[(-1)^{\mathsf{hp}^{\pi}_{\theta_i \# \theta_j}}\ \mathsf{G}^{(\mathbf{a}^{ij})', (\mathbf{b}^{ij})'}_{(\theta_i \# \theta_j)'}(\mathbf{x}_m; \mathbf{U}, \mathbf{V})\right]_{1 \le i,j \le k}. 
    \end{equation*}
    By definitiion, 
    $
        \theta'_j \# \theta'_i = (\theta_i \# \theta_j)'. 
    $ 
    Now we need to show that: 
    \begin{equation*}
        (\mathbf{a}^{ij})' = \mathbf{a}'^{ji} \qquad (\mathbf{b}^{ij})' = \mathbf{b}'^{ji} \qquad (\mathbf{a}^{\varepsilon})' = \mathbf{a}'^{\varepsilon} \qquad (\mathbf{b}^{\varepsilon})' = \mathbf{b}'^{\varepsilon}. 
    \end{equation*}
    Let us prove the first equation, and others can be proved similarly. 
    Let $\alpha_r$ be the rightmost cell in $r$-th row of $\mathcal{L}(\theta_i \# \theta_j)$, and $M_r = \mathrm{max}(\lambda/\mu, c(\alpha_r))$. Let $\beta_r$ be the bottommost cell in $r$-th column of $\mathcal{A}(\theta'_j \# \theta'_i)$, and $M'_r = \mathrm{max}(\lambda'/\mu', c(\beta))$. Note that $M'_r$ is conjugate to $M_r$; thus if $M_r$ is $\mathrm{undefined}$ then $M'_r$ is $\mathrm{undefined}$ as well, and $(a^{ij}_r)' = a'^{ji}_r = -\infty$. If $M_r$ is not undefined, by definition we have: 
    \begin{align*}
        (a^{ij}_r)' &= a^{ij}_r + c(\alpha_r) \\ 
        &= a_{\mathrm{row}(M_r)} + c(M_r) \\ 
        &= a'_{\mathrm{row}(M_r)} - \mu_{\mathrm{row}(M_r)} + \mathrm{row}(M_r) + \mathrm{col}(M_r) - \mathrm{row}(M_r) \\
        &= a'_{\mathrm{col}(M'_r)} - \mu_{\mathrm{col}(M'_r)}  + \mathrm{row}(M'_r) \\
        &= a'^{ji}_r. 
    \end{align*}
    Thus, by transposing the matrix the result follows. 
\end{proof}

\subsection{Specialization to canonical refined stable Grothendieck polynomials}
\begin{proposition}
    Let $\lambda, \mu \in \mathsf{Par}_n$. Define column $\mathsf{G}$-flags $\mathbf{a}' = (\mu_i), \mathbf{b}' = (i - 1)$, whose conjugate flags are $\mathbf{a} = (i), \mathbf{b} = (\lambda_i)$. By taking the $G$-specialization: $u_{i,j} = \alpha_i - \beta_j$, $v_{i,j} = \beta_i - \alpha_j$ and setting $\alpha_i = \beta_i = 0$ for $i \le 0$ we have 
    \begin{align*} 
        G_{\lambda'/\mu'}(\mathbf{x}_m; \boldsymbol{\alpha}, \boldsymbol{\beta}) = D \cdot \mathsf{G}^{\mathbf{a}', \mathbf{b}'}_{\lambda'/\mu'}(\mathbf{x}_m; \mathbf{U}, \mathbf{V})_{G\text{-spec.}}, \quad D = \prod_{i = 1}^n \prod_{j = 1}^m (1 - \alpha_i x_j)^{-1}.
    \end{align*}
\end{proposition}

\begin{proof}
It follows from \cite[Thm.~6.4]{hwang1} that $G_{\lambda/\mu}(\mathbf{x}_m ; \boldsymbol{\alpha}, \boldsymbol{\beta})/C$ can be expressed as weighted sums over trio tableaux with the given concrete flags. 
\end{proof}

\section{Some special decompositions}\label{sec:spec}

Recall that the elementary and complete homogeneous symmetric functions are given by 
$$
e_n(\mathbf{x}) = \sum_{i_1 < \cdots < i_n} x_{i_1} \cdots x_{i_n}, \quad h_n(\mathbf{x}) = \sum_{i_1 \le \cdots \le i_n} x_{i_1} \cdots x_{i_n}, \quad n \ge 1, 
$$
with $e_0 = h_0 = 1$ and $e_n = h_n = 0$ for $n < 0$. 

For $a \in \mathbb{Z}$ we denote $\mathbf{U}_a = (u_{i,j})_{i \le a, j \in \mathbb{Z}}$, $\mathbf{V}_a = (v_{i,j})_{i \le a, j \in \mathbb{Z}}$. 
Then  for $n \ge 1, c \in \mathbb{Z}$ we define the following extensions of elementary and complete homogeneous symmetric functions:  
    \begin{align*}
        h_{n, c}(\mathbf{V}_a) & := \sum_{\substack{
            i_n \le \ldots \le i_1 \le a \\
        }}\ \prod_{r = 1}^n v_{i_r, i_r + c - r + 1}, \\
        e_{n, c}(\mathbf{V}_a) & := \sum_{\substack{
            i_n < \ldots < i_1 \le a \\
        }}\ \prod_{r = 1}^n v_{i_r, i_r + c + r - 1},
    \end{align*}
    with $h_{0, c} = e_{0, c} = 1$ and $h_{n, c} = e_{n, c} = 0$ for $n < 0$.

    For $n \in \mathbb{Z}$ we also define the following functions 
    \begin{align*}
        h_{n}(\mathbf{x} \,|\, \tau^{c_1}\mathbf{U}_b, \tau^{c_0}\mathbf{V}_a) &:= \sum_{k \ge 0}h_{n + k}(\mathbf{x}) \sum_{p + q = k} h_{p, c_0}(\mathbf{V}_a) h_{q, c_1}(\mathbf{U}_b), \\
        e_{n}(\mathbf{x} \,|\, \tau^{c_1}\mathbf{U}_b, \tau^{c_0}\mathbf{V}_a) &:= \sum_{k \ge 0}e_{n + k}(\mathbf{x}) \sum_{p + q = k} e_{p, c_0}(\mathbf{V}_a) e_{q, c_1}(\mathbf{U}_b).
    \end{align*}

\begin{lemma}
    \label{lem:minor-enumerator-new}
    Let $(\lambda_1)/(\mu_1)$ be a single row shape with $n$ cells. Let $c_0, c_1$ be contents of the rightmost cell in $\mathcal{L}((\lambda_1)/(\mu_1))$, and the leftmost cell in $\mathcal{R}((\lambda_1)/(\mu_1))$ (resp.), i.e. $c_0 = \mu_1 - 1, c_1 = (\lambda_1 + 1) - 1 = \lambda_1$. Let $\mathbf{a} = (a), \mathbf{b} = (b)$ be its row $\mathsf{G}$-flags. 
    Then we have:
    \begin{equation*}
        \mathsf{G}^{\mathbf{a}, \mathbf{b}}_{(\lambda_1) / (\mu_1)}(\mathbf{x}; \mathbf{U}, \mathbf{V}) = 
        h_{n}(\mathbf{x} \,|\, \tau^{-c_1}\mathbf{U}_b, \tau^{c_0}\mathbf{V}_a).
    \end{equation*}
\end{lemma}
\begin{proof}
    Let $a, b \in \mathbb{Z}$. For simplicity we will write $a / b$ instead of $(a) / (b)$ for a single row shape. 
    
    By definition, $\mathsf{G}^{\mathbf{a}, \mathbf{b}}_{\lambda_1 / \mu_1}(\mathbf{x}; \mathbf{U}, \mathbf{V})$ is the enumerator for tableaux $\mathbf{T} = (T^{\circ}, \widetilde{T}, T^{\bullet}) \in TT_{\lambda_1 / \mu_1}(\mathbf{a}, \mathbf{b})$ where tableaux $T^{\circ}, \widetilde{T}, T^{\bullet}$ from the sets $TT^{\circ}_{\nu_1 / \mu_1}, TT_{\tau_1 / \nu_1}, TT^{\bullet}_{\lambda_1 / \tau_1}$ (resp.) where $\nu_1 \le \mu_1,\ \lambda_1 \le \tau_1,\ \nu_1 \le \tau_1$.
    
    For $k, p, q \in \mathbb{Z}$ with $p + q = k$ we have:
    \begin{align*}
        h_{p, c_0}(\mathbf{V}_a) &= \sum_{T^{\circ} \in TT^{\circ}_{\nu_1 / \mu_1}(\mathbf{a})} \mathrm{wt}(T^{\circ}), \\
        h_{n + k}(\mathbf{x}) &= \sum_{\widetilde{T} \in TT_{\tau_1 / \nu_1}} \mathrm{wt}(\widetilde{T}), \\ 
        h_{q, -c_1}(\mathbf{U}_b) &= \sum_{T^{\bullet} \in TT^{\bullet}_{\lambda_1 / \tau_1}(\mathbf{b})} \mathrm{wt}(T^{\bullet}), 
    \end{align*}
    where $\mu_1 - \nu_1 = p,\ \tau_1 - \nu_1 = n + k,\ \tau_1 - \lambda_1 = q$; which completes the proof.
\end{proof}

\subsection{Jacobi--Trudi-type formulas}
Let $\lambda, \mu \in \mathsf{Sign}_n$. 
Let us consider the decomposition of the shape $\lambda/\mu$ into horizontal ribbons, i.e. $\pi(c) = (0, 0)$ for all $c \in \mathbb{Z}$. Then we have 
\[
\theta_j \# \theta_i = (\lambda_i - i + 1) / (\mu_j - j + 1), \quad i, j \in [n].
\] 
The induced $\mathsf{G}$-flags are $\mathbf{a}^{ji} = (a_j), \mathbf{b}^{ji} = (b_i)$ (flags with a single element).

\begin{corollary}[Jacobi--Trudi-type formula]
    \label{cor:main-enumerator}
    Let $\lambda/\mu$ be skew signature shape with $\lambda, \mu \in \mathsf{Sign}_n$ and row $\mathsf{G}$-flags $\mathbf{a} = (a_i), \mathbf{b} = (b_i)$. We have: 
    \begin{align*}
        \mathsf{G}^{\mathbf{a}, \mathbf{b}}_{\lambda/\mu}(\mathbf{x}; \mathbf{U}, \mathbf{V}) = 
        \det \left[ h_{\lambda_i - \mu_j - i + j}(\mathbf{x}_m \,|\, \tau^{i -\lambda_i - 1}\mathbf{U}_{b_i}, \tau^{\mu_j - j}\mathbf{V}_{a_j}) \right]_{1 \le i, j \le n}.
    \end{align*}
\end{corollary}

\begin{proof}
    From Lemma~\ref{lem:minor-enumerator-new} we have 
    \begin{equation*}
        \mathsf{G}^{\mathbf{a}^{ji}, \mathbf{b}^{ji}}_{\theta_j \# \theta_i}(\mathbf{x}_m; \mathbf{U}, \mathbf{V}) = 
        h_{\lambda_i - \mu_j - i + j}(\mathbf{x}_m \,|\, \tau^{i -\lambda_i - 1}\mathbf{U}_{b_i}, \tau^{\mu_j - j}\mathbf{V}_{a_j}). 
    \end{equation*}
    Note that there are no cycles in the lattice $TL_{\lambda/\mu}(\pi, m, \mathbf{a}, \mathbf{b})$ (as $\pi(c) = (0, 0)$, for all $c \in \mathbb{Z}$), and then the statement follows from Theorem~\ref{thm:HG-sfG}.
\end{proof}

\begin{example}
    Let $\theta_1 \# \theta_2$ be the shape and $\varepsilon$ be the corresponding empty shape from Example~\ref{ex:hg-ex-1}. Recall that $\theta_1 \# \theta = (2,2) / (1,-1)$ and $\varepsilon = (2,2,0) / (2,2,0)$. The corresponding $\mathsf{G}$-flags are: 
    \begin{align*}
        &\mathbf{a}^{12} = (a_1, a_3),  \quad \mathbf{b}^{12} = (b_2, b_3), \\ 
        &\mathsf{a}^{\varepsilon} = (-\infty, a_1, a_3), \quad \mathsf{b}^{\varepsilon} = (b_2, b_3, -\infty),
    \end{align*}
    whose conjugate flags are: 
    \begin{align*}
        &\mathbf{a}'^{12} = (a_1, a_3 - 3), \quad \mathbf{b}'^{12} = (b_2 - 2, b_3 - 1), \\
        &\mathsf{a}'^{\varepsilon} = (-\infty, a_1, a_3 - 3), \quad \mathsf{b}'^{\varepsilon} = (b_2 - 2, b_3 - 1, -\infty).
    \end{align*}
    Then using Corollary~\ref{cor:main-enumerator} we can express the functions $\mathsf{G}^{\mathbf{a}^{1, 2}, \mathbf{b}^{1, 2}}_{\theta_1 \# \theta_2}, \mathsf{G}^{\mathbf{a}^{\varepsilon}, \mathbf{b}^{\varepsilon}}_{\varepsilon}$ in terms of $h$. 
\end{example}

\begin{lemma}
    \label{lem:omega-on-h}
    Let $n, a, b, c_0, c_1 \in \mathbb{Z}$ 
    and $a' = a + c_0, b' = b - c_1$. Then we have  
    \begin{equation*}
        \omega(h_{n}(\mathbf{x} \,|\, \tau^{c_1}\mathbf{U}_b, \tau^{c_0}\mathbf{V}_a)) = 
        e_{n}(\mathbf{x} \,|\, \tau^{-c_1}\mathbf{U}^T_{b'}, \tau^{-c_0}\mathbf{V}^T_{a'}).
    \end{equation*}
\end{lemma}
\begin{proof}
    We have  
    \begin{align*}
        h_{n, c}(\mathbf{V}_a) &= \sum_{\substack{
            i_n \le \ldots \le i_1 \le a \\
        }}\ \prod_{r = 1}^n v_{i_r, i_r + c - r + 1} \\
        &= \sum_{\substack{
            j_n < \ldots < j_1 \le a' \\
        }}\ \prod_{r = 1}^n v_{j_r, j_r + (-c) + r - 1} \big|_{
            v_{i, j} \rightarrow v_{j, i}
        } \\
        &= e_{n, -c}(\mathbf{V}_{a'}) \big|_{
            v_{i, j} \rightarrow v_{j, i}
        } = e_{n, -c}(\mathbf{V}^T_{a'}). 
    \end{align*}
    Therefore,  
    \begin{align*}
        \omega(h_{n}(\mathbf{x} \,|\, \tau^{c_1}\mathbf{U}_b, \tau^{c_0}\mathbf{V}_a)) &= \sum_{k \ge 0}\omega(h_{n + k}(\mathbf{x})) \sum_{p + q = k} h_{p, c_0}(\mathbf{V}_a) h_{q, c_1}(\mathbf{U}_b) \\
        &= \sum_{k \ge 0}e_{n + k}(\mathbf{x}) \sum_{p + q = k} e_{p, -c_0}(\mathbf{V}_{a'}) e_{q, -c_1}(\mathbf{U}_{b'}) \big|_{\substack{
            u_{i, j} \rightarrow u_{j, i} \\
            v_{i, j} \rightarrow v_{j, i}
        }} \\
        &= e_{n}(\mathbf{x} \,|\, \tau^{-c_1}\mathbf{U}_{b'}, \tau^{-c_0}\mathbf{V}_{a'})\ \big|_{\substack{
            u_{i, j} \rightarrow u_{j, i} \\
            v_{i, j} \rightarrow v_{j, i}
        }}
        = e_{n}(\mathbf{x} \,|\, \tau^{-c_1}\mathbf{U}^T_{b'}, \tau^{-c_0}\mathbf{V}^T_{a'})
    \end{align*}
    as needed. 
\end{proof}
        
\begin{corollary}[Dual-Jacobi--Trudi-type formula]
    \label{cor:main-dual-enumerator}
    Let $\lambda/\mu$ be skew shape with $\lambda, \mu \in \mathsf{Sign}_n$ and column $\mathsf{G}$-flags $\mathbf{a}' = (a'_i), \mathbf{b}' = (b'_i)$, and row $\mathsf{G}$-flags $\mathbf{a} = (a_i), \mathbf{b} = (b_i)$ be their conjugates. We have: 
    \begin{align*}
        \mathsf{G}^{\mathbf{a}', \mathbf{b}'}_{\lambda'/\mu'}(\mathbf{x}; \mathbf{U}, \mathbf{V}) = 
        \det \left[ e_{\lambda_i - \mu_j - i + j}(\mathbf{x} \,|\, \tau^{\lambda_i - i + 1}\mathbf{U}_{b'_i}, \tau^{j - \mu_j}\mathbf{V}_{a'_j}) \right]_{1 \le i, j \le n}. 
    \end{align*}
\end{corollary}
\begin{proof}
    Apply the involution $\omega$ to Corollary~\ref{cor:main-enumerator}, then using Lemma~\ref{lem:omega-on-h} and changing variables the result follows. 
\end{proof}

\begin{corollary}[Dual Jacobi--Trudi-type formula for $\mathsf{G}^{\mathbf{a}, \mathbf{b}}_{\lambda/\mu}$]
    Let $\lambda/\mu$ be skew shape with $\lambda, \mu \in \mathsf{Sign}_n$ and $\mathsf{G}$-flags $\mathbf{a}, \mathbf{b}$. 
    Let $\Theta = (\theta_1, \ldots, \theta_r)$ be decomposition of $\lambda/\mu$ into vertical ribbons. 
    Let $p_i, q_i$ be the contents of the head and tail of $\theta_i$. 
Define contents $a = \mathrm{min}(\lambda_n, \mu_n) + 1 - n, b = \mathrm{max}(\lambda_1, \mu_1) - 1$.     
The pipe function is given by: 
    \begin{equation*}
        \pi(c)  = 
        \begin{cases}
            (1, 1), \quad c \in [a, b] \\
            (0, 1), \quad c = a - 1 \\
            (1, 0), \quad c = b + 1 \\
            (0, 0), \quad \text{otherwise}
        \end{cases}
    \end{equation*}
    Let $\mathbf{a}^{ij}, \mathbf{b}^{ij}$ be induced ribbon $\mathsf{G}$-flags. Let $\mathbf{a}^{\varepsilon} = (a^{\varepsilon}_i), \mathbf{b}^{\varepsilon} = (b^{\varepsilon}_i)$ be flags for the empty shape with conjugates $\mathbf{a}'^{\varepsilon} = (a'^{\varepsilon}_i), \mathbf{b}'^{\varepsilon} = (b'^{\varepsilon}_i)$. Denote
    \[
    1^{(a, b, p, q)} := 
    \begin{cases}
        ((b + 1)^{b - p + 1},\ b^{p - a}) / ((b + 1)^{b - q},\ b^{q + 1 - a}), &\quad \lambda/\mu \text{ is positive}, \\
        ((b + 1)^{b - p},\ b^{p + 1 - a}) / ((b + 1)^{b - q + 1},\ b^{q - a}), &\quad \text{otherwise},
    \end{cases}
    \]
    where $p, q \in [a, b]$.
    Denote $n' = b - a + 2$. Then we have 
    \begin{equation*}
        \mathsf{G}^{\mathbf{a}, \mathbf{b}}_{\lambda/\mu}(\mathbf{x}; \mathbf{U}, \mathbf{V}) = S^{-r+1} \cdot \det \left[\mathsf{G}^{\mathbf{a}^{ij}, \mathbf{b}^{ij}}_{1^{(a, b, p_i, q_j)}}(\mathbf{x}; \mathbf{U}, \mathbf{V})\right]_{1 \le i,j \le r}, 
    \end{equation*}
    where 
    \begin{align*}
        S &= 
        \det\left[ h_{j - i}(\mathbf{x} \,|\, \tau^{i - b - 2}\mathbf{U}_{b^{\varepsilon}_i}, \tau^{b - j + 1}\mathbf{V}_{a^{\varepsilon}_j}) \right]_{1 \le i, j \le n'}
    \end{align*}
\end{corollary} 
\begin{proof}
    The infinite ribbon $\hat{\theta}$ is $(\infty, (b + 2)^{n'})\ /\ ((b + 1)^{n'}, -\infty)$, and hence  $\theta_i \# \theta_j = 1^{(a, b, p_i, q_j)}$. The empty ribbon $\varepsilon$ is $\kappa/\kappa$ where $\kappa = ((b + 1)^{n'})$.
    Thus, using Theorem~\ref{thm:HG-sfG} and Corollary~\ref{cor:main-enumerator} the result follows. 
\end{proof}

\begin{corollary}[Jacobi--Trudi-type formula for $\mathsf{G}^{\mathbf{a}', \mathbf{b}'}_{\lambda'/\mu'}$]
    Let $\lambda'/\mu'$ be conjugate skew shape with $\lambda, \mu \in \mathsf{Sign}_n$. 
    Let $\Theta' = (\theta'_1, \ldots, \theta'_r)$ be outer decomposition of $\lambda'/\mu'$ into horizontal ribbons. 
        Let $p_i, q_i$ be the contents of head and tail of $\theta'_i$. 
    Define contents $b = -\mathrm{min}(\lambda_n, \mu_n) - 1 + n, a = -\mathrm{max}(\lambda_1, \mu_1) + 1$. 
The pipe function $\pi$ for $\lambda/\mu$ is given by: 
        \begin{equation*}
        \pi(c)  = 
        \begin{cases}
            (1, 1), \quad c \in [-b, -a] \\
            (0, 1), \quad c = -b - 1 \\
            (1, 0), \quad c = -a + 1 \\
            (0, 0), \quad \text{otherwise}. 
        \end{cases}
    \end{equation*}
    Let $\mathbf{a}'^{ij}, \mathbf{b}'^{ij}$ be induced column ribbon $\mathsf{G}$-flags. Let $\mathbf{a}'^{\varepsilon} = (a'^{\varepsilon}_i), \mathbf{b}'^{\varepsilon} = (b'^{\varepsilon}_i)$ be flags for the empty shape with conjugates $\mathbf{a}^{\varepsilon} = (a^{\varepsilon}_i), \mathbf{b}^{\varepsilon} = (b^{\varepsilon}_i)$. Denote 
    \[
    (a, b, p, q) := 
    \begin{cases}
        ((-a + 1)^{-a + q + 1},\ b^{-q + a})' / ((-a + 1)^{-a + p},\ (-a)^{-p + 1 + b})', \quad \lambda/\mu \text{ is positive}, \\
        ((-a + 1)^{-a + q},\ b^{-q + 1 + a})' / ((-a + 1)^{-a + p + 1},\ (-a)^{-p + b})', \quad \text{otherwise},
    \end{cases}
    \]
    where $p, q \in [a, b]$. 
    Denote $n' = b - a + 2$. Then we have 
    \begin{equation*}
        \mathsf{G}^{\mathbf{a}', \mathbf{b}'}_{\lambda'/\mu'}(\mathbf{x}; \mathbf{U}, \mathbf{V}) = S^{-r+1} \cdot \det \left[\mathsf{G}^{\mathbf{a}'^{ij}, \mathbf{b}'^{ij}}_{(a, b, p_i, q_j)}(\mathbf{x}; \mathbf{U}, \mathbf{V})\right]_{1 \le i,j \le r}, 
    \end{equation*}
    where 
    \begin{align*}
        S &= 
        \det \left[ e_{j - i}(\mathbf{x} \,|\, \tau^{2 - a - i}\mathbf{U}_{b'^{\varepsilon}_i}, \tau^{a + j - 1}\mathbf{V}_{a'^{\varepsilon}_j}) \right]_{1 \le i, j \le n'}
    \end{align*}
\end{corollary} 
\begin{proof}
    Note that $-a = \mathrm{max}(\lambda_1, \mu_1) - 1,\ -b = \mathrm{min}(\lambda_n, \mu_n) + 1 -n$ and $-p_i, -q_j$ are tail and head of $\theta_i, \theta_j$ (resp.), thus we have 
    \begin{align*}
        &\theta'_i \# \theta'_j = (\theta_j \# \theta_i)' \\
        &= 
        \begin{cases}
        ((-a + 1)^{-a + q + 1},\ b^{-q + a})' / ((-a + 1)^{-a + p},\ (-a)^{-p + 1 + b})', \quad \lambda/\mu \text{ is positive}, \\
        ((-a + 1)^{-a + q},\ b^{-q + 1 + a})' / ((-a + 1)^{-a + p + 1},\ (-a)^{-p + b})', \quad \text{otherwise}.
        \end{cases}
    \end{align*}
    Thus, using Corollaries~\ref{cor:dual-HG-sfG}, \ref{cor:main-dual-enumerator} the result follows. 
\end{proof}

\subsection{Skew Giambelli-type formulas} 
Let $\lambda/\mu$ be skew shape 
where $\lambda, \mu \in \mathsf{Sign}_n$. We define Frobenius notation for shape $\lambda/\mu$. Define the content $f = (\mu_n + 1) - 1 = \mu_n$. Let $k$ be the size of {\it Durfee square} top-left cell of which is $(1, f)$ and fit inside $\mathcal{L}(\lambda)$; similarly let $\ell$ be the size of Durfee square with the same top-left cell, but fit inside $\mathcal{L}(\mu)$. 
Let $u_i, v_i$ be the number of cells to the right and below (resp.) of the cell $(i, f + i)$ in $\mathcal{L}(\lambda)$ for $i \in [k]$. Then we can write the signature $\lambda$ in {\it Frobenius notation} $\lambda = (u_1, \ldots, u_k \mid v_1, \ldots, v_k)_f$. 
Similarly, let $\mu = (c_1, \ldots, c_{\ell} \mid d_1, \ldots, d_{\ell})_f$. 
Let $a = \mu_n - n$. 
The corresponding pipe function is defined as follows:
\begin{equation*}
    \pi(c) = 
    \begin{cases}
        (0, 0), &\quad c < a \\
        (0, 1), &\quad c = a \\
        (1, 1), &\quad a < c < f \\
        (1, 0), &\quad c = f \\
        (0, 0), &\quad c > f. \\
    \end{cases}
\end{equation*}

\begin{figure}
    \centering
    \begin{tikzpicture}[scale=0.4, font=\small]

    \foreach \x in {4,5,6,7} {
    \draw[thick] (\x,7) rectangle (\x+1,8);
    }

    \foreach \x in {4,5,6,7} {
    \draw[thick] (\x,6) rectangle (\x+1,7);
    }

    \foreach \x in {3,4,5,6} {
    \draw[thick] (\x,5) rectangle (\x+1,6);
    }

    \foreach \x in {1,2,3,4,5,6} {
    \draw[thick] (\x,4) rectangle (\x+1,5);
    }

    \foreach \x in {0,1,2,3,4,5} {
    \draw[thick] (\x,3) rectangle (\x+1,4);
    }

    \foreach \x in {0,1,2,3,4,5} {
    \draw[thick] (\x,2) rectangle (\x+1,3);
    }

    \foreach \x in {0,1,2,3,4} {
    \draw[thick] (\x,1) rectangle (\x+1,2);
    }

    \foreach \x in {0,1,2} {
    \draw[thick] (\x,0) rectangle (\x+1,1);
    }

    \draw [red, thick] (0.5, 0) -- (0.5, 4);
    \draw [red, thick] (1.5, 0) -- (1.5, 5);
    \draw [red, thick] (2.5, 0) -- (2.5, 5);
    \draw [red, thick] (3, 5.5) -- (7, 5.5);
    \draw [red, thick] (4, 6.5) -- (8, 6.5);
    \draw [red, thick] (4, 7.5) -- (8, 7.5);
    \draw [red, thick] (3.5, 1) -- (3.5, 4.5) -- (7, 4.5);
    \draw [red, thick] (4.5, 1) -- (4.5, 3.5) -- (6, 3.5);
    \draw [red, thick] (5.5, 2) -- (5.5, 2.5) -- (6, 2.5);

    \draw [red, densely dashed, thick] (0.5, 4) -- (0.5, 7.5) -- (4, 7.5);
    \draw [red, densely dashed, thick] (1.5, 5) -- (1.5, 6.5) -- (4, 6.5);
    \draw [red, densely dashed, thick] (2.5, 5) -- (2.5, 5.5) -- (3, 5.5);
    \draw [red, densely dashed, thick] (-5.3, 0.5) -- (-0.5, 0.5) -- (-0.5, 8);
    \draw [red, densely dashed, thick] (-5.3, 1.5) -- (-1.5, 1.5) -- (-1.5, 8);
    \draw [red, densely dashed, thick] (-5.3, 2.5) -- (-2.5, 2.5) -- (-2.5, 8);
    \draw [red, densely dashed, thick] (-5.3, 3.5) -- (-3.5, 3.5) -- (-3.5, 8);
    \draw [red, densely dashed, thick] (-5.3, 4.5) -- (-4.5, 4.5) -- (-4.5, 8);

    \draw [red, densely dashed, thick] (3.5, 0) -- (3.5, 1);
    \draw [red, densely dashed, thick] (4.5, 0) -- (4.5, 1);
    \draw [red, densely dashed, thick] (5.5, 0) -- (5.5, 2);
    \draw [red, densely dashed, thick] (6.5, 0) -- (6.5, 1.5) -- (13.3, 1.5);
    \draw [red, densely dashed, thick] (7.5, 0) -- (7.5, 0.5) -- (13.3, 0.5);
    
    \draw [red, densely dashed, thick] (8, 7.5) -- (13.3, 7.5);
    \draw [red, densely dashed, thick] (8, 6.5) -- (13.3, 6.5);
    \draw [red, densely dashed, thick] (7, 5.5) -- (13.3, 5.5);
    \draw [red, densely dashed, thick] (7, 4.5) -- (13.3, 4.5);
    \draw [red, densely dashed, thick] (6, 3.5) -- (13.3, 3.5);
    \draw [red, densely dashed, thick] (6, 2.5) -- (13.3, 2.5);

    \foreach \x in {-5, ..., 12} {
        \foreach \y in {0, ..., 7} {
            \draw[gray, dashed] (\x, \y) rectangle ++(1, 1);
        }
    }
    
    \fill[black] (0.5, -0.5) circle (0pt) node {$\theta_1$};
    \fill[black] (1.5, -0.5) circle (0pt) node {$\theta_2$};
    \fill[black] (2.5, -0.5) circle (0pt) node {$\theta_3$};
    \fill[black] (3.5, 0.5) circle (0pt) node {$\theta_7$};
    \fill[black] (4.5, 0.5) circle (0pt) node {$\theta_8$};
    \fill[black] (5.5, 0.5) circle (0pt) node {$\theta_9$};

    \fill[black] (9, 7.5) circle (0pt) node {$\theta_4$};
    \fill[black] (9, 6.5) circle (0pt) node {$\theta_5$};
    \fill[black] (9, 5.5) circle (0pt) node {$\theta_6$};

    \foreach \x in {0,...,8} {
        \draw[gray] (-5, \x) -- (-5.3, \x);
    }
    \foreach \x in {0,...,8} {
        \draw[gray] (13, \x) -- (13.3, \x);
    }

    \fill (0.5, 7.5) circle (4pt) node[above=0.2cm] {$f$};

    \end{tikzpicture}
\caption{Ribbon decomposition into hooks and strips. The cell with content $f$ is shown by the dot.}
\label{fig:hook-decomposition}
\end{figure}
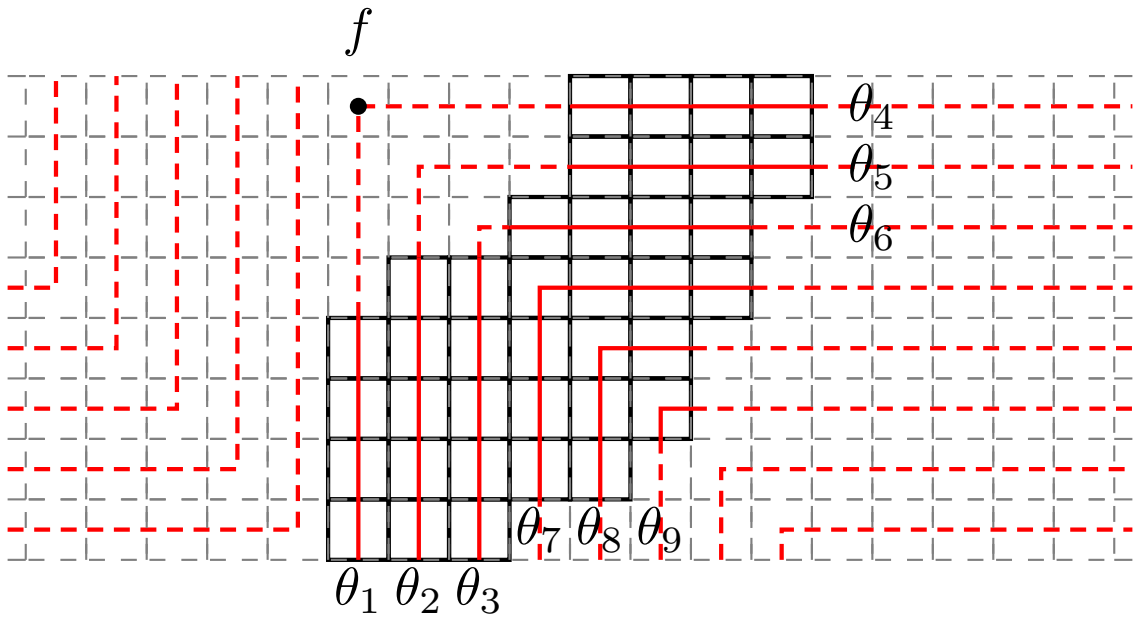

Let $\mathbf{a}, \mathbf{b}$ be row $\mathsf{G}$-flags for the shape $\lambda/\mu$, and $\Theta = (\theta_1, \ldots, \theta_{k + \ell})$ be the outer decomposition of the shape $\lambda/\mu$ into hooks and strips, as shown in Fig.~\ref{fig:hook-decomposition}. 

First let us assume that $\lambda/\mu$ is positive shape, i.e. $\mu \subseteq \lambda$. 
We specify the diagrams (ribbons) in the outer decomposition (Fig.~\ref{fig:hook-decomposition}) as follows:
$$
\theta_i = 
\left\{
\begin{alignedat}{3}
    &[f - v_{i}, f - d_{i} - 1]_{\hat{\theta}}, &&\quad i \in [\ell], && \quad \text{vertical ribbon}, \\
    &[f + c_{i - \ell} + 1, f + u_{i - \ell}]_{\hat{\theta}}, &&\quad i \in [\ell + 1, 2\ell], && \quad \text{horizontal ribbon}, \\
    &[f - v_{i - \ell}, f + u_{i - \ell}]_{\hat{\theta}}, &&\quad i \in [2\ell+1, k + \ell], && \quad \text{hook ribbon}. 
\end{alignedat}
\right.
$$
Define the following notations:
\begin{align*}
    (u_j \mid v_i)_f &:= 
    \begin{cases}
        \theta_i \# \theta_{j + \ell}, &\quad i \in [\ell], j \in [k], \\
        \theta_{\ell + i} \# \theta_{j + \ell}, &\quad i \in [\ell + 1, k], j \in [k], 
    \end{cases} \\[0.4cm]
    1^{(v_i, d_j)_f} &:= \begin{cases}
        \theta_i \# \theta_j, &\quad i, j \in [\ell], \\
        \theta_{\ell + i} \# \theta_j, &\quad i \in [\ell + 1, k], j \in [\ell],
    \end{cases} \\[0.4cm] 
    (c_i \mid d_j)_f &:= \theta_{i + \ell} \# \theta_j, \quad i,j \in [\ell], \\[0.4cm]
    (u_j, c_i)_f &:= \theta_{\ell + i} \# \theta_{\ell + j}, \quad i \in [\ell], j \in [k].
\end{align*}
The new notations are motivated by the diagrams of the corresponding shapes. For example, the diagram of $(u_j \mid v_i)_f$ is similar to hook $(u_j \mid v_i)$, where they differ only by the contents of cells. 


Define the following induced column flags:
\begin{itemize}
    \item $\mathbf{a}^{ij} = (a^{ij}_r)_r, \mathbf{b}^{ij} = (b^{ij}_r)_r$ 
    for hooks shapes $(u_j \mid v_i)_f$, $i,j\in [k]$
    \item $\mathbf{\tilde{a}}^{ij} = (\tilde{a}^{ij}_r)_r, \mathbf{\tilde{b}}^{ij} = (\tilde{b}^{ij}_r)_r$ 
    for negative hook shapes $(c_i \mid d_j)_f$, $i,j\in [\ell]$
    \item $\mathbf{p}^{ij} = (p^{ij}_r)_r,\mathbf{q}^{ij} = (q^{ij}_r)_r$
    for vertical shapes $1^{(v_i, d_j)_f}$, $i \in [k]\ , j\in [\ell]$
    \item $\mathbf{\tilde{p}}^{ij} = (\tilde{p}^{ij}_r)_r,\mathbf{\tilde{q}}^{ij} = (\tilde{q}^{ij}_r)_r$ 
    for horizontal shapes $(u_j, c_i)_f$, $i \in [\ell]\ , j\in [k]$
\end{itemize}

\begin{corollary}[Skew Giambelli formulas for positive shapes]
Let $\lambda/\mu$ be positive shape as above.  
Let $\varepsilon$ be the empty shape induced by $\pi$, and $\mathbf{a}^{\varepsilon} = (a^{\varepsilon}_i), \mathbf{b}^{\varepsilon} = (b^{\varepsilon}_i)$ be flags for the $\varepsilon$ with conjugates $\mathbf{a}'^{\varepsilon} = (a'^{\varepsilon}_i), \mathbf{b}'^{\varepsilon} = (b'^{\varepsilon}_i)$. Let $n' = f - a$. Then the following formulas hold: 
    \label{cor:Giambelli-main-shape}
    \[
    \mathsf{G}^{\mathbf{a}, \mathbf{b}}_{\lambda/\mu}(\mathbf{x}; \mathbf{U}, \mathbf{V})
    = (-1)^{\ell}\ S^{1 - k - \ell}\ \det \begin{bmatrix}
    \bigl[\mathsf{G}^{\mathbf{a}^{ij}, \mathbf{b}^{ij}}_{(u_j \mid v_i)_f} \bigr]_{\begin{subarray}{l}
            i,j \in [k] 
        \end{subarray}} 
    & \bigl[\mathsf{G}^{\mathbf{p}^{ij},\mathbf{q}^{ij}}_{1^{(v_i, d_j)_f}}\bigr]_{\begin{subarray}{l}
            i \in [k] \\ j \in [\ell] 
        \end{subarray}} \\\\
    \bigl[\mathsf{G}^{\mathbf{\tilde{p}}^{ij},\mathbf{\tilde{q}}^{ij}}_{(u_j, c_i)_f}\bigr]_{\begin{subarray}{l}
            i \in [\ell] \\ j \in [k] 
        \end{subarray}} 
    & \bigl[-\mathsf{G}^{\mathbf{\tilde{a}}^{ij},\mathbf{\tilde{b}}^{ij}}_{(c_i \mid d_j)_f}\bigr]_{\begin{subarray}{l}
            i,j \in [l] 
        \end{subarray}} 
\end{bmatrix},
    \]
    where 
    \begin{align*}
        S &= \det \left[
            h_{j-i} (\mathbf{x} \,|\, 
            \tau^{i - f - 1} \mathbf{U}_{b^{\varepsilon}_{i}} , \tau^{f - j} \mathbf{V}_{a^{\varepsilon}_j} )
        \right]_{1 \le i, j \le n'}.
    \end{align*}
\end{corollary}

\begin{proof}
    Let us specify the shapes $\theta_i \# \theta_j, i, j \in [k + \ell]$: 
    \[
    \theta_i \# \theta_j = 
    \begin{cases}
        1^{(v_i, d_j)_f},\ & i,j \in [\ell], \\[0.5em]
        (u_{j - \ell} \mid v_i)_f ,\ & i \in [\ell], j\in [\ell + 1, k + \ell], \\[0.5em]
        (u_{j - \ell}, c_{i - \ell})_f,\ & i \in [\ell + 1, 2\ell], j\in [\ell + 1, k + \ell], \\[0.5em]
        (c_{i - \ell} \mid d_j)_f  ,\ & i \in [\ell + 1, 2\ell]\ , j\in [\ell], \\[0.5em]
        1^{(v_{i - \ell}, d_j)_f} ,\ & i \in [2\ell + 1, k + \ell]\ , j\in [\ell], \\[0.5em]
        (u_{j - \ell} \mid v_{i - \ell})_f,\ & i \in [2\ell + 1, k + \ell], j\in [\ell + 1, k + \ell].
    \end{cases}
    \]
    Now, using our Hamel--Goulden-type formula we get the determinant of some 
    matrix which can be split into block matrices. 
    To simplify notation, we shall use only shapes (i.e. indices) of corresponding polynomials:
    $$
    \begin{bmatrix}
        \bigl[ 1^{(v_i, d_j)_f} \bigr]_{i,j \in [\ell]} 
        & \bigl[ (u_{j - \ell} \mid v_i)_f \bigr]_{\begin{subarray}{l}
            i \in [\ell] \\ j \in [\ell + 1, k + \ell] 
        \end{subarray}}\\\\
        \bigl[ (c_{i - \ell} \mid d_j)_f \bigr]_{\begin{subarray}{l}
            i \in [\ell + 1, 2\ell] \\ j \in [\ell] 
        \end{subarray}} 
        &  
        \bigl[ (u_{j-l}, c_{i - \ell})_f \bigr]_{\begin{subarray}{l}
            i \in [\ell + 1, 2\ell] \\ j \in [\ell + 1, k + \ell] 
        \end{subarray}} \\\\
        \bigl[ 1^{(v_{i - \ell}, d_j)_f} \bigr]_{\begin{subarray}{l}
            i \in [2\ell + 1, k + \ell] \\ j \in [\ell] 
        \end{subarray}} 
        &
        \bigl[ (u_{j - \ell} \mid v_{i - \ell})_f \bigr]_{\begin{subarray}{l}
            i \in [2\ell + 1, k + \ell] \\ j \in [\ell + 1, k + \ell] 
        \end{subarray}} 
    \end{bmatrix}.
    $$
    which can be rewritten as:
    $$
    \begin{bmatrix}
        \bigl[ 1^{(v_i, d_j)_f} \bigr]_{i,j \in [\ell]} 
        & \bigl[ u_{j} \mid v_i \bigr]_{\begin{subarray}{l}
            i \in [\ell] \\ j \in [k] 
        \end{subarray}}\\\\
        \bigl[ (c_i \mid d_j)_f \bigr]_{\begin{subarray}{l}
            i, j \in [\ell] 
        \end{subarray}} 
        &  
        \bigl[ (u_{j}, c_{i})_f \bigr]_{\begin{subarray}{l}
            i \in [\ell] \\ j \in [k] 
        \end{subarray}} \\\\
        \bigl[ 1^{(v_i - d_j)_f} \bigr]_{\begin{subarray}{l}
            i \in [\ell + 1, k] \\ j \in [\ell] 
        \end{subarray}} 
        &
        \bigl[ (u_j \mid v_i)_f \bigr]_{\begin{subarray}{l}
            i \in [\ell + 1, k] \\ j \in [k] 
        \end{subarray}} 
    \end{bmatrix}.
    $$
    Permuting the second and third row blocks, and  
    the first and second column blocks gives 
    the desired matrix. The number of overall 
    permutations of rows and columns is  
    $2k\ell - \ell^2$ , and so we need to multiply the 
    determinant by $(-1)^{2k\ell - \ell^2} = (-1)^{\ell}$.
    
        Also, we have the infinite ribbon $\hat{\theta} = (\infty, (f + 1)^{n'}) / (f^{n'}, -\infty)$, and hence using the rule for obtaining the empty shape we deduce that $\epsilon = (f^{n'}) / (f^{n'})$. 
\end{proof}

Now we give the Giambelli-type formulas for the negative shape $\mu/\lambda$. The ribbons $\theta_i \in \Theta$ are the same, except that they are now negative (as the shape $\mu/\lambda$ is negative).


Define the following notations:
\begin{align*}
    (u_i \mid v_j)_f &:= 
    \begin{cases}
        \theta_{i + \ell} \# \theta_j, &\quad i \in [k], j \in [\ell], \\
        \theta_{i + \ell} \# \theta_{j + \ell}, &\quad i \in [k], j \in [\ell + 1, k], 
    \end{cases} \\[0.4cm]
    (c_j \mid d_i)_f &:= \theta_i \# \theta_{j + \ell}, \quad i,j \in [\ell], \\[0.4cm]
    1^{(v_j, d_i)_f} &:= \begin{cases}
        \theta_i \# \theta_j, &\quad i, j \in [\ell], \\
        \theta_i \# \theta_{j + \ell}, &\quad i \in [\ell], j \in [\ell + 1, k],
    \end{cases} \\[0.4cm] 
    (u_i, c_j)_f &:= \theta_{i + \ell} \# \theta_{j + \ell}, \quad i \in [k], j \in [\ell].
\end{align*}
Define the following induced column flags:
\begin{itemize}
    \item $\mathbf{a}^{ij} = (a^{ij}_r)_r, \mathbf{b}^{ij} = (b^{ij}_r)_r$ 
    for hook shapes $(u_i \mid v_j)_f$, $i,j\in [k]$
    \item $\mathbf{\tilde{a}}^{ij} = (\tilde{a}^{ij}_r)_r, \mathbf{\tilde{b}}^{ij} = (\tilde{b}^{ij}_r)_r$ 
    for positive hook shapes $(c_j \mid d_i)_f$, $i,j\in [\ell]$
    \item $\mathbf{p}^{ij} = (p^{ij}_r)_r,\mathbf{q}^{ij} = (q^{ij}_r)_r$
    for vertical shapes $1^{(v_j, d_i)_f}$, $i \in [\ell]\ , j\in [k]$
    \item $\mathbf{\tilde{p}}^{ij} = (\tilde{p}^{ij}_r)_r,\mathbf{\tilde{q}}^{ij} = (\tilde{q}^{ij}_r)_r$ 
    for horizontal shapes $(u_i, c_j)_f$, $i \in [k]\ , j \in [\ell]$
\end{itemize}

\begin{corollary}[Skew Giambelli formulas for negative shapes]
Let $\mu/\lambda$ be negative shape as above.  
Let $\varepsilon$ be the empty shape induced by $\pi$, and $\mathbf{a}^{\varepsilon} = (a^{\varepsilon}_i), \mathbf{b}^{\varepsilon} = (b^{\varepsilon}_i)$ be flags for the $\varepsilon$ with conjugates $\mathbf{a}'^{\varepsilon} = (a'^{\varepsilon}_i), \mathbf{b}'^{\varepsilon} = (b'^{\varepsilon}_i)$. Let $n' = f - a$. The following formulas hold: 
    \label{cor:Giambelli-neg-main}
    \[
    \mathsf{G}^{\mathbf{a}, \mathbf{b}}_{\mu/\lambda}(\mathbf{x}; \mathbf{U}, \mathbf{V})
    = (-1)^k\ S^{1 - k - \ell}\ \det \begin{bmatrix}
    \bigl[-\mathsf{G}^{\mathbf{a}^{ij}, \mathbf{b}^{ij}}_{(u_i \mid v_j)_f}\bigr]_{\begin{subarray}{l}
            i,j \in [k] 
        \end{subarray}} 
    & \bigl[\mathsf{G}^{\mathbf{p}^{ij},\mathbf{q}^{ij}}_{1^{(v_j, d_i)_f}}\bigr]_{\begin{subarray}{l}
            i \in [k] \\ j \in [\ell] 
        \end{subarray}} \\\\
    \bigl[\mathsf{G}^{\mathbf{\tilde{p}}^{ij},\mathbf{\tilde{q}}^{ij}}_{(u_i, c_j)_f}\bigr]_{\begin{subarray}{l}
            i \in [\ell] \\ j \in [k] 
        \end{subarray}} 
    & \bigl[\mathsf{G}^{\mathbf{\tilde{a}}^{ij},\mathbf{\tilde{b}}^{ij}}_{(c_j \mid d_i)_f}\bigr]_{\begin{subarray}{l}
            i,j \in [l] 
        \end{subarray}} 
\end{bmatrix},
    \]
    where 
    \begin{align*}
        S &= \det \left[
            h_{j-i} (\mathbf{x} \,|\, 
            \tau^{i - f - 1}\mathbf{U}_{b^{\epsilon}_i}, \tau^{f - j}\mathbf{V}_{a^{\epsilon}_j}) 
        \right]_{1 \le i, j \le n'}.
    \end{align*}
\end{corollary}

\begin{proof}
    Applying Hamel--Goulden-type formula we get the transpose of the matrix from the proof of Corollary~\ref{cor:Giambelli-main-shape}, thus we need to multiply the determinant by the same $(-1)^{\ell}$. As the shape $\mu/\lambda$ is negative and the number of hook pipes in cells of $\lambda/\mu$ is $k - \ell$, and the $\mathrm{sgn}_{\pi}(\mu/\lambda) = (-1)^{k - \ell}$; thus the result follows. 
\end{proof}

\section{Some specializations}\label{sec:spec1}
In this section we list several specializations of flagged Grothendieck enumerators which we show by verifying the known Jacobi--Trudi-type formulas or tableaux formulas for the corresponding functions. 

\subsection{Some useful functions}
For $n \in \mathbb{Z}$ and variable sets $\mathbf{x}, \mathbf{y}$ let us also define the following {\it supersymmetric} versions of elementary and complete homogenous symmetric functions:
    \[
        h_n(\mathbf{x} / \mathbf{y}) := \sum_{p + q = n} h_p(\mathbf{x})\,(-1)^q\,e_q(\mathbf{y}), \qquad e_n(\mathbf{x} / \mathbf{y}) := \sum_{p + q = n} e_p(\mathbf{x})\,(-1)^q\,h_q(\mathbf{y}). 
    \]
    We also define the following functions: 
    \[
        h_n(\mathbf{x} \,|\, \mathbf{y} / \mathbf{z}) := \sum_{k \ge 0} h_{n + k}(\mathbf{x}) \cdot h_k(\mathbf{y} / \mathbf{z}), \qquad 
         e_n(\mathbf{x} \,|\, \mathbf{y} / \mathbf{z}) := \sum_{k \ge 0} e_{n + k}(\mathbf{x}) \cdot e_k(\mathbf{y} / \mathbf{z}). 
    \]
    with an extra variable set $\mathbf{z}$.  
    
    \begin{remark}
     We note that the functions $h_n(\mathbf{x} \,|\, \mathbf{y} / \mathbf{z})$ were denoted using $\ominus$ operation as $h_n(\mathbf{x} \ominus (\mathbf{y} / \mathbf{z}))$ in \cite{hwang1}.
     \end{remark}

\begin{lemma}[\cite{kim2}]
    \label{lem:plethystic-identities}
    Let $n \in \mathbb{Z}_{\ge 0}$, $a, b \in \mathbb{Z}$. The following identities hold:
    \begin{align*}
        h_{n}(\mathbf{z}_a / \mathbf{y}_b) &= h_n(\mathbf{z}_{a - 1} / \mathbf{y}_b) + z_a \cdot h_{n - 1}(\mathbf{z}_a / \mathbf{y}_b), \\
        h_{n}(\mathbf{z}_a / \mathbf{y}_b) &= h_n(\mathbf{z}_a / \mathbf{y}_{b - 1}) - y_b \cdot h_{n - 1}(\mathbf{z}_a / \mathbf{y}_{b - 1}).
    \end{align*}
\end{lemma}

We shall use the following variable sets:  
\[
\mathbf{y} = (y_i)_{i \in \mathbb{Z}}, \quad \mathbf{z} = (z_i)_{i \in \mathbb{Z}}, \quad \mathbf{y'} = (y'_i)_{i \in \mathbb{Z}}, \quad \mathbf{z'} = (z'_i)_{i \in \mathbb{Z}}. 
\]

\begin{lemma}
    \label{lem:spec-to-supersym}
    Let $n, c_0, c_1, a, b \in \mathbb{Z}$ and $a' = a + c_0, b' = b - c_1$. Then 
    \begin{equation}
        \label{eq:5}
        h_{n}(\mathbf{x} \,|\, \tau^{-c_1}\mathbf{U}_b, \tau^{c_0}\mathbf{V}_a) \big|_{\substack{
            u_{i, j} \rightarrow y'_i - z'_j \\
            v_{i, j} \rightarrow z_i - y_j
        }} = 
        h_{n}(\mathbf{x} \,|\, \mathbf{y}'_b, \mathbf{z}_a / \mathbf{z}'_{b'}, \mathbf{y}_{a'}). 
    \end{equation}
\end{lemma}
\begin{proof}
By definition, we have:
\[
h_{n}(\mathbf{x} \,|\, \tau^{-c_1}\mathbf{U}_b, \tau^{c_0}\mathbf{V}_a) = \sum_{k \ge 0}h_{n + k}(\mathbf{x}) \sum_{p + q = k} h_{p, c_0}(\mathbf{V}_a) h_{q, -c_1}(\mathbf{U}_b),
\]
It is enough to show that for all $p \in \mathbb{Z}$ the following equation holds: 
    \[
        h_{p, c_0}(\mathbf{V}_a) \big|_{v_{i, j} \rightarrow z_i - y_j} = h_p(\mathbf{z}_a / \mathbf{y}_{a'}).
    \]
    For $p \le 0$ it follows from definitions of $h_{n, c}, h_n$. Suppose that $p \ge 1$, then we proceed by induction on $p$. The statement is obvious for $p = 1$ and suppose it holds for all $p \le r$ where $r \ge 1$, let us show that it is also true for $p = r + 1$. 
    Using the induction hypothesis we have 
    \begin{align}
        h_{r + 1, c_0}(\mathbf{V}_a) \big|_{v_{i, j} \rightarrow z_i - y_j} = \sum_{d = 0}^{\infty}h_r(\mathbf{z}_{a - d} / \mathbf{y}_{a' - 1 - d})(z_{a - d} - y_{a' - d}). \label{eq:4}
    \end{align}
    Now from the Lemma~\ref{lem:plethystic-identities} it follows that:
    \begin{equation*}
        h_{r + 1}(\mathbf{z}_{a - d} / \mathbf{y}_{a' - d}) = h_{r + 1}(\mathbf{z}_{a - d - 1} / \mathbf{y}_{a' - d - 1}) + h_r(\mathbf{z}_{a - d} / \mathbf{y}_{a' - d - 1}) (z_{a - d} - y_{a' - d}).
    \end{equation*}
    Using the last equation we can rewrite the sum in (\ref{eq:4}) as follows:
    \begin{align*}
        h_{r + 1, c_0}(\mathbf{V}_a) \big|_{v_{i, j} \rightarrow z_i - y_j} = \sum_{d = 0}^{\infty} h_{r + 1}(\mathbf{z}_{a - d} / \mathbf{y}_{a' - d}) - h_{r + 1}(\mathbf{z}_{a - d - 1} / \mathbf{y}_{a' - d - 1}) 
        = h_{r + 1}(\mathbf{z}_a / \mathbf{y}_{a'}),
    \end{align*}
    as needed. 
\end{proof}

\begin{lemma}
    \label{lem:spec-to-supersym-dual}
    Let $n, c'_0, c'_1, a', b' \in \mathbb{Z}$ and $a = a' + c'_0, b = b' - c'_1$. Then 
    \begin{equation*}
        e_{n}(\mathbf{x} \,|\, \tau^{-c'_1}\mathbf{U}_{b'}, \tau^{c'_0}\mathbf{V}_{a'})
        \big|_{\substack{
            u_{i, j} \rightarrow y'_i - z'_j \\
            v_{i, j} \rightarrow z_i - y_j
        }} = 
        e_{n}(\mathbf{x} \,|\, \mathbf{y}'_{b'}, \mathbf{z}_{a'} / \mathbf{z}'_b, \mathbf{y}_a). 
    \end{equation*}
\end{lemma}
\begin{proof}
    Let $c_0 = -c'_0, c_1 = -c_1$, then $a' = a + c_0, b' = b - c_0$. 
    Note that 
    \[
    h_n(\mathbf{x} / \mathbf{y}) = e_n(\overline{\mathbf{y}} / \overline{\mathbf{x}}), 
    \]
    where we denote $\overline{\mathbf{x}} = (-x_i)$. 
    Using the above equation and by applying the involution $\omega$ to both sides of $\eqref{eq:5}$, we have 
    \[
        e_{n}(\mathbf{x} \,|\, (\tau^{c_1}\mathbf{U}_{b'}, \tau^{-c_0}\mathbf{V}_{a'})) \big|_{\substack{
            u_{i, j} \rightarrow y'_j - z'_i \\
            v_{i, j} \rightarrow z_j - y_i
        }} = 
        e_{n}(\mathbf{x} \,|\, \overline{\mathbf{z}}'_{b'}, \overline{\mathbf{y}}_{a'} / \overline{\mathbf{y}}'_b, \overline{\mathbf{z}}_a). 
    \]
    Then  the result follows after changing the variables $z_i = -y_i,\, y_i = -z_i$ and $z'_i = -y'_i,\, y'_i = -z'_i$ for all $i \in \mathbb{Z}$. 
\end{proof}
\begin{definition}
    It will be useful to introduce the following functions: 
    \begin{align*}
        \mathsf{G}^{\mathbf{a}, \mathbf{b}}_{\lambda/\mu}(\mathbf{x}; \mathbf{y}, \mathbf{z}, \mathbf{y}', \mathbf{z}') &:= 
        \mathsf{G}^{\mathbf{a}, \mathbf{b}}_{\lambda/\mu}(\mathbf{x}; \mathbf{U}, \mathbf{V})\big|_{\substack{
            u_{i, j} \rightarrow y'_i - z'_j \\
            v_{i, j} \rightarrow z_i - y_j
        }}, \\ 
        \mathsf{G}^{\mathbf{a}', \mathbf{b}'}_{\lambda'/\mu'}(\mathbf{x}; \mathbf{y}, \mathbf{z}, \mathbf{y}', \mathbf{z}') &:= \mathsf{G}^{\mathbf{a}', \mathbf{b}'}_{\lambda'/\mu'}(\mathbf{x}; \mathbf{U}, \mathbf{V})\big|_{\substack{
            u_{i, j} \rightarrow y'_i - z'_j \\
            v_{i, j} \rightarrow z_i - y_j
        }}.
    \end{align*}
\end{definition}

\begin{lemma}
    \label{lem:jacobi-trudi-for-spec-1}
    Let $\lambda, \mu \in \mathsf{Sign}_n$ with row $\mathsf{G}$-flags $\mathbf{a} = (a_i), \mathbf{b} = (b_i)$ and conjugate flags $\mathbf{a}' = (a'_i), \mathbf{b}' = (b'_i)$. Then we have   
    \begin{align*}
        \mathsf{G}^{\mathbf{a}, \mathbf{b}}_{\lambda/\mu}(\mathbf{x}; \mathbf{y}, \mathbf{z}, \mathbf{y}', \mathbf{z}') 
        &= \det \left[ h_{\lambda_i - \mu_j - i + j}(\mathbf{x} \,|\, \mathbf{y}'_{b_i}, \mathbf{z}_{a_j} / \mathbf{y}_{a'_j}, \mathbf{z}'_{b'_i}) \right]_{1 \le i, j \le n},  \\ 
        \mathsf{G}^{\mathbf{a}', \mathbf{b}'}_{\lambda'/\mu'}(\mathbf{x}; \mathbf{y}, \mathbf{z}, \mathbf{y}', \mathbf{z}') 
        &= \det \left[ e_{\lambda_i - \mu_j - i + j}(\mathbf{x} \,|\, \mathbf{y}'_{b'_i}, \mathbf{z}_{a'_j} / \mathbf{y}_{a_j}, \mathbf{z}'_{b_i}) \right]_{1 \le i, j \le n}. 
    \end{align*}
\end{lemma}
\begin{proof}
    Directly follows from Corollaries~\ref{cor:main-enumerator}, \ref{cor:main-dual-enumerator} and Lemmas~\ref{lem:spec-to-supersym}, \ref{lem:spec-to-supersym-dual}. 
\end{proof}

\subsection{Canonical refined stable Grothendieck polynomials}
Using Lemma~\ref{lem:jacobi-trudi-for-spec-1}, it is also not hard to see the following special case. 

\begin{proposition} 
Let $\lambda, \mu \in \mathsf{Par}_n$ with $\lambda \supseteq \mu$ and $\mathsf{G}$-flags $\mathbf{a} = (i), \mathbf{b} = (\lambda_i)$ whose conjugates are $\mathbf{a}' = (\mu_i), \mathbf{b}' = (i - 1)$. 
Define $G$-specialization by substituting: 
  $$u_{i,j} = \alpha_i - \beta_j, \quad  
  v_{i,j} = \beta_i - \alpha_j,$$ 
  and $\alpha_i = \beta_i = 0$ for $i \le 0$.
Then flagged Grothendieck enumerators specialize to the canonical refined stable Grothendieck polynomials $G_{\lambda/\mu}(\mathbf{x}; \boldsymbol{\alpha}, \boldsymbol{\beta})$ as follows:   
\begin{align*}
G_{\lambda/\mu}(\mathbf{x}; \boldsymbol{\alpha}, \boldsymbol{\beta}) &= 
C \cdot \mathsf{G}^{\mathbf{a}, \mathbf{b}}_{\lambda/\mu}(\mathbf{x}; \mathbf{U}, \mathbf{V}) |_{G\text{-spec}}, \quad C = \prod_{i = 1}^n \prod_{j} (1 - \beta_i x_j), \\
G_{\lambda'/\mu'}(\mathbf{x}; \boldsymbol{\alpha}, \boldsymbol{\beta}) &= 
D \cdot \mathsf{G}^{\mathbf{a}', \mathbf{b}'}_{\lambda'/\mu'}(\mathbf{x}; \mathbf{U}, \mathbf{V}) |_{G\text{-spec}}, \quad D = \prod_{i = 1}^n \prod_{j} (1 - \alpha_i x_j)^{-1}.
\end{align*}
\end{proposition}

\subsection{Dual canonical refined stable Grothendieck polynomials}
\begin{proposition}
Let $\lambda, \mu \in \mathsf{Par}_n$ with $\lambda \supseteq \mu$. Consider the $\mathsf{G}$-flags $\mathbf{a} = (m+i-1)$, $\mathbf{b} = (\mu_i)$ with conjugates $\mathbf{a}' = (m+ \lambda_i - 1)$, $\mathbf{b}' = (i-1)$ for negative shape $\mu/\lambda$. Define {\it $g$-specialization} as the following substitution of variables:
    $$
    u_{i,j} = \alpha_i - \beta_j, \quad v_{i,j} = z_i - y_j, 
    $$
    where $z_i = x_i,\,y_i = 0$ for $i \in [m]$, and $z_i = \beta_{i -m},\,y_i = \alpha_{i - m}$ for $i > m$, and $\alpha_i = \beta_i  = z_i = y_i = 0$ for $i \le 0$.
Then flagged Grothendieck enumerators specialize to dual canonical refined stable Grothendieck polynomials $g_{\lambda/\mu}(\mathbf{x}; \boldsymbol{\alpha}, \boldsymbol{\beta})$ as follows: 
$$
g_{\lambda/\mu}(\mathbf{x}_m; \boldsymbol{\alpha}, \boldsymbol{\beta}) = \mathsf{G}^{\mathbf{a}, \mathbf{b}}_{\mu/\lambda}(\varnothing; \mathbf{U}, \mathbf{V}) |_{g\text{-spec}}.
$$
\end{proposition}
\begin{proof}
    Note that for $k \ge 1,\, \ell \ge 0$ we have  
    \[
        h_{k} (\varnothing \,|\, \mathbf{x} / \mathbf{y}) = h_{-k}(\mathbf{x} / \mathbf{y}) = 0, \qquad h_{-\ell} (\varnothing \,|\, \mathbf{x} / \mathbf{y}) = h_{\ell}(\mathbf{x} / \mathbf{y}). 
    \]
Using Lemma~\ref{lem:jacobi-trudi-for-spec-1}, we obtain  
\begin{align*}
    \mathsf{G}^{\mathbf{a}, \mathbf{b}}_{\mu/\lambda}(\varnothing; \mathbf{U}, \mathbf{V}) |_{g\text{-spec}} 
    &= \det \left[ h_{\mu_i - \lambda_j - i + j} (\varnothing \,|\, \mathbf{x}_m, \boldsymbol{\alpha}_{\mu_i}, {\boldsymbol{\beta}_{j - 1}} / \boldsymbol{\alpha}_{\lambda_j - 1}, \boldsymbol{\beta}_{i - 1}) \right]_{1 \le i, j \le n}  \\
    &= \det \left[ h_{\lambda_i - \mu_j - i + j} (\mathbf{x}_m, \boldsymbol{\alpha}_{\mu_j}, \boldsymbol{\beta}_{i - 1} / \boldsymbol{\alpha}_{\lambda_i - 1}, \boldsymbol{\beta}_{j - 1} ) \right]_{1 \le i, j \le n}, \\
    &= g_{\lambda/\mu}(\mathbf{x}_m; \boldsymbol{\alpha}, \boldsymbol{\beta}),  
\end{align*}
where we refered to the known Jacobi--Trudi formula for $g_{\lambda/\mu}$ \cite{hwang1}. 
\end{proof}



\begin{remark} 
As flagged Grothendieck enumerators generalize both stable Grothendieck polynomials and their duals, it is interesting to compare them with hybrid Grothendieck polynomials studied in \cite{guo}, and self-dual inhomogeneous $q$-Whittaker polynomials studied in \cite{gunna}. 
\end{remark}

\subsection{Factorial supersymmetric Schur functions }
We refer to \cite{mol98} for these functions.  

\begin{proposition}
Let $\lambda, \mu \in \mathsf{Par}_n$ with $\lambda \supseteq \mu$. Consider the $\mathsf{G}$-flags $\mathbf{a} = (m)$, $\mathbf{b} = (m' + \mu_i + 1 - i),\ m, m' \ge 1$ for negative shape $\mu/\lambda$. Define {\it $f$-specialization} as the following substitution of variables:
    $$
    u_{i,j} = y_j + z_i, \, v_{i,j} = x_i - z_j, 
    $$
where $y_i = x_i = 0$ for $i \le 0$.
Then flagged Grothendieck enumerators specialize to factorial supersymmetric Schur funcions $s_{\lambda/\mu}(\mathbf{x}_m / \mathbf{y}_{m'}\,|\,\mathbf{z})$ as follows: 
$$
s_{\lambda/\mu}(\mathbf{x}_m / \mathbf{y}_{m'}\,|\,\mathbf{z}) = \mathsf{G}^{\mathbf{a}, \mathbf{b}}_{\mu/\lambda}(\varnothing; \mathbf{U}, \mathbf{V}) |_{f\text{-spec}}.
$$
\end{proposition}
\begin{proof}
It follows that the function $\mathsf{G}^{\mathbf{a}, \mathbf{b}}_{\mu/\lambda}(\varnothing; \mathbf{U}, \mathbf{V}) |_{f\text{-spec}}$ enumerates trio tableaux $(T^{\circ}, \varnothing, T^{\bullet})$ such that $T^{\circ} \in TT^{\circ}_{\lambda/\nu}, T^{\bullet} \in TT^{\bullet}_{\nu/\mu}$ where $\nu \in \mathrm{Par}_n$ with $\mu \subseteq \nu \subseteq \lambda$. On the other hand, $s_{\lambda/\mu}(\mathbf{x}_m / \mathbf{y}_{m'}\,|\,\mathbf{z})$ enumerates tableaux $(T, T')$ of shapes $\lambda/\nu,\, \nu/\mu$ (resp.) where $\nu \in \mathrm{Par}_n$ with $\mu \subseteq \nu \subseteq \lambda$; the tableau $T$ is semistandard and filled with elements from $[m]$, and $T'$ is filled by elements from $[m']$ such that they strictly decrease along rows, and weakly decrease along columns, and the weight is defined as follows (see \cite{mol98}): 
\[
    \mathrm{wt}(T) = \prod_{i,j} (x_{T_{i, j}} - z_{T_{i, j} + c(i, j)}), \qquad \mathrm{wt}(T') = \prod_{i,j} (y_{T'_{i, j}} + z_{T'_{i, j} + c(i, j)}). 
\]
We define the map $\phi' : (T^{\circ}, T^{\bullet}) \mapsto (T, T')$, where $T = T^{\circ}$, and $T'$ is defined by $T'_{i, j} = T^{\bullet}_{i, j} - c(i, j)$. It is not hard to see that $\phi'$ is a weight-preserving bijection, and the result follows.     
\end{proof}

\

\subsection{Macdonald's ninth variation of Schur functions.} 
Ninth variation of Schur functions are due to Macdonald \cite{macdonald1}, and we refer to \cite{fkingb} for the tableau version used here.  

\begin{proposition}
Let $\lambda, \mu \in \mathsf{Par}_n$ with $\lambda \supseteq \mu$. Consider the $\mathsf{G}$-flags $\mathbf{a} = (m)$, $\mathbf{b} = (-\infty), m \ge 1$ for negative shape $\mu/\lambda$. Define variable set of matrix $\mathbf{X} = (x_{i,j})$ for $i \in [m], j \in \mathbb{Z}$. Define {\it $s$-specialization} as the following substitution of variables:
    $$
    v_{i,j} = x_{i, j - i},  
    $$
where $u_{i, j} = 0$ for $i, j \in \mathbb{Z}$ and $v_{i, j} = 0$ for $i \in \mathbb{Z} \setminus [m], j \in \mathbb{Z}$.
Then flagged Grothendieck enumerators specialize to tableau version of ninth variation of Schur functions $s_{\lambda/\mu}(\mathbf{X})$ as follows: 
$$
s_{\lambda/\mu}(\mathbf{X}) = \mathsf{G}^{\mathbf{a}, \mathbf{b}}_{\mu/\lambda}(\varnothing; \mathbf{U}, \mathbf{V}) |_{s\text{-spec}}.
$$
\end{proposition}
\begin{proof}
    It follows that the function $\mathsf{G}^{\mathbf{a}, \mathbf{b}}_{\mu/\lambda}(\varnothing; \mathbf{U}, \mathbf{V}) |_{s\text{-spec}}$ enumerates trio tableaux $(T^{\circ}, \varnothing, \varnothing)$ such that $T^{\circ} \in TT^{\circ}_{\lambda/\mu}$, and moreover 
    \[
    \mathrm{wt}(T^{\circ}) = \prod_{i,j} x_{T^{\circ}_{i, j}, c(i, j)}.
    \]
   We note that $s_{\lambda/\mu}(\mathbf{X})$ enumerates the same tableaux (see \cite{fkingb}), and the result follows. 
\end{proof}



\newpage
\begin{thebibliography}{abcdefghij}

\bibitem[AY25]{ay25}
A. Adilzhan and D. Yeliussizov, Flagged Hamel--Goulden formulas, arXiv:2512.13265, 2025. 


\bibitem[Buc02]{buch}
A. Buch, A Littlewood Richardson rule for the K-theory of Grassmannians, Acta Math. {\bf 189} (2002), 37--78.

\bibitem[CP21]{cp}
M. Chan, and N. Pflueger, Combinatorial relations on skew Schur and skew stable Grothendieck polynomials. Algebraic Combin. {\bf 4} (2021), 175--188.

\bibitem[CLL02]{cll}
W. Y. C. Chen, B. Li, and J. D. Louck, The flagged double Schur function, J. Algebraic Combin. {\bf 15} (2002), 7--26.

\bibitem[CYY05]{cyy}
W. Y. C. Chen, G.-G. Yan, A. L. B. Yang, Transformations of border strips and Schur function determinants, J. Algebraic Combin. {\bf 21} (2005), 379--394.


\bibitem[FK21]{fkingb}
A. M. Foley and R. C. King, Determinantal and Pfaffian identities for ninth variation skew Schur functions and $Q$-functions. European J. Combin. {\bf 93} (2021), 103271.

\bibitem[FK94]{fk1}
S. Fomin, A. N. Kirillov, Grothendieck polynomials and the Yang-Baxter equation, Proc. 6th Intern. Conf. on Formal Power Series and Algebraic Combinatorics, DIMACS, (1994), 183--190.

\bibitem[FK96]{fk2}
S. Fomin, and A. N. Kirillov, The Yang-Baxter equation, symmetric functions, and Schubert polynomials, Discrete Math. {\bf 153} (1996), 123--143.

\bibitem[GGL16]{ggl}
P. Galashin, D. Grinberg, and G. Liu, Refined dual stable Grothendieck polynomials and generalized Bender-Knuth involutions, Electronic J. Combin. {\bf 23} (2016): 3--14.

\bibitem[Gal17]{gal}
P. Galashin, A Littlewood–Richardson rule for dual stable Grothendieck polynomials. J. Combin. Theory Ser. A, {\bf 151} (2017), 23--35.

\bibitem[GP24]{gp}
S. Gavrilova and L. Petrov, Tilted biorthogonal ensembles, Grothendieck random partitions, and determinantal tests, Selecta Math. {\bf 30} (2024), 56.

\bibitem[GV89]{lgv}
I. Gessel and X. Viennot, Determinants, paths, and plane partitions, preprint (1989).

\bibitem[GWZ25]{gunna}
A. Gunna, M. Wheeler, and P. Zinn-Justin, Inhomogeneous $q$-Whittaker Polynomials I: Duality and Expansions, arXiv:2512.04468, 2025.

\bibitem[GKL25]{guo}
P. L. Guo, M. Kang, and J. Liu, Hybrid Grothendieck polynomials. arXiv:2505.19072, 2025.

\bibitem[HG95]{hg}
A. M. Hamel and I. P. Goulden, Planar decompositions of tableaux and Schur function determinants, European J. Combin. {\bf 16} (1995), 461--477.

\bibitem[HS20]{hs}
G. Hawkes and T. Scrimshaw, Crystal structures for canonical Grothendieck functions. Algebraic Combin. {\bf 3} (2020), 727--755.

\bibitem[HJK+24]{hwang1}
B. H. Hwang, J. Jang, J. S. Kim, M. Song, and U. K.  Song, Refined canonical stable Grothendieck polynomials and their duals, Part 1,  Adv. Math. {\bf 446} (2024), article 109670.

\bibitem[HJK+25]{hwang2}
B. H. Hwang, J. Jang, J. S. Kim, M. Song, and U. K.  Song,  Refined canonical stable Grothendieck polynomials and their duals, Part 2, European J. Combin. {\bf 127} (2025), article 104166. 

\bibitem[IMS24]{ims}
S. Iwao, K. Motegi, and T. Scrimshaw, Free fermions and canonical Grothendieck polynomials, Algebraic Combin. {\bf 7} (2024), 245--274.

\bibitem[IMS25]{ims2}
S. Iwao, K. Motegi, and T. Scrimshaw, Free fermionic probability theory and K-theoretic Schubert calculus. Forum Math. Sigma {\bf 13} (2025), e197. 
 
\bibitem[Kim22]{kim2}
J. S. Kim, Jacobi--Trudi formulas for flagged refined dual stable Grothendieck polynomials, Algebraic Combin. {\bf 5} (2022), 121--148.

\bibitem[LP07]{lp}
T.~Lam and P.~Pylyavskyy, 
Combinatorial Hopf algebras and K-homology of Grassmannians,
{Int. Math. Res. Not.} Vol. 2007, (2007), rnm 125.
      
\bibitem[LP88]{lasp}
A. Lascoux and P. Pragacz, Ribbon Schur functions, European J. Combin. {\bf 9} (1988), 561--574.

\bibitem[LS82]{ls}
A. Lascoux and M.-P. Sch\"utzenberger, Polyn\^omes de Schubert, C. R. Acad. Sci. Paris S\'er. I Math. {\bf 294} (1982), 447--450.

\bibitem[LS83]{ls0}
A. Lascoux and M.-P. Sch\"utzenberger, Symmetry and flag manifolds, Lecture Notes in Mathematics Vol. 996, (1983), 118--144.

\bibitem[Len00]{lenart}
C. Lenart, Combinatorial aspects of the K-theory of Grassmannians, Ann. Comb. {\bf 4} (2000), 67--82.

  
\bibitem[Mac92]{macdonald1}
I. G. Macdonald, Schur functions: theme and variations, S\'eminaire Lotharingien de Combinatoire {\bf 28} (1992), B28a-39.

\bibitem[Mol98]{mol98}
A. Molev, Factorial supersymmetric Schur functions and super Capelli identities. Amer. Math. Soc. Transl. (1998), 109--138.

\bibitem[MPS21]{pech}
C. Monical, O. Pechenik, T. Scrimshaw, Crystal structures for symmetric Grothendieck polynomials. Transformation Groups, {\bf 26} (2021), 1025--1075.

\bibitem[MS13]{motegi}
K. Motegi and K. Sakai, Vertex models, TASEP and Grothendieck polynomials, J. Phys. A, 46 (2013), no. 35, p. 26.
  
\bibitem[NP26]{np}
S. Nguyen, and P. Pylyavskyy, Temperley-Lieb Immanants of Ribbon Decomposition Matrices, arXiv:2605.12880, 2026.


\bibitem[Tal12]{tal}
K. Talaska, Determinants of weighted path matrices, arXiv:1202.3128, 2012.


\bibitem[Wac85]{wachs}
M. L. Wachs, Flagged Schur functions, Schubert polynomials, and symmetrizing operators, J. Combin. Theory Ser. A {\bf 40} (1985), 276--289.
  
\bibitem[Yel17]{yel17}
D. Yeliussizov, Duality and deformations of stable Grothendieck polynomials, J. Algebraic Combin. {\bf 45} (2017), 295--344.

\bibitem[Yel19]{yel19}
D. Yeliussizov, Symmetric Grothendieck polynomials, skew Cauchy identities, and dual filtered Young graphs, 
J. Combin. Theory Ser. A {\bf 161} (2019), 453--485.

\bibitem[Yel20]{yel20}
D. Yeliussizov, Positive specializations of symmetric Grothendieck polynomials. Adv. Math.
{\bf 363} (2020), article 107000.

 
  
  


\end{thebibliography}
\end{document}